\documentclass[letterpaper,11pt]{amsart}
\usepackage{etex}
\usepackage[all]{xy}
\usepackage{amsmath}
\usepackage{amssymb,amsthm}
\usepackage{amscd,amsfonts}
\usepackage{stmaryrd,mathabx,tikz,tikz-cd,mathrsfs,float,enumitem,accents}
\usepackage{ulem}
\usepackage{graphpap,color}
\usepackage{pstricks}
\usepackage{cancel}
\usepackage[verbose]{hyperref}

\numberwithin{equation}{section}

\def\wh#1{\widehat{#1}}
\def\wc#1{\widecheck{#1}}
\def\ti#1{\widetilde{#1}}
\def\ov#1{\overline{#1}}
\def\ud#1{\underline{#1}}
\def\udt#1{\underaccent{\dot}{#1}}
\def\jia#1{\left({#1}\right)}
\def\ut#1{\undertilde{#1}}

\def\tn{\textnormal}
\def\ts{\textsf}

\def\inn{\!\in\!}
\def\eq{\!=\!}
\def\bsl{\backslash}
\def\lra{\longrightarrow}

\def\lr#1{\langle{#1}\rangle}
\def\blr#1{\big\langle{#1}\big\rangle}

\def\fk#1{\mathfrak{#1}}

\def\mwt{\fM^{\tn{wt}}}

\def\mdv{\mathfrak D}

\def\wt{\tn{wt}}

\def\mn{\tn{mn}}
\def\mc{\tn{mc}}
\def\sm{\tn{sm}}

\def\der{\tn{dr}}

\def\ny{\mu}
\def\tw{\eta}

\def\A{\mathbb A}
\def\fB{\mathfrak B}

\def\cC{\mathcal C}

\def\fD{\mathfrak D}

\def\fE{\mathfrak E}
\def\bE{\mathbf E}

\def\cE{\mathcal E}
\def\sfE{\mathsf C}

\def\bF{\mathbf F}
\def\cF{\mathcal F}
\def\fF{\mathfrak F}

\def\bH{\mathbf H}
\def\bbI{\mathbb I}
\def\fI{\mathfrak I}

\def\cK{\mathcal K}

\def\fM{\mathfrak M}

\def\cM{\mathcal M}

\def\fN{\mathfrak N}

\def\sO{\mathscr O}

\def\P{\mathbb P}

\def\bT{\mathbf T}
\def\fT{\mathbf T}
\def\cU{\mathcal U}
\def\fU{\mathfrak U}
\def\cV{\mathcal V}
\def\fV{\mathfrak V}

\def\Z{\mathbb Z}

\def\fe{\mathfrak e}
\def\ff{\mathfrak f}

\def\fp{\mathfrak p}
\def\fr{\mathfrak r}

\def\bu{\mathbf u}

\def\bv{\mathbf v}

\def\bfw{\mathbf w}

\def\De{\Delta}
\def\Up{\Upsilon}
\def\Ga{\Gamma}
\def\Th{\Theta}
\def\Si{\Sigma}
\def\La{\Lambda}
\def\Om{\Omega}

\def\al{\alpha}
\def\ga{\gamma}
\def\de{\delta}

\def\ka{\kappa}

\def\th{\theta}
\def\ve{\varepsilon}

\def\vs{\varsigma}
\def\u{\upsilon}
\def\ze{\zeta}
\def\vr{\varrho}

\def\cht{{\mathbf{m}}}

\def\lt{{\mathfrak t}}

\def\ls{\mathfrak s}

\def\tf{\ov{\tn{tf}}}
\def\rtf{\underline{\tn{tf}}}
\def\wt{\tn{wt}}

\def\ex{\tn{gr}}
\def\dg#1{\rho_{#1}}

\def\core{\tn{cor}}

\def\set{locally Euclidean}

\def\PT{\tn{pt}}
\def\dt{\tn{d}}
\def\Dm{\tn{dom}}
\def\ND{\tn{ndm}}
\def\ms{0}
\def\mfi{\mathfrak i}
\def\score{\tn{El}}
\def\be{\jmath}
\def\cj{\tn{cnj}}

\usetikzlibrary{arrows,decorations.pathmorphing,backgrounds,positioning,fit,petri}
\usetikzlibrary{decorations.pathreplacing}
\usetikzlibrary{patterns}
\tikzcdset{arrow style=tikz, diagrams={>= stealth }}

\newtheorem{thm}{Theorem}[section]
\newtheorem{prp}[thm]{Proposition} 
\newtheorem{lmm}[thm]{Lemma}  
\newtheorem{crl}[thm]{Corollary}
\newtheorem{cnj}[thm]{Conjecture}
\newtheorem{dfn}[thm]{Definition}
\newtheorem{cnv}[thm]{Convention}

\theoremstyle{remark}
\newtheorem{rmk}[thm]{Remark}
\newtheorem{eg}[thm]{Example}

\title
[STF theory and resolution of genus two stable map moduli]
{A theory of stacks with twisted fields and resolution of moduli of genus two stable maps}
\date{}
\author{Yi Hu}
\address{Department of Mathematics, University of Arizona, USA.}
\email{yhu@math.arizona.edu}
\author{Jingchen Niu}
\address{Department of Mathematics, University of Arizona, USA.
\newline\indent
Current address:
Santakatu 12 G 109, Helsinki 00180, Finland}
\email{jingchen.niu@mail.huji.ac.il}

\begin{document}

\begin{abstract} 
 We construct a smooth algebraic stack of tuples consisting of
  genus two nodal curves, simple effective divisors away from the nodes, and twisted fields.
  It provides a desingularization of the moduli of genus two stable maps to projective spaces. 
 The construction is based on  systematic application of the theory of stacks with twisted fields (STF),
 which has its prototype appeared in~\cite{HLN, g1modular} and
 is fully developed in this article.
 As a byproduct of the STF theory,
 we also obtain a novel desingularization of the moduli of genus one stable maps to projective spaces, which is isomorphic to the blowup that reverses the order used by Vakil-Zinger and Hu-Li.
 The results of this article are the second step of a program toward the resolutions of the moduli of stable maps of higher genera. 
\end{abstract}

\maketitle 
\setcounter{tocdepth}{2}
\tableofcontents

\section{Introduction}\label{Sec:Intro}
This paper is  sequel to \cite{g1modular}.
The series aims to resolve the singularities of the  moduli $\ov M_g(\P^n,d)$,
which possess arbitrary singularities when varying $g$, $n$, and $d$ (c.f.~\cite{V}).
The problem of resolution of singularities is arguably one of the hardest problems in algebraic geometry (\cite{Hironaka64a, Hironaka64b, deJong96, K07}).

Concretely,
we aim to construct a Deligne-Mumford moduli stack $\ti M_g(\P^n,d)$, along with a proper morphism $$\varpi:\ti M_g(\P^n,d)\lra\ov M_g(\P^n,d),$$ enjoying the following properties.
\begin{enumerate}
[leftmargin=*,label=($\mathsf{mr}_{\arabic*}$)]
\item\label{Cond:MainSmooth} $\ti M_g(\P^n,d)$ has smooth irreducible components and admits at worst 
normal crossing singularities.

\item \label{Cond:MainMainComponent}
Let $\ov M_g(\P^n,d)^{\mc}\!\subset\!\ov M_g(\P^n,d)$ be the component whose general points are the stable maps with smooth domain curves,
and $\ov M_g(\P^n,d)^{\mc}_{\sm}$ be its smooth locus (if it is nonempty).
Denote by $\ti M_g(\P^n,d)^{\mc}$ and $\ti M_g(\P^n,d)^{\mc}_{\sm}$ their preimages in $\ti M_g(\P^n,d)^{\mc}$, respectively.
Then,
$$
\ti M_g(\P^n,d)^{\mc}\,\big\bsl 
\ti M_g(\P^n,d)^{\mc}_{\sm}
$$
is a simple normal crossing divisor of $\ti M_g(\P^n,d)^{\mc}$.

\item\label{Cond:MainBirational}
$\ti M_g(\P^n,d)^{\mc}_{\sm}$ is isomorphic onto $\ov M_g(\P^n,d)^{\mc}_{\sm}$, 
hence the restriction of $\varpi$ to $\ti M_g(\P^n,d)^{\mc}$ is birational.

\item\label{Cond:MainLocallyFree}
Let
$$
\ti\pi:\ti\cC\lra\ti M_g(\P^n,d)^{\mc},\qquad
\ti\ff:\ti\cC\lra\P^n,
$$
be the pullback of the universal family $$\pi:\cC\lra\ov M_g(\P^n,d),\qquad\ff:\cC\lra\P^n$$ of $\ov M_g(\P^n,d)$.
Then, for every integer $k\!\ge\!1$,
the derived object $\mathbf{R}\ti\pi_*\ti\ff^*\sO_{\P^n}(k)$ is locally diagonalizable in the sense of \cite[Definitions~3.8 \& 3.9]{HL11},
which implies the direct image sheaf $\ti\pi_*\ti\ff^*\sO_{\P^n}(k)$ is locally free.
\end{enumerate} 

Given $g$, $n$ and $d$, if such $\ti M_g(\P^n,d)\big/\ov M_g(\P^n,d)$ exists, we call it a \ts{modular resolution} of $\ov M_g(\P^n,d)$,
and call $\ti M_g(\P^n,d)^{\mc}$ and $\ov M_g(\P^n,d)^{\mc}$ the \ts{main components} of $\ti M_g(\P^n,d)$ and $\ov M_g(\P^n,d)$, respectively. 
Here,
Conditions\ref{Cond:MainSmooth}-\ref{Cond:MainBirational} are natural  in view of the theory of resolution,
while Condition~\ref{Cond:MainLocallyFree} is desirable for applications in computational Gromov-Witten 	theory.

\subsection
{Local defining equations of the moduli of genus 2 stable maps: motivation}

Let  $\fD_g$ be the smooth Artin stack consisting  of pre-stable pairs
$\big(C, D\big)$ where~$C$  are connected genus $g$ nodal curves 
and $D\!\subset\!C$ are simple divisors away from the nodes of $C$. A pair $(C,D)$ is \ts{pre-stable} if any smooth rational
component of $C$ missing the divisor $D$ has at least three nodes.

For every $[C,u]\inn \ov M_g(\P^n,d)$, we have
$\big(C,u^{-1}(H)\big)\inn\fD_g$ for a generic hyperplane $H\!\subset\!\P^n$.
For every neighborhood $U\!\subset\!\ov M_g(\P^n,d)$ of $[C,u]$, we set
\begin{align}\label{Eqn:H}
	\bH_U:=
	\big\{\,H\inn\tn{Gr}(n,n\!+\!1)\,:\,
	\big(C', (u')^{-1}(H)\big)\inn\fD_g~
	\forall~[C',u']\inn U\,
	\big\}\,,
\end{align}
which is nonempty as long as $U$ is small.
Each $H\inn \bH_U$ gives rise to a morphism 
\begin{align}\label{Eqn:Mg_tauto}
	f_{U,H}:\,
	U\lra \fD_g, \qquad
	[C,\bu]\mapsto \big(C, \bu^{-1}(H)\big).
\end{align}

By (the proof of)  \cite[Theorem 2.17]{HL10} (cf.~also  \cite[Theorem 2.5]{HLN}),
we can cover $\fD_g$ by open charts $\{\cV\}$ such that 
$\{U\!\times_{\fD_g}\!\cV\}$ covers $\ov M_g(\P^n,d)$; further, there exists a trivial bundle $\cE_\cV$ 
of rank $(d+1)n$
and a trivial bundle $\cF_{\cE_\cV}$ of rank $gn$ over  $\cE_\cV$ equipped with a section $\Phi_{\cE_\cV}$ such that 
$\{U\!\times_{\fD_g}\!\cV\}$ is canonically embedded in $\cE_\cV$ as the zero locus  $(\Phi_{\cE_\cV}\eq 0)$.
When $g\eq 1$ (\cite{HL10}) and $g\eq 2$ (\cite{ HLN}), we can choose trivializations of the bundles such that  $(\Phi_{\cE_\cV}\eq 0)$
is equivalent to a matrix equation,  
\begin{align}\label{Eqn:str_hom_intro}
	{\varphi_\cV}\cdot[\bfw^1,\cdots, \bfw^n]=0,\qquad
\tn{where}\quad {\varphi_\cV}=[\chi_{ij}]_{g \times d}.
\end{align}
Here, the entries (and $2\!\times\! 2$ minors if $g\eq 2$) of ${\varphi_\cV}$ are products of regular functions  whose
vanishing loci are irreducible components of nodal curves (and/or the loci of Weierstrass and conjugate points in $\cV\!\subset\!\fD_2$ if $g\eq 2$),
and $\bfw^i$ are vectors of free variables.

Suppose we have a sequence of blowups of $\fD_g$ along {\it smooth closed centers,}
\begin{equation}\label{base blowups}
	\widetilde\fD_g \lra \cdots \lra \fD_g.
\end{equation}
Then $\widetilde\fD_g$ can be covered by some charts $\{\tilde\cV/\cV\}$; over the pullback bundle
$\tilde\cE_{\tilde\cV}$, the  section $\Phi_{\cE_\cV}$ pulls back to a section $\Phi_{\tilde\cE_{\tilde\cV}}$, and its equivalent matrix equation takes the pullback form
\begin{align*}
	\ti\varphi_\cV\cdot[\bfw^1,\cdots, \bfw^n]=0;
\end{align*}
further, $U\!\times_{f_{U,H};\fD_g} \widetilde\fD_g$ is covered by 
$\{U\! \times_{\fD_g}\!\tilde\cV\}$, and each  $\cU \!\times_{\fD_g}\!\tilde\cV$ admits the induced embedding 
in~$\tilde\cE_{\tilde\cV}$ as the zero locus  $(\Phi_{\tilde\cE_{\tilde\cV}}\!\eq 0)$.
We will show the natural morphisms
\begin{align*}
	\ti U_H:=
	U\!\times_{f_{U,H};\fD_g} \widetilde\fD_g
	\lra U,\qquad
	H\in \bH_U,
\end{align*}
are isomorphic to each other;
in addition,
these $\ti U_H/U$ glue to form a morphism $$\ti M_g(\P^n,d)\lra\ov M_g(\P^n,d).$$
If $\ti\varphi_\cV$ can be diagonalized in the sense of~\cite[Definition 3.2]{HL11} such that
the diagonal entries $z_{ii}$ are all products of regular functions whose vanishing loci are smooth, 
then  $\ti M_g(\P^n,d)\!\lra\!\ov M_g(\P^n,d)$ is a modular resolution of $\ov M_g(\P^n,d)$.

The existence of such  $\widetilde\fD_1$ is a consequence of the main result of \cite{HL10}, and
the existence of such~$\widetilde\fD_2$ is established in \cite{HLN}.

To extend the above to a genus higher than 2, the global smooth centers in the sequence of \eqref{base blowups}
quickly become formidable to describe. Fortunately, the local
matrix equations seem relatively  more manageable. 
This motivates us to seek a local-to-global approach,
avoiding the difficulty of explicitly prescribing  the global smooth blowup centers in (\ref{base blowups}).
This strategy can be summarized as follows:
first, we classify the local matrix equations by their combinatorial types;  then, we perform local blowups
in terms of adding \ts{twisted fields}, according to the type of the matrix equation, 
to the data $(C, D)$. The way that these twisted fields are added is coherently designed 
so  that these local constructions can
be naturally glued together to form a smooth stack $\ti\fD_g^{\rm tf}$, which is isomorphic to the desired  $\widetilde\fD_g$. Needless to say,
the gluing implies that there exist corresponding global blowups along ideal sheaves, 
albeit intentionally not described explicitly.
Here, each type of the matrix equation  also determines a stratum, leading to the notion of
\ts{locally Euclidean stratification (LES)}.

The existence of such $\ti\fD_1^{\rm tf}$ under this local-to-global  approach
is a consequence of the main results of \cite{g1modular}.
The existence of such $\ti\fD_2^{\rm tf}$ under the similar idea, albeit much more complex, 
will be  constructed in the current article.
Furthermore, we lay a theoretical foundation in this work, aiming to apply it
to the case of genera higher than 2, in particular the case of $g=3$.

\subsection
{Introduction to {\set} stratifications and treelike structures} 

We take $\fM$ to be either $\fD_g$ or any of the  intermediate smooth moduli in 
\eqref{base blowups}.  It can be covered by smooth charts; over
every chart $\cV$, 
by abuse of notation, we write the pullback of the matrix of (\ref{Eqn:str_hom_intro}) still as $\varphi_\cV$.
Assume each entry of $\varphi_\cV$ contains a product
of regular functions whose vanishing loci are smooth (i.e.~a product of {\it local parameters}) as a factor.
An LES of $\fM$ is a stratification whose strata over any chart can be
described in terms of the intersections of these vanishing loci in a coherent way.
To each stratum, we assign a rooted tree such that on each chart $\cV$ centered at the given stratum, the edges of the tree are labeled by
distinct local parameters aforementioned.
This assignment of rooted trees to the strata of $\fM$, 
termed as a \ts{treelike structure} on $\fM$, is compatible with the LES.  

Coming with the treelike structure are the products of local parameters on the linear paths 
from the root to the leaves of the rooted trees, called the \ts{modular monomials}.
On a chart $\cV$ centered at a given stratum of $\fM$,	
each  entry of $\varphi_\cV$ contains a modular monomial as a factor.
The set of all the modular monomials on $\cV$ is called a \ts{tautological set of monomials}; it lives on the chart $\cV$
and varies when varying $\cV$.
(The pullback of) a tautological set of monomials
is said to admit a \ts{divisibly minimal element}  if there exists one member $m_0$ such that it divides all other members $m$,  that is, $m/m_0$ is a product of the regular functions.	

\subsection{Introduction to stacks with twisted fields} 
Given a smooth stack $\fM$ equipped with an LES,
after endowing it with a treelike structure $\La$,
we obtain two new smooth stacks 
$$\fM_{\La}^{\tf}\qquad\tn{and}\qquad \fM_{\La}^{\rtf}\,;$$
the construction is referred to as 
the theory of \ts{stacks with twisted fields (STF)}.
The following key statement of the STF theory is deduced from Theorem~\ref{Thm:tf_smooth}~\ref{Cond:smooth_tf} and Theorem~\ref{Thm:tf_smooth_revert}~\ref{Cond:smooth_tf_revert}.

\begin{thm}\label{Thm:TSTF}
Let $\fM\eq\bigsqcup_{\al\in A}\fM_\al$ be a smooth algebraic stack with an LES  as in Definition~\ref{Dfn:G-adim_fixture}, 
endowed with a treelike structure $\La$ as in Definition~\ref{Dfn:Treelike_structure} that assigns a rooted tree $\tau_\al$ to each stratum $\fM_\al$.

Then, there exists a smooth algebraic stack $\fM^{\tf}_\La$ parameterizing the objects of $\fM$ along with the root-to-leaf (RL) twisted fields with respect to $\La$, as well as a proper and birational morphism $$\ov\varpi:\fM^{\tf}_\La\lra\fM,$$
whose restriction to the preimage of $\fM^\mn$ $(\subset\!\fM)$ as in~(\ref{Eqn:boundary}) is an isomorphism.
Furthermore, $\fM^{\tf}_\La$  admits an induced LES.

Likewise, there exists a smooth algebraic stack $\fM^{\rtf}_\La$ parameterizing the objects of $\fM$ along with the leaf-to-root (LR) twisted fields with respect to $\La$, as well as a proper and birational morphism $$\ud\varpi:\fM^{\rtf}_\La\lra\fM,$$
whose restriction to the preimage of $\fM^\mn$ is also an isomorphism.
Furthermore, $\fM^{\rtf}_\La$  admits an induced LES.
\end{thm}

We call $\fM^{\tf}_\La$ and $\fM^{\rtf}_\La$ the \ts{stacks of $\fM$ with RL- and LR-twisted fields}, respectively.
These two stacks are both referred to as stacks with twisted fields (STFs) when the context is clear.

The main benefit of the STF theory is the following.

\begin{thm}
	\label{Thm:Tauto_monom}
Let $\fM$ and $\La$ be as in Theorem~\ref{Thm:TSTF}.
On each of $\fM^{\tf}_\La$ and $\fM^{\rtf}_\La$,
the pullback of every tautological set of monomials admits a divisibly minimal element in the sense of Definition~\ref{Dfn:Tauto_monom}.
\end{thm}

The precise statements are given in Corollaries~\ref{Crl:prod_RL} and~\ref{Crl:prod_LR}.
Theorem~\ref{Thm:Tauto_monom} brings us closer to the final stack of (\ref{base blowups}).

\subsection{Applications to genus one and two}
\subsubsection{The case of $g\eq 2$}
Starting from $\fM\eq\fD_2$,
progressively,
we construct nine STFs:
in each step, we construct a new treelike structure on the stack obtained in the previous step, then apply Theorem~\ref{Thm:TSTF} to obtain a new STF.
The final stack obtained is denoted by
$\widetilde \fD_2^{\tn{tf}}$,
which is isomorphic to $\ti \fD_2$ of (\ref{base blowups}). 

\begin{thm}
\label{Thm:Main}  There exist a smooth algebraic stack $\ti \fD_2^{\tn{tf}}$  parameterizing stable
tuples of nodal curves of genus two, simple effective divisors away from the nodes, 
and twisted fields,
as well as 
a proper and birational
forgetful morphism $\ti\fD^{\tn{tf}}_2\!\lra\!\fD_2$,
whose restriction to the preimage of  $\{(C,D)\inn\fD_2:C~\tn{is~smooth}\}$ is isomorphic.

Moreover, 
there exists a modular resolution $\ti M_2^{\tn{tf}}(\P^n,d)\big/\ov M_2(\P^n,d)$ (i.e.~a proper morphism $\ti M_2^{\tn{tf}}(\P^n,d)\!\lra\!\ov M_2(\P^n,d)$ of a Deligne-Mumford stack $\ti M_2^{\tn{tf}}(\P^n,d)$, satisfying~\ref{Cond:MainSmooth}-\ref{Cond:MainLocallyFree}) 
such that for every small open $U\!\subset\!\ov M_2(\P^n,d)$ and every $H\inn \bH_U$ as in (\ref{Eqn:H}),
with $f_{U,H}\!:U\!\lra\!\fD_2$ as in (\ref{Eqn:Mg_tauto}), 
there exists a commutative diagram
\[
\begin{tikzcd}[column sep=.5em,row sep=1.5em]
	U\times_{f_{U,H};\fD_2} \ti{\fD}^{\tn{tf}}_2  \arrow[dr,""] \arrow[rr,"\sim"] && U\times_{\ov M_2(\P^n,d)}\ti M_2^{\tn{tf}}(\P^n,d)
	\arrow[dl,""]
	\\
	& U 
\end{tikzcd}
\]
where the horizontal arrow is an isomorphism (of stacks over $U$),
and the other two arrows are natural morphisms.
\end{thm}

The stack $\ti{\fD}^{\tn{tf}}_2$ should not be isomorphic to the blowup stack $\ti\fD_2$ constructed in~\cite{HLN}.
In fact,
being a typical phenomenon of a resolution problem,
the order of the blowups performed in~\cite{HLN}  versus this article are not identical,
which should lead to different resolutions of $\ov M_2(\P^n,d)$ in general.

In addition, with $(\ti\pi,\ti\ff)$ as in Condition~\ref{Cond:MainLocallyFree},
\begin{align*}
\blr{\,e\big(\ti\pi_*\ti\ff^*\sO_{\P^n}(k)\big)\,,\,[\ti M_2^{\tn{tf}}(\P^n,d)^{\mc}]\,}
\end{align*}
should equal the expected reduced genus 2 Gromov-Witten (GW) invariants of the corresponding complete intersection, parallel to~\cite[(1.4)]{VZ} and~\cite[(1.7)]{LZ}.
The reduced genus 1 GW-invariants,
as well as its comparison with the standard genus 1 GW-invariants, 
are introduced in~\cite{Zi2,Zi1} and further studied in~\cite{LZ,VZ,CL,CM,LO1,LO2,LL},
and lead to important results such as A.~Zinger's proof~\cite{Zi3} of the prediction of~\cite{BCOV}
for genus 1 GW-invariants of a quintic 3-fold.

\subsubsection{The case of $g\eq 1$}
The STF theory also gives a novel  application when $g\eq 1$.

Recall $\fM_1^\wt$ denotes the moduli stack of stable weighted genus 1 curves as in 
\cite{HL10} and \cite{g1modular}. It comes equipped with the closed substacks $\Theta_k$
whose general points are smooth weight-0 genus~1 curves attached with $k$ rational tails
for all $k \inn {\mathbb Z_{>0}}$.
The stack $\widetilde{\fM}_1^\wt$ constructed in \cite{HL10} is 
the iterated blowup of 
$\fM_1^\wt$ along the proper transforms of $\Theta_2,\Theta_3, \cdots$,
which is finite on each connected component of $\fM^\wt_1$.
As revealed in~\cite{HL10},
the restriction of the morphism
\begin{align}\label{Eqn:HL10}
\ti M_1(\P^n,d):=
{\ov M}_1(\P^n,d)\times_{\fM^\wt_1}\ti\fM_1^{\wt}\;
\lra\,{\ov M}_1(\P^n,d)
\end{align}
to the main component ${\ti M}_1(\P^n,d)^\mc$  gives an algebro-geometric approach to Vakil-Zinger's virtual blowup in~\cite{VZ},
where the morphism ${\ov M}_1(\P^n,d)\!\lra\!\fM^\wt_1$ is given by
\begin{align}\label{Eqn:M1wt_tauto}
	\ov M_1(\P^n,d)\lra \fM_1^\wt, \qquad
	[C,\bu]\mapsto\big(\,C,\, c_1\big(\bu^*\sO_{\P^n}(1)\big)\,\big)\,.
\end{align}

On the other hand,
as shown in~\cite{g1modular},
$\fM^\wt_1$ naturally carries an LES 
determined by the topology of the nodal curves and the distribution of the weights, 
as well as a treelike structure that assigns to each stratum the corresponding weighted dual tree (in the sense of \cite[\S2.2]{g1modular}; also known as the terminally weighted tree in \cite{HL10}). 
We apply
Theorem~\ref{Thm:TSTF} to $\fM^\wt_1$ with respect to this treelike structure only {\it once} and obtain two STFs
$$ (\fM_1^{\wt})^{\tf}\qquad\tn{and}\qquad(\fM_1^{\wt})^{\rtf}.$$

\begin{thm} \label{Thm:g1} 
	With notation as above, we have the following.
	\begin{enumerate}[leftmargin=*,label=(\arabic*)]
		\item
		\label{Cond:g1_tf_wt}
		$(\fM_1^{\wt})^{\tf}\big/\fM_1^{\wt}$ is isomorphic to 
		the blowup stack $\widetilde{\fM}_1^\wt\big/\fM_1^{\wt}$ constructed in \cite{HL10}. 
		
		\item 
		\label{Cond:g1_tf_stable_map}
		The morphism
		\begin{align*}
			\ov M_1^{\tf}(\P^n,d):=
			{\ov M}_1(\P^n,d)\times_{\fM^\wt_1}(\fM_1^{\wt})^{\tf}\;
			\lra\,{\ov M}_1(\P^n,d)
		\end{align*}
		is a modular resolution of $\ov M_1(\P^n,d)$,
		and is isomorphic to the morphism (\ref{Eqn:HL10}).
		Its restriction to the main component thus gives an algebro-geometric approach to Vakil-Zinger's virtual blowup in~\cite{VZ}.
		
		\item 
		\label{Cond:g1_rtf_wt} $(\fM_1^{\wt})^{\rtf}\big/\fM_1^{\wt}$ is  isomorphic to the iterated blowup of $\fM_1^\wt$ along the proper transforms of
		$$\cdots,\ \Th_3\,,\ \Theta_2\,,$$
		which is finite on each connected component of $\fM_1^\wt$. 
		Moreover, there does not exist any morphism between  $(\fM_1^{\wt})^{\tf}\big/\fM_1^{\wt}$ and $(\fM_1^{\wt})^{\rtf}\big/\fM_1^{\wt}$.
		
		\item 
		\label{Cond:g1_rtf_stable_map}
		Similarly, 
		$$
		\ov M_1^{\rtf}(\P^n,d):=
		{\ov M}_1(\P^n,d)\times_{\fM^\wt_1}(\fM_1^{\wt})^{\rtf}\;
		\lra\,{\ov M}_1(\P^n,d)
		$$ 
		is a modular resolution of $\ov M_1(\P^n,d)$.
		Moreover, 
		for $d\!>\!2$,
		there does not exist any morphism between $$
		\ov M_1^{\tf}(\P^n,d)^\mc\big/ {\ov M}_1(\P^n,d)^\mc
		\qquad\tn{and}\qquad\ov M_1^{\rtf}(\P^n,d)^\mc\big/{\ov M}_1(\P^n,d)^\mc,$$
		where $\ov M_1^{\rtf}(\P^n,d)^\mc$ denotes the main component of $\ov M_1^{\rtf}(\P^n,d)$.
	\end{enumerate}
\end{thm}

The first two statements of Theorem~\ref{Thm:g1} are \cite[Theorem~1.1]{g1modular} and \cite[Corollary~1.2]{g1modular}, respectively, restated using the current notation.	
The last two statements of Theorem~\ref{Thm:g1} are proved in~\S\ref{Sec:Proof_g1}.

\subsection
{Further and related works}

\begin{cnj} Staring from $\fD_3$, there exist finitely many stacks 
carrying  treelike structures or their analogues so that the STF theory can be progressively applied to
obtain $\widetilde\fD_3$ as in \eqref{base blowups}.	
\end{cnj}

In \cite{g1log} D.~Ranganathan, K.~Santos-Parker, and J.~Wise
provide a logarithmic geometry theoretical approach to Vakil-Zinger's desingularization \cite{VZ}.
Using logarithmic geometry and Gorenstein curves,  L.~Battistella and F.~Carocci  in \cite{g2log},
construct a desingularization of ${\ov M}_2(\P^n,d)^\mc$.
 A related work on desingularizations of sheaves and reduced GW invariants can be found in~\cite{rmmp} by
	A.~Rabano, E.~Mann, C.~Manolache, and R.~Picciotto.
We also anticipate our twisted field stack $ \ti \fD_2^{\tn{tf}}$ in Theorem \ref{Thm:Main}
admits interpretation and construction in terms of logarithmic algebraic geometry.

Based on the application of the blowup construction of~\cite{HLN} to the moduli of quasi-maps (see~\cite{LLO} by S.~Lee, M.-L. Li, and J.~Oh),
the STF theory should be applied to the moduli of quasi-maps equally well to obtain their modular resolutions. Rigorous formulations of these results will appear in the forthcoming work.

\smallskip
\textbf{Acknowledgments.} We thank Dawei Chen, Qile Chen, Jack Hall, and YP Lee for the valuable discussions. 
We thank the anonymous referee for the thoughtful comments, suggestions, and corrections.
During the final revision of this paper, the first named author was a visiting professor at Great Bay University, whose generous hospitality and support are gratefully acknowledged.

\section{Stacks with twisted fields (STF)}
\label{Sec:tf}

\subsection{Locally Euclidean stratifications (LES)}
\label{Subsec:STF}
We begin a class of stacks to which the theory of stacks with twisted fields (STF) will be developed.

\begin{dfn}\label{Dfn:G-adim_fixture}
	Let $\fM$ be a smooth algebraic stack and $B$ be a finite set.
	A stratification of $\fM$:
	$$
	\fM\,=\,\fM_{\ms}\sqcup\bigsqcup_{\al\in B}\!\fM_\al\,=\, \bigsqcup_{\al\in A}\fM_\al,\qquad
	\tn{where}\quad A:=\{\ms\}\sqcup B\,,
	$$ 
	by an open substack $\fM_{\ms}$ and  connected locally closed substacks $\fM_\al$, $\al\inn B$,
	is called a \ts{{\set} stratification},
	abbreviated as an \ts{LES}, if
	\begin{enumerate}
	[leftmargin=*,label=$\bullet$]
	\item
	(\emph{atlas}) for every $\al\inn B$, 
	there exists a set $\fV_\al\eq\{\cV\}$ of affine smooth charts of $\fM$ that covers $\fM_\al$ (i.e.~$\fM_\al\!
	\subset\!
	\bigcup_{\cV\in\fV_\al}\!\!\cV$),
	known as the \ts{modular charts}, and
	
	\item
	(\emph{local parameters}) for every $\al\inn B$ and $\cV\inn\fV_\al$, there exists a subset $\{\ze_e^\cV\}_{e\in S_\al}$ of a system of local parameters on $\cV$ indexed by a nonempty finite set $S_\al$,
	known as the \ts{modular parameters},	satisfying the following compatibility conditions:
	
	\item (\emph{compatibility conditions})
	for every $\al,\al'\inn B$, $\cV\inn\fV_{\al}$, and $\cV'\inn\fV_{\al'}$,
	if	$\fM_{\al'}\!\cap\!\cV\!\cap\!\cV'\!\ne\!\emptyset,$
	then
	there exists an injection $ \iota_{\cV,\cV'}\!:S_{\al'}\!\hookrightarrow\!S_{\al}$ such that
	\begin{align}
		&\ze^{\cV'}_{e'}\big/\ze^\cV_{\iota_{\cV,\cV'}(e')}
		\,\in\,
		\Ga(\sO^*_{\cV\cap\cV'})\quad\forall~e'\in S_{\al'}\,,\qquad\tn{and}
		\label{Eqn:loc_Euc_transition}
		\\
		&
		\fM_{\al'}\cap\cV\cap\cV'\,=\,
		\big\{~\ze_e^{\cV}\!\!\eq 0\ \;
		\forall\,e\in \iota_{\cV,\cV'}(S_{\al'})\,;\ \;
		\ze_e^{\cV}\!\!\!\neq\! 0\ \;
		\forall\,e\in S_\al\big\bsl\iota_{\cV,\cV'}(S_{\al'})~\big\} \, \cap\,\cV'
		\,.
		\label{Eqn:loc_Euc_equation}
	\end{align}
	\end{enumerate}
\end{dfn}

For conciseness, the open stratum $\fM_{\ms}$ and its complement in $\fM$ are denoted by 
\begin{align}\label{Eqn:boundary}
	\fM^{\mn}=\fM_{\ms}\qquad\tn{and}\qquad
	\De=\fM\bsl\fM^{\mn}=\bigsqcup_{\al\in B}\fM_\al,
\end{align}
respectively.
Such $\De$ is called the \ts{boundary} of $\fM$,
which is closed in $\fM$.	

In Definition~\ref{Dfn:G-adim_fixture},
the requirement 
(\ref{Eqn:loc_Euc_transition}) describes a transition of the local parameters of $\cV'$ into $\cV$,
whereas
(\ref{Eqn:loc_Euc_equation}) describes the equations for the locally closed subset $\fM_{\al'}\!\cap\!\cV\!\cap\!\cV'$ in the chart~$\cV$.	
Particularly,
by setting $\al'\eq\al$ and $\cV'\eq\cV$,
we observe that $\ze^{\cV}_e$ and $\ze^\cV_{\iota_{\cV,\cV}(e)}$ differ by a unit on $\cV$ for all $e\inn S_\al$,
hence $\iota_{\cV,\cV}$ is the identity map, and
\begin{align}\label{Eqn:M_strata_local'}
	\fM_\al\cap\cV=
	\{\,\ze_e^{\cV}\!\eq 0\ \ 
	\forall\,e\inn S_\al\,\}
\end{align}
holds for all $\al\inn B$ and $\cV\inn\fV_\al$.
Therefore,
if both $\fM_{\al'}\!\cap\!\cV\!\cap\!\cV'$ and $\fM_{\al}\!\cap\!\cV\!\cap\!\cV'$ are nonempty for some modular charts $\cV\inn\fV_\al$ and $\cV'\inn\fV_{\al'}$,
then they must be given by the same ideal of the restriction of (\ref{Eqn:M_strata_local'}) to $\cV\!\cap\!\cV'$,
hence
$\al\eq\al'$.

By (\ref{Eqn:M_strata_local'}), 
we have $|S_\al|\eq\tn{codim}\,\fM_\al$, so the nonemptyness of $S_\al$ implies $\tn{codim}\,\fM_\al\!>\!0$ for all $\al\inn B$.
Therefore, $\fM^{\mn}$ is nonempty.
The set $B$, however, is allowed to be empty, so $\De$ is possibly empty.

Observe that (\ref{Eqn:M_strata_local'}) and
(\ref{Eqn:loc_Euc_transition}) imply that
for every $\al\inn B$ and $\cV,\cV',\cV''\inn\fV_\al$,
we have
\begin{align*}
	\iota_{\cV,\cV'}\circ\iota_{\cV',\cV''}=
	\iota_{\cV,\cV''}\qquad\tn{whenever}\qquad
	\fM_{\al}\cap\cV\cap\cV'\cap\cV''\ne\emptyset\,.
\end{align*}	
Taking the connectedness of each $\fM_\al$, $\al\inn B$, into consideration,
we can thus re-label
the modular parameters associated with $\fM_\al$ if necessary so that
\begin{align}\label{Eqn:loc_Euc_id}
	\iota_{\cV,\cV'}=\tn{Id}:\,
	S_\al\lra S_\al\qquad\forall\quad
	\al\in B\ \ \tn{and}\ \ 
	\cV,\cV'\in\fV_\al\  \ \tn{s.t.}\ \ 
	\fM_{\al}\cap\cV\cap\cV'\ne\emptyset\,.
\end{align}	

\begin{cnv}
	\label{Convention:loc_Euc_id}
	When an LES of a stack is considered,
	we always assume (\ref{Eqn:loc_Euc_id}) holds.
\end{cnv}

\begin{eg}\label{Eg:G-adim_fixture}
	Consider the moduli stack $\fM_g^{\rm wt}$ of stable weighted genus $g$ curves introduced in~\cite{HL10}.
	It consists of the pairs $(C,\bfw)$ of genus $g$ nodal curves $C$ and weights $\bfw\inn H^2(C,\Z)$,
	satisfying $\bfw(\Si)\!\ge\!0$ for all irreducible $\Si\!\subset\! C$,
	and every rational $\Si\!\subset\! C$ with $\bfw(\Si)\eq 0$ must contain at least three nodal points. 
	Then, $\mwt_g$ is equipped with an LES given by the dual graphs (in the usual sense; c.f.~\cite[\S23.4]{MirSym}) and the distribution of the weights.
	The boundary of $\mwt_g$ is comprised of the strata whose dual graphs have at least one edge (i.e.~whose underlying curves are nodal). 
	The modular parameters correspond to the smoothing of the nodes.
\end{eg}

Next, we delve a bit deeper into the local structure of an LES.
For every $\al,\al'\inn B$ and $\cV\inn\fV_\al$, it is possible that
$\fM_{\al'}\!\cap\!\cV$, if not empty, is disconnected even if $\fM_{\al'}$ is connected.
The connected components of $\fM_{\al'}\!\cap\!\cV$ are analyzed in Lemma~\ref{Lm:M_strata_local'} and Proposition~\ref{Prp:M_strata_local} below.

For every $\al\inn B$, $\cV\inn\fV_\al$, and $E\!\subset\!S_\al$,
let
\begin{align}\label{Eqn:M_strata_local}
	Z_{(E)}^\cV=\big\{\,
	\ze_e^\cV\eq 0~\forall~e\inn S_\al\bsl E,~~
	\ze_e^\cV\!\ne\! 0~\forall~e\inn E\,
	\big\}\,.
\end{align}
Since $\cV$ is an affine smooth chart and $\ze_e^\cV$'s are local parameters on $\cV$,
we see  $Z_{(E)}^\cV$ is connected.	

\begin{lmm}\label{Lm:M_strata_local'}
	Let $\fM$ be endowed with an LES as in Definition~\ref{Dfn:G-adim_fixture}.
	For every $\al,\al'\inn B$  and $\cV\inn\fV_\al$,
	if $\fM_{\al'}\!\cap\cV\!\ne\!\emptyset$,
	then for every connected component~$X$ of $\fM_{\al'}\!\cap\!\cV$,
	there exists $E\!\subset\!S_\al$ such that
	$$
	X=Z^\cV_{(E)}\,.
	$$
	Moreover,
	for every $\cV'\inn\fV_{\al'}$ such that $X\!\cap\!\cV'\!\ne\!\emptyset$,
	we have 
	$$
	 E\eq S_\al\big\bsl\big(\iota_{\cV,\cV'}(S_{\al'})\big)\,,
	$$
	where $\iota_{\cV,\cV'}$ is the injection in Definition~\ref{Dfn:G-adim_fixture}.	
\end{lmm}

\begin{proof}
	Fix an arbitrary connected component $X$ of $\fM_{\al'}\!\cap\!\cV$ and an arbitrary $\cV'\inn\fV_{\al'}$ such that $X\!\cap\!\cV'\!\ne\!\emptyset$.
	Take $E\eq S_\al\big\bsl\big(\iota_{\cV,\cV'}(S_{\al'})\big)$.
	Then, by Definition~\ref{Dfn:G-adim_fixture}, $X\!\cap\!\cV'$ is nonempty and open in  $Z^\cV_{(E)}$.
	Since $X$ is connected,
	for any other $\ti\cV'\inn\fV_{\al'}$ with $X\!\cap\!\ti\cV'\!\ne\emptyset$,
	we have $\iota_{\cV,\ti\cV'}(S_{\al'})\eq\iota_{\cV,\cV'}(S_{\al'})$.
	Hence $E$ is independent of the choice of $\cV'\inn\fV_{\al'}$ as long as $X\!\cap\!\cV'\!\ne\!\emptyset$, and moreover, $X$ is a nonempty open subset of $Z^\cV_{(E)}$.
	
	Suppose $X\!\subsetneq\!Z^\cV_{(E)}$.
	For any $y\inn Z^\cV_{(E)}\bsl X$,
	if $y\inn \fM^{\mn}$,
	then $\fM^{\mn}\!\cap\!Z^\cV_{(E)}$ is an open neighborhood of $y$ in $Z^\cV_{(E)}\bsl X$;
	if $y\inn \fM_{\al''}$ for some $\al''\inn B$,
	then we can apply the preceding paragraph to the connected component $Y$ of $\fM_{\al''}\!\cap\!\cV$ containing $y$
	and conclude that $Y$ is an open neighborhood of $y$ in $Z^\cV_{(E)}\bsl X$.
	In sum, $Z^\cV_{(E)}\bsl X$ is also open and nonempty, which contradicts the connectedness of $Z^\cV_{(E)}$.
	Therefore,
	we have $X\eq Z^\cV_{(E)}$.
\end{proof}

For every subset $Y\!\subset\!X$, 
we denote by 
$$
\tn{Cl}_X(Y)
~(\,\subset X\,)
$$ the closure of $Y$ in $X$.
Throughout this paper, we avoid the $\bar{~}$ notation for closures so that no confusion with the RL sequences (c.f.~Definition~\ref{Dfn:R2L_seq}) should occur.
Recall that for every topological space $X$,
the set of all connected components of $X$ is denoted by $\pi_0(X)$.

\begin{prp}\label{Prp:M_strata_local} 
	Let $\fM$ be as in Lemma~\ref{Lm:M_strata_local'}.
	Then
	for every $\al\inn B$ and $E\!\subset\!S_\al$, 
	there exists 
	$$ 
	\al_{(E)}\in A\qquad\tn{s.t.}\qquad
	Z^\cV_{(E)}\,\subset\,\fM_{\al_{(E)}}\!\cap\cV\quad\forall\ \cV\inn\fV_{\al}\,.
	$$
	Moreover, 
	one of the following two statements must hold:
	\begin{itemize}
		[leftmargin=*]
		\item either
		$\al_{(E)}\eq\ms$ (i.e.~
		$Z_{(E)}^\cV\!\subset\!\fM^{\mn}$  for all $\cV\inn\fV_{\al}$); 
		
		\item or
		$$ 
		\al_{(E)}\in B\qquad\tn{and}\qquad
		Z_{(E)}^\cV\in\pi_0\big(\fM_{\al_{(E)}}\!\cap\!\cV\big)
		\quad\forall\ 
		\cV\inn\fV_{\al}
		\,,
		$$
		and furthermore,
		for every $E'\!\subset\!E$, we have
		$$
		\al_{(E')}\in B\qquad\tn{and}\qquad
		Z_{(E')}^\cV\in\pi_0\big(\fM_{\al_{(E')}}\!\cap\!\cV\big)\quad\forall\ 
		\cV\inn\fV_{\al}\,.
		$$
	\end{itemize} 
\end{prp}

For every $e\inn S_\al$,
we write $\al_{(\{e\})}$ as $\al_{(e)}$ for conciseness.

\begin{proof}[Proof of Proposition~\ref{Prp:M_strata_local}]
	First, fix an arbitrary $\cV\inn\fV_\al$.
	Assume $Z_{(E)}^\cV\!\not\subset\!\fM^{\mn}$,
	i.e.~$Z_{(E)}^\cV\!\cap\!\De\!\ne\!\emptyset$.
	
	Then, there exists some $\al_{(E)}\inn B$ so that $\fM_{\al_{(E)}}\!\cap\!Z_{(E)}^\cV\!\ne\!\emptyset$, which implies $\fM_{\al_{(E)}}\!\cap\cV\!\ne\!\emptyset$.
	By  Lemma~\ref{Lm:M_strata_local'},
	there then exist $E_1,\ldots,E_\ell\!\subset\!S_\al$ such that 
	$$
	\fM_{\al_{(E)}}\!\cap\!\cV=\bigsqcup_{1\le i\le\ell}\!\!Z^\cV_{(E_i)}.
	$$
	Now take an arbitrary $x\inn\fM_{\al_{(E)}}\!\cap\!Z_{(E)}^\cV$.
	Since $x\inn Z_{(E)}^\cV$,
	we have 
	$$
	\ze_e^\cV(x)\eq 0\quad\forall~e\inn S_\al\bsl E;\qquad
	\ze_e^\cV(x)\!\neq\! 0\quad\forall~e\inn E.
	$$
	However, $x\inn\fM_{\al_{(E)}}$ implies there exists $1\!\le\!i\!\le\!\ell$ such that $x\inn Z_{(E_i)}^\cV$, hence 
	$$
	\ze_e^\cV(x)\eq 0\quad\forall~e\inn S_\al\bsl E_i;\qquad
	\ze_e^\cV(x)\!\neq\! 0\quad\forall~e\inn E_i.
	$$
	Combining the above two displays,
	we see $E\eq E_i$, so $Z_{(E)}^\cV\eq Z_{(E_i)}^\cV\inn\pi_0(\fM_{\al_{(E)}}\!\cap\!\cV).$ 
	
	Next, for every $E'\!\subset\!E$, we aim to show $Z_{(E')}^\cV\!\not\subset\fM^\mn$ so that the preceding paragraph can be applied to $Z_{(E')}^\cV$.
	In fact,
	we have 
	$$ 
	Z_{(E')}^\cV\!\cap\!\De~\supset~
	Z_{(E')}^\cV\!\cap\!\tn{Cl}_\fM\big(Z_{(E)}^\cV\big)~\supset~
	Z_{(E')}^\cV\!\cap\!\tn{Cl}_\cV\big(Z_{(E)}^\cV\big)
	\;=\,Z_{(E')}^\cV\,\ne\,\emptyset,
	$$
	where the first inclusion follows from the closedness of $\De$.
	
	Finally, we show $\al_{(E)}$ is independent of the choice of $\cV\inn\fV_\al$.
	To this end, we show that for every $\cV_1,\cV_2\inn\fV_\al$ satisfying $\cV_1\!\cap\!\cV_2\!\cap\!\fM_\al\!\ne\!\emptyset,$
	we have
	\begin{align}\label{Eqn:al_(E)}
		Z^{\cV_1}_{(E)}\cap\cV_1\cap\cV_2\,=\, Z^{\cV_2}_{(E)}\cap\cV_1\cap\cV_2\,\ne\,\emptyset.
	\end{align}
	This implies $Z^{\cV_1}_{(E)}$ and $Z^{\cV_2}_{(E)}$ are contained in the same $\fM_{\al_{(E)}}$.
	The connectedness of $\fM_{\al}$ then implies $\al_{(E)}$ is independent of $\cV$.
	
	It remains to justify (\ref{Eqn:al_(E)}).
	Indeed, the equality in (\ref{Eqn:al_(E)}) follows from Convention~\ref{Convention:loc_Euc_id},
	while the inequality there holds because $\cV_1\!\cap\!\cV_2$ is an open neighborhood of some $z_0\inn\fM_\al$,
	and such $z_0$ lies in the closure of $Z^{\cV_2}_{(E)}$.
\end{proof}

\begin{crl}
	\label{Crl:M_strata_local}
Let $\fM$, $\al$ and $E$ be as in Proposition~\ref{Prp:M_strata_local}.
Then, there exists a bijection 
\begin{align*}
	\iota_{\al;E}:\,
	S_{\al_{(E)}}\lra 
	S_\al\bsl E\quad\big(\hookrightarrow S_\al\big)
\end{align*}
such that  for every $\cV\inn\fV_\al$ and $\cV'\inn\fV_{\al_{(E)}}$ satisfying $Z^\cV_{(E)}\!\cap\!\cV'\!\ne\!\emptyset$, 
we have $\iota_{\al;E}\eq\iota_{\cV,\cV'}$.
Moreover,
for every pair of disjoint subsets $E,E'\!\subset\!S_{\al}$,
we have
\begin{align*}
	\big(\al_{(E)}\big)_{\big(\iota_{\al;E}^{-1}(E')\big)}=\al_{(E\sqcup E')}\qquad\tn{and}\qquad
	\iota_{\al;E}\circ\iota_{\al_{(E)};\,\iota^{-1}_{\al;E}(E')}=
	\iota_{\al;\,E\sqcup E'}.
\end{align*}
\end{crl}

For every $e\inn S_\al$,
we write $\iota_{\al;\{e\}}$ as $\iota_{\al;e}$ for conciseness.

\begin{proof}[Proof of Corollary~\ref{Crl:M_strata_local}]	
The former statement follows from (\ref{Eqn:loc_Euc_transition}) and Convention~\ref{Convention:loc_Euc_id}.

To see the latter statement,
take arbitrary $\cV\inn\fV_\al$ and $\cV'\inn\fV_{\al_{(E)}}$ satisfying $Z^\cV_{(E)}\!\cap\!\cV'\!\ne\!\emptyset$.
Shrinking $\cV'$ if necessary,
we assume that $\cV'\!\cap\!\cV\!\subset\!\big\{\ze_{e}^\cV\!\neq\! 0:
e\inn E\big\}$,
which is consistent with (\ref{Eqn:loc_Euc_equation}).
Then, by (\ref{Eqn:loc_Euc_transition}),
we have
\begin{align*}
	Z^{\cV'}_{\big(\iota_{\al;E}^{-1}(E')\big)}\cap\cV
	=
	Z^{\cV}_{(E'\sqcup E)}\cap \cV',
\end{align*}
which is not empty because neither is $Z^\cV_{(E)}\!\cap\!\cV'$.
Since the LHS and RHS of the above equality are respectively contained in 
$\fM_{\big(\al_{(E)}\big)_{\big(\iota_{\al;E}^{-1}(E')\big)}}$ and $\fM_{\al_{(E\sqcup E')}}$,
we establish the first identity of the latter statement of Corollary~\ref{Crl:M_strata_local}.

The last identity of the latter statement of Corollary~\ref{Crl:M_strata_local} follows from (\ref{Eqn:loc_Euc_transition}).
\end{proof}

\begin{eg}
	We continue with the setting of Example~\ref{Eg:G-adim_fixture}.
	Consider two boundary strata $\fM_\al$ and $\fM_{\al'}$ of $\fM^\wt_g$ as follows.
	\begin{itemize}[leftmargin=*]
		\item 
		$\fM_\al$ consists of $(C,\bfw)$ such that $C\eq C_g\!\cup\!C_a\!\cup\!C_b$, $C_g$ is smooth and of genus $g$, $C_a$ and $C_b$ are both smooth and rational (each meeting $C_g$ at a node),
		$\bfw(C_g)\eq k$ for some $k\inn\Z_{\ge 0}$, and $\bfw(C_a)\eq\bfw(C_b)\eq\ell$ for some $\ell\inn\Z_{>0}$;
		\item 
		$\fM_{\al'}$ consists of $(C,\bfw)$ such that $C\eq C_g\!\cup\!C_c$, $C_g$ is smooth and of genus $g$, $C_c$ is smooth and rational (meeting $C_g$ at a node),
		$\bfw(C_0)\eq k\!+\!\ell$, and $\bfw(C_c)\eq \ell$.
	\end{itemize}
	We set $S_\al\eq\{a,b\}$ so that the modular parameter $\ze_a$ (resp.~$\ze_b$) corresponds to the smoothing of the node between $C_g$ and $C_a$ (resp.~$C_b$),
	and set $S_{\al'}\eq\{c\}$ likewise.
	
	It is a direct check that $\al_{(a)}\eq\al_{(b)}\eq\al'$,
	and $\iota_{\al;a}$ (resp.~$\iota_{\al;b}$) maps $c$ to $b$ (resp.~to $a$).
\end{eg}

\subsection{Rooted trees}
\label{Subsec:graphs_and_levels}

In the STF theory, rooted trees are assigned to the strata of the given LES in a coherent way,
so we introduce necessary terminology and notation about rooted trees.

In the usual language of the graph theory, a \ts{tree} refers to a connected graph whose first Betti number is 0, and a \ts{rooted tree} is a tree along with a selected vertex known as the \ts{root}.
On the set $E$ of the edges of a rooted tree $\tau$,
consider the partial order $\preceq$ given by $e\!\preceq\!e'$  if any path containing the root and $e$ must contain $e'$ as well.
We call $\preceq$ the \ts{tree order} of $\tau$.
Notice that any elements $e'$ and $e''$ of $E$ are comparable if and only if there exists $e\inn E$ such that $e\!\preceq\!e'$ and $e\!\preceq\!e''$ hold simultaneously.

Conversely, any poset $(E,\preceq)$ satisfying the aforementioned property, up to isomorphisms of posets, uniquely determines a rooted tree $\tau$ whose set of the edges is $E$ and the corresponding tree order is $\preceq$.
The maximal and minimal elements of $E$ are connected to the root and the \ts{leaves}
of $\tau$, respectively.

The poset description turns out to be more convenient for the STF theory,
hence we utilize the following definition of rooted trees throughout the paper.

\begin{dfn}\label{Dfn:Rooted_tree}
A \ts{rooted tree} is a finite (possibly empty) partially ordered set (or poset) $$\tau=(E,\preceq),$$ satisfying for any $e',e''\inn E$,
\begin{align}\label{Eqn:tree_order}
	\big\lgroup
	e'\ \tn{and}\ e''\ \tn{are~comparable}\big\rgroup\quad
	\Longleftrightarrow\quad
	\big\lgroup \,
	\exists\, e\in E\ \ \tn{s.t.}\ \ e\preceq e'\ \tn{and}\ e\preceq e''\big\rgroup.
\end{align}
We call the elements of $E$ the \ts{edges} of $\tau$, and the order $\preceq$ the \ts{tree order} of $\tau$. 
When the context is clear,
we identify $\tau$ with $E$ and simply write $e\inn \tau$ to refer to an edge of $\tau$.

The set $\max(\tau)$ of the maximal elements of $\tau$ is called the \ts{root} of $\tau$. 
Each minimal element of $\tau$ is called a \ts{leaf} of $\tau$.
The set of the leaves of $\tau$ is denoted by $\min(\tau)$.

A \ts{root-to-leaf} path is a linearly ordered subset of $\tau$ that is maximal with respect to inclusion.

The edge-less rooted tree $(\emptyset,-)$ is said to be \ts{trivial} and denoted by $\tau_\bullet$.

The set of all the rooted trees is denoted by $\bT$.

Two rooted trees are said to be \ts{isomorphic} if they are isomorphic as posets.
\end{dfn}

\begin{rmk}\label{Rmk:rooted_tree}
	Throughout this paper,
	although we consider the rooted trees as posets instead of graphs,
	it is convenient to illustrate them in the usual sense of the graph theory;
	c.f.~Figure~\ref{Fig:rooted_tree}.
	The root as per Definition~\ref{Dfn:Rooted_tree} is exactly the set of the incident edges of the  root vertex,
	which is labeled as $o$ (if needed)
	and always placed no lower than the other vertices.
	For every edge,
	we call its endpoint that is closer to the root the \ts{upper endpoint} and the other the \ts{lower endpoint}.
	
	We emphasize the trivial rooted tree $\tau_\bullet$, when visualized as a graph, is the edge-less graph that has exactly one vertex, instead of the empty graph.
\end{rmk}

\begin{figure}[htp]
	\begin{center}
		\begin{tikzpicture}
			\filldraw
			(0,0) circle (1.2pt)
			(0,.5) circle (1.2pt)
			(-.3,0) circle (1.2pt)
			(-.6,0) circle (1.2pt)
			(.3,0) circle (1.2pt)
			(.9,0) circle (1.2pt)
			(.3,1) circle (1.2pt)
			(.6,.5) circle (1.2pt);
			\draw
			(0,0)--(0,.5)
			(-.3,0)--(.3,1)--(.9,0)
			(-.6,0)--(0,.5)--(.3,0);
			\draw
			(0,.75) node {\tiny{$a$}}
			(-.31,.16) node {\tiny{$c$}}
			(-.5,.25) node {\tiny{$b$}}
			(-.09,.16) node {\tiny{$d$}}
			(.3,.3) node {\tiny{$f$}}
			(.92,.25) node {\tiny{$z$}}
			(.62,.75) node {\tiny{$g$}}
			(.3,1.15) node {\tiny{$o$}};
			\draw
			(2,1) node[right] {\tiny{$\tau=\{a,b,c,d,f,g,z\}$}}
			(5,1) node[right] {\tiny{satisfying $b,c,d,f\prec a;\ \ z\prec g.$}}
			(2,.5) node[right] {\tiny{$\Xi(\tau)=\{\fE,\fF,\mathfrak G,\mathfrak H\},$}}
			(5,.5) node[right] {\tiny{where $\fE:=\{a,g\},\ \fF:=\{a,z\},\ \mathfrak G:=\{b,c,d,f,g\},\ \mathfrak H:=\{b,c,d,f,z\}$,}}
			(5,0) node[right] {\tiny{satisfying $\mathfrak H\prec \fF\prec\fE;$\ \ $\mathfrak H\prec \mathfrak G\prec\fE.$}}
			;
		\end{tikzpicture}		
	\end{center}
	\caption{A rooted tree $\tau$}\label{Fig:rooted_tree}		
\end{figure}

For every $E'\!\subset\!\tau$,
let
\begin{align}\label{Eqn:E_R}
	(E')^{\tn R}=\bigcup_{e'\in E'}\{\,e\inn\tau\,:\,e\,\tn R\, e'\,\}\qquad
	\tn{with}\quad \tn R=\,\prec,\,\preceq,\,\succ,\,\tn{and}\,\succeq.
\end{align}
For instance,
$(E')^\prec$ consists of $e\inn\tau$ satisfying $e\!\prec\!e'$ for {\it some} $e'\inn E'$ (not necessarily all $e'\inn E'$).
In Figure~\ref{Fig:rooted_tree},
for example,
$\{a,g\}^\prec\eq\{b,c,d,f,z\}$ and $\{b,z\}^\succeq\eq\{a,b,g,z\}$.

For every $\tau\eq(E,\preceq)\inn\bT$ and $E'\!\subset\!E$,
it is direct to verify that
the complement $E\bsl E'$, along with the restriction of $\preceq$, satisfies (\ref{Eqn:tree_order}) and thus is a rooted tree.
\begin{dfn}\label{Dfn:Edge_contraction_tree}
	The rooted tree above is denoted by $$\tau\bsl E':=(E\bsl E',\preceq).$$
	Such an operation is called \ts{contracting (the edges of) $E'$ from $\tau$}, or simply an \ts{edge contraction}.
\end{dfn}

Intuitively,
contracting an edge means deleting the edge from the tree and identify the two endpoints of the deleted edge.

Next, we recall the notion of transverse sections introduced in~\cite{HLN}.
\begin{dfn}
\label{Dfn:transverse_sections}
Let $\tau\inn\bT$.
A nonempty subset $\fE\!\subset\!\tau$ is called a \ts{transverse section} of $\tau$ if
it meets any of the following three equivalent conditions: 
\begin{itemize}[leftmargin=*]
	\item $\fE$ meets every root-to-leaf path of $\tau$ at exactly one edge;
	\item $\fE$ is a maximal subset of incomparable edges of $\tau$;
	\item every pair of distinct edges of $\fE$ are incomparable, and every edge of $\tau$ is comparable with an edge of  $\fE$.
\end{itemize}
The set of all the transverse sections of $\tau$ is denoted by $\Xi(\tau)$.
\end{dfn}


The tree order on $\tau$ induces a partial order, still denoted by $\preceq$,
on $\Xi(\tau)$ such that for every $\fE,\fE'\inn\Xi(\tau)$,
\begin{align}\label{Eqn:transverse_sections_order}
	\lgroup\,\fE\prec \fE'\,\rgroup\quad\Longleftrightarrow\quad
	\lgroup\, \fE\ne \fE'\,\rgroup\ \ \tn{and}\ \ 
	\lgroup\,\fE\subset(\fE')^\preceq\,\rgroup\,.
\end{align}

For example,
in Figure~\ref{Fig:rooted_tree},
we have
\begin{align*}
	&\Xi(\tau)=\big\{\fE,\fF,\mathfrak G,
	\mathfrak H\big\},\qquad
	\tn{where}\\
	&\fE:=\{a,g\},\ \ 
	\fF:=\{a,z\},\ \ 
	\mathfrak G:=\{b,c,d,f,g\},\ \ 
	\mathfrak H:=\{b,c,d,f,z\}.
\end{align*}
The partial order (\ref{Eqn:transverse_sections_order}) is given by $\mathfrak H\!\prec\! \fF\!\prec\!\fE$ and $\mathfrak H\!\prec\! \mathfrak G\!\prec\!\fE$ in this case.

Every $E\!\subset\!\tau$ is a poset with the induced partial order.
As in Definition~\ref{Dfn:Rooted_tree},
we write the set of the maximal (resp.~minimal) elements of $E$ as $\max(E)$ (resp.~$\min(E)$).

\begin{lmm}\label{Lm:Transverse_sections}
For every nontrivial $\tau\inn\bT$ (i.e.~$\tau\!\ne\!\tau_\bullet$), the following holds.
\begin{enumerate}[leftmargin=*,label=(\arabic*)]
\item \label{Cond:TS_max_min}
The subset $\max(\tau)$ (resp.~$\min(\tau)$) is the greatest (resp.~least) element of the poset $\Xi(\tau)$.

\item \label{Cond:TS_max} For every $E\!\subset\!\tau$ and $\fE\inn\Xi(\tau)$,
we have $\max (E\!\cup\!\fE)\in\Xi(\tau)$.

\item \label{Cond:TS_min} For every nonempty $S\!\subset\!\Xi(\tau)$,
we have $\min\big(\bigcup_{\fE\in S}\fE\big)\in\Xi(\tau)$.
\end{enumerate}
\end{lmm}

\begin{proof}
Notice that a rooted tree is required to be a finite poset.
The statement of Part~\ref{Cond:TS_max_min} is thus straightforward.

To show Part~\ref{Cond:TS_max},
observe that $E\!\cup\!\fE$ is nonempty because $\fE$, a transverse section, is so.
In addition, distinct maximal elements of $E\!\cup\!\fE$ are incomparable.
Finally,
for every $e'\inn\tau$,
there exists $e''\inn\fE$  comparable with $e'$.
In addition, there exists $e\inn\max(E\!\cup\!\fE)$ satisfying $e''\!\preceq\!e$ because $E\!\cup\!\fE$ is finite.
Therefore,
if $e'\!\preceq\!e''$, then $e'\!\preceq\!e$;
if $e'\!\succ\!e''$,
then by (\ref{Eqn:tree_order}),
$e'$ and $e$ are still comparable.
In sum, $\max(E\!\cup\!\fE)$ satisfies the last  criterion (hence all the criteria) of Definition~\ref{Dfn:transverse_sections}.

In Part~\ref{Cond:TS_min},
we write $\fF\!:=\!\min\big(\bigcup_{\fE\in S}\fE\big)$ for conciseness.
Mimicking the proof of Part~\ref{Cond:TS_max}, we see $\fF\!\ne\!\emptyset$, and distinct elements of $\fF$ are incomparable.
Moreover,
for every $e'\inn\tau$ and $\fE\inn S$,
there exists $e_\fE\inn\fE$ comparable with $e'$.
If $e'\!\succeq\!e_\fE$ for some $\fE\inn S$,
then there exists $e\inn\fF$ such that $e\!\preceq\!e_\fE$ (because $\bigcup_{\fE\in S}\fE$ is finite), hence $e\!\preceq\!e'$.
If $e'\!\prec\!e_\fE$ for all $\fE\inn S$,
then by (\ref{Eqn:tree_order}),
the edges $e_\fE$, $\fE\inn S$, are pairwise comparable,
so there exists $\fE^\flat\inn S$ such that $e_{\fE^\flat}\!\preceq\!e_\fE$ for all $\fE\inn S$,
and $e'\!\prec\!e_{\fE^\flat}$. 
It remains to show $e_{\fE^\flat}\inn\fF$.
Indeed,
suppose there exist $\fE''\inn S$ and $e'' \inn\fE''$ such that $e''\!\prec\!e_{\fE^\flat}$.
Then, $e''\!\prec\!e_{\fE''}$,
so $\fE''$ contains two elements that are comparable,
which contradicts the fact  $\fE''\inn\Xi(\tau)$.
\end{proof}

\begin{crl}
	\label{Crl:Transverse_sections}
For every $e\inn\tau$, there exists $\fE\inn\Xi(\tau)$ containing $e$.
\end{crl}

\begin{proof}
By Lemma~\ref{Lm:Transverse_sections},
we have
$e\in\max\big(\{e\}\!\cup\!\min(\tau)\big)\in\Xi(\tau)$.
\end{proof}

We conclude this subsection with the notion of adjoint sets that will play a key role in the description of the treelike structures in \S\ref{Subsec:treelike} (c.f.~Definition~\ref{Dfn:Treelike_structure}).

\begin{dfn}\label{Dfn:E^}
For every $E\!\subset\!\tau$,
the \ts{adjoint set} $E^\wedge$ of $E$ is the complement in $\tau$ of the union of the transverse sections that do not meet $E$.
In other words,
\begin{align*}
E^\wedge
:=\,
\tau\Big\bsl\Big(\!\bigcup_{\fE\in\Xi^{\tn c}(\tau;E)}\hspace{-.15in}\fE\,\Big)\,,
\qquad
\tn{where}\quad
\Xi^{\tn c}(\tau;E)\,:=\,
\big\{\,
\fE\in\Xi(\tau):\,
\fE\cap E=\emptyset\,
\big\}\,.
\end{align*}
For every $e\inn\tau$, $\{e\}^\wedge$ is simply written as $e^\wedge$.
\end{dfn}


On some occasions,
contracting the adjoint set of certain subset of $\tau$ is involved.
For every $E,F\!\subset\!\tau$,
we emphasize the notation
$E\bsl\, F^\wedge$ means the set difference of $E$ and $F^\wedge$.
Similarly, 
$\tau\bsl\,E^\wedge$
refers to the contraction of $E^\wedge$ from $\tau$.

\begin{lmm}\label{Lm:E^_simple}
For every $\tau\inn\bT$ and $F\!\subset\!E\!\subset\!\tau$, we have 
\begin{align*}
	E\subset E^\wedge,\qquad 
	F^\wedge\subset E^\wedge,
	\qquad\tn{and}\qquad
	(E^\wedge)^\wedge=E^\wedge.
\end{align*}
\end{lmm}
\begin{proof}
The first statement follows directly from Definition~\ref{Dfn:E^}.
The middle statement is true because $\Xi^{\tn c}(\tau;E)\!\subset\!\Xi^{\tn c}(\tau;F)$,
which also follows from Definition~\ref{Dfn:E^}.

The first two statements imply $E^\wedge\!\subset\!(E^\wedge)^\wedge$.
Moreover, for every $\fE\inn\Xi(\tau)\bsl\Xi^{\tn c}(\tau;E^\wedge)$, we have $\fE\!\cap\!E^\wedge\!\ne\!\emptyset$, so
by Definition~\ref{Dfn:E^}, we have $\fE\!\not\in\!\Xi^{\tn c}(\tau;E)$.
Therefore, 
$\Xi^{\tn c}(\tau;E)\!\subset\! \Xi^{\tn c}(\tau;E^\wedge)$, which is equivalent to
$(E^\wedge)^\wedge\!\subset\!E^\wedge$.
In this way, the last statement has been established.
\end{proof}

Given $e\inn\tau$, if $e\!\not\in\!\min(\tau)$, then $e^\wedge\eq\{e\}$. 
Indeed, for each $e'\inn\tau\bsl\{e\}$, 
by Lemma~\ref{Lm:Transverse_sections}, we have
\begin{align*}
	e'\in \max\big(\{e'\}\!\cup\!\min(\tau)\big) 
	\in\Xi(\tau),\qquad
	e\not\in \max\big(\{e'\}\!\cup\!\min(\tau)\big).
\end{align*}
Hence $e'\!\not\in\!e^\wedge$.
If $e\inn\min(\tau)$,
however, it is possible that $e^\wedge$ is strictly larger than $\{e\}$, which follows from the statement below.

\begin{lmm}
	\label{Lm:path^}
For every $\tau\inn\bT$, $e\inn\tau$, and $e^\flat\inn \{e\}^\preceq\!\cap\!\min(\tau)$,
let $[e^\flat,e]\!:=\!\{e'\inn\tau: e^\flat\!\preceq\!e'\!\preceq\!e\}.$
Then, one of the following holds.
\begin{itemize}[leftmargin=*]
\item 
If $e\!\not\in\!\max(\tau)$,
then with $e^\sharp$ denoting the least element $\{e\}^\succ$,
we have $[e^\flat,e]^\wedge\eq\{e^\sharp\}^\prec$.

\item 
If $e\!\in\!\max(\tau)$,
then $[e^\flat,e]^\wedge\eq\tau$.
\end{itemize}
\end{lmm}

\begin{proof}
First, assume $e\!\not\in\max(\tau)$.

For each $e'\!\not\in\!\{e^\sharp\}^\prec$,
Corollary~\ref{Crl:Transverse_sections} guarantees there exists $\fE'\inn\Xi(\tau)$ containing $e'$.
From Lemma~\ref{Lm:Transverse_sections}~\ref{Cond:TS_max}, we then conclude that $\max\big(\{e^\sharp\}\!\cup\!\fE'\big)\inn\Xi(\tau)$,
which contains $e'$ while being disjoint from  $[e^\flat,e]$,
hence $e'\!\not\in\![e^\flat,e]^\wedge$.

For each $e'\!\in\!\{e^\sharp\}^\prec$,
suppose $e'\!\not\in\![e^\flat,e]^\wedge$.
Then by Definition~\ref{Dfn:E^},
there exists $\fE'\inn\Xi(\tau)$ that contains $e'$ while being disjoint from $[e^\flat,e]$.
As a transverse section, $\fE'$ must contain an element $e''$ comparable with $e$.
Since $\fE'\!\cap\![e^\flat,e]\eq\emptyset$ and $e^\flat$ is minimal, we have $e''\!\succ\!e$.
Consequently, $e''\!\succeq\!e^\sharp\!\succ\!e'$, yet $e''$ and $e'$ both lie in the transverse section $\fE'$, which contradicts Definition~\ref{Dfn:transverse_sections}.
Therefore, $e\inn [e^\flat,e]^\wedge$.

To summarize, the statement for $e\!\not\in\max(\tau)$ has been proved.
The statement for $e\!\in\max(\tau)$ follows from an argument similar to but simpler than the last paragraph of the $e\!\not\in\!\max(\tau)$ case.
\end{proof}

The following results will play important roles in~\S\ref{Subsec:treelike} and~\S\ref{Subsec:levels}.

\begin{lmm}\label{Lm:E^}
For every $\tau\inn\bT$ and $F\!\subset\! E\!\subset\!\tau$, we have the following.

\begin{enumerate}[leftmargin=*,label=(\arabic*)]

\item \label{Cond:adj_TS_edge_contr}
$\Xi(\tau\bsl\, F^\wedge)=\Xi^{\tn c} (\tau; F)$.

\item \label{Cond:adj_edge_contrs}
The adjoint set of $E\bsl\,F^{\wedge}$ in $\tau\bsl\,F^{\wedge}$ is $E^\wedge\bsl\,F^{\wedge}$.		
\end{enumerate}
\end{lmm}

\begin{proof}
We begin with~\ref{Cond:adj_TS_edge_contr}.
First, notice that every $\fE\inn\Xi^{\tn c}(\tau;F)$ is a subset of $\tau\bsl F^\wedge$.
Moreover, it is a direct check that $\fE\inn\Xi(\tau\bsl F^\wedge)$, using the fact that $\fE\inn\Xi(\tau)$.
Therefore, $\Xi^{\tn c}(\tau;F)\!\subset\!\Xi(\tau\bsl F^\wedge)$.

To show $\Xi(\tau\bsl F^\wedge)\!\subset\!\Xi^{\tn c}(\tau;F)$,
consider the subsets
\begin{align*}
	\ti\fF:=\min(\tau\bsl\,F^\wedge)\qquad\tn{and}\qquad
	\fF:=\max\big(\,\ti\fF\cup\min(\tau)\,\big)
\end{align*}
of $\tau$.
By Parts~\ref{Cond:TS_max_min} and~\ref{Cond:TS_max} of Lemma~\ref{Lm:Transverse_sections},
we have $\ti\fF\inn\Xi(\tau\bsl F^\wedge)$ and $\fF\inn\Xi(\tau)$.
Moreover,
given $e\inn\ti\fF$,
any other element of $\ti\fF$ cannot be in $\{e\}^\succ$ because $\ti\fF$ is a transverse section, while any element of $\min(\tau)$, by definition, cannot be in $\{e\}^\succ$ either.
Hence $e\inn\fF$,
i.e.~$\ti\fF\!\subset\!\fF$.

\begin{itemize}[leftmargin=*]
\item 
If $\ti\fF\!\subsetneq\!\fF$,
there then exists $e\inn\min(\tau)$ such that $\ti\fF\!\cap\!\{e\}^\succeq\eq\emptyset$.
This implies 
\begin{align*}
	\{e\}^\succeq \cap (\tau\bsl\, F^\wedge)=\emptyset\,,\qquad\tn{i.e.}\qquad
	\{e\}^\succeq\subset F^\wedge.
\end{align*}
By Lemmas~\ref{Lm:path^} and~\ref{Lm:E^_simple},
we have $F^\wedge\eq\tau$,
hence $\Xi(\tau\bsl\, F^\wedge)\eq \emptyset\!\subset\! \Xi^{\tn c} (\tau; F)$.

\item 
If $\ti\fF\!=\!\fF$,
then $\ti\fF\inn\Xi(\tau)$,
hence for every $\fE\inn \Xi(\tau\bsl F^\wedge)$, by Lemma~\ref{Lm:Transverse_sections}~\ref{Cond:TS_max},
we have
\begin{align*}
	\fE=\max(\fE\cup\ti\fF)\in\Xi(\tau).
\end{align*}
In addition, $\fE\!\cap\!F\eq\emptyset$ because $\fE\!\subset\!\tau\bsl F^\wedge$.
Consequently, $\fE\inn\Xi^{\tn c}(\tau;F)$,
i.e.~$\Xi(\tau\bsl\, F^\wedge)\!\subset\! \Xi^{\tn c} (\tau; F)$.
\end{itemize}
To summarize, Part~\ref{Cond:adj_TS_edge_contr} has been established.

By Definition~\ref{Dfn:E^}, Part \ref{Cond:adj_TS_edge_contr}, and the assumption $F\!\subset\!E$, we have
\begin{align*}
	\Xi^{\tn c}(\tau\bsl F^\wedge;E\bsl F^\wedge)
	&=\big\{\,\fE\inn\Xi^{\tn c}(\tau; F):\,
	\fE\!\cap\! E\eq\emptyset\,\big\}\\
	&=\big\{\,\fE\inn\Xi(\tau):\,\fE\!\cap\! F\eq\emptyset,\ 
	\fE\!\cap\! E\eq\emptyset\,\big\}
	=\Xi^{\tn c}(\tau; E)\,.
\end{align*}
Therefore,
the adjoint set of $E\bsl F^\wedge$ in $\tau\bsl F^\wedge$ can be written as
\begin{align*}
	\big(\bigcup_{\fE\in \Xi^{\tn c}(\tau;F)}\hspace{-.15in}\fE\big)
	\,\Big\bsl\,
	\big(\bigcup_{\fE\in\Xi^{\tn c}(\tau; E)}\hspace{-.15in}\fE\big),
\end{align*}
which, by Definition~\ref{Dfn:E^}, is the same as $E^\wedge\bsl F^\wedge$.
This completes the proof of Part~\ref{Cond:adj_edge_contrs}.
\end{proof}

In Figure~\ref{Fig:E^}, we provide examples in illustration of Lemma~\ref{Lm:E^}.
Here, for $e\eq c,g,x,f$, the outer adjoint in $\{E\bsl\,e^\wedge\}^\wedge$  is taken in $\tau\bsl\,e^\wedge$ instead of  $\tau$.

\begin{figure}[htp]
	\begin{center}
		\begin{tikzpicture}
			\filldraw
			(0,.5) circle (1.2pt)
			(.3,0) circle (1.2pt)
			(-.6,-.5) circle (1.2pt)
			(-.3,0) circle (1.2pt)
			(.9,0) circle (1.2pt)
			(1.2,-.5) circle (1.2pt)
			(.3,1) circle (1.2pt)
			(.6,.5) circle (1.2pt)
			(.6,-.5) circle (1.2pt)
			(0,-.5) circle (1.2pt)
			;
			\draw
			(.6,-.5)--(0,.5)
			(0,.5)--(.3,1)--(1.2,-.5)
			(0,-.5)--(.3,0)
			(-.6,-.5)--(0,.5);
			\draw
			(.3,1) node[above] {\tiny{$o$}}
			(0,.75) node {\tiny{$a$}}
			(-.35,.25) node {\tiny{$c$}}
			(-.65,-.25) node {\tiny{$g$}}
			(.6,.75) node {\tiny{$b$}}
			(.95,.25) node {\tiny{$f$}}
			(1.22,-.25) node {\tiny{$z$}}
			(-.05,-.25) node {\tiny{$x$}}
			(.65,-.25) node {\tiny{$y$}}
			(.3,.25) node {\tiny{$d$}}
			(1.8,-.4) node {\small{$\tau$}};
			
			\draw[xshift=3cm, yshift=-.6cm]
			(0,1.8) node[right] {\tiny{$E\eq\{c,g,x,f\}$}}
			(0,1.4) node[right] {\tiny{${c}^\wedge\eq\{c\}$}}
			(0,1) node[right] {\tiny{${g}^\wedge\eq\{g\}$}}
			(0,.6) node[right] {\tiny{$x^\wedge\eq\{x,y\}$}}
			(0,.2) node[right] {\tiny{$f^\wedge\eq\{f\}$}}
			;
			
			\draw[xshift=5.5cm, yshift=-.6cm]
			(0,1.8) node[right] {\tiny{$E^\wedge\eq\{c,g,d,x,y,f\}$}}
			(0,1.4) node[right] {\tiny{$(E\bsl\,c^\wedge)^{\wedge}\eq\{g,d,x,y,f\}$}}
			(0,1) node[right] {\tiny{$(E\bsl\,g^\wedge)^{\wedge}\eq\{c,d,x,y,f\}$}} 
			(0,.6) node[right] {\tiny{$(E\bsl\,x^\wedge)^{\wedge}\eq\{c,g,d,f\}$}}
			(0,.2) node[right] {\tiny{$(E\bsl\,f^\wedge)^{\wedge}\eq\{c,g,d,x,y\}$}}
			;
			
			\draw[xshift=9cm, yshift=-.6cm]
			(0,1.8) node[right] {\tiny{$\Xi(\tau\bsl\, E^\wedge)\eq\{\,\{a,b\},\,\{a,z\}\,\}\eq\Xi^{\tn c}(\tau;E)$}}
			(0,1.4) node[right] {\tiny{$\Xi(\tau\bsl\,x^\wedge)\eq\{\,\{a,b\},\,\{a,f\},\,\{a,z\},$}}
			(1.545,1) node[right] {\tiny{$\{c,d,b\},\,\{c,d,f\},\,\{c,d,z\},$}} 
			(1.545,.6) node[right] {\tiny{$\{g,d,b\},\,\{g,d,f\},\,\{g,d,z\}\,\}$}}
			(1.12,.2) node[right] {\tiny{$=\!\Xi^{\tn c}(\tau;\{x\})$}}
			;
		\end{tikzpicture}
	\end{center}
	\caption{Examples of the adjoint sets}\label{Fig:E^}
\end{figure} 

\subsection{Treelike structures}
	\label{Subsec:treelike}
In this subsection,
we introduce the notion of treelike structures on {\set} stratified stacks. 
They provide necessary platforms on which the twisted fields (to be formally introduced in \S\ref{Subsec:STF_main_statement}) will be added.	

Consider the LES of $\fM_g^{\rm wt}$ in Example~\ref{Eg:G-adim_fixture}.
Fix a dual graph as well as the weights $\bfw$.  By contracting the smallest subgraph of genus $g$ from the graph and pruning some rational tails according to the weights, one obtains a rooted tree. 
The notion of treelike structures below is modeled on similar consideration.

Recall that for $\tau\inn\bT$ and $E\!\subset\!\tau$,
$E^\wedge$ is the adjoint set of $E$ as per Definition~\ref{Dfn:E^}.
Also recall for any $\fM$ with an LES as in Definition~\ref{Dfn:G-adim_fixture},
a nonempty index set $S_\al$ for the modular parameters are assigned to each $\al\inn B\eq A\bsl\{\ms\}$.
Since no modular parameters are associated with the stratum $\fM^{\mn}\eq\fM_{\ms}$,
we set $S_{\ms}\!:=\!\emptyset$.

\begin{dfn}\label{Dfn:Treelike_structure}
	Let $\fM$ be endowed with an LES as in Definition~\ref{Dfn:G-adim_fixture}. 
	We call an indexed family
	\begin{equation}\label{Eqn:Treelike_structure}
		\La=
		\big(\,
		\tau_{\al},\,
		\be_{\al}\,
		\big)_{\al\in A}
	\end{equation}
	a \ts{treelike structure} on $\fM$ if 
	\begin{itemize}[leftmargin=*]
		\item 
		(\emph{assignment}) it assigns to each $\al\inn A$ a rooted tree and an injection between sets:
		$$\tau_{\al}\in\bT\qquad
		\tn{and}\qquad
		\be_{\al}:
		\tau_\al\hookrightarrow S_{\al},$$
		satisfying 
		\item
		(\emph{compatibility}) for every $\al\inn B$ and  $e\inn S_\al$,
		with $\al_{(e)}\inn A$ as in Proposition~\ref{Prp:M_strata_local},
		there exists an isomorphism
		$$
		\phi_{\al;e}:\,\tau_{\al_{(e)}}\lra\tau_{\al}\big\bsl\,\big(\be_\al^{-1}\{e\}\big)^{\!\wedge}
		$$ 
		such that
		\begin{align}\label{Eqn:Treelike}
			\iota_{\al;e}\circ\be_{\al_{(e)}}=\be_\al\circ\phi_{\al;e}\,,
		\end{align}
		where $\iota_{\al;e}$ is the injection in Corollary~\ref{Crl:M_strata_local}.
	\end{itemize}
\end{dfn}

\begin{rmk}\label{Rmk:Treelike}
In Definition~\ref{Dfn:Treelike_structure},
the injection $\be_\al$ assigns to each edge of $\tau_\al$ a distinct modular parameter,
and the condition (\ref{Eqn:Treelike}) guarantees the compatibility of the assignment.
Particularly, we have $\tau_{\ms}\eq\tau_\bullet$, i.e.~the trivial rooted tree is assigned to $\fM^\mn$.

Notice that the adjoint set $(\be_\al^{-1}\{e\})^\wedge$ is empty if $e$ is not in the image of $\be_\al$, and just consists of $\be_\al^{-1}(e)$ itself if $\be_\al^{-1}(e)$ is not a minimal edge in $\tau_\al$.
If $\be_\al^{-1}(e)$ is a minimal edge, intuitively, $\big(\be_\al^{-1}(e)\big)^\wedge$ contains all the edges of $\tau_\al$ {\it below} the  ``higher'' vertex $v'$ of $\be_\al^{-1}(e)$ (recall we always put the root at the top).
Contracting $\big(\be_\al^{-1}(e)\big)^\wedge$ is equivalent to {\it pruning} $v'$ in the sense of \cite[\S3.2]{HL10}.
In Figure~\ref{Fig:Treelike} we provides some examples of the resulting rooted tree.
\end{rmk}

\begin{dfn}
	\label{Dfn:Tauto_monom}
	Let $\fM$ and $\La$ be as in Definition~\ref{Dfn:Treelike_structure}.
	For every $\al\inn A$, $\cV\inn\fV_\al$, and root-to-leaf path $\wp\!\subset\!\tau_\al$, the product 
	$
	\prod_{e\in\wp}\ze_e^\cV
	$
	is called a \ts{modular monomial} on $\cV$.
	The set of all the modular monomials on $\cV$ is called a \ts{tautological set of monomials}.
	
	A set of regular functions is said to have a \ts{divisibly minimal element} if there exists one member $m_0$ that divides any  member $m$ of the set,
	i.e.~$m/m_0$ is a regular function. 
\end{dfn}

Since every root-to-leaf path is a subset of $\tau$, it does not allow duplicated elements, hence every modular monomial is square-free.

\begin{figure}[htp]
	\begin{center}
		\begin{tikzpicture} 
			\filldraw
			(.35,-.66) circle (1pt)
			(-.35,-.66) circle (1pt)
			(0,0) circle (1pt)
			(.35,.66) circle (1pt)
			(.7,0) circle (1pt);
			\draw
			(-.35,-.66)--(.35,.66)--(.7,0)
			(0,0)--(.35,-.66);
			\draw
			(.35,.66) node[above] {\tiny{$o_\al$}}
			(.175,.33) node[left] {\tiny{$a$}}
			(.525,.33) node[right] {\tiny{$d$}}
			(-.175,-.33) node[left] {\tiny{$b$}}
			(.175,-.33) node[right] {\tiny{$c$}}
			(1.05,-.66) node {\tiny{$\tau_{\al}$}}
			;
			
			\draw[xshift=2.8cm]
			(0,.5) node [left] {\tiny{$\{a,b,c,d\}$}}
			(.5,.5) node [right] {\tiny{$S_\al\eq\{e_a,e_b,e_c,e_d,e'\}$}}
			(.25,.5) node [above] {\tiny{$\be_\al$}}
			(.18,.25) node [below] {\tiny{$\be_\al(j)\eq e_j,\ \ j\eq a,b,c,d$}}
			;
			
			\draw[xshift=2.8cm,->,>=stealth']
			(0,.5)--(.5,.5);
			
			\filldraw[xshift=7cm]
			(0,0) circle (1pt)
			(.35,.66) circle (1pt)
			(.7,0) circle (1pt);
			\draw[xshift=7cm]
			(0,0)--(.35,.66)--(.7,0)
			;
			\draw[xshift=7cm]
			(.35,.66) node[above] {\tiny{$o_\al$}}
			(.175,.33) node[left] {\tiny{$a$}}
			(.525,.33) node[right] {\tiny{$d$}}
			(.35,-.5) node {\tiny{$\tau_{\al}\bsl\,(\be_\al^{-1}\{e_b\})^\wedge$}}
			;
			
			\filldraw[xshift=9.5cm]
			(0,0) circle (1pt)
			(.5,0) circle (1pt)
			(.5,.66) circle (1pt)
			(1,0) circle (1pt);
			\draw[xshift=9.5cm]
			(0,0)--(.5,.66)--(1,0)
			(.5,.66)--(.5,0)
			;
			\draw[xshift=9.5cm]
			(.5,.66) node[above] {\tiny{$o_\al$}}
			(.25,.33) node[left] {\tiny{$b$}}
			(.4,.23) node[right] {\tiny{$c$}}
			(.75,.33) node[right] {\tiny{$d$}}
			(.5,-.5) node {\tiny{$\tau_{\al}\bsl\,(\be_\al^{-1}\{e_a,e'\})^\wedge$}}
			;
			
			\filldraw[xshift=11.25cm]
			(1.5,.66) circle (1pt)
			;
			\draw[xshift=11.25cm]
			(1.5,.66) node[above] {\tiny{$o_\al$}}
			(2.88,.25) node[left] {\tiny{$\tau_{\al}\bsl\,(\be_\al^{-1}\{e_a,e_b\})^\wedge$}}
			(2.88,-.2) node[left] {\tiny{$\eq \tau_{\al}\bsl\,(\be_\al^{-1}\{e_d\})^\wedge$}}
			;
		\end{tikzpicture}
	\end{center}
	\caption{Examples of $\tau_{\al}\big\bsl\,\big(\be_\al^{-1}(E)\big)^{\!\wedge}$}\label{Fig:Treelike}
\end{figure}

\begin{prp}
	\label{Prp:Treelike_contraction}
	Let $\fM$ and $\La$ be as in Definition~\ref{Dfn:Treelike_structure}.
	For every $\al\inn B$ and $E\!\subset\!S_\al$,
	let $\al_{(E)}\inn A$ be as in
	Proposition~\ref{Prp:M_strata_local}.
	Then,
	there exists an isomorphism
	$$
	\phi_{\al;E}:\,\tau_{\al_{(E)}} \lra \tau_{\al}\big\bsl\,\big(\be_\al^{-1}(E)\big)^{\!\wedge}
	$$ 
	such that 
	$$
	\iota_{\al;E}\circ \be_{\al_{(E)}}=\be_\al\circ\phi_{\al;E}\,,
	$$
	where $\iota_{\al;E}$ is the injection in Corollary~\ref{Crl:M_strata_local}.
\end{prp}


\begin{proof}
	The statements follow from applying Definition~\ref{Dfn:Treelike_structure} and Lemma~\ref{Lm:E^}~\ref{Cond:adj_edge_contrs} repeatedly.
\end{proof}

\begin{eg}\label{Eg:genus_1}
	We continue with the $g\eq 1$ case of the LES of $\fM_g^{\wt}$ in Example~\ref{Eg:G-adim_fixture}.
	Then, the index set $A$ consists of all the pairs $(\ga,\bfw)$,
	where $\ga$ are the dual graphs  of arithmetic genus 1, and $\bfw$ assign the weights of the irreducible components to the corresponding vertices.
	
	To each $\al\eq(\ga,\bfw)\inn A$,
	we assign a rooted tree $\tau_{\al}$ that is obtained from the underlying graph of $\ga$ by 
	\begin{itemize}
		[leftmargin=*]
		\item contracting all the edges of the smallest connected genus 1 subgraph $\ga_{\core}$ of the underlying graph of $\ga$ into one vertex~$o$ and obtaining a rooted tree $\ti\tau_\al$ whose root is given by $o$, 
		\item then setting the weight of $o$ to be the total weight of $\ga_{\core}$, and
		\item finally, with $v_e^+$ denoting the incident vertex of $e$ in $\ti\tau_\al$ that either is $o$ or lies between $o$ and~$e$, setting 
		\begin{align*}
			\tau_\al=
			\ti\tau_\al\big\bsl
			\big(\{e\inn\ti\tau_\al:\,
			\bfw(v_e^+)\!>\!0\}^\preceq\big).
		\end{align*}
	\end{itemize}
	Particularly,
	$\tau_\al\eq\tau_\bullet$ if $\bfw(o)\!>\!0$.
	Such $\tau_{\al}$ is exactly the {\it weighted dual tree} of $\fM_\al$ in~\cite[\S2.2]{g1modular}.
	
	For every $\al\inn A$,
	we take $S_\al$ to be
	the set of the edges of $\ga$,
	and $\be_{\al}$ to be the natural inclusion.
	
	To see such $\tau_{\al}$ and $\be_{\al}$ give a treelike structure on $\fM^\wt_1$,
	we take an arbitrary $e\inn S_\al$. 
	Then,
	there exists a unique $\al_{(e)}\eq(\ga_{(e)},\bfw')\inn A$ satisfying Proposition~\ref{Prp:M_strata_local}.
	(More concretely, $\ga_{(e)}$ is obtained by smoothing out the node corresponding to $e$, and the only difference between $\bfw'$ and $\bfw$ is if $e$ has distinct incident vertices, then $\bfw'$ assigns the sum of their weights to the new vertex.)
	
	If $e$ is not in the image of $\be_\al$ (i.e.~$e\!\not\in\!\tau_\al$),
	then either $e$ is in $\ga_{\core}$, or there exists a positively weighted vertex between $e$ and $\ga_{\core}$.
	Either way, it is easy to see that $\tau_{\al_{(e)}}$ and $\tau_{\al}$ are isomorphic in accordance with (\ref{Eqn:Treelike}).
	
	If $e\inn \tau_{\al}$ is not minimal,
	then $e^\wedge\eq\{e\}$,
	and both of the endpoints of~$e$ are of weight 0, 
	hence their images in~$\al_{(e)}$ are of weight 0,
	and thus $\tau_{\al_{(e)}}$ and $\tau_{\al}\bsl e$ are isomorphic in accordance with (\ref{Eqn:Treelike}).
	
	If $e\inn \tau_{\al}$ is minimal,
	the above construction of $\tau_{\al}$ implies the endpoint of $e$ that is different from $v_e^+$ is positively weighted and
	so is its image $\ov v$ in $\al_{(e)}$.
	Therefore, $\ov v$ is a minimal vertex of $\tau_{\al_{(e)}}$, and the image of any element of $e^\wedge\bsl\{e\}$ in $\ga_{(e)}$ is thus not in~$\tau_{\al_{(e)}}$.
	Consequently, $\tau_{\al_{(e)}}$ and $\tau_{\al}\bsl\,e^\wedge$ are isomorphic in accordance with (\ref{Eqn:Treelike}).
	
	In sum, we obtain a treelike structure $(\tau_\al,\be_\al)_{\al\in A}$ on $\fM^\wt_1$ as per Definition~\ref{Dfn:Treelike_structure},
	which in turn gives rise to the stack with twisted fields $\fM^{\tf}_1$ by Theorem~\ref{Thm:tf_smooth}.
	Such $\fM^{\tf}_1$ is identical to the one constructed in~\cite[(2.13)]{g1modular}.
\end{eg}

\subsection{Root-to-leaf and leaf-to-root sequences with dominant edges}
\label{Subsec:levels}

As motivated in the introduction, adding twisted fields to the data $(C,D)$ corresponds to local blowups.
Precisely, for a stack $\fM$ endowed with an LES and a treelike structure,
each local blowup center is given by the common vanishing locus of the modular parameters corresponding to a transverse section of the associated rooted tree.
The order of these blowups comes in two natural ways: one from the root  towards its leaves, while the other in the opposite direction, from the leaves to the root. This brings up two notions:  \ts{root-to-leaf sequence (RLS)}
and \ts{leaf-to-root sequence (LRS)}. 
In either way,
the resulting stack with the twisted fields still possesses an LES covered by modular
charts. 
These local charts contain certain modular parameters corresponding to the exceptional divisors (thus labeled by the transverse sections);
such parameters become common factors of the pullbacks of all the modular monomials after local blowups.
The remaining edges of the rooted tree are called \ts{non-dominate edges}. 
Among others, these notions will be used to establish an LES on the  stack with the twisted fields.

\begin{dfn}\label{Dfn:R2L_seq}
	A \ts{root-to-leaf sequence (RLS)}  is a tuple $(\tau,\ov\bE)$ comprised of a rooted tree $\tau\!\in\!\bT$ and a finite sequence of transverse sections of $\tau$ as in Definition~\ref{Dfn:transverse_sections}
	$$\ov \bE:= \fE_1,\cdots, \fE_\cht\,,
	\qquad
	\fE_i\in\Xi(\tau)\quad\forall~1\!\le\!i\!\le\!\cht$$  so that exactly one of the following two conditions holds:
	\begin{itemize}[leftmargin=*]
		\item either $(\tau,\ov\bE)\eq(\tau_\bullet,\emptyset)$,
		\item or the following criteria are satisfied:
		$$\max(\tau)\eq\fE_1\succ \fE_2\succ\cdots\succ \fE_\cht,
		\qquad
		\fE_\cht\cap\min(\tau)\ne\emptyset,\qquad\tn{and}\qquad
		\big(\!\!\bigcup_{1\le i\le \cht}\!\!\!\!\fE_i\big)\sqcup \fE_\cht^\prec=\tau.
		$$ 
	\end{itemize}
	Similarly,
	a \ts{leaf-to-root sequence (LRS)} is a tuple $(\tau,\ud\bE)$ comprised of a rooted tree $\tau\inn\bT$ and a finite sequence of transverse sections of $\tau$
	$$\ud \bE:= \fE_1,\cdots, \fE_\cht\,,\qquad
	\fE_i\in\Xi(\tau)\quad\forall~1\!\le\!i\!\le\!\cht$$  so that exactly one of the following two conditions holds:
	\begin{itemize}[leftmargin=*]
		\item either $(\tau,\ud\bE)\eq(\tau_\bullet,\emptyset)$,
		\item or the following criteria are satisfied:
		$$\min(\tau)\eq\fE_1\prec \fE_2\prec\cdots\prec \fE_\cht,
		\qquad
		\fE_\cht\cap\max(\tau)\ne\emptyset,\qquad\tn{and}\qquad
		\big(\!\!\bigcup_{1\le i\le \cht}\!\!\!\!\fE_i\big)\sqcup \fE_\cht^\succ=\tau.
		$$
	\end{itemize}
\end{dfn}

\begin{rmk}
	Throughout this paper, we always use $\fE_\cht$ to denote the {\it final term} of an RLS or LRS.
\end{rmk}

\begin{cnv}\label{Convention:R2L}
	For 
	every RLS $(\tau,\ov\bE)$, we identify the finite sequence $\ov\bE$  with the set of all the terms of $\ov \bE$.
	Such convention also applies to the LRS's.
	
	In this way,
	we treat $\ov \bE$ and $\ud\bE$ as  subsets of $\Xi(\tau)$ (although not every subset of $\Xi(\tau)$ corresponds to an RLS or LRS). 
\end{cnv}

Convention~\ref{Convention:R2L} is introduced because each term of $\ov\bE$ (resp.~$\ud\bE$) corresponds to a special local parameter in Theorem~\ref{Thm:tf_smooth} (resp.~Theorem~\ref{Thm:tf_smooth_revert}). 
It is thus convenient to consider $\ov\bE$ and $\ud\bE$ themselves as index sets,
rather than introducing additional indices for the transverse sections and the corresponding local parameters.
In this way,
we may simply write $\fE\inn\ov\bE$ (resp.~$\fE\inn\ud\bE$) to refer to a term of $\ov\bE$ (resp.~$\ud\bE$).
The final term $\fE_\cht$ satisfies
\begin{align*}
	\fE_\cht=\min(\ov\bE)\quad\forall~\tn{RLS}~(\tau,\ov\bE)\qquad
	\tn{and}\qquad
	\fE_\cht=\max(\ud\bE)\quad\forall~\tn{LRS}~(\tau,\ud\bE).
\end{align*}
We further write
\begin{align}\label{Eqn:fE_+-}
	\fE_+:=
	\begin{cases}
		\fE_{i+1} &\tn{if}~\fE\eq \fE_{i},~1\!\le\!i\!<\!\cht,\\
		\emptyset &\tn{if}~\fE\eq \fE_{\cht};
	\end{cases}
	\qquad
	\fE_-:=
	\begin{cases}
		\fE_{i-1} &\tn{if}~\fE\eq \fE_{i},~1\!<\!i\!\le\!\cht,\\
		\emptyset &\tn{if}~\fE\eq \fE_1;
	\end{cases}
\end{align}
for all $\fE\inn\ov\bE$ and $\fE\inn\ud\bE$.

It is possible that an RLS and an LRS determine the same subset of $\Xi(\tau)$.
For instance,
in Figure~\ref{Fig:rooted_tree},
$\{\fE,\fF\}$ can be viewed as an RLS but not an LRS, 
whereas $\{\fE,\fF,\mathfrak H\}$ can be viewed either as an RLS or an LRS.
Nonetheless, this does not cause any ambiguity as long as the direction (RL or LR) of sequence is specified;
c.f.~Theorems~\ref{Thm:tf_smooth} and~\ref{Thm:tf_smooth_revert} for example.

\begin{dfn}\label{Dfn:Enhanced_rtl_seq}
	A \ts{root-to-leaf sequence (RLS) with dominant edges}
	$$
	\ov\lt=(\,\tau\,,\,\ov \bE\,,\,
	\Dm(\ov\lt)\,)
	$$
	is a tuple consisting of an RLS $\big(\tau,\,\ov \bE\eq \fE_1,\ldots, \fE_\cht\big)$  as in Definition~\ref{Dfn:R2L_seq}
	and a subset $\Dm(\ov\lt)\!\subset\!\tau$	such that exactly one of the following two conditions holds:
	\begin{enumerate}[leftmargin=*,label=$\ov{\alph*}.$]
	\item either $$(\tau,\ov\bE,\Dm(\ov\lt))=(\tau_\bullet,\emptyset,\emptyset)=:\ov\lt_\bullet\,,$$
	\item\label{Cond:RLS_dom} or $\tau\!\ne\!\tau_\bullet$, and $\Dm(\ov\lt)$ is a union of nonempty subsets of $\tau$:
	$$
	\Dm(\ov\lt)=\dot\fE_1\cup\cdots\cup\dot\fE_\cht\ \big(\subset\tau\big),\qquad\tn{where}\quad
	\emptyset\ne\dot \fE_i\subset\fE_i\ \  \forall\ 1\!\le\!i\!\le\!\cht,
	$$
	satisfying
	\begin{enumerate}[leftmargin=*,label=({\ref{Cond:RLS_dom}}\arabic*)]
		\item 
		$\dot\fE_\cht\!\cap\!\min(\tau)\!\ne\!\emptyset$, and
		\item for every $1\!\le\!i\!<\!\cht$,
		$\dot \fE_i\eq\fE_i\bsl \fE_{i+1}$.
	\end{enumerate}
	\end{enumerate}
	The edges in $\dot \fE_i$, $1\!\le\!i\!\le\!\cht$, are said to be \ts{dominant} in $\fE_i$.
	The complement of $\Dm(\ov\lt)$ in $\tau$ is denoted by $$\ND(\ov\lt)\;,$$
	whose elements are called \ts{non-dominant} edges of $\ov\lt$.
	The underlying RLS of $\ov\lt$ is denoted by $(\tau_{\ov\lt},\ov\bE_{\ov\lt})$.
	
	Similarly, a \ts{leaf-to-root sequence (LRS) with dominant edges}
	$$
	\ud\lt=(\,\tau\,,\,\ud \bE\,,\,
	\Dm(\ud\lt)\,)
	$$
	is a tuple consisting of an LRS $\big(\tau,\,\ud
	\bE\eq \fE_1,\ldots, \fE_\cht\big)$  as in Definition~\ref{Dfn:R2L_seq}
	and 
	a subset $\Dm(\ud\lt)\!\subset\!\tau$	such that exactly one of the following two conditions holds:
	\begin{enumerate}[leftmargin=*,label=\ud{\alph*}.]
		\item either $$(\tau,\ud\bE,\Dm(\ud\lt))=(\tau_\bullet,\emptyset,\emptyset)=:\ud\lt_\bullet\,,$$
		\item\label{Cond:LRS_dom} or $\tau\!\ne\!\tau_\bullet$, and $\Dm(\ud\lt)$ is a union of nonempty subsets of $\tau$:
		$$
		\Dm(\ud\lt)=\udt\fE_1\cup\cdots\cup\udt\fE_\cht\ \big(\subset\tau\big),\qquad\tn{where}\quad
		\emptyset\ne\udt\fE_i\subset\fE_i\ \  \forall\ 1\!\le\!i\!\le\!\cht,
		$$
	satisfying
	\begin{enumerate}[leftmargin=*,label=(\ref{Cond:LRS_dom}\arabic*)]
		\item \label{Cond:LRS}
		$\udt\fE_\cht\!\cap\!\max(\tau)\!\ne\!\emptyset$,
		\item for every $1\!\le\!i\!<\!\cht$, $\udt\fE_i\!\subset\!\fE_i\bsl \fE_{i+1}$, and
		\item\label{Cond:ud_dom_3} for every  $1\!\le\!i\!<\!\cht$ and $e'\inn \fE_{i+1}\bsl \fE_i$,
		there exists $e\inn \udt\fE_i$ such that $e\!\prec\!e'$.
	\end{enumerate}
	\end{enumerate}
	The edges in $\udt\fE_i$ are said to be \ts{dominant} in $\fE_i$.
	The complement of $\Dm(\ud\lt)$ in $\tau$ is denoted by $$\ND(\ud\lt)\;,$$
	whose elements are called \ts{non-dominant} edges of $\ud\lt$.
	The underlying LRS of $\ud\lt$ is denoted by $(\tau_{\ud\lt},\ud\bE_{\ud\lt})$.
	
	The sets of all RLS's with dominant edges and LRS's with dominant edges are denoted by $$\ov\fT\qquad\tn{and}\qquad\ud\fT,$$ respectively.
\end{dfn}

\begin{rmk}
It is straightforward that in $\ov\lt$ (resp.~$\ud\lt$), 
the subsets $\dot\fE_i$'s (resp.~$\udt\fE_i$'s) are pairwise disjoint,
hence the union $\Dm(\ov\lt)$ (resp.~$\Dm(\ud\lt)$) is a disjoint union.

We also remark that for any $\ov\lt\inn\ov\fT$ and $1\!\le\!i\!<\!\cht$,
requiring
$\dot \fE_i\eq \fE_i\bsl \fE_{i+1}$ is equivalent to requiring the following two conditions hold simultaneously:
\begin{enumerate}[leftmargin=*]
	\item $\dot \fE_i\!\subset\!\fE_i\bsl \fE_{i+1}$ and
	\item for  every $e\inn \fE_{i+1}\bsl \fE_i$,
	there exists $e'\inn \dot\fE_i$ such that
	$e\!\prec\!e'$.
\end{enumerate}
For $\ud\lt\inn\ud\fT$, however, it is possible that $\udt \fE_i\subsetneq\fE_i\bsl \fE_{i+1}$ holds for any $1\!\le\!i\!<\!\cht$; c.f.~the following example.
\end{rmk}

\begin{rmk}
The RLS's with dominant edges can be illustrated by adding levels to the underlying rooted trees.	
Such approach has its prototype in \cite[Definition~2.1]{g1modular},
which is closely related to the notion of level graphs in~\cite{BCGGM}.
Precisely, for every $\ov\lt\eq\big(\tau,\ov\bE,\Dm(\ov\lt)\big)\inn\ov\bT$,
we consider the rooted tree $\tau$ in the usual sense of graph theory (i.e.~$\tau$ is a connected graph whose first Betti number is 0, along with a special root vertex $o$).
We begin with $o$ and put it on the top level,
and the lower endpoints of the elements of $\dot\fE_1$ one level lower than $o$ (c.f.~Remark~\ref{Rmk:rooted_tree} for terminology).
Likewise, the lower endpoints of the elements of each $\dot\fE_i$, $1\!<\!i\!\le\!\cht$, are put one level lower than those of $\dot\fE_{i-1}$.
The lower endpoints of the non-dominant edges are placed below the $\cht$-th level,
counted downward from the top;
they need not to be on any specific levels.
In this way, each
$\fE_i$ of $\ov\bE$ consists of every edge $e$ of $\tau$ that $e$ itself, or a portion of $e$, lies between the $(i-1)$-th and $i$-th levels below $o$.
An example is given in Figure~\ref{Fig:ov_t}.

\begin{figure}[htp]
\begin{center}
\begin{tikzpicture}
	\filldraw
	(0,0) circle (1.2pt)
	(0,.5) circle (1.2pt)
	(-.3,0) circle (1.2pt)
	(-.6,0) circle (1.2pt)
	(.9,0) circle (1.2pt)
	(.3,1) circle (1.2pt)
	(.6,-.5) circle (1.2pt)
	(1.2,-.5) circle (1.2pt);
	\draw
	(0,0)--(0,.5)
	(-.3,0)--(.3,1)--(1.2,-.5)
	(-.6,0)--(0,.5)--(.6,-.5);
	\draw[dotted]
	(-.9,1)--(1.5,1)
	(-.9,.5)--(1.5,.5)
	(-.9,0)--(1.5,0);
	\draw[decorate,decoration=brace,thin]
	(1.5,.95)--(1.5,.55);
	\draw[decorate,decoration=brace,thin]
	(1.5,.45)--(1.5,.05);
	\draw
	(0,.75) node {\tiny{$a$}}
	(-.31,.16) node {\tiny{$c$}}
	(-.5,.25) node {\tiny{$b$}}
	(-.09,.16) node {\tiny{$d$}}
	(.6,-.2) node {\tiny{$f$}}
	(.92,.25) node {\tiny{$g$}}
	(1.2,-.25) node {\tiny{$z$}}
	(1.5,.75) node[right] {\tiny{$\fE_1$}}
	(1.5,.25) node[right] {\tiny{$\fE_2$}};
	\draw
	(5,1) node[right] {\tiny{$\ov\lt=(\,\tau\,,\,\{\fE_1,\fE_2\}
			\,,\,\{a,b,c,d,g\}\,)$}}
	(5,.5) node[right] {\tiny{$\fE_{1}=\{a,g\},$}}
	(7.5,.5) node[right] {\tiny{$\fE_{2}=\{b,c,d,f,g\},$}}
	(5,0) node[right] {\tiny{$\dot\fE_1=\{a\},$}}
	(7.5,0) node[right] {\tiny{$\dot\fE_2=\{b,c,d,g\},$}}
	;
\end{tikzpicture}	
\end{center}
\caption{Illustration of an RLS with dominant edges}\label{Fig:ov_t}	
\end{figure}  

Illustrating an LRS with dominant edges is somewhat different,
mainly because  $\udt\fE_i$ is possibly a strict subset of $\fE_i\bsl\fE_{i+1}$.
It thus may lose some data on the dominant edges if we simply put vertices on various levels as in the RL case.
So for every $\ud\lt\eq\big(\tau,\ud\bE,\Dm(\ud\lt)\big)\inn\ud\bT$,
we still consider $\tau$ in the usual sense,
but begin with $\fE_1$, the set of the leaves, instead.
We put the lower endpoints of the elements of $\fE_1$ on the bottom level,
and the upper endpoints of the elements of $\udt\fE_1$ one level higher.
Likewise,
the upper endpoints of the elements of each $\udt\fE_i$,
$1\!<\!i\!\le\!\cht$, are put one level higher than those of $\udt\fE_{i-1}$;
in particular, the root of $\tau$ is on the top level due to Condition~\ref{Cond:LRS} of Definition~\ref{Dfn:Enhanced_rtl_seq}.
For every $e\inn\ND(\ud\lt)$,
we write
\begin{align*}
	\fE':=\max\big\{\,\fE_i\,:\;
	1\le i\le\cht\,,\;
	\fE_i\cap\{e\}^\preceq\ne\emptyset\,\big\}\,,
\end{align*}
where the set on the right-hand side is non-empty because it contains $\fE_1$, and the maximum is taken with respect to (\ref{Eqn:transverse_sections_order}).
We then place the upper endpoint of $e$ on the same level as those of $\udt\fE'$.
In this way, each 
$\fE_i$ of $\ud\bE$ consists of every edge $e$ of $\tau$ that $e$ itself, or a portion of $e$, lies between the $i$-th and $(i+1)$-th levels, counted upward from the bottom;
moreover, for every edge $e$ that is not contained in any $\fE_i$ of $\ud\bE$,
the endpoints of $e$ must be on the same level,
so we put such $e$ entirely on that level.
Finally, we thicken
the dominant edges for distinction.
An example is given in Figure~\ref{Fig:ud_t}.

\begin{figure}[htp]
\begin{center}
\begin{tikzpicture}
	\filldraw
	(.2,.5) circle (1.2pt)
	(-.6,.9) circle (1.2pt)
	(-.8,.5) circle (1.2pt)
	(-.4,1.3) circle (1.2pt)
	(.4,1.3) circle (1.2pt)
	(-.2,.1) circle (1.2pt)
	(.2,.1) circle (1.2pt)
	(-1,.1) circle (1.2pt)
	(-.6,.1) circle (1.2pt)
	(.6,.1) circle (1.2pt)
	(1,.1) circle (1.2pt)
	;
	
	\draw
	(.4,1.3)--(-.4,1.3)--(-1,.1)
	(.2,.5)--(.4,1.3)--(1,.1)
	(-.6,.1)--(-.6,.9)--(-.2,.1)
	(.6,.1)--(.2,.5)--(.2,.1)
	;
	
	\draw[ultra thick]
	(-.8,.5)--(-1,.1)
	(-.4,1.3)--(-.6,.9)--(-.2,.1)
	(.2,.1)--(.2,.5)
	;
	
	\draw[dotted]
	(-1.1,.1)--(1.1,.1)
	(-1.1,.5)--(1.1,.5)
	(-1.1,.9)--(1.1,.9)
	(-1.1,1.3)--(1.1,1.3)
	;
	
	\draw[decorate,decoration=brace,thin]
	(1.1,1.25)--(1.1,.95);
	\draw[decorate,decoration=brace,thin]
	(1.1,.85)--(1.1,.55);
	\draw[decorate,decoration=brace,thin]
	(1.1,.45)--(1.1,.15);
	
	\draw
	(.12,.7) node {\tiny{$w$}}
	(.37,.2) node {\tiny{$v$}}
	(.8,.25) node {\tiny{$z$}}
	(.07,.25) node {\tiny{$u$}}
	(-.41,.25) node {\tiny{$y$}}
	(-.65,1.1) node {\tiny{$a$}}
	(-1.05,.26) node {\tiny{$c$}}
	(-.85,.7) node {\tiny{$b$}}
	(-.7,.25) node {\tiny{$d$}}
	(0,1.18) node {\tiny{$x$}}
	(-.4,1.45) node {\tiny{$o$}}
	(1.15,1.1) node[right] {\tiny{$\fE_3$}}
	(1.15,.7) node[right] {\tiny{$\fE_2$}}
	(1.15,.3) node[right] {\tiny{$\fE_1$}}
	;
	
	\draw[xshift=-2cm,yshift=.3cm]
	(5,1) node[right] {\tiny{$\ud\lt=(\,\tau\,,\,\{\fE_1,\fE_2,\fE_3\}\,,\,\{a,c,y,u\}\,)$}}
	(5,.5) node[right] {\tiny{$\fE_{1}=\{c,d,y,u,v,z\},$}}
	(8,.5) node[right] {\tiny{$\fE_{2}=\{b,d,y,w,z\},$}}
	(11,.5) node[right] {\tiny{$\fE_{3}=\{a,w,z\},$}}
	(5,0) node[right] {\tiny{$\udt\fE_1=\{c,u\},$}}
	(8,0) node[right] {\tiny{$\udt\fE_2=\{y\},$}}
	(11,0) node[right] {\tiny{$\udt\fE_3=\{a\},$}}
	;
\end{tikzpicture}		
\end{center}
\caption{Illustration of an LRS with dominant edges}\label{Fig:ud_t}		
\end{figure}
\end{rmk}

\begin{cnv}\label{Convention:R2L_dom}
	For every $\ov\lt\eq\big(\tau,\ov\bE,\Dm(\ov\lt)\big)\inn\ov\fT$ and every $\fE\inn\ov\bE$,
	the corresponding set of dominant edges is written as $\dot\fE$ when the context is clear.
	
	Similarly, for every $\ud\lt\eq\big(\tau,\ud\bE,\Dm(\ud\lt)\big)\inn\ud\fT$ and every $\fE\inn\ud\bE$,
	the corresponding set of dominant edges is written as $\udt\fE$ when the context is clear.
\end{cnv}

By Convention~\ref{Convention:R2L_dom},
for $\ov\lt\inn\ov\fT$ and $\ud\lt\inn\ud\fT$, we can write
\begin{align*}
	\Dm(\ov\lt)=\bigsqcup_{\fE\in\ov\bE}\dot\fE\qquad\tn{and}\qquad
	\Dm(\ud\lt):=\bigsqcup_{\fE\in\ud\bE}\udt\fE,
\end{align*}
respectively.

\begin{lmm}\label{Lm:dominant}
Let $\ov\lt\inn\ov\fT$ and $\ud\lt\inn\ud\fT$.
Then,
\begin{enumerate}[leftmargin=*,label=(\arabic*)]
	\item \label{Cond:dom_complement}
	we have
	\begin{align*}
	\tau_{\ov\lt}&=\Dm(\ov\lt)\sqcup\ND(\ov\lt)=\Dm(\ov\lt)\sqcup(\fE_\cht\bsl\dot\fE_\cht)\sqcup\fE_\cht^\prec\,,\\
	\tau_{\ud\lt}&=\Dm(\ud\lt)\sqcup\ND(\ud\lt)=\Dm(\ud\lt)\sqcup\bigsqcup_{\fE\in\ud\bE}\!\big((\fE\bsl\fE_+)\big\bsl\udt\fE\big) \sqcup\fE_\cht^\succ\,;
	\end{align*}
	
	\item \label{Cond:dom_above}
	for any $e\inn\Dm(\ov\lt)$ and any $e'\!\succeq\!e$,
	we have $e'\inn\Dm(\ov\lt)$; and
	
	\item \label{Cond:dom_below}
	for any $e\inn\Dm(\ud\lt)$,
	there exists $e'\inn\Dm(\ud\lt)\!\cap\!\min(\tau_{\ud\lt})$ such that $e''\inn\Dm(\ud\lt)$ for all $e'\!\preceq\!e''\!\preceq\!e$.
\end{enumerate}
\end{lmm}

\begin{proof}
Part~\ref{Cond:dom_complement} follows from Definitions~\ref{Dfn:R2L_seq} and~\ref{Dfn:Enhanced_rtl_seq},
while Part~\ref{Cond:dom_above} follows from Part~\ref{Cond:dom_complement}.

To establish Part~\ref{Cond:dom_below},
let $\fE\inn\ud\bE$ be such that $e\inn\udt\fE$.
If $e$ is minimal, then Part~\ref{Cond:dom_below} trivially holds.
If $e$ is not minimal,
let $\fE'\eq\max\{\fE''\inn\ud\bE:\fE''\!\prec\!\fE,\,e\!\not\in\!\fE''\}$,
then $e\inn\fE'_+\bsl\fE'$,
so by Criterion~\ref{Cond:ud_dom_3} of Definition~\ref{Dfn:Enhanced_rtl_seq},
we can find $e'\inn\udt\fE'$ satisfying $e'\!\prec\!e$.
Applying the above procedure repeatedly if necessary, we obtain Part~\ref{Cond:dom_below}.
\end{proof}

Lemma~\ref{Lm:dominant}~\ref{Cond:dom_below} implies that given $e\inn\Dm(\ud\lt)$, not every edge in $\{e\}^\preceq$ is necessarily dominant.
We thus write
\begin{align}\label{Eqn:Omega_e}
\Om(e):=
\big\{\,e'\in\{e\}^\preceq\!\cap\min(\tau_{\ud\lt})\,:~
e''\in\Dm(\ud\lt)\;\forall\,e'\!\preceq\!e''\!\preceq\!e\,\big\}\,.
\end{align}
As long as $\ud\lt\!\ne\!\ud\lt_\bullet$, 
we always have  
$$
 \Om(e)\ne\emptyset\qquad\tn{and}\qquad
 \Om(e)^\succeq\subset\Dm(\ud\lt)
$$
because of Definition~\ref{Dfn:Enhanced_rtl_seq}~\ref{Cond:ud_dom_3} and (\ref{Eqn:Omega_e}).
If $e$ itself is not minimal, then for every $e'\inn\Om(e)$,
we further set
\begin{align}\label{Eqn:fp_e}
\fp_e(e'):=\max\{e'':\,e'\!\preceq\!e''\prec\!e\}\,.
\end{align}
It is direct to check that $\fp_e(e')\inn\Dm(\ud\lt)$ and $e'\inn\Om\big(\fp_e(e')\big).$
Such $\Om(e)$ and $\fp_e(e')$ will be needed for constructing the stacks with LR twisted fields later on.

Every isomorphism of rooted trees induces a natural way of identifying the RLS's (resp.~LRS's) with dominant edges over them:
\begin{lmm}
\label{Lm:RLS_LRS_isom}
Let $\phi:\tau\!\lra\!\tau'$ be an isomorphism of rooted trees $\tau$ and $\tau'$.
Then, for every $\ov\lt=\big(\tau,\ov\bE,\Dm(\ov\lt)\big)\inn\ov\bT$, we have
\begin{align*}
	\phi(\ov\lt):=
	\Big(\,\tau'\,,\ \phi(\ov\bE)\!:=\{\phi(\fE):\fE\inn\ov\bE\}\,,\ \phi\big(\Dm(\ov\lt)\big)\,\Big)\in\ov\bT.
\end{align*}
Similarly,
for every $\ud\lt=\big(\tau,\ud\bE,\Dm(\ud\lt)\big)\inn\ud\bT$, we have
\begin{align*}
	\phi(\ud\lt):=
	\Big(\,\tau'\,,\ \phi(\ud\bE)\!:=\{\phi(\fE):\fE\inn\ud\bE\}\,,\ \phi\big(\Dm(\ud\lt)\big)\,\Big)\in\ud\bT.
\end{align*}
\end{lmm} 	
\begin{proof}
Both of the statements follow from direct verification. We omit the details.
\end{proof}

In the STF theory,
an STF is described via stratification,
and each stratum corresponds to some $\ov\lt\inn\ov\bT$ or $\ud\lt\inn\ud\bT$.
Near such a stratum,
there are local parameters (called twisted parameters) which serve as the modular parameters on the STF in Theorem~\ref{Thm:Main}.
Part of these local parameters are indexed by the following sets:
\begin{equation}\begin{split}\label{Eqn:bbI}
	&\bbI(\ov\lt):=
	\ov\bE\sqcup\ND(\ov\lt)=
	\ov\bE\sqcup \big(\fE_\cht\bsl\dot\fE_\cht\big) \sqcup\fE_\cht^\prec,\\
	&
	\bbI(\ud\lt):=
	\ud\bE\sqcup\ND(\ud\lt)=
	\ud\bE\sqcup \bigsqcup_{\fE\in\ud\bE}\!\big((\fE\bsl\fE_+)\big\bsl\udt\fE\big) \sqcup\fE_\cht^\succ.
\end{split}\end{equation}
Notice that for $\ov\lt\inn\ov\fT$,
the union $\bigsqcup_{\fE\in\ov\bE}\big((\fE\bsl\fE_+)\big\bsl\dot\fE\big)$ is the same as $\fE_\cht\bsl\dot\fE_\cht$.

In order to describe the nearby strata of a given stratum of an STF,
we observe that
every subset $J\!\subset\!\bbI(\ov\lt)$ determines an RLS with dominant edges
as follows:
\begin{itemize}[leftmargin=*]
	\item
	The underlying rooted tree is given by $\tau_{(J)}=\tau\big\bsl\,\ov\sfE_{J}^{\,\wedge}$, where
	\begin{align}\label{Eqn:C_J}
	\ov\sfE_{J}:=
	\bigg(\Big(\Dm(\ov\lt)\sqcup\big(J\!\cap\!(\fE_\cht\bsl\dot\fE_\cht)\big)\Big)
	\Big\bsl 
	\Big(\bigcup_{\fE\in\ov\bE\bsl J}\!\!\!\fE\Big)\bigg)
	\sqcup\,
	\big(J\!\cap\!\fE_\cht^\prec\big).
	\end{align}
	
	\item 
	The RLS is given by
	$\ov\bF\eq\ov\bE\bsl J$.
	If $\ov\bF\!\ne\!\emptyset$,
	then the terms of $\ov\bF$ are written as $\fF_1\!\succ\!\ldots\!\succ\!\fF_r$, corresponding to $\fE_{i_1}\!\succ\!\ldots\!\succ\!\fE_{i_r}\inn\ov\bE$.
	
	\item 
	If $\ov\bF\!\ne\!\emptyset$,
	the set of the dominant edges is given by
	\begin{align}\label{Eqn:t_J_dom'}
		\Big(\Dm(\ov\lt)\sqcup 
		\big(J\!\cap\!(\fE_\cht\bsl\dot\fE_\cht)\big)\Big)\cap 
		\big(\bigcup_{\fF\in\ov\bF}\fF\,\big).
	\end{align} 
	In other words, for every $1\!\le\!k\!\le\!r$,
	\begin{align}
	\label{Eqn:t_J_dom}
		\dot\fF_k=
		\Big(\Dm(\ov\lt)\!\sqcup
		\! 
		\big(J\!\cap\!(\fE_\cht\bsl\dot\fE_\cht)\big)\Big)\cap \fE_{i_k}=
		\begin{cases}
		\bigcup_{i_k\le j<i_{k+1}}\!\!\big(\dot\fE_j\!\cap\!\fE_{i_k}
		\big)
		&
		\tn{if}~k\!<\!r,
		\\
		\bigcup_{i_r\le j\le\cht}\!\big(\dot\fE_j\!\cap\!\fE_{i_r}
		\big)\sqcup\big(
		 J\!\cap\!\fE_{i_r}\big)
		&
		\tn{if}~k\eq r.	
		\end{cases}
	\end{align}
\end{itemize}

\begin{lmm}
	\label{Lm:t_J}
With notation as above,
the tuple
\begin{align*}
	\ov\lt_{(J)}:=\big(\,\tau_{(J)}\,,\,
	\ov\bF\,,\,
	\Dm(\ov\lt_{(J)})\big)
\end{align*}
is an RLS with dominant edges as per Definition~\ref{Dfn:Enhanced_rtl_seq};
i.e.~$\ov\lt_{(J)}\inn\ov\bT$.
\end{lmm}

\begin{proof}
First,
observe that if $\ov\bE\!\subset\! J$,
then by Definition~\ref{Dfn:Enhanced_rtl_seq} and Lemma~\ref{Lm:dominant},
there exists $e\inn\min(\tau)$ such that $\{e\}^\succeq\!\subset\!\Dm(\ov\lt)\!\subset\!\ov\sfE_J$.
Lemmas~\ref{Lm:E^_simple} and~\ref{Lm:path^} then imply $\ov\sfE_J^{\,\wedge}\eq(\{e\}^\succeq)^\wedge\eq\tau$.
Hence $\ov\lt_{(J)}\eq\ov\lt_\bullet\inn\ov\bT$.

Hereafter, we assume $\ov\bE\!\not\subset\! J$, i.e.~$\ov\bF\!\ne\!\emptyset$.
We shall prove that
\begin{equation}\begin{split}
	\label{Eqn:t_J_check}
&\fF\in\Xi(\tau_{(J)})\quad\forall~\fF\inn\ov\bF,\qquad
\fF_1=\max(\tau_{(J)}),\qquad
\tau_{(J)}= \fF_r^\prec\sqcup\bigcup_{\fF\in\ov\bF}\fF,\qquad\\
&\dot\fF_r\cap\min(\tau_{(J)})\ne\emptyset,\qquad\tn{and}\qquad
\dot\fF_k=\fF_k\bsl\fF_{k+1}\quad\forall~1\!\le\!k\!<\!r,
\end{split}\end{equation}
which, along with the fact that $\fF_1\!\succ\!\ldots\!\succ\!\fF_r$, will imply $\ov\lt_{(J)}\inn\ov\bT$.

In fact, 
by (\ref{Eqn:C_J}),
each $\fF\inn\ov\bF\eq\ov\bE\bsl J$ satisfies $\fF\!\cap\!\ov\sfE_J\eq\emptyset$,
so $\fF\inn\Xi^{\tn c}(\tau;\ov\sfE_J)$. 
The first equality of (\ref{Eqn:t_J_check}) then follows from Lemma~\ref{Lm:E^}~\ref{Cond:adj_TS_edge_contr}.

By Lemma~\ref{Lm:dominant},
we have
\begin{align*}
	\fE_{i_1}^\succ\ \subset\ \Dm(\ov\lt)
	\Big\bsl 
	\Big(\bigcup_{\fE\in\ov\bE\bsl J}\!\!\!\fE\Big) \ 
	\subset\ \ov\sfE_J.
\end{align*}
This justifies the second equality of (\ref{Eqn:t_J_check}).

To see the third equality of (\ref{Eqn:t_J_check}),
notice that for every $e'\inn\tau_{(J)}\!\subset\!\tau$,
either $e'\inn\fE_\cht^\prec$, or $e'\inn\bigcup_{\fE\in\ov\bE}\fE$.
In the former case,
since $\fE_\cht\!\preceq\!\fE_{i_r}\eq\fF_r$,
we have $e'\inn\fF^\prec$.
In the latter case,
either $e\inn\fF$ for some $\fF\inn\ov\bF$,
hence $e\inn \bigcup_{\fF\in\ov\bF}\fF$;
or $e'\!\not\in\!\fF$ for all $\fF\inn\ov\bF$,
which implies $e'\inn(\fE_\cht\bsl\dot\fE_\cht)\bsl J$ (for otherwise $e'$ would be in $\ov\sfE_J$ and thus not in $\tau_{(J)}$),
hence $e\inn\fF_r^\preceq$.
In sum, the third equality of (\ref{Eqn:t_J_check}) holds.

Regarding $\dot\fF_r\!\cap\!\min(\tau_{(J)})$, since
$\ov\lt\inn\ov\bT$, there exists $e\inn\dot\fE_\cht\!\cap\!\min(\tau)$.
In addition,
by (\ref{Eqn:C_J}) and Definition~\ref{Dfn:E^}, for every $e\inn\Dm(\ov\lt)$, we have
\begin{align}
	\label{Eqn:C_J^}
	e\notin \ov\sfE_J^{\,\wedge}\quad\big(\,\tn{i.e.}~e\in\tau_{(J)}\,\big)\qquad\tn{whenever}\qquad
	e\notin\ov\sfE_J\,.
\end{align}
\begin{itemize}[leftmargin=*]
\item 
If $e\!\notin\!\ov\sfE_J$,
from (\ref{Eqn:C_J^}) we see that $e\inn\tau_{(J)}$, hence $e\inn(\dot\fE_\cht\!\cap\!\fE_{i_r})\!\subset\!\dot\fF_r$.
In addition, $e\inn\min(\tau)$ implies $e\inn\min(\tau_{(J)})$.

\item
If $e\!\in\!\ov\sfE_J$,
since $\fE_{i_r}$ is a transverse section,
we see that $\{e\}^\succeq\!\cap\!\fE_{i_r}$ contains exactly one element, which is denoted by  $e^\sharp$.
Then, $e^\sharp\!\notin\!\ov\sfE_J$, and by (\ref{Eqn:C_J^}),
we have $e^\sharp\inn\tau_{(J)}$.
Taking Lemmas~\ref{Lm:path^} and~\ref{Lm:E^_simple} into consideration,
we conclude that $e^\sharp\inn\min(\tau_{(J)})$.
In addition,
Lemma~\ref{Lm:dominant} and the $k\eq r$ case of (\ref{Eqn:t_J_dom}) imply $e^\sharp\inn\dot\fF_r$.
\end{itemize}
In sum, $\dot\fF_r\!\cap\!\min(\tau_{(J)})\ne\emptyset$.

Finally, the last equality of (\ref{Eqn:t_J_check}) follows directly from the $k\!<\!r$ case of (\ref{Eqn:t_J_dom}).
\end{proof}

Similarly,
every subset $ J\!\subset\!\bbI(\ud\lt)$ determines an LRS with dominant edges 
\begin{align}
	\label{Eqn:t_I_LRS} \ud\lt_{(J)}=\big(\,\tau_{( J)}\,,\,
\ud\bF\,,\,
\Dm(\ud\lt_{( J)})\big)\,\in\ud\bT
\end{align} as follows:
\begin{itemize}[leftmargin=*]
	\item
	The underlying rooted tree is given by $\tau_{( J)}=\tau\big\bsl\,\ud\sfE_{ J}^{\;\wedge}$, where
	\begin{align}\label{Eqn:C_J_revert}
	\ud\sfE_{ J}:=
	\Bigg(\bigg(\Dm(\ud\lt)\sqcup
	\bigsqcup_{\fE\in\ud\bE}
	\Big( J\cap\big((\fE\bsl\fE_+)\big\bsl\udt\fE\big)\Big)
	\bigg)
	\bigg\bsl 
	\Big(\bigcup_{\fE\in\ud\bE\bsl J}\!\!\!\fE\Big)\Bigg)
	\sqcup\,
	\big( J\!\cap\!\fE_\cht^\succ\big).
	\end{align}
	\item 
	The LRS is given by $\ud\bF\eq\ud\bE\bsl J$.
	The terms of $\ud\bF$ are written as 
	$\fF_1\!\prec\!\ldots\!\prec\!\fF_r$,
	corresponding to $\fE_{i_1}\!\prec\!\ldots\!\prec\!\fE_{i_r}\inn\ud\bE$.
	\item 
	If $\ud\bF\!\ne\!\emptyset$,
	the set of the dominant edges is given by
	\begin{align*}
	\bigg(\Dm(\ud\lt)\sqcup
	\bigsqcup_{\fE\in\ud\bE}
	\Big( J\cap\big((\fE\bsl\fE_+)\big\bsl\udt\fE\big)\Big)
	\bigg)\cap 
	\big(\bigcup_{\fF\in\ud\bF}\fF\,\big).
	\end{align*} 
	In other words, for every $1\!\le\!k\!\le\!r$,
	\begin{align*}
	\udt\fF_k&=
	\bigg(\Dm(\ud\lt)\sqcup
	\bigsqcup_{\fE\in\ud\bE}\!
	\Big( J\!\cap\!\big((\fE\bsl\fE_+)\big\bsl\udt\fE\big)\Big)\!
	\bigg)\cap \fE_{i_k}
	=
	\begin{cases}
		\bigcup_{i_k\le j<i_{k+1}}\!\Big(\big(\udt\fE_j\!\sqcup\!( J\!\cap\!\fE_j)\big)\cap\!\fE_{i_k}\Big)
		&
		\tn{if}~k\!<\!r,
		\\
		\bigcup_{i_k\le j\le\cht}\!\Big(\big(\udt\fE_j\!\sqcup\!( J\!\cap\!\fE_j)\big)\cap\!\fE_{i_k}\Big)
		&
		\tn{if}~k\!=\!r.
	\end{cases}	
	\end{align*}
\end{itemize}

The reason that $\ud\lt_{(J)}\inn\ud\bT$ is parallel to its RL counterpart; c.f.~the proof of Lemma~\ref{Lm:t_J}.
We omit further details.

In Figure~\ref{Fig:t(I)},
we provide some examples of $\ov\lt_{(J)}$ and $\ud\lt_{(J)}$.

\begin{figure}[htp]
\begin{center}
\begin{tikzpicture}
	\filldraw
	(0,0) circle (1.2pt)
	(0,.5) circle (1.2pt)
	(-.3,0) circle (1.2pt)
	(-.6,0) circle (1.2pt)
	(.9,0) circle (1.2pt)
	(.3,1) circle (1.2pt)
	(.6,-.5) circle (1.2pt)
	(1.2,-.5) circle (1.2pt);
	\draw
	(0,0)--(0,.5)
	(-.3,0)--(.3,1)--(1.2,-.5)
	(-.6,0)--(0,.5)--(.6,-.5);
	\draw[dotted]
	(-.9,0)--(1.5,0)
	(-.9,1)--(1.5,1)
	(-.9,.5)--(1.5,.5);
	\draw[thin, decorate, decoration=brace]
	(1.5,.95)--(1.5,.55);
	\draw[thin, decorate, decoration=brace]
	(1.5,.45)--(1.5,.05);
	\draw
	(0,.75) node {\tiny{$a$}}
	(-.31,.16) node {\tiny{$c$}}
	(-.5,.25) node {\tiny{$b$}}
	(-.09,.16) node {\tiny{$d$}}
	(.6,-.2) node {\tiny{$f$}}
	(.92,.25) node {\tiny{$g$}}
	(1.2,-.25) node {\tiny{$z$}}
	(1.55,.75) node[right] {\tiny{$\fE_1$}}
	(1.55,.25) node[right] {\tiny{$\fE_2$}}
	(.3,-1) node {\tiny{$\ov\lt$}}
	;
	
	\filldraw[xshift=4cm]
	(0,0) circle (1.2pt)
	(0,.5) circle (1.2pt)
	(-.3,0) circle (1.2pt)
	(-.6,0) circle (1.2pt)
	(.9,0) circle (1.2pt)
	(.3,1) circle (1.2pt)
	(.6,-.5) circle (1.2pt)
	;
	\draw[xshift=4cm]
	(0,0)--(0,.5)
	(-.3,0)--(.3,1)--(.9,0)
	(-.6,0)--(0,.5)--(.6,-.5);
	\draw[dotted,xshift=4cm]
	(-.9,0)--(1.2,0)
	(-.9,1)--(1.2,1)
	(-.9,.5)--(1.2,.5);
	\draw[thin, decorate, decoration=brace, xshift=4cm]
	(1.2,.95)--(1.2,.55);
	\draw[thin, decorate, decoration=brace, xshift=4cm]
	(1.2,.45)--(1.2,.05);
	\draw[xshift=4cm]
	(0,.75) node {\tiny{$a$}}
	(-.31,.16) node {\tiny{$c$}}
	(-.5,.25) node {\tiny{$b$}}
	(-.09,.16) node {\tiny{$d$}}
	(.6,-.2) node {\tiny{$f$}}
	(.92,.25) node {\tiny{$g$}}
	(1.25,.75) node[right] {\tiny{$\fE_1$}}
	(1.25,.25) node[right] {\tiny{$\fE_2$}}
	(.5,-1) node {\tiny{$\ov\lt_{(J_1)}~\tn{where}~J_1\!:=\!\{z\}$}}
	;
	
	\filldraw[xshift=7.7cm]
	(0,.5) circle (1.2pt)
	(0,1) circle (1.2pt)
	(-.3,.5) circle (1.2pt)
	(-.6,.5) circle (1.2pt)
	(.3,.5) circle (1.2pt)
	(.6,.5) circle (1.2pt)
	(1.2,0) circle (1.2pt)
	;
	\draw[xshift=7.7cm]
	(0,1)--(0,.5)
	(-.3,.5)--(0,1)--(1.2,0)
	(-.6,.5)--(0,1)--(.3,.5);
	\draw[dotted,xshift=7.7cm]
	(-.9,1)--(1.5,1)
	(-.9,.5)--(1.5,.5)
	;
	\draw[thin, decorate, decoration=brace, xshift=7.7cm]
	(1.5,.95)--(1.5,.55);
	\draw[xshift=7.7cm]
	(-.31,.66) node {\tiny{$c$}}
	(-.5,.75) node {\tiny{$b$}}
	(-.09,.66) node {\tiny{$d$}}
	(.1,.66) node {\tiny{$f$}}
	(.52,.75) node {\tiny{$g$}}
	(1.1,.25) node {\tiny{$z$}}
	(1.55,.75) node[right] {\tiny{$\fE_2$}}
	(.7,-1) node {\tiny{$\ov\lt_{(J_2)}~\tn{where}~J_2\!:=\!\{\fE_1,f\}$}}
	;
	
	\filldraw[xshift=11.4cm]
	(0,1) circle (1.2pt)
	(-.3,.5) circle (1.2pt)
	(.3,.5) circle (1.2pt)
	(.6,0) circle (1.2pt)
	;
	\draw[xshift=11.4cm]
	(-.3,.5)--(0,1)--(.6,0)
	;
	\draw[dotted,xshift=11.4cm]
	(-.6,1)--(.9,1)
	(-.6,.5)--(.9,.5)
	;
	\draw[thin, decorate, decoration=brace, xshift=11.4cm]
	(.9,.95)--(.9,.55);
	\draw[xshift=11.4cm]
	(-.3,.75) node {\tiny{$a$}}
	(.3,.75) node {\tiny{$g$}}
	(.6,.25) node {\tiny{$z$}}
	(.95,.75) node[right] {\tiny{$\fE_1$}}
	(.8,-1) node {\tiny{$\ov\lt_{(J_3)}~\tn{where}~J_3\!:=\!\{\fE_2\}$}}
	;
	
	\filldraw[yshift=-4cm,xshift=.4cm]
	(.2,.5) circle (1.2pt)
	(-.6,.9) circle (1.2pt)
	(-.8,.5) circle (1.2pt)
	(-.4,1.3) circle (1.2pt)
	(.4,1.3) circle (1.2pt)
	(-.2,.1) circle (1.2pt)
	(.2,.1) circle (1.2pt)
	(-1,.1) circle (1.2pt)
	(-.6,.1) circle (1.2pt)
	(.6,.1) circle (1.2pt)
	(1,.1) circle (1.2pt)
	;
	
	\draw[yshift=-4cm,xshift=.4cm]
	(.4,1.3)--(-.4,1.3)--(-1,.1)
	(.2,.5)--(.4,1.3)--(1,.1)
	(-.6,.1)--(-.6,.9)--(-.2,.1)
	(.6,.1)--(.2,.5)--(.2,.1)
	;
	
	\draw[ultra thick,yshift=-4cm,xshift=.4cm]
	(-.8,.5)--(-1,.1)
	(-.4,1.3)--(-.6,.9)--(-.2,.1)
	(.2,.1)--(.2,.5)
	;
	
	\draw[dotted,yshift=-4cm,xshift=.4cm]
	(-1.1,.1)--(1.1,.1)
	(-1.1,.5)--(1.1,.5)
	(-1.1,.9)--(1.1,.9)
	(-1.1,1.3)--(1.1,1.3)
	;
	
	\draw[decorate,decoration=brace,thin,yshift=-4cm,xshift=.4cm]
	(1.1,1.25)--(1.1,.95);
	\draw[decorate,decoration=brace,thin,yshift=-4cm,xshift=.4cm]
	(1.1,.85)--(1.1,.55);
	\draw[decorate,decoration=brace,thin,yshift=-4cm,xshift=.4cm]
	(1.1,.45)--(1.1,.15);
	
	\draw[yshift=-4cm,xshift=.4cm]
	(.12,.7) node {\tiny{$w$}}
	(.37,.2) node {\tiny{$v$}}
	(.8,.25) node {\tiny{$z$}}
	(.07,.25) node {\tiny{$u$}}
	(-.41,.25) node {\tiny{$y$}}
	(-.65,1.1) node {\tiny{$a$}}
	(-1.05,.26) node {\tiny{$c$}}
	(-.85,.7) node {\tiny{$b$}}
	(-.7,.25) node {\tiny{$d$}}
	(0,1.18) node {\tiny{$x$}}
	(-.4,1.45) node {\tiny{$o$}}
	(1.15,1.1) node[right] {\tiny{$\fE_3$}}
	(1.15,.7) node[right] {\tiny{$\fE_2$}}
	(1.15,.3) node[right] {\tiny{$\fE_1$}}
	(-.1,-.4) node {\tiny{$\ud\lt$}}
	;
	
	\filldraw[yshift=-4cm,xshift=4.2cm]
	(.2,.5) circle (1.2pt)
	(-.6,.9) circle (1.2pt)
	(-.8,.5) circle (1.2pt)
	(-.4,1.3) circle (1.2pt)
	(-.2,.1) circle (1.2pt)
	(.2,.1) circle (1.2pt)
	(-1,.1) circle (1.2pt)
	(-.6,.1) circle (1.2pt)
	(.6,.1) circle (1.2pt)
	(1,.1) circle (1.2pt)
	;
	
	\draw[yshift=-4cm,xshift=4.2cm]
	(.2,.5)--(-.4,1.3)--(-1,.1)
	(-.6,.1)--(-.6,.9)--(-.2,.1)
	(.6,.1)--(.2,.5)--(.2,.1)
	;
	
	\draw[ultra thick,yshift=-4cm,xshift=4.2cm]
	(-.8,.5)--(-1,.1)
	(-.4,1.3)--(-.6,.9)--(-.2,.1)
	(.2,.1)--(.2,.5)
	(-.4,1.3)--(1,.1)
	;
	
	\draw[dotted,yshift=-4cm,xshift=4.2cm]
	(-1.1,.1)--(1.1,.1)
	(-1.1,.5)--(1.1,.5)
	(-1.1,.9)--(1.1,.9)
	(-1.1,1.3)--(1.1,1.3)
	;
	
	\draw[decorate,decoration=brace,thin,yshift=-4cm,xshift=4.2cm]
	(1.1,1.25)--(1.1,.95);
	\draw[decorate,decoration=brace,thin,yshift=-4cm,xshift=4.2cm]
	(1.1,.85)--(1.1,.55);
	\draw[decorate,decoration=brace,thin,yshift=-4cm,xshift=4.2cm]
	(1.1,.45)--(1.1,.15);
	
	\draw[yshift=-4cm,xshift=4.2cm]
	(-.12,.7) node {\tiny{$w$}}
	(.37,.2) node {\tiny{$v$}}
	(.5,.7) node {\tiny{$z$}}
	(.07,.25) node {\tiny{$u$}}
	(-.41,.25) node {\tiny{$y$}}
	(-.65,1.1) node {\tiny{$a$}}
	(-1.05,.26) node {\tiny{$c$}}
	(-.85,.7) node {\tiny{$b$}}
	(-.7,.25) node {\tiny{$d$}}
	(-.4,1.45) node {\tiny{$o$}}
	(1.15,1.1) node[right] {\tiny{$\fE_3$}}
	(1.15,.7) node[right] {\tiny{$\fE_2$}}
	(1.15,.3) node[right] {\tiny{$\fE_1$}}
	(.3,-.4) node {\tiny{$\ud\lt_{(J_1)}~\tn{where}~J_1\!:=\!\{x,z\}$}}
	;
	
	\filldraw[yshift=-4cm,xshift=8cm]
	(.2,.5) circle (1.2pt)
	(-.6,.9) circle (1.2pt)
	(-.4,1.3) circle (1.2pt)
	(.4,1.3) circle (1.2pt)
	(-.2,.5) circle (1.2pt)
	(-1,.5) circle (1.2pt)
	(-.6,.5) circle (1.2pt)
	(1,.5) circle (1.2pt)
	;
	
	\draw[yshift=-4cm,xshift=8cm]
	(.4,1.3)--(-.4,1.3)--(-.6,.9)--(-1,.5)
	(.2,.5)--(.4,1.3)--(1,.5)
	(-.6,.5)--(-.6,.9)
	;
	
	\draw[ultra thick,yshift=-4cm,xshift=8cm]
	(-.4,1.3)--(-.6,.9)--(-.2,.5)
	;
	
	\draw[dotted,yshift=-4cm,xshift=8cm]
	(-1.1,.5)--(1.1,.5)
	(-1.1,.9)--(1.1,.9)
	(-1.1,1.3)--(1.1,1.3)
	;
	
	\draw[decorate,decoration=brace,thin,yshift=-4cm,xshift=8cm]
	(1.1,1.25)--(1.1,.95);
	\draw[decorate,decoration=brace,thin,yshift=-4cm,xshift=8cm]
	(1.1,.85)--(1.1,.55);
	
	\draw[yshift=-4cm,xshift=8cm]
	(.12,.7) node {\tiny{$w$}}
	(.73,.7) node {\tiny{$z$}}
	(0,1.18) node {\tiny{$x$}}
	(-.23,.7) node {\tiny{$y$}}
	(-.65,1.1) node {\tiny{$a$}}
	(-.95,.7) node {\tiny{$b$}}
	(-.7,.65) node {\tiny{$d$}}
	(-.4,1.45) node {\tiny{$o$}}
	(1.15,1.1) node[right] {\tiny{$\fE_3$}}
	(1.15,.7) node[right] {\tiny{$\fE_2$}}
	(.3,-.4) node {\tiny{$\ud\lt_{(J_2)}~\tn{where}~J_2\!:=\!\{\fE_1\}$}}
	;
	
	\filldraw[yshift=-4cm,xshift=12cm]
	(.2,.5) circle (1.2pt)
	(.6,.5) circle (1.2pt)
	(-.6,.9) circle (1.2pt)
	(.6,.9) circle (1.2pt)
	(-.6,1.3) circle (1.2pt)
	(.6,1.3) circle (1.2pt)
	(-.2,.5) circle (1.2pt)
	(-1,.5) circle (1.2pt)
	(-.6,.5) circle (1.2pt)
	(1,.5) circle (1.2pt)
	;
	
	\draw[yshift=-4cm,xshift=12cm]
	(-1,.9)--(-.6,.9)--(-.6,.5)
	(-.6,1.3)--(.6,1.3)--(1,.5)
	(.6,.5)--(.6,1.3)
	;
	
	\draw[ultra thick,yshift=-4cm,xshift=12cm]
	(-.6,1.3)--(-.6,.9)--(-.2,.5)
	(-1,.5)--(-1,.9)
	(.2,.5)--(.6,.9)
	;
	
	\draw[dotted,yshift=-4cm,xshift=12cm]
	(-1.1,.5)--(1.1,.5)
	(-1.1,.9)--(1.1,.9)
	(-1.1,1.3)--(1.1,1.3)
	;
	
	\draw[decorate,decoration=brace,thin,yshift=-4cm,xshift=12cm]
	(1.1,1.25)--(1.1,.95);
	\draw[decorate,decoration=brace,thin,yshift=-4cm,xshift=12cm]
	(1.1,.85)--(1.1,.55);
	
	\draw[yshift=-4cm,xshift=12cm]
	(.48,1.1) node {\tiny{$w$}}
	(.83,.62) node {\tiny{$z$}}
	(0,1.18) node {\tiny{$x$}}
	(.52,.62)  node {\tiny{$v$}}
	(.23,.7)  node {\tiny{$u$}}
	(-.23,.7) node {\tiny{$y$}}
	(-.46,1.1) node {\tiny{$a$}}
	(-.85,1.02) node {\tiny{$b$}}
	(-1.1,.65) node {\tiny{$c$}}
	(-.7,.65) node {\tiny{$d$}}
	(-.6,1.45) node {\tiny{$o$}}
	(1.15,1.1) node[right] {\tiny{$\fE_3$}}
	(1.15,.7) node[right] {\tiny{$\fE_1$}}
	(.3,-.4) node {\tiny{$\ud\lt_{(J_3)}~\tn{where}~J_3\!:=\!\{\fE_2\}$}}
	;
\end{tikzpicture}	
\end{center}
\caption{Examples of $\ov\lt_{(J)}$ and $\ud\lt_{(J)}$}\label{Fig:t(I)}	
\end{figure}  

For $J\!\subset\!\bbI(\ov\lt)$ and $K\!\subset\!\bbI(\ud\lt)$,
by (\ref{Eqn:bbI}) we have
\begin{align*}
	\bbI(\ov\lt_{(J)})=
	(\ov\bE\bsl J)\sqcup
	\ND(\ov\lt_{(J)})\,,\qquad
	\bbI(\ud\lt_{(K)})=
	(\ud\bE\bsl K)\sqcup
	\ND(\ud\lt_{(K)})\,.
\end{align*}

\begin{lmm}
\label{Lm:I(t_J)}
Let $\ov\lt_{(J)}$ be as in Lemma~\ref{Lm:t_J}, and $\ud\lt_{(K)}$ be as in~(\ref{Eqn:t_I_LRS}).
Then,
\begin{alignat*}{2}
	&\Dm(\ov\lt_{(J)}) \sqcup \ov\sfE_J
	= 
	\Dm(\ov\lt)\sqcup\big(\ND(\ov\lt)\cap J\big)
	&&
	\ND(\ov\lt_{(J)}) 
	= 
	\ND(\ov\lt)\big\bsl\,\big(J\cup\ov \sfE_J^{\,\wedge}\big)\,,
	\\
	&\Dm(\ud\lt_{(K)}) \sqcup \ud\sfE_K
	= 
	\Dm(\ud\lt)\sqcup\big(\ND(\ud\lt)\cap K\big)
	&
	\qquad
	&
	\ND(\ud\lt_{(K)})
	=
	\ND(\ud\lt)\big\bsl\,\big(K\cup\ud \sfE_K^{\;\wedge}\big)\,.
\end{alignat*}	
\end{lmm}

\begin{proof}
On the one hand, (\ref{Eqn:t_J_dom'}) implies
\begin{align*}
	\Dm(\ov\lt_{(J)})\sqcup 
	\ov\sfE_J
	=
	\Dm(\ov\lt)\sqcup \big(J\cap(\fE_\cht\bsl\dot\fE_\cht)\big)\sqcup (J\cap\fE_\cht^\prec)
	=\Dm(\ov\lt)\sqcup\big(\ND(\ov\lt)\cap J\big),
\end{align*}
hence
\begin{align*}
	\tau=\Dm(\ov\lt)\sqcup\ND(\ov\lt)=
	\Dm(\ov\lt_{(J)})\sqcup 
	\ov\sfE_J\sqcup\big(\ND(\ov\lt)\bsl J\big).
\end{align*}
On the other hand, considering the underlying sets of the posets $\tau$ and $\tau\bsl\,\ov\sfE_J^{\,\wedge}$, we have
\begin{align*}
	\tau=\big(\tau\bsl\,\ov\sfE_J^{\,\wedge}\big)\sqcup  \ov\sfE_J^{\,\wedge}
	=
	\Dm(\ov\lt_{(J)})\sqcup 
	\ND(\ov\lt_{(J)})\sqcup 
	\ov\sfE_J^{\,\wedge}.
\end{align*}
Combining the above two expressions of $\tau$, we have
\begin{align*}
	\ND(\ov\lt_{(J)})= 
	\Big(\ov\sfE_J\sqcup\big(\ND(\ov\lt)\bsl J\big)\Big)\Big\bsl \,\ov\sfE_J^{\,\wedge}
	=\ND(\ov\lt)\big \bsl \, (J\cup\ov\sfE_J^{\,\wedge})\,.
\end{align*} 
The proof for $\bbI(\ud\lt_{(K)})$ is parallel, hence is omitted.	
\end{proof}	

\begin{rmk}\label{Rmk:t(I)}
There is a subtle difference between the construction of $\ov\lt_{(J)}$ here and in~\cite[(2.5)]{g1modular}:
in this article,
we take $\tau_{(J)}\eq(\tau_{\ov\lt})_{(\ov\sfE_{J}^\wedge)}$,
whereas in~\cite[(2.6)]{g1modular},
we take $\tau_{(J)}\eq(\tau_{\ov\lt})_{(\ov\sfE_{J})}$.
The difference is caused by the ``cutoff'' level $\cht$ in~\cite{g1modular} that utilizes the additional data of the weights,
which would be unnatural for the more general setting of this article.
Nonetheless,
comparing~\cite[(2.13)]{g1modular}
with Theorem~\ref{Thm:tf_smooth}~\ref{Cond:smooth_tf} of \S\ref{Subsec:STF_main_statement},
we see the twisted fields can only be affected by {\it dominant} edges,
hence the two approaches to $\ov\lt_{(J)}$ yield the same $\fM^{\tf}$.
\end{rmk} 

\subsection{Main theorems of STF}
\label{Subsec:STF_main_statement}
We are ready to present the main statement of the STF theory.

Let
$\fM$ and $\La$ be as in Definition~\ref{Dfn:Treelike_structure}.
To each $\al\inn A$, the treelike structure $\La$ assigns a rooted tree $\tau_{\al}$,
which in turn determines a subset $\ov\fT_\al$ of $\ov\fT$ and a subset $\ud\fT_\al$ of $\ud\fT$, respectively:
\begin{equation}\label{Eqn:Treelike_Level}
 \ov\fT_{\al}=\ov\fT_\al(\La):=
 \big\{\,
 \ov \lt\inn\ov\fT\,:\,
 \tau_{\ov\lt}\eq\tau_{\al}
 \,\big\}\,,\qquad
 \ud\fT_{\al}=\ud\fT_{\al}(\La):=
 \big\{\,
 \ud \lt\inn\ud\fT\,:\,
 \tau_{\ud\lt}\eq\tau_{\al}
 \,\big\}\,.
\end{equation}
These subsets are used for indexing the strata of $\fM^{\tf}$ and $\fM^{\rtf}$ in Theorem~\ref{Thm:tf_smooth} below.

\begin{figure}[htp]
\begin{center}
\begin{tikzpicture}
	\filldraw
	(0,0) circle (1pt)
	(.35,.66) circle (1pt)
	(.7,0) circle (1pt);
	\draw
	(0,0)--(.35,.66)--(.7,0);
	\draw
	(.35,.66) node[above] {\tiny{$o$}}
	(.175,.33) node[left] {\tiny{$b$}}
	(.525,.33) node[right] {\tiny{$c$}}
	(0,0) node[left] {\tiny{$v_b$}}
	(.7,0) node[right] {\tiny{$v_c$}}
	;
	
	\filldraw[xshift=3cm]
	(.117,.22) circle (1pt)
	(.35,.66) circle (1pt)
	(.7,0) circle (1pt);
	\draw[xshift=3cm]
	(.117,.22)--(.35,.66)--(.7,0);
	\draw[xshift=3cm, dotted]
	(-.2,.22)--(.9,.22)
	(-.2,.66)--(.9,.66);
	\draw[xshift=3cm, thin, decorate, decoration=brace]
	(.95,.61)--(.95,.27);
	\draw[xshift=3cm]
	(.233,.44) node[left] {\tiny{$b$}}
	(.467,.44) node[right] {\tiny{$c$}}
	(.95,.44) node[right] {\tiny{$\fE_1$}}
	;
	
	\filldraw[xshift=6cm]
	(0,0) circle (1pt)
	(.35,.66) circle (1pt)
	(.583,.22) circle (1pt);
	\draw[xshift=6cm]
	(0,0)--(.35,.66)--(.583,.22);
	\draw[xshift=6cm, dotted]
	(-.2,.22)--(.9,.22)
	(-.2,.66)--(.9,.66);
	\draw[xshift=6cm, thin, decorate, decoration=brace]
	(.95,.61)--(.95,.27);
	\draw[xshift=6cm]
	(.233,.44) node[left] {\tiny{$b$}}
	(.467,.44) node[right] {\tiny{$c$}}
	(.95,.44) node[right] {\tiny{$\fE_1$}}
	;
	
	\filldraw[xshift=9cm]
	(0,0) circle (1pt)
	(.35,.66) circle (1pt)
	(.7,0) circle (1pt);
	\draw[xshift=9cm]
	(0,0)--(.35,.66)--(.7,0);
	\draw[xshift=9cm, dotted]
	(-.2,0)--(.9,0)
	(-.2,.66)--(.9,.66);
	\draw[xshift=9cm, thin, decorate, decoration=brace]
	(.95,.61)--(.95,.05);
	\draw[xshift=9cm]
	(.175,.33) node[left] {\tiny{$b$}}
	(.525,.33) node[right] {\tiny{$c$}}
	(.95,.33) node[right] {\tiny{$\fE_1$}}
	;
	
	\draw[dashed]
	(2.5,-.5) rectangle (10.8,1);
	
	\draw[yshift=-.5cm]
	(.35,-.25) node {\tiny{$\ga=\tau_{\al}$}}
	(3.35,0.25) node {\tiny{$\ov\lt_{\al,1}$}}
	(6.35,0.25) node {\tiny{$\ov\lt_{\al,2}$}}
	(9.35,0.25) node {\tiny{$\ov\lt_{\al,3}$}}
	(6.65,-.25) node {\tiny{$\ov\fT_{\al}$}};	
\end{tikzpicture}	
\end{center}
\caption{An example of $\ov\fT_{\al}$}\label{Fig:Treelik_Level}	
\end{figure}

\begin{eg}\label{Eg:Treelike_Level}
We continue with the setting of Example~\ref{Eg:genus_1}.
Consider the graph $\ga$ given by the leftmost graph of Figure~\ref{Fig:Treelik_Level}, and
the stratum $\fM_\al$ of $\fM^\wt_1$ comprised of all $(C,\bfw)$ whose dual graph is $\ga$, satisfying $g(C_o)\eq 1$, $\bfw(C_o)\eq 0$, $\bfw(C_{v_b})\eq 2$, and $\bfw(C_{v_c})\eq 3$.
Following Example~\ref{Eg:genus_1}, we see $\tau_{\al}\eq\ga$ with the root $o$.
Then $$\ov\fT_{\al}=\big\{\,\ov\lt_{\al,1}\eq\big(\tau_\al,\{\fE_1\},\{b\}\big)\,,\,\ov\lt_{\al,2}\eq\big(\tau_\al,\{\fE_1\},\{c\}\big)\,,\,\ov\lt_{\al,3}\eq\big(\tau_\al,\{\fE_1\},\{b,c\}\big)\,\big\},$$
where $\fE_1\!:=\!\{b,c\}$, and $\ov\lt_{\al,i}$'s are given in Figure~\ref{Fig:Treelik_Level}.

On a small neighborhood $\cV$ of a point of $\fM_\al$,
the sequential blowup of $\mwt_1$ in~\cite{HL10} locally just blows up along the locus $\fM_\al\!\cap\!\cV$.
The strata $\jia{\fM_\al}_{\ov\lt_{\al,1}}$, $\jia{\fM_\al}_{\ov\lt_{\al,2}}$, and $\jia{\fM_\al}_{\ov\lt_{\al,3}}$ in Theorem~\ref{Thm:tf_smooth}\ref{Cond:smooth_tf} below correspond to the $[1\!:\!0]$-, $[0\!:\!1]$-, and the remaining directions of the exceptional divisor, respectively.
\end{eg}

Recall $\De$ denotes the boundary of $\fM$; c.f.~(\ref{Eqn:boundary}).
For every $\fM_\al\!\subset\!\De$ (i.e.~$\al\inn B$) and $e\inn S_{\al}$,
we denote by
\begin{equation}\label{Eqn:L_e}
 L_e\lra\fM_\al
\end{equation}
the line bundle such that on each chart $\cV\inn\fV_\al$,
\begin{align*}
	L_e|_{\fM_\al\cap\cV}=
	\sO_\cV(X_e^\cV)|_{\fM_\al\cap\cV},\qquad
	\tn{where}\quad
	X_e^\cV=(\ze_e^\cV)\subset\cV.
\end{align*}
The existence of $L_e$ is assured by 	Convention~\ref{Convention:loc_Euc_id}.
For every $\fe\inn \tau_{\al}$, every pair $\fe\!\preceq\!\fe'\inn \tau_{\al}$,
and every pair $\fe\!\prec\!\fe'\inn \tau_{\al}$,
let
\begin{equation}\label{Eqn:Twisted_line_bundle}
 L_{[\fe,o)}
 :=\bigotimes_{\fe''\succeq\fe}
 L_{\be_{\al}(\fe'')}\big/\fM_\al\,,\quad
 L_{[\fe,\fe']}
 :=\bigotimes_{\fe\preceq\fe''\preceq\fe'}\!\!\!\!\!
 L_{\be_{\al}(\fe'')}\big/\fM_\al\,,
 \quad
 L_{[\fe,\fe')}
 :=\bigotimes_{\fe\preceq\fe''\prec\fe'}\!\!\!\!\!
 L_{\be_{\al}(\fe'')}\big/\fM_\al\,,
\end{equation}
respectively,
where $\preceq$ is the tree order on $\tau_{\al}$.

Given a finite set $I$ and a direct sum of line bundles $V\eq\oplus_{i\in I}L_i'$ over an arbitrary base,
let $V_{(i)}$ be the sub-bundle of $V$ given by $V_{(i)}\eq\oplus_{j\in I\bsl\{i\}}L_j'$.
Then, we write
$$
 \mathring\P(V)=\P(V)\big\bsl \big(\,\bigcup_{i\in I}\P(V_{(i)})\big).
$$
In particular, when $V\eq L$ is a line bundle, then $\mathring\P(V)\eq \P(V)$.
Given $k\inn\Z_{>0}$ and morphisms $M_i\!\lra\! S$ with $1\!\le\!i\!\le\!k$,
we write
$$
 \Big\lgroup\!\!\prod_{1\le i\le k}\!\!\Big\rgroup_{\!\!S}\,
 M_i \,
 :=
 M_1\times_S M_2\times_S\cdots
 \times_S M_k.
$$

\begin{thm}\label{Thm:tf_smooth}
Let $\fM$ be a smooth algebraic stack  endowed with an LES and  a treelike structure $\La$ as in Definition~\ref{Dfn:Treelike_structure}.
Then,
there exists a smooth algebraic stack $$\fM^{\tf}=\fM^{\tf}_{\La},$$
known as the \ts{stack with RL twisted fields over $\fM$ with respect to $\La$},
admitting an induced LES described below.

\begin{enumerate}
[leftmargin=*,label=$(\mathsf{\ov{p}_\arabic*})$]

\item 
(\emph{strata of STF}).
\label{Cond:smooth_tf}
$\fM^{\tf}$ has an LES  given by
\begin{align*}
 &\fM^{\tf}
 =\bigsqcup_{(\al,
 \ov\lt)\in \ov\La}
 \hspace{-.075in}
 \fM^{\tf}_{(\al,\ov\lt)}\,,\qquad
 \tn{where}\\
 & \ov\La
 :=\big\{\,(\al,\ov\lt):\;
 \al\inn A,\ \ov\lt\inn\ov\bT_\al(\La)\,\big\}\,,\qquad
 \fM^{\tf}_{(\al,\ov\lt)}
 :=
 \Big\lgroup\!\!
 \prod_{\fE\in \ov \bE_{\ov\lt}}\!\!
 \Big\rgroup_{\!\!\fM_\al}
 \mathring{\P}\Big(\!
 \bigoplus_{\;
 	\begin{subarray}{c}
 		\fe\in\dot \fE
 	\end{subarray}
 }
 L_{[\fe,o)}
 \Big) \,.
\end{align*}
Moreover,
the projections $\fM^{\tf}_{(\al,\ov\lt)}\!\lra\!\fM_\al$ altogether determine a proper and birational morphism $$\ov\varpi:\fM^{\tf}\lra\fM$$ known as the \ts{forgetful morphism},
satisfying $$ \ov\varpi|_{(\fM^{\tf})^{\mn}}:\, \big(\fM^{\tf}\big)^\mn=\fM^{\tf}_{(\ms,\ov\lt_\bullet)}\,
\lra\,\fM^\mn
$$ is an isomorphism.

\item \label{Cond:smooth_parameters}
(\emph{modular charts and parameters}).
For every 
$$
 \Big(\,\al\,,\,\ov\lt\eq\big(\tau_\al,\ov\bE,\Dm(\ov\lt)\big)\,\Big)\ \in\  
 \ov\La\,
 \big\bsl \big\{(\ms,\ov\lt_\bullet)\big\}\,,
$$
let
\begin{align*}
	 S_{(\al,\ov\lt)}:=
	\ov\bE\sqcup \be_\al\big(\ND(\ov\lt)\big)\sqcup\big(S_{\al}\big\bsl\,\be_\al(\tau_\al)\big)
	=\ov\bE\sqcup \Big(S_\al\!\big\bsl\,\be_\al\big(\Dm(\ov\lt)\big)\Big) .
\end{align*}
Then, for every $(\al,\ov\lt)\inn\ov\La\,
\big\bsl \big\{(\ms,\ov\lt_\bullet)\big\}$
and  $\cV\inn\fV_\al$,
there exist a smooth chart called a \ts{twisted chart}: 
$$\cU=\cU_{\cV,\ov{\lt}} \lra \fM^{\tf}
\qquad\tn{satisfying}\qquad
\ov\varpi(\cU)\subset\cV,
$$  
and a subset of a system of local parameters called the 
\ts{twisted parameters}:
\begin{align*}
 \xi_s^{\cU}\,,\quad s\in
 S_{(\al,\ov\lt)}\,,
\end{align*}
satisfying the following.
\begin{itemize}[leftmargin=*]
\item 
For every $\fe\inn\tau_{\al}$, we have
\begin{align}\label{Eqn:mod_par_pullback}
\ov\varpi^*\big(\ze_{\be_{\al}(\fe)}^\cV\,\big)
= u_{\fe}\cdot\!\prod_{
	\begin{subarray}{c}
		\fe\in\fE\in\ov\bE
\end{subarray}}\!\!\xi_\fE^{\cU},
\qquad\tn{where}\quad
\begin{cases}
	u_{\fe}\in\Ga(\sO^*_{\cU})
	&
	\tn{if}~\fe\in\Dm(\ov\lt),\\
	u_\fe/\xi_{\be_\al(\fe)}^{\cU}\in\Ga(\sO^*_{\cU})
	&
	\tn{if}~\fe\in \ND(\ov\lt),
\end{cases}
\end{align}
and for every $e\inn S_\al\bsl\be_\al(\tau_\al)$,
we have
\begin{align}\label{Eqn:mod_par_pullback'}
 \ov\varpi^*\big(\ze_{e}^\cV\big)
 \,\big/\,\xi_e^\cU\,\in\,
 \Ga(\sO^*_{\cU}).
\end{align}

\item 
The twisted charts in
$$
\fV_{(\al,\ov\lt)}:=\big\{\,\cU_{\cV,\ov\lt}:\,\cV\inn\fV_\al\,\big\},\qquad(\al,\ov\lt)\in \ov\La\,
\big\bsl \big\{(\ms,\ov\lt_\bullet)\big\},
$$
and the local parameters in
$$
\big\{\xi^\cU_s\big\}_{s\in S_{(\al,\ov\lt)}},\qquad
\cU\inn\fV_{(\al,\ov\lt)},\ \ 
(\al,\ov\lt)\in \ov\La\,
\big\bsl \big\{(\ms,\ov\lt_\bullet)\big\},
$$
can respectively be used as modular charts and modular parameters of the LES of $\fM^{\tf}$ in \ref{Cond:smooth_tf}, in the sense of Definition~\ref{Dfn:G-adim_fixture}.

\item
For every system of local parameters on $\cV$:
$$
\{\ze_e^\cV\}_{e\in S_\al}\sqcup
\{v_r\}_{r\in I'}
$$ 
containing $\{\ze_{e}^\cV\}_{e\in S_\al}$ as a subset,
every subset $\{e_\fE\}_{\fE\in\ov\bE}\!\subset\!\Dm(\ov\lt)$ satisfying $e_\fE\inn\fE$ for all $\fE\inn\ov\bE$,
and every point $y\inn\fM^{\tf}_{(\al,\ov\lt)}\!\cap\!\cU$,
the following set
\begin{align*}
&	\{\xi_s^{\cU}\}_{s\in S_{(\al,\ov\lt)}}\sqcup
	\big\{
	g^{\cU}_e-\big(g^{\cU}_e(y)\big)
	\big\}_{e\in\Dm(\ov\lt)\bsl\{e_\fE:\fE\in\ov\bE\}}
	\sqcup
	\{\varpi^*v_r\}_{r\in I'},\\
	&
	\tn{where}\quad
	g^{\cU}_e:=
	\frac{u_e}{\prod_{e\in\fE\in\ov\bE}u_{e_\fE}}\,,\qquad
	e\inn\Dm(\ov\lt)\bsl\{e_\fE:\fE\inn\ov\bE\},
\end{align*}
forms a system of local parameters on $\cU$.
Here, $u_e$, $e\inn\Dm(\ov\lt)$ are as in (\ref{Eqn:mod_par_pullback}).
\end{itemize}

\item\label{Cond:smooth_Z}
(\emph{nearby strata}).
Let $(\al,\ov\lt)$, $\cV$ and $\cU$ be as in \ref{Cond:smooth_parameters}.
For every $S'\!\subset\!S_{(\al,\ov\lt)}$,
we set
\begin{alignat*}{2}
J&:=(S'\!\cap\!\ov\bE)\sqcup\be_\al^{-1}\Big(S'\!\cap\!\be_\al\big(\ND(\ov\lt)\big)\Big)
&\qquad&\big(\subset\bbI(\ov\lt)= \ov\bE\!\sqcup\!\ND(\ov\lt)\,\big)\,,\\
E&:=\be_\al(\ov\sfE_{J})\sqcup \Big(S'\cap\big(S_\al\big\bsl\be_\al(\tau_\al)\big)\Big)&\qquad&
\big(\subset S_\al\big)\,,
\end{alignat*}
where $\bbI(\ov\lt)$ is as in (\ref{Eqn:bbI}), and
$\ov\sfE_{J}$ $(\subset\!\tau_\al)$ is as in (\ref{Eqn:C_J}). 
Let $\ov\lt_{(J)}$ be as in Lemma~\ref{Lm:t_J}, $\al_{(E)}$ be as in Proposition~\ref{Prp:M_strata_local}, $$\phi_{\al;E}:\,\tau_{\al_{(E)}}\lra\tau_\al\big\bsl\,\big(\be_\al^{-1}(E)\big)^{\!\wedge}= \tau_\al\bsl\,\ov\sfE_J^{\,\wedge}=\tau_{\ov\lt_{(J)}}$$ 
be as in Proposition~\ref{Prp:Treelike_contraction},
and $\phi_{\al;E}^{-1}\big(\ov\lt_{(J)}\big)\inn\ov\bT_{\al_{(E)}}$ be as in Lemma~\ref{Lm:RLS_LRS_isom}.

Then, 
we have
$$
 (\al,\ov\lt)_{(S')}=\big(\,\al_{(E)}\,,\,\phi_{\al;E}^{-1}(\ov\lt_{(J)})\,\big),\qquad\tn{i.e.}\qquad
 Z_{(S')}^{\cU}\subset \fM^{\tf}_{\big(\al_{(E)},\ \phi_{\al;E}^{-1}(\ov\lt_{(J)})\big)}\!\cap\cU\,.
$$ 
Moreover,
the injection $\iota_{(\al,\ov\lt);\,S'}\!:S_{(\al,\ov\lt)_{(S')}}\!\hookrightarrow\!S_{(\al,\ov\lt)}$ as in Corollary~\ref{Crl:M_strata_local} satisfies
\begin{alignat*}{2}
&\iota_{(\al,\ov\lt);\,S'}(\fE)=\phi_{\al;E}(\fE)&\qquad&\forall\ \fE\in\phi^{-1}_{\al;E}(\ov\bE\bsl J)\,,\\
&\iota_{(\al,\ov\lt);\,S'}(e)=\iota_{\al;E}(e)&\qquad&
\forall\ e\in S_{\al_{(E)}}\!\!\Big\bsl\; \be_{\al_{(E)}}\!\Big(\phi^{-1}_{\al;E}\big(\Dm(\ov\lt_{(J)})\big)\Big)\,.
\end{alignat*}

\item \label{Cond:smooth_local_blowup}
(\emph{RL local blowups}).
For every $\al\inn B$ and $\cV\inn\fV_\al$,
let  $\ti\cV\!\lra\!\cV$ be a sequence of blowups of $\cV$ successively along the proper transforms of 
\begin{align*}
	X^\cV_{\be_\al(\fE)}:=
	\big\{\,\ze_{\be_\al(e)}^\cV\!\eq 0\,:\,
	e\inn\fE\,\big\},\qquad
	\fE\in\Xi(\tau_\al)\,,
\end{align*}
with respect to any linear order extending the partial order (\ref{Eqn:transverse_sections_order}) on $\Xi(\tau_{\al})$,
starting from $\fE\eq\max(\tau_\al)$.
Then, $\ti\cV/\cV$ is independent of the choice of the linear order extending (\ref{Eqn:transverse_sections_order}),
and is isomorphic to $\big(\cV\!\times_{\fM}\!\fM^{\tf}\big)\big/\cV$ 
(i.e.~$\cV^{\tf}\big/\cV$).
\end{enumerate}
\end{thm}

Theorem~\ref{Thm:tf_smooth_revert} below is
the leaf-to-root counterpart of Theorem~\ref{Thm:tf_smooth}.
For every $\ud\lt\inn\ud\bT$, let $\Om(\fe)$, $\fe\inn\Dm(\ud\lt),$ be as in (\ref{Eqn:Omega_e}), and $\fp_{\fe}(\fe')$, $\fe'\inn\Om(\fe)$, be as in (\ref{Eqn:fp_e}).

\begin{thm}\label{Thm:tf_smooth_revert}
Let $\fM$ be a smooth algebraic stack  endowed with an LES and  a treelike structure $\La$ as in Definition~\ref{Dfn:Treelike_structure}.
Then,
there exists a smooth algebraic stack $$\fM^{\rtf}=\fM^{\rtf}_{\La},$$
known as the \ts{stack with LR twisted fields over $\fM$ with respect to $\La$},
enjoying the following properties.
	
\begin{enumerate}
	[leftmargin=*,label=$(\mathsf{\ud{p}_\arabic*})$]
	
	\item  \label{Cond:smooth_tf_revert}
	(\emph{strata of STF}).
	$\fM^{\rtf}$ has an LES  given by 
	\begin{align*}
	\fM^{\rtf}
	=\bigsqcup_{(\al,
		\ud\lt)\in \ud\La}
	\hspace{-.075in}
	\fM^{\rtf}_{(\al,\ud\lt)}\,,\qquad&
	\tn{where}\qquad
	\ud\La:=\big\{\,(\al,\ud\lt):\;
	\al\inn A,\ \ud\lt\inn\ud\bT_\al\,\big\}\,,
	\\
	\fM^{\rtf}_{(\al,\ud\lt)}
	:=
	\Big\{\,
	&\big(\,[v^\fE_{\fe,\fe'}]_{\fe\in\udt\fE,\fe'\in\Om(\fe)}\big)_{\fE\in\ud\bE_{\ud\lt}}
	\in
	\Big\lgroup\!\!
	\prod_{\fE\in \ud \bE_{\ud\lt}}\!\!
	\Big\rgroup_{\!\!\fM_\al}
	\mathring{\P}\Big(\!
	\bigoplus_{\;
		\begin{subarray}{c}
			\fe\in\udt \fE
		\end{subarray}
	}\;
	\bigoplus_{
		\begin{subarray}{c}
			\fe'\in\Om(\fe)
		\end{subarray}
	}\!\!
	L_{[\fe',\fe]}
	\Big):\\
	&
	v^{\fE}_{\fe,\fe'}\!\otimes
	v^{\fE_-}_{p_{\fe}(\fe''),\fe''}=
	v^{\fE_-}_{p_{\fe}(\fe'),\fe'}\!\otimes
	v^{\fE}_{\fe,\fe''}\ \ \big(
	\in 
	L_{[\fe',\fe)}\!\otimes_{\,\fM_\al}\!L_{[\fe'',\fe)}\!\otimes_{\,\fM_\al}\!L_{\fe}\,\big)
	\\
	&\forall~
	\fE\in\ud\bE_{\ud\lt}\bsl\{\fE_1\eq\min(\tau_{\ud\lt})\},~
	\fe\inn \udt\fE,~
	\fe',\fe''\inn \Om(\fe)\,
	\Big\}\,.
	\end{align*}
	Moreover,
	the projections $\jia{\fM_\al}_{\ud\lt}\!\lra\!\fM_\al$ altogether determine a proper and birational morphism $$\ud\varpi:\fM^{\rtf}\lra\fM$$ known as the \ts{forgetful morphism},
	satisfying $$ \ud\varpi|_{(\fM^{\rtf})^{\mn}}:\, \big(\fM^{\rtf}\big)^\mn=\jia{\fM^\mn}_{\ud\lt_\bullet}\,
	\lra\,\fM^\mn
	$$ is an isomorphism.
	
	\item \label{Cond:smooth_parameters_revert}
	(\emph{modular charts and parameters}).
	For every 
	$$
	\Big(\,\al\,,\;\ud\lt\eq\big(\tau_\al,\ud\bE,\Dm(\ud\lt)\big)\,\Big)\ \in\  
	\ud\La
	\big\bsl \big\{(\ms,\ud\lt_\bullet)\big\}\,,
	$$
	let 
	$$
	S_{(\al,\ud\lt)}:=
	\ud\bE\sqcup \be_\al\big(\ND(\ud\lt)\big)
	\sqcup\big(S_{\al}\big\bsl\,\be_\al(\tau_\al)\big)
	=\ud\bE\sqcup \Big(S_\al\!\big\bsl\,\be_\al\big(\Dm(\ud\lt)\big)\Big).
	$$
	Then, for every $(\al,\ud\lt)\in\ud\La
	\big\bsl \big\{(\ms,\ud\lt_\bullet)\big\}$
	and  $\cV\inn\fV_\al$,
	there exist a smooth chart called a \ts{twisted chart}: 
	$$\cU=\cU_{\cV,\ud{\lt}} \lra \fM^{\rtf}
	\qquad\tn{satisfying}\qquad
	\ud\varpi(\cU)\subset\cV,
	$$  
	and a subset of a system of local parameters called the 
	\ts{twisted parameters}:
	\begin{align*}
		\xi_s^{\cU}\,,\quad s\in
		S_{(\al,\ud\lt)}\,,
	\end{align*}
	satisfying the following.
\begin{itemize}[leftmargin=*]
	\item For every $\fe\inn\tau_{\al}$, we have
	\begin{align}\label{Eqn:mod_par_pullback_revert}
		\ud\varpi^*\big(\ze_{\be_{\al}(\fe)}^\cV\,\big)
		= u_{\fe}\cdot\!\prod_{
			\begin{subarray}{c}
				\fE\in\ud\bE~\tn{s.t.}~\fE\ni\fe
		\end{subarray}}\!\!\!\!\!\xi_\fE^{\cU},
		\qquad\tn{where}\quad
		\begin{cases}
			u_{\fe}\in\Ga(\sO^*_{\cU})
			&
			\tn{if}~\fe\in\Dm(\ud\lt),\\
			u_\fe/\xi_{\be_\al(\fe)}^{\cU}\in\Ga(\sO^*_{\cU})
			&
			\tn{if}~\fe\in \ND(\ud\lt),
		\end{cases}
	\end{align}
	and for every $e\inn S_\al\bsl\be_\al(\tau_\al)$,
	we have
	\begin{align}\label{Eqn:mod_par_pullback'_revert}
		\ud\varpi^*\big(\ze_{e}^\cV\big)
		\,\big/\,\xi_e^\cU\,\in\,
		\Ga(\sO^*_{\cU}).
	\end{align}
	
	\item The twisted charts in
	$$
	\fV_{(\al,\ud\lt)}:=\big\{\,\cU_{\cV,\ud\lt}:\,\cV\inn\fV_\al\,\big\},\qquad(\al,\ud\lt)\in \ud\La
	\big\bsl \big\{(\ms,\ud\lt_\bullet)\big\},
	$$
	and the local parameters in
	$$
	\big\{\xi^\cU_s\big\}_{s\in S_{(\al,\ud\lt)}},\qquad
	\cU\inn\fV_{(\al,\ud\lt)},\ \ 
	(\al,\ud\lt)\in \ud\La
	\big\bsl \big\{(\ms,\ud\lt_\bullet)\big\},
	$$
	can respectively be used as modular charts and modular parameters of the LES of $\fM^{\rtf}$ in \ref{Cond:smooth_tf_revert}, in the sense of Definition~\ref{Dfn:G-adim_fixture}.
	
	\item For every system of local parameters on $\cV$:
	$$
	\{\ze_e^\cV\}_{e\in S_\al}\sqcup
	\{v_r\}_{r\in I'}
	$$ 
	containing $\{\ze_{e}^\cV\}$ as a subset,
	every subset $\{e_\fE\}_{\fE\in\ud\bE}\!\subset\!\Dm(\ud\lt)$ satisfying $e_\fE\inn\fE$ for all $\fE\inn\ud\bE$,
	and every point $y\inn\fM^{\rtf}_{(\al,\ud\lt)}\!\cap\!\cU$,
	the following set
	\begin{align*}
		&	\{\xi_s^{\cU}\}_{s\in S_{(\al,\ud\lt)}}\sqcup
		\big\{
		g^{\cU}_e-\big(g^{\cU}_e(y)\big)
		\big\}_{e\in\Dm(\ud\lt)\bsl\{e_\fE:\fE\in\ud\bE\}}
		\sqcup
		\{\varpi^*v_r\}_{r\in I'},\\
		&
		\tn{where}\quad
		g^{\cU}_e:=
		\frac{u_e}{\prod_{e\in\fE\in\ud\bE}u_{e_\fE}}\,,\qquad
		e\inn\Dm(\ud\lt)\bsl\{e_\fE:\fE\inn\ud\bE\},
	\end{align*}
	forms a system of local parameters on $\cU$.
	Here, $u_e$, $e\inn\Dm(\ud\lt)$ are as in (\ref{Eqn:mod_par_pullback_revert}).
\end{itemize}	
	
	\item\label{Cond:smooth_Z_revert}
	(\emph{nearby strata}).
	Let $(\al,\ud\lt)$, $\cV$ and $\cU$ be as in \ref{Cond:smooth_parameters_revert}.
	For every $S'\!\subset\!S_{(\al,\ud\lt)}$,
	we set
	\begin{alignat*}{2}
		J&:=(S'\!\cap\!\ud\bE)\sqcup\be_\al^{-1}\Big(S'\!\cap\!\be_\al\big(\ND(\ud\lt)\big)\Big)
		&\qquad&\big(\subset\bbI(\ud\lt)= \ud\bE\!\sqcup\!\ND(\ud\lt)\,\big)\,,\\
		E&:=\be_\al(\ud\sfE_{J})\sqcup \Big(S'\cap\big(S_\al\big\bsl\be_\al(\tau_\al)\big)\Big)&\qquad&
		\big(\subset S_\al\big)\,,
	\end{alignat*}
	where $\bbI(\ud\lt)$ is as in (\ref{Eqn:bbI}), and
	$\ud\sfE_{J}$ $(\subset\!\tau_\al)$ is as in (\ref{Eqn:C_J_revert}). 
	Let $\ud\lt_{(J)}$ be as in (\ref{Eqn:t_I_LRS}), $\al_{(E)}$ be as in Proposition~\ref{Prp:M_strata_local}, $$\phi_{\al;E}:\,\tau_{\al_{(E)}}\lra\tau_\al\big\bsl\,\big(\be_\al^{-1}(E)\big)^{\!\wedge}= \tau_\al\bsl\,\ud\sfE_J^{\;\wedge}=\tau_{\ud\lt_{(J)}}$$ be as in Proposition~\ref{Prp:Treelike_contraction},
	and $\phi_{\al;E}^{-1}\big(\ud\lt_{(J)}\big)\inn\ud\bT_{\al_{(E)}}$ be as in Lemma~\ref{Lm:RLS_LRS_isom}.
	
	Then, 
	we have
	$$
	(\al,\ud\lt)_{(S')}=\big(\,\al_{(E)}\,,\,\phi_{\al;E}^{-1}(\ud\lt_{(J)})\,\big),\qquad\tn{i.e.}\qquad
	Z_{(S')}^{\cU}\subset \fM^{\rtf}_{\big(\al_{(E)},\ \phi_{\al;E}^{-1}(\ud\lt_{(J)})\big)}\!\cap\cU\,.
	$$ 
	Moreover,
	the injection $\iota_{(\al,\ud\lt);\,S'}\!:S_{(\al,\ud\lt)_{(S')}}\!\hookrightarrow\!S_{(\al,\ud\lt)}$ as in Corollary~\ref{Crl:M_strata_local} satisfies
	\begin{alignat*}{2}
		&\iota_{(\al,\ud\lt);\,S'}(\fE)=\phi_{\al;E}(\fE)&\qquad&\forall\ \fE\in\phi^{-1}_{\al;E}(\ud\bE\bsl J)\,,\\
		&\iota_{(\al,\ud\lt);\,S'}(e)=\iota_{\al;E}(e)&\qquad&
		\forall\ e\in S_{\al_{(E)}}\!\!\Big\bsl\, \be_{\al_{(E)}}\!\Big(\phi^{-1}_{\al;E}\big(\Dm(\ud\lt_{(J)})\big)\Big)\,.
	\end{alignat*}
		
	\item \label{Cond:smooth_local_blowup_revert}
	(\emph{LR local blowups}).
	For every $\al\inn B$ and $\cV\inn\fV_\al$,
	let $\undertilde\cV\!\lra\!\cV$ be a sequence of blowups of $\cV$ successively along the proper transforms of 
	\begin{align*}
		X^\cV_{\be_\al(\fE)}:=
		\big\{\,\ze_{\be_\al(e)}^\cV\!\eq 0\,:\,
		e\inn\fE\,\big\},\qquad
		\fE\in\Xi(\tau_\al)\,,
	\end{align*}
	with respect to any linear order extending the reverse order of the partial order (\ref{Eqn:transverse_sections_order}) on $\Xi(\tau_{\al})$,
	starting from $\fE\eq\min(\tau_\al)$.
	Then, $\undertilde\cV/\cV$ is independent of the choice of the linear order extending the reverse order of (\ref{Eqn:transverse_sections_order}),
	and is isomorphic to $\big(\cV\!\times_{\fM}\!\fM^{\rtf}\big)\big/\cV$ 
		(i.e.~$\cV^{\rtf}\big/\cV$).
\end{enumerate}	
\end{thm}

In Theorem~\ref{Thm:tf_smooth}~\ref{Cond:smooth_parameters},
since
the twisted parameters $\xi^\cU_s$, $s\inn S_{(\al,\ov\lt)}$ are modular parameters on $\cU$, by Definition~\ref{Dfn:G-adim_fixture}, (\ref{Eqn:mod_par_pullback}) and (\ref{Eqn:mod_par_pullback'}), we have
\begin{align*}
 &\jia{\fM_{\al}}_{\ov\lt}\cap\cU=
 \big\{\,
 \xi_s^{\cU}\!\eq 0\ \ 
 \forall\,s\inn S_{(\al,\ov\lt)}\,\big\}\\
 = \,&
 \big\{\,
 \xi_{\fE}^{\cU}\!\eq 0\ 
 \forall\,\fE\inn\ov\bE\ ;\ \ 
 \xi_{\be_\al(\fe)}^{\cU}\!\eq 0\ 
 \forall\,\fe\inn\fE_\cht\bsl\dot\fE_\cht\ ;\ \ 
 \ov\varpi^*\ze_e^\cV\eq 0\ 
 \forall\,e\inn\be_\al(\fE_\cht^\prec)\!\sqcup\! \big(S_\al\bsl\be_\al(\tau_\al)\big)\,\big\}.
\end{align*}
Similarly, in
Theorem~\ref{Thm:tf_smooth_revert}~\ref{Cond:smooth_parameters_revert},
we have
\begin{align*}
&\jia{\fM_{\al}}_{\ud\lt}\cap\cU=
\big\{\,
\xi_s^{\cU}\!\eq 0\ \ 
\forall\,s\inn S_{(\al,\ud\lt)}\,\big\}\\
= \,&
\big\{\,
\xi_{\fE}^{\cU}\!\eq 0\ 
\forall\,\fE\inn\ud\bE\ ;\ \ 
\xi_{\be_\al(\fe)}^{\cU}\!\eq 0\ 
\forall\,\fe\inn\bigsqcup_{\fE\in\ud\bE}\!\big((\fE\bsl\fE_+)\big\bsl\udt\fE\big)\ ;\ \ 
\ud\varpi^*\ze_e^\cV\eq 0\ 
\forall\,e\inn\be_\al(\fE_\cht^\succ)\!\sqcup\! \big(S_\al\bsl\be_\al(\tau_\al)\big)\,\big\}.
\end{align*}
These are the analogues of (\ref{Eqn:M_strata_local'}) for $\fM^{\tf}$ and $\fM^{\rtf}$, respectively.

\begin{proof}[Proof of Theorem~\ref{Thm:tf_smooth}]

The statements \ref{Cond:smooth_tf}-\ref{Cond:smooth_Z} follow from the same argument as in~\cite[\S3.1~\&~3.2]{g1modular}.
The key reason is as follows.

Given $(\al,\ov\lt)\inn\ov\La$ and $J\!\subset\!\bbI(\ov\lt)$,
let $E\!\subset\!\be_\al(\tau_\al)$ be determined by $S':=(J\!\cap\!\ov\bE)\sqcup\be_\al\big(J\!\cap\!\ND(\ov\lt)\big)$ as in~\ref{Cond:smooth_Z}.
Given $\cV\inn\fV_\al$,
let $\{\ze_e^\cV\}_{e\in S_\al}\sqcup
\{v_r\}_{r\in I'}$ be an extension of the set of the modular parameters into a system of local parameters on $\cV$.

We fix a subset $\{e_\fE:\fE\in\ov\bE\}\!\subset\!\Dm(\ov\lt)$ satisfying $e_\fE\inn\fE$ for all $\fE\inn\ov\bE$,
and mimic the paragraphs containing~\cite[(3.7)-(3.10)]{g1modular} and define the open subset and the locus 
$$
 \cU\subset\mathbb A^{S_{(\al,\ov\lt)}\sqcup (\Dm(\ov\lt)\bsl \{e_\fE:\fE\in\ov\bE\})\sqcup I'}\qquad
 \tn{and}\qquad
 \cU_{(J)}
 :=
 \bigsqcup_{E_1	\subset S_\al\bsl\be_\al(\tau_\al)}\hspace{-.2in}
 Z_{(S'\sqcup E_1)}^{\cU}\quad
 \subset \cU
$$ 
as in the displays above~\cite[(3.8) \& (3.10)]{g1modular}, respectively (where $\cU$ and $\cU_{(J)}$ are respectively denoted by $\fU$ and $\fU_{x;[\lt_{(\bbI')}]}$).
In addition, we take the following locus of $\cV$
\begin{equation*}\begin{split}
 \cV_{(J)}
 :=
 \big\{\,
 \ze_{e}^\cV\eq 0\ \ \forall\ e\inn\be_\al(\tau_{\al})\bsl E
 \,;\ \,\ze_{e}^\cV\!\ne\! 0\ \ \forall\,e\inn E\,
 \big\} =
 \bigsqcup_{E_1\subset S_\al\bsl\be_\al(\tau_\al)}\hspace{-.22in}
 Z_{(E\sqcup E_1)}^\cV
 \lra
 \bigsqcup_{E_1\subset S_\al\bsl\be_\al(\tau_\al)}\hspace{-.22in}
 \fM_{\al_{(E\sqcup E_1)}}\,,
\end{split}\end{equation*}
which is the analogue of~\cite[(3.2)]{g1modular}.
We thus have the morphism
\begin{align}\label{Eqn:Phi_I'}
 \Phi_{\cU;(J)}:\ \cU_{(J)}\lra
 \cV_{(J)}\times_{{}_{\fM}}
 \Big(\bigsqcup_{E_1\subset S_\al\bsl\be_\al(\tau_\al)}\hspace{-.2in}
 \fM_{\big(\al_{(E\sqcup E_1)},\,\phi^{-1}_{\al;E\sqcup E_1}(\ov\lt_{(J)})\big)}^{\tf}\Big)
\end{align}
as in the paragraph containing~\cite[(3.14)]{g1modular}.
By~\cite[Lemma~3.3]{g1modular},
$\Phi_{\cU;(J)}$ is an isomorphism to an open subset of the right hand side of (\ref{Eqn:Phi_I'}).
The equality~\cite[(3.14)]{g1modular} and the last equality in~\cite[(3.12)]{g1modular}
together imply
\begin{align}\label{Eqn:Phi_I'_strata}
 \Phi_{\cU;(J)}\big(Z_{(S'\sqcup E_1)}^{\cU}\big)
 \ \subset\ 
 Z_{(E\sqcup E_1)}^\cV
 \times_{{}_{\fM}}
 \fM_{\big(\al_{(E\sqcup E_1)},\,\phi^{-1}_{\al;E\sqcup E_1}(\ov\lt_{(J)})\big)}^{\tf}
\end{align}
for all $E_1\!\subset\!S_\al\bsl\be_\al(\tau_\al)$.

As shown in~\cite[(3.16)]{g1modular}, 
the morphisms  $\Phi_{\cU;(J)}$,
$J\!\subset\!\bbI(\ov\lt)$, together give rise to a map
\begin{align*}
	\Phi_\cU:\cU\lra\fM^{\tf}\,.
\end{align*}
As proved in~\cite[Corollary~3.8]{g1modular},
these $\Phi_\cU$ form a smooth atlas of the algebraic stack $\fM^{\tf}$.

The former and latter statements of \ref{Cond:smooth_Z} follow from (\ref{Eqn:Phi_I'_strata}) and from the two displays above \cite[(3.32)]{g1modular},
respectively.

The twisted parameters of~\ref{Cond:smooth_parameters}, along with the regular functions $g_e^\cU$, $e\inn\Dm(\ov\lt\bsl\{e_\fE:\fE\inn\ov\bE\})$, are written as $\ve_i$, $u_e$, and $z_e$ in~\cite[(3.9)]{g1modular}.
The equalities (\ref{Eqn:mod_par_pullback}) and (\ref{Eqn:mod_par_pullback'}) follow from~\cite[(3.12)]{g1modular}.

Using the twisted parameters in \ref{Cond:smooth_parameters} and the injections in \ref{Cond:smooth_Z},
we conclude from (\ref{Eqn:Phi_I'_strata}) that the disjoint union of the strata of $\fM^{\tf}$ in \ref{Cond:smooth_tf} is an LES.
To establish the birationality of $\varpi$, recall the trivial rooted tree is assigned to $\fM^\mn$ (c.f.~Remark~\ref{Rmk:Treelike}),
hence $\varpi$ restricts to the identity map Since $\De$ is closed in $\fM$, so is $\varpi^{-1}(\De)$ in $\fM^{\tf}$. Therefore, $\varpi$ is birational.
The properness of $\varpi$ follows directly from~\ref{Cond:smooth_local_blowup}.

It remains to prove~\ref{Cond:smooth_local_blowup},
which follows immediately from Proposition~\ref{Prp:RL/LR-cptb}; see \S\ref{Subsec:local_blowups} below.
\end{proof}

\begin{proof}[Proof of Theorem~\ref{Thm:tf_smooth_revert}]
	The proof is parallel to that of Theorem~\ref{Thm:tf_smooth},
	hence is omitted.
\end{proof}

The following proposition is a restatement of~\cite[Proposition~3.12]{g1modular} under the current setup.
Let $\fM$ be as in Theorems~\ref{Thm:tf_smooth} and~\ref{Thm:tf_smooth_revert}.
Assume that $\pi\!:\cC\!\lra\!\fM$ is the universal family of $\fM$.

\begin{prp}
\label{Prp:moduli}
With notation as above, the induced stratification $$\cC\eq\bigsqcup_{\al\in A}\cC_\al,\qquad
\tn{where}\quad\cC_\al:=\pi^{-1}(\fM_\al),
$$ is an LES of $\cC$.

Moreover, with notation as in Theorem~\ref{Thm:tf_smooth},
the union 
\begin{align*}
 \cC^{\tf}
 =\bigsqcup_{(\al,
 \ov\lt)\in A^{\ov\der}}
 \hspace{-.15in}
 \jia{\cC_{\al}}_{\ov\lt}\ ,\qquad
 \tn{where}
 \qquad
 \jia{\cC_{\al}}_{\ov\lt}:=
 \Big\lgroup\!\!
 \prod_{\fE\in \ov \bE_{\ov\lt}}\!\!
 \Big\rgroup_{\!\!\cC_\al}
 \mathring{\P}\Big(\!
 \bigoplus_{\;
 	\begin{subarray}{c}
 		\fe\in\dot \fE
 	\end{subarray}
 }
 \pi^*L_{[\fe,o)}
 \Big)\,,
\end{align*}
has a smooth algebraic stack structure,
and the projection 
$$
 \pi^{\tf} :\cC^{\tf}\lra\fM^{\tf}
$$ induced by $\pi$ gives the universal family of $\fM^{\tf}$.

Similarly, with notation as in Theorem~\ref{Thm:tf_smooth_revert}, the union
\begin{align*}
	&\hspace{-.66in}\cC^{\rtf}
	=\bigsqcup_{(\al,
		\ud\lt)\in A^{\ud\der}}
	\hspace{-.15in}
	\jia{\cC_\al}_{\ud\lt}\ ,\qquad
	\tn{where}\\
	\jia{\cC_\al}_{\ud\lt}
	:=
	\Big\{\,
	&\big(\,[c^\fE_{\fe,\fe'}]_{\fe\in\udt\fE,\fe'\in\Om(\fe)}\big)_{\fE\in\ud\bE_{\ud\lt}}
	\in
	\Big\lgroup\!\!
	\prod_{\fE\in \ud \bE_{\ud\lt}}\!\!
	\Big\rgroup_{\!\!\cC_\al}
	\mathring{\P}\Big(\!
	\bigoplus_{\;
		\begin{subarray}{c}
			\fe\in\udt \fE
		\end{subarray}
	}\;
	\bigoplus_{
		\begin{subarray}{c}
			\fe'\in\Om(\fe)
		\end{subarray}
	}\!\!
	\pi^*L_{[\fe',\fe]}
	\Big):\\
	&
	c^{\fE}_{\fe,\fe'}\!\otimes
	c^{\fE_-}_{p_{\fe}(\fe''),\fe''}=
	c^{\fE_-}_{p_{\fe}(\fe'),\fe'}\!\otimes
	c^{\fE}_{\fe,\fe''}\ \ \big(
	\in 
	\pi^*L_{[\fe',\fe)}\!\otimes_{\,\cC_\al}\!\pi^*L_{[\fe'',\fe)}\!\otimes_{\,\cC_\al}\!\pi^*L_{\fe}\,\big)
	\\
	&\forall~
	\fE\in\ud\bE_{\ud\lt}\bsl\{\fE_1\eq\min(\tau_{\ud\lt})\},~
	\fe\inn \udt\fE,~
	\fe',\fe''\inn \Om(\fe)\,
	\Big\}\,,
\end{align*}
has a smooth algebraic stack structure,
and the projection 
$$
\pi^{\rtf} :\cC^{\rtf}\lra\fM^{\rtf}
$$ induced by $\pi$ gives the universal family of $\fM^{\rtf}$.
\end{prp}

\subsection{Comparison with local blowups}
\label{Subsec:local_blowups}

In this subsection,
we establish the local blowup properties~\ref{Cond:smooth_local_blowup} and~\ref{Cond:smooth_local_blowup_revert}.
This will complete the proofs of Theorems~\ref{Thm:tf_smooth} and~\ref{Thm:tf_smooth_revert}.

\begin{dfn}\label{Dfn:Local_RL_Comptbl}
	A morphism $\ti\fM\!\lra\!\fM$ is called a \ts{root-to-leaf (RL-) compatible local blowup} (resp.~a \ts{leaf-to-root (LR-) compatible local blowup}) if there exists an atlas $\{\cV\}$ of $\fM$ so that for every $\cV\inn\{\cV\}$,
	there exist
	a rooted tree~$\tau$, a partition of $\Xi(\tau)$:
	$$
	\Xi(\tau)
	=\bigsqcup_{k\ge 1}
	\Xi_k(\tau),
	$$
	a subset $\{z_e\}_{e\in\tau}$ of a system of local parameters on $\cV$,
	and a sequence of closed loci 
	$$Z_1^\cV,\, Z_2^\cV,\, \cdots\subset\cV,$$ satisfying the following conditions:
	\begin{itemize}[leftmargin=*]
		\item
		for every $k\ge 1$,
		$$
		Z_{k}^\cV= 
		\bigcup_{\fE\in\Xi_k(\tau)} \!\!\!\!
		X_\fE^\cV,\qquad\tn{where}\quad
		X_E^\cV:=\big\{\,z_e\eq 0:\,
		e\inn E\,\big\}\ \ \forall~E\subset\tau\,;
		$$
		
		\item
		if $\fE'\inn\Xi_{k'}(\tau)$, $\fE''\inn\Xi_{k''}(\tau)$, and $\fE'\!\succ\!\fE''$,
		then $k'\!<\!k''$ (resp.~$k'\!>\!k''$); and
		
		\item 
		$\ti\fM/\cV$ is isomorphic to the blowup of $\ti\cV/\cV$ successively along the proper transforms of $Z_1^\cV,$ $Z_2^\cV,$ $\cdots$.
	\end{itemize}
\end{dfn}

Notice that for RL- or LR-compatible local blowups, distinct elements of each $\Xi_k(\tau)$ are incomparable. 
In addition, the blowup procedure is finite on each chart $\cV$ (i.e.~only finitely many $Z_k^\cV$ are nonempty), because so is $\tau$ and hence the set $\Xi(\tau)$ of its transverse sections.


\begin{lmm}
	\label{Lm:RL/LR_cptb_blowups}
If the morphism $\ti\fM\!\lra\!\fM$  is an RL- or LR-compatible local blowup,
then for every $\cV\inn\{\cV\}$ and $k\!>\!1$,
after the $(k\!-\!1)$-th step of the blowup $\ti\cV/\cV$,
the proper transforms of the irreducible components $X_\fE^\cV$ of $Z_k^\cV$ are all smooth and pairwise disjoint. 
\end{lmm}

\begin{proof}
First, we assume that $\ti\fM/\fM$ is LR-compatible.
For every $\fE,\fE'\inn\Xi_k(\tau)$ with $k\!>\!1$,
by Part~\ref{Cond:TS_min} of  Lemma~\ref{Lm:Transverse_sections}, we have
$$
\fF:=\min\big(\fE\cup\fE'\big)\in\Xi(\tau),\qquad\tn{satisfying}\qquad
\fF\prec\fE,\fE'.
$$
In fact, $\fF$ is the greatest common predecessor (w.r.t.~$\prec$) of $\fE$ and $\fE'$.
Notice that $\fF\!\subset\!\fE\!\cup\!\fE'$, so
\begin{align*}
X_\fE^\cV\cap X_{\fE'}^\cV=X_{\fE\cup\fE'}^\cV\subset X_\fF^\cV.
\end{align*}
In addition,
with $L_e/X^\cV_{\{e\}}$, $e\inn\tau$, denoting the line bundles as in (\ref{Eqn:L_e}),
we see the restriction of the normal bundle of $X^\cV_\fF$ in $\cV$ to $X^\cV_\fE\!\cap\!X^\cV_{\fE'}$ can be written as
$$
 \Big(
 \bigoplus_{e\in\fF\cap\fE\cap\fE'}\!\!\!\!\!\!\big(L_e\big/(X^\cV_\fE\!\cap\!X^\cV_{\fE'})\big)
 \Big)
 \oplus
 \Big(
 \bigoplus_{e\in\fF\cap(\fE\bsl\fE')}\!\!\!\!\!\!\big(L_e\big/(X^\cV_\fE\!\cap\!X^\cV_{\fE'})\big)
 \Big)
 \oplus
 \Big(
 \bigoplus_{e\in\fF\cap(\fE'\bsl\fE)}\!\!\!\!\!\!\big(L_e\big/(X^\cV_\fE\!\cap\!X^\cV_{\fE'})\big)
 \Big),
$$
where the leftmost direct sum is in the normal direction of both $X^\cV_{\fE}$ and $X^\cV_{\fE'}$, while
the middle (resp.~rightmost) sum is in the normal direction of $X^\cV_{\fE}$ (resp.~$X^\cV_{\fE'}$) and tangent direction of $X^\cV_{\fE'}$ (resp.~$X^\cV_{\fE}$).

Let $h\!\ge\!1$ be such that $X_\fF^\cV\!\subset\! Z_{h}^\cV$.
Then,  $h\!<\!k$.
The last paragraph implies the proper transforms of $X_\fE^\cV$ and $X_{\fE'}^\cV$ becomes disjoint after $X_\fF^\cV$ is blown up in the $h$-th step.
The smoothness of the proper transforms of $X_\fE^\cV$ and $X_{\fE'}^\cV$ follows from direct computation.

The RL-compatible case follows from a parallel argument,
with $\fF$ replaced by 
$
\max \big(\fE\!\cup\!\fE'\big),
$
which is also a transverse section by Part~\ref{Cond:TS_max} of Lemma~\ref{Lm:Transverse_sections}. 
(In fact, it is the least common successor of $\fE$ and $\fE'$ in $\Xi(\tau)$).
The rest of the argument holds verbatim.
\end{proof}

The RL-compatible case of Lemma~\ref{Lm:RL/LR_cptb_blowups} generalizes a similar statement of~\cite{HLN} (where $\ti\fM/\fM$ is assumed to be a {\it global} blowup).
	
\begin{crl}\label{Crl:RL/LR_cptb_blowups}
If the morphism $\ti\fM\!\lra\!\fM$ is an RL- or LR-compatible local blowup,
then the blowup $\ti\cV'\!\lra\!\cV$ successively along the total transforms of
$$
 Y_1^\cV=Z_1^\cV,\quad
 Y_2^\cV=Z_1^\cV\cup Z_2^\cV,\quad 
 Y_3^\cV=Z_1^\cV\cup Z_2^\cV\cup Z_3^\cV,\quad\ldots
$$
yields the same blowup as $\ti\cV\!\lra\!\cV$.
\end{crl}

\begin{proof}
	The proof of the RL-compatible case is parallel to that of~\cite[Lemma~2.40]{g1modular}.
	Below we verify the LR-compatible case in a similar way.
	
Assume that $\ti\fM/\fM$ is an LR-compatible local blowup.
We aim to show by induction that after each step $k$, the blowup $\ti\cV_{(k)}$ of $\cV$ along the proper transforms of $Z_1^\cV,\ldots,Z_k^\cV$ is the same as the blowup $\ti\cV'_{(k)}$ of $\cV$ along the total transforms of $Y_1^\cV,\ldots,Y_k^\cV$.

The base case is trivial.
For $k\!>\!1$, assume that $\ti\cV_{(k-1)}/\cV\eq\ti\cV'_{(k-1)}/\cV$.
Fix an arbitrary $x\inn\cV$ as well as a lift $\ti x$ of $x$ in $\ti\cV_{(k-1)}$.
It suffices to verify that blowing up along the proper transform $\wc Z_k^\cV$ of $Z_k^\cV$ and along the total transform $\ti Y_k^\cV$ of $Y_k^\cV$ have the same effect on a neighborhood of $\ti x$.
Since such verification is local, w.l.o.g.~we also assume that $x\inn\bigcap_{1\le h\le k} Z_h^\cV$ (otherwise we simply omit the loci $Z_h^\cV$ not
containing $x$ and change the indices of $Z_h^\cV$ and $Y_h^\cV$ accordingly).

In fact, by Lemma~\ref{Lm:RL/LR_cptb_blowups},
we have 
$$H\subset\{1,\cdots,k\!-\!1\}\qquad\tn{and}\qquad
\fE_{h}\in\Xi_h(\tau)\quad\forall~h\inn H
$$ 
such that $\fE_h\!\succ\!\fE_{h'}$ for all $h,h'\inn H$ with $h\!<\!h'$, and $\ti x$ is contained in the exceptional divisor corresponding to $X_{\fE_{h}}^\cV$ for all $h\inn H$, and not contained in any other exception divisor.
Set $\hbar\!:=\!\max( H)$.
Applying Lemma~\ref{Lm:RL/LR_cptb_blowups} to $Z_k^\cV$ and mimicking the proof of~\cite[Lemma~4.3]{g1modular},
we obtain an affine smooth chart $\ti\cV_{\ti x}\!\lra\!\ti\cV_{(k-1)}$ and a subset $\{\ve_h\}_{h\in H}$ of a system of local parameters
on $\ti\cV_{\ti x}$ such that
the total transform of each $Z_h^\cV$, $h\inn H$ (resp.~$h\inn\{1,\ldots,k\!-\!1\}\bsl H$), is locally given by $\{\ve_h=0\}$ (resp.~$\emptyset$).
Therefore,
we have
$$
\ti Y_k^\cV\!\cap\!\ti \cV_{\ti x}=
\big(\wc Z_k^\cV\!\cap\!\ti \cV_{\ti x}\big)\cup
\big\{\prod_{h\in H}\ve_h\eq 0\,\big\}.
$$

With $\wc z_e$ denoting the proper transforms of $z_e$ in Definition~\ref{Dfn:Local_RL_Comptbl},
consider the following subsets of $\tau$:
$$
E':=
\big\{\,\min\{e'\succ e\}:~ e\inn\fE_{\hbar},~\wc z_e(\ti x)\!\ne\!0\,\big\},\qquad
\fF:=E'\sqcup \{e\inn\fE_{\hbar}:~e\!\not\prec\!e'~\forall~e'\inn E'\}.
$$
It is direct to check that $\fF\inn\Xi(\tau)$,
and $\{\wc z_e\}_{e\in\fF}\cup\{\ve_h\}_{h\in H}$ is a subset of a system of local parameters on $\ti\cV_{\ti x}$.

If $\fF\!\not\in\!\Xi_k(\tau)$,
then 
$$
 \wc Z_k^\cV\!\cap\!\ti \cV_{\ti x}=\emptyset\qquad\tn{and}\qquad
 \ti Y_k^\cV\!\cap\!\ti \cV_{\ti x}=
 \big\{\prod_{h\in H}\ve_h\eq 0\,\big\},
$$
so blowing up along $\wc Z_k^\cV$ and along $\ti Y_k^\cV$ does not affect $\ti\cV_{\ti x}$.

If $\fF\!\in\!\Xi_k(\tau)$, then
$$
\wc Z_k^\cV\!\cap\!\ti \cV_{\ti x}=\{\wc z_e:e\inn\fF\}\qquad\tn{and}\qquad
\ti Y_k^\cV\!\cap\!\ti \cV_{\ti x}=\{\wc z_e:e\inn\fF\}\cup
\big\{\prod_{h\in H}\ve_h\eq 0\,\big\}.
$$
In other words, on $\ti\cV_{\ti x}$,
$\wc Z_k^\cV$ and $\ti Y_k^\cV$ are defined by the ideals
$$
\mathscr I_{\wc Z_k^\cV}:=
\lr{\wc z_e: e\inn\fF}\qquad\tn{and}\qquad
\mathscr I_{\wc Z_k^\cV}\,\big(\prod_{h\in H}\!\ve_h\big),
$$
respectively.
Therefore, blowing up along $\wc Z_k^\cV$ has the same effect on $\ti\cV_{\ti x}$ as that along $\ti Y_k^\cV$.
\end{proof}

\begin{prp}\label{Prp:RL/LR-cptb}
Let $\fM^{\tf}/\fM$ (resp.~$\fM^{\rtf}/\fM$) be a stack with RL (resp.~LR) twisted fields over $\fM$ as in Theorem~\ref{Thm:tf_smooth} (resp.~Theorem~\ref{Thm:tf_smooth_revert}).
Then, $\fM^{\tf}/\fM$ (resp.~$\fM^{\rtf}/\fM$) is an RL-compatible (resp.~LR-compatible) local blowup.
\end{prp}

\begin{proof}
In the RL circumstance, the proof is the same as that of \cite[Proposition~4.4]{g1modular}.
The proof in the LR circumstance follows the same idea,
which we elaborate on below.

Let $\fM^{\rtf}/\fM$ be as in Theorem~\ref{Thm:tf_smooth_revert}.
Fix $\al\inn A$.
If $\tau_\al\eq\tau_\bullet$, then on any chart $\cV\inn\fV_\al$ if $\al\inn B$ and any $\cV\!\lra\!\fM_\ms\eq\fM^\mn$ if $\al\eq\ms$,
no twisted field is added, and no blowup is exerted, so $\fM^{\rtf}/\cV$ is the same as the trivial blowup of $\cV$.
Hereafter we assume $\tau_\al\!\ne\!\tau_\bullet$ and fix $\cV\inn\fV_\al$.
Let
$$
\ell=|\min(\tau_\al)|\ (\ge 1).
$$

With $\ze_e$, $e\inn S_\al$ denoting the modular parameters on $\cV$,
we set
\begin{align}
&X_{\fE}=X_{\fE}^\cV=\{\ze_e^\cV\eq 0:\,e\inn\fE\}&&\forall\quad\fE\in\Xi(\tau_\al)\,,\nonumber\\
&Z_k=Z^\cV_k=\bigcup_{\fE\in\Xi(\tau_\al)~\tn{s.t.}~|\fE|=\ell-k+1}\hspace{-.22in}X_{\fE}&&
\forall\quad 1\le k\le \ell.
\label{Eqn:Z_k}
\end{align}
Let $\pi:\ti\cV\!\lra\!\cV$ be the blowup of $\cV$ successively along the proper transforms of $Z_1, Z_2, \ldots, Z_\ell$.
To show $\fM^{\rtf}/\fM$ is an RL-compatible local blowup,
it suffices to show the restriction of the forgetful morphism $\ud\varpi:\fM^{\rtf}\!\lra\!\fM$ to $$\cV^{\rtf}:=\ud\varpi^{-1}(\cV)$$ is isomorphic to $\ti\cV/\cV$.

To this end,
it suffices to construct morphisms $\psi_1:\cV^{\rtf}\!\lra\!\ti\cV$ and $\psi_2:\ti\cV\!\lra\!\cV^{\rtf}$ such that
\begin{align*}
	\pi\circ\psi_1=\ud\varpi\qquad\tn{and}\qquad
	\ud\varpi\circ\psi_2=\pi\,;
\end{align*}
c.f.~the first paragraph of the proof of \cite[Proposition~4.4]{g1modular}.

First, we construct $\psi_1$.
By Corollary~\ref{Crl:RL/LR_cptb_blowups},
$\ti\cV/\cV$ can be obtained by blowing up $\cV$ successively along the total transforms of
$$
Y_1:=Z_1=\bigcup_{\fE\in\Xi(\tau_\al)~\tn{s.t.}~|\fE|\ge \ell}\hspace{-.22in}X_{\fE}\,,\qquad
\cdots,\qquad 
Y_k:=Z_1\cup\cdots\cup Z_k =
\bigcup_{\fE\in\Xi(\tau_\al)~\tn{s.t.}~|\fE|\ge \ell-k+1}\hspace{-.22in}X_{\fE}\,,\qquad\cdots
$$
We will show for every $k\!\ge\!1$, 
$\ud\varpi^{-1}(Y_k)$ is a Cartier divisor in $\cV^{\rtf}$.
The universality of blowups will then give rise to a morphism $\psi_1$ of $\cV^{\rtf}/\cV$ to $\ti\cV/\cV$.

The proof of $\ud\varpi^{-1}(Y_k)$ being a Cartier divisor is slightly more complicated than that of \cite[Proposition~4.4]{g1modular},
so we elaborate on it here.
Fix arbitrary $$\ud\lt=\big(\,\tau_\al,\ \ud\bE\eq\fE_1,\ldots,\fE_\cht,\ \Dm(\ud\lt)\,\big)\in\ud\bT_\al.$$
Let $\cU\eq\cU_{\cV,\ud\lt}\!\lra\fM^{\rtf}$ be the twisted chart as in Theorem~\ref{Thm:tf_smooth_revert}~\ref{Cond:smooth_parameters_revert}.
Since $\varpi(\cU)\!\subset\!\cV$,
$\cU$ can be viewed as a chart of $\cV^{\rtf}$.
Let 
$$
\fF=
\max\{\,\fE\in\ud\bE:\,
|\fE|\ge\ell-k+1\,
\}
$$
Here the maximum is taken with respect to the order (\ref{Eqn:transverse_sections_order}).
The set in the above equality is not empty because $\min(\tau_\al)\inn\ud\bE$ and $|\min(\tau_\al)|\eq\ell$.

On the one hand,
for each $\fE'\inn\ud\bE$ with $\fE'\!\preceq\!\fF$,
we have $|\fE'|\!\ge\!|\fF|\!\ge\!\ell\!-\!k\!+\!1$,
because for any $\fE,\fE'\inn\Xi(\tau_\al)$,
\begin{align}\label{Eqn:|trans_sec|_order}
	\fE\preceq\fE'\quad\Longrightarrow\quad
	|\fE|\ge|\fE'|.
\end{align} 
Therefore, $X_{\fE'}$ is an irreducible component of $Y_k$.
Notice that by~(\ref{Eqn:mod_par_pullback}),  we have
$$
\ud\varpi\,\big\{\xi^\cU_{\fE'}\eq 0\big\}\subset
X_{\fE'}.
$$
Consequently,
\begin{align}
	\label{Eqn:Y_k}
\big\{\prod_{\fE'\in\ud\bE~\tn{s.t.}~\fE'\preceq\fF}\!\!\!\!\!\!\!\!\xi^\cU_{\fE'}= 0\,\big\}
\,\subset\, \ud\varpi^{-1}(Y_k)\,.
\end{align}

On the other hand,
for every irreducible component $X_\fE$ of $Y_k$,
consider the intersection
$$\fI:=\fE\cap\Big(\big(\fF^\preceq\!\cap\Dm(\ud\lt)\big)^\wedge\Big).$$
\begin{itemize}[leftmargin=*]
\item If $\fF^\preceq\!\cap\!\Dm(\ud\lt)\!\cap\!\max(\tau_\al)\!\ne\!\emptyset$,
then by Lemma~\ref{Lm:dominant}~\ref{Cond:dom_below} and Lemmas~\ref{Lm:path^} and~\ref{Lm:E^_simple},
we see  $\big(\fF^\preceq\!\cap\Dm(\ud\lt)\big)^\wedge\eq\tau_\al$,
hence $\fI\!\ne\!\emptyset$.

\item 
If $\fF^\preceq\!\cap\!\Dm(\ud\lt)\!\cap\!\max(\tau_\al)\!=\!\emptyset$,
then by Definition~\ref{Dfn:Enhanced_rtl_seq}, this implies
$\fF\!\ne\!\fE_\cht$.
With $\fF_+\inn\ud\bE$  as in (\ref{Eqn:fE_+-}),
we then have
$$
 \big(\fF^\preceq\!\cap\Dm(\ud\lt)\big)^\wedge
 =
 \fF^\preceq\bsl\fF_+\,.
$$
To see the above equality,
just observe the same equality holds when $\fF$ is replaced by any $\fE'\inn\ud\bE\bsl\{\fE_\cht\}$.
This can be shown by applying Part~\ref{Cond:adj_edge_contrs} of Lemma~\ref{Lm:E^} repeatedly, starting with $\fE'\eq\min(\tau_\al)$.

We then claim that 
$$
\fE\cap(\fF^\preceq\bsl\fF_+)\ne\emptyset.
$$
If this was not true,
then $\fE\!\subset\!\fF_+^\succeq$,
hence by (\ref{Eqn:|trans_sec|_order}),
we would have 
$|\fF_+|\!\ge\!|\fE|\!\ge\!\ell\!-\!k\!+\!1$,
which would contradict the description of $\fF$.

In sum, $\fI\!\ne\!\emptyset$ still holds in this case.
\end{itemize}

Since the intersection $\fI$ is always nonempty, and $\fE$ is a transverse section,
we have
$$
\fE\cap \fF^\preceq\!\cap\Dm(\ud\lt)\ne\emptyset\,.
$$
Therefore, there exists $\fE''\inn\ud\bE$ with $\fE''\!\preceq\!\fF$ such that $\fE\!\cap\!\udt\fE''\!\ne\!\emptyset$.
Fix $e''\inn \fE\!\cap\!\udt\fE''\!\ne\!\emptyset$.
Then,
$$
X_{\fE}\subset\{\ze_{e''}^\cV\eq 0\},
$$
hence
$$
\ud\varpi^{-1}(X_\fE)\,\subset\,
\{\varpi^*\ze_{e''}^\cV\eq 0\}\,\subset\,
\big\{\prod_{\fE'\in\ud\bE~\tn{s.t.}~\fE'\preceq\fE''}\!\!\!\!\!\!\!\!\xi^\cU_{\fE'}= 0\,\big\}
\,\subset\,
\big\{\prod_{\fE'\in\ud\bE~\tn{s.t.}~\fE'\preceq\fF}\!\!\!\!\!\!\!\!\xi^\cU_{\fE'}= 0\,\big\}.
$$
The middle inclusion above follows from 
(\ref{Eqn:mod_par_pullback}) and the fact that $e''\inn\udt\fE''$.
Since $X_\fE$ is an arbitrary irreducible component of $Y_k$,
we conclude that
\begin{align}
	\label{Eqn:Y_k'}
\ud\varpi^{-1}(Y_k)
\,\subset\, \big\{\prod_{\fE'\in\ud\bE~\tn{s.t.}~\fE'\preceq\fF}\!\!\!\!\!\!\!\!\xi^\cU_{\fE'}= 0\,\big\}\,.
\end{align}

Combining (\ref{Eqn:Y_k}) and (\ref{Eqn:Y_k'}),
we have 
$$
\ud\varpi^{-1}(Y_k)
\,=\, \big\{\prod_{\fE'\in\ud\bE~\tn{s.t.}~\fE'\preceq\fF}\!\!\!\!\!\!\!\!\xi^\cU_{\fE'}= 0\,\big\}\,,
$$
which is a Cartier divisor of $\cV^{\rtf}$.
In this way, we obtain the morphism $\psi_1$ .

The morphism $\psi_2$ is constructed using the same direct approach as in the proof of~\cite[Proposition~4.4]{g1modular}.
Given $x\inn\fM_\al\cap\cV$ and $\ti x\inn\ti\cV$ overlying $x$,
$\ti x$ is contained in a unique maximal sequence of exceptional divisors,
each obtained from blowing up along the proper transform of some $X_\fE$.
Thus, we obtain an LR-sequence $(\tau_\al,\ud\bE)$.
The dominant edges are determined by the non-zero components of the exceptional divisors,
which then gives rise to an LR-sequence with dominant edges
$\ud\lt\inn\ud\bT_\al$.
In this way,
the stratum $\jia{\fM_\al}_{\ud\lt}$ containing the proposed image of $\ti x$ under $\psi_2$ is determined.

The explicit expression of $\psi_2(\ti x)$, as well as the reason that such point-wise defined map $\psi_2$ is indeed a morphism,
are parallel to that of the proof of~\cite[Proposition~4.4]{g1modular}.
We omit further details.
\end{proof}

\subsection{Properties of STF}
Now that the stacks with twisted fields have been constructed in Theorems~\ref{Thm:tf_smooth} and~\ref{Thm:tf_smooth_revert},
we explore some key properties of them.
Throughout this subsection,
we fix an LES $\fM\eq\bigsqcup_{\al\in A}\fM_\al$ and a treelike structure $\La$ on $\fM$.

We begin with a direct conclusion of Theorem~\ref{Thm:tf_smooth}~\ref{Cond:smooth_parameters}.
\begin{crl}
	\label{Crl:prod_RL}
	Let $\al$,  $\ov\lt$, $\cV\inn\fV_{\al}$ and $\cU\inn\fU_{(\al,\ov\lt)}$ be as in Theorem~\ref{Thm:tf_smooth}~\ref{Cond:smooth_parameters}.
	Assume $\tau_\al\!\ne\!\tau_\bullet$.
	Then, for every $\fe\inn\tau_{\al}$, there exists $u_{\fe}\inn\Ga(\sO^*_\cU)$ such that
	$$
	\ov\varpi^*\big(\prod_{\fe'\succeq\fe}\ze_{\be_{\al}(\fe')}^\cV\,\big)
	\,=\; 
	u_{\fe}\cdot
	\prod_{
		\fe'\in\,\ND(\ov\lt)~\tn{s.t.}\;
		\fe'\succeq \fe
	}
	\hspace{-.25in}\xi_{\be_\al(\fe')}^{\cU}\hspace{.25in}
	\cdot\hspace{.1in}
	\prod_{
		\begin{subarray}{c}
			\fE\in\ov\bE~\tn{s.t.}\;
			\fe\in\fE^\preceq
	\end{subarray}}\!\!\!\!\!\xi_\fE^{\cU}\;.
	$$	
	Moreover,
	the tautological set of monomials $\big\{\ov\varpi^*\big(\prod_{\fe'\succeq\fe}\ze_{\be_{\al}(\fe')}^\cV\,\big):
	\fe\inn \min(\tau_\al)\big\}$ admits a divisibly minimal element.
	Precisely,
	there exists $\fe_0\inn\min(\tau_\al)$ so that $$
	\ov\varpi^*\big(\prod_{\fe'\succeq{\fe}_0}\ze_{\be_{\al}(\fe')}^\cV\,\big)~
	\Big|~
	\ov\varpi^*\big(\prod_{\fe'\succeq\fe}\ze_{\be_{\al}(\fe')}^\cV\,\big)\qquad
	\forall~\fe\in\min(\tau_\al).
	$$
\end{crl}

\begin{proof}
	The first statement follows immediately from~(\ref{Eqn:mod_par_pullback}).
	
	To see the second statement, consider an arbitrary $\fe\inn\min(\tau_\al)$.
	By Lemma~\ref{Lm:dominant}~\ref{Cond:dom_above},
	$$
	\prod_{
		\fe'\in\,\ND(\ov\lt)~\tn{s.t.}\;
		\fe'\succeq \fe
	}
	\hspace{-.25in}\xi_{\be_\al(\fe')}^{\cU}\hspace{.15in}
	\in \Ga(\sO^*_\cU)
	\qquad\Longleftrightarrow\qquad
	\fe\in\dot\fE_\cht\cap\min(\tau_\al)\,.
	$$
	In addition,
	we have
	$$
	\prod_{
		\begin{subarray}{c}
			\fE\in\ov\bE~\tn{s.t.}\;
			\fe\in\fE^\preceq	
	\end{subarray}}\!\!\!\!\!\!\!\xi_\fE^{\cU}\ 
	=\ \prod_{\fE\in\ov\bE}\xi^\cU_\fE\,,
	$$
	which is a common factor for all $\fe\inn\min(\tau_\al)$.
	Since $\dot\fE_\cht\!\cap\!\min(\tau_\al)$ is nonempty by Definition~\ref{Dfn:Enhanced_rtl_seq},
	we choose $\fe_0\inn\min(\tau_\al)$ so that the second statement is established.
\end{proof}

The same argument also applies to the LR circumstances below,
where Part~\ref{Cond:dom_below} of Lemma~\ref{Lm:dominant} is needed instead of Part~\ref{Cond:dom_above}.
We omit further details.

\begin{crl}
	\label{Crl:prod_LR}
	Let $\al$,  $\ud\lt$, $\cV\inn\fV_{\al}$ and $\cU\inn\fU_{(\al,\ud\lt)}$ be as in Theorem~\ref{Thm:tf_smooth_revert}~\ref{Cond:smooth_parameters_revert}.
	Assume $\tau_\al\!\ne\!\tau_\bullet$.
	Then, for every $\fe\inn\tau_{\al}$, there exists $u_{\fe}\inn\Ga(\sO^*_\cU)$ such that
	$$
	\ud\varpi^*\big(\prod_{\fe'\succeq\fe}\ze_{\be_{\al}(\fe')}^\cV\,\big)
	\,=\;
	u_{\fe}\cdot
	\prod_{
		\fe'\in\,\ND(\ud\lt)~\tn{s.t.}\;
		\fe'\succeq \fe
	}
	\hspace{-.25in}\xi_{\be_\al(\fe')}^{\cU}\hspace{.25in}
	\cdot\hspace{.1in}
	\prod_{
		\begin{subarray}{c}
			\fE\in\ud\bE~\tn{s.t.}\;
			\fe\in\fE^\preceq
	\end{subarray}}\!\!\!\!\!\xi_\fE^{\cU}\;.
	$$	
	Moreover,
	the tautological set of monomials $\big\{\ov\varpi^*\big(\prod_{\fe'\succeq\fe}\ze_{\be_{\al}(\fe')}^\cV\,\big):
	\fe\inn \min(\tau_\al)\big\}$ admits a divisibly minimal element.
	Precisely,
	there exists $\fe_0\inn\min(\tau_\al)$ so that $$
	\ud\varpi^*\big(\prod_{\fe'\succeq{\fe}_0}\ze_{\be_{\al}(\fe')}^\cV\,\big)~
	\Big|~
	\ud\varpi^*\big(\prod_{\fe'\succeq\fe}\ze_{\be_{\al}(\fe')}^\cV\,\big)\qquad
	\forall~\fe\inn\min(\tau_\al).
	$$
\end{crl}	

Next,
we study the universality of the STF, based on Theorem~\ref{Thm:tf_smooth}~\ref{Cond:smooth_local_blowup} and Theorem~\ref{Thm:tf_smooth_revert}~\ref{Cond:smooth_local_blowup_revert}.
For every $\al\inn B$ with $\tau_\al\!\ne\!\tau_\bullet$,
$\cV\inn\fV_\al$, and $\fE\inn\Xi(\tau_\al)$,
observe that the locus of $\cV$
\begin{align*}
	X^\cV_{\be_\al(\fE)}:=
	\{\,\ze_{\be_\al(e)}^\cV\eq 0:\,
	e\inn\fE\,\}
\end{align*}
is an irreducible component of 
\begin{align*}
	\tn{Cl}_\fM\big(\fM_{\al_{(S_\al\bsl\be_\al(\fE))}}\big)\cap\cV
	=
	\tn{Cl}_\cV\big(\fM_{\al_{(S_\al\bsl\be_\al(\fE))}}\cap\cV\big)\,;
\end{align*}
the equality above holds because $\cV\!\lra\!\fM$ is a chart.
We shall analyze the pullbacks of
\begin{align}
	\label{Eqn:univ}
 \ov X^\cV_{\be_\al(\fE)}:=
 \bigcup_{\fE'\in\Xi(\tau_\al)~\tn{s.t.}\;\fE'\succeq\fE}\hspace{-.35in}
 X^\cV_{\be_\al(\fE')}\qquad\tn{and}\qquad
 \ud X^\cV_{\be_\al(\fE)}:= \bigcup_{\fE'\in\Xi(\tau_\al)~\tn{s.t.}\;\fE'\preceq\fE}\hspace{-.35in}
 X^\cV_{\be_\al(\fE')}\qquad
 (\subset\cV)\,,
\end{align}
where $\succeq$ and $\preceq$ refer to the order (\ref{Eqn:transverse_sections_order}) on $\Xi(\tau_\al)$,
in the following universality statement.

\begin{prp}
	\label{Prp:univ}
Let $f\!:\fN\!\lra\! \fM$ be a morphism of an irreducible stack $\fN$ to $\fM$
satisfying
\begin{align*}
	f(\fN)\!\cap\!\cV\not\subset \bigcup_{\fE\in\Xi(\tau_\al)}\!\!\!\!
	X^\cV_{\be_\al(\fE)}
	\qquad\forall~
	\cV\inn\fV_\al,\ 
	\al\inn B~\tn{with}~\tau_\al\!\ne\!\tau_\bullet\,,
\end{align*}
and $\ov\varpi\!:\fM^{\tf}\!\lra\!\fM$ and $\ud\varpi\!:\fM^{\rtf}\!\lra\!\fM$ be respectively the forgetful morphisms.
Then, $f$ uniquely factors through $\ov\varpi$,
if and only if for every
$\al\inn B$ with $\tau_\al\!\ne\!\tau_\bullet$,
$\cV\inn\fV_\al$, and $\fE\inn\Xi(\tau_\al)$,
the locus
$\ov X^\cV_{\be_\al(\fE)}$ as in (\ref{Eqn:univ}) pulls back to a Cartier divisor of $f^{-1}(\cV)$.

Similarly,
$f$ uniquely factors through $\ud\varpi$, if and only if for every
$\al\inn B$ with $\tau_\al\!\ne\!\tau_\bullet$,
$\cV\inn\fV_\al$, and $\fE\inn\Xi(\tau_\al)$,
the locus
$\ud X^\cV_{\be_\al(\fE)}$ as in (\ref{Eqn:univ}) pulls back to a Cartier divisor  of $f^{-1}(\cV)$.
\end{prp}

\begin{proof}
We begin with the RL case.

Consider arbitrary $\al\inn B$ with $\tau_\al\!\ne\!\tau_\bullet$,  $\fE\inn\Xi(\tau_\al)$, 
$\cV\inn\fV_\al$,   $\ov\lt\eq\big(\tau_\al,\{\fE_1,\ldots,\fE_\cht\},\Dm(\ov\lt)\big)\inn\ov\bT_\al$, 
and twisted chart $\cU\eq\cU_{\cV,\ov\lt}$ of $\fM^{\tf}$ over $\cV$.
Let
\begin{align*}
	h:=\min\big\{\,i\inn\{1,2,\cdots,\cht\}\,:\,
	\fE_i\!\cap\!\fE\!\ne\!\emptyset\,\big\}.
\end{align*}
The set on the right-hand side is nonempty because $\fE_\cht\!\cap\!\min(\tau_{\al})\!\ne\!\emptyset$ and $\fE$ is a transverse section.
Notice that
\begin{align*}
	\fE_i\succeq\fE\quad\forall~1\!\le\!i\!\le\!h\,,\qquad\tn{and}\qquad
	\fE'\cap\big(\!\bigcup_{1\le i\le h}\!\!\!\fE_i\big)\ne\emptyset\quad
	\forall~\fE'\inn\Xi(\tau_\al)~\tn{s.t.}~
	\fE'\!\succeq\!\fE\,.
\end{align*}
It is then a direct check 
\begin{align*}
	\ov\varpi^{-1}\big(\ov X^\cV_{\be_\al(\fE)}\big)=
	\big\{\,\prod_{1\le i\le h}\!\!\xi^\cV_{\fE_i}=0\,\big\}\quad\subset
	\cU.
\end{align*}
Particularly, it is a Cartier divisor of $\cU$. 
Since $\al$, $\fE$, $\cV$, and $\cU$ are arbitrary,
we conclude
if~$f$ factors through $\ov\varpi$,
then $\ov X^\cV_{\be_\al(\fE)}$ pulls back to a Cartier divisor of $f^{-1}(\cV)$.

Conversely, assume that
for every
$\al\inn B$ with $\tau_\al\!\ne\!\tau_\bullet$,
$\cV\inn\fV_\al$, and $\fE\inn\Xi(\tau_\al)$,
$\ov X^\cV_{\be_\al(\fE)}$ pulls back to a Cartier divisor of  $f^{-1}(\cV).$
Let $\le$ be an arbitrary linear order on $\Xi(\tau_\al)$ that extends the partial order (\ref{Eqn:transverse_sections_order}).
By Theorem~\ref{Thm:tf_smooth}~\ref{Cond:smooth_local_blowup} and Corollary~\ref{Crl:RL/LR_cptb_blowups},
$\ov\varpi|_{\cV^{\tf}}\!:\cV^{\tf}\!\lra\!\cV$ is isomorphic to the blowup $\ti\cV$ of $\cV$ successively along the total transforms of 
\begin{align*}
	\wh X^\cV_{\be_\al(\fE)}:=
	\bigcup_{\fE'\in\Xi(\tau_\al)~\tn{s.t.}\;\fE'\ge\fE}\hspace{-.2in} X^\cV_{\be_\al(\fE)}\,,\qquad
	\fE\in\Xi(\tau_\al)\,,
\end{align*}
with respect to the aforementioned linear order $\le$, starting from $\fE\eq\max(\tau_\al)$.

We then show by induction over $
\fE\inn\Xi(\tau_\al)$ that each $\wh X^\cV_{\be_\al(\fE)}$ pulls back to a Cartier divisor in $f^{-1}(\cV)$.
The base case $\fE\eq\max(\tau_\al)$ is true because $\wh X^\cV_{\be_\al(\max(\tau_\al))}\eq 
\ov X^\cV_{\be_\al(\max(\tau_\al))}$.
Now assume for some $\fF\inn\Xi(\tau_\al)$, $\wh X^\cV_{\be_\al(\fF)}$ pulls back to a Cartier divisor of $f^{-1}(\cV).$
We denote by $\fE$ the immediate predecessor of $\fF$ (w.r.t.~$\le$) in $\Xi(\tau_\al)$.
Observe that
\begin{align*}
	\wh X^\cV_{\be_\al(\fE)}=
	\wh X^\cV_{\be_\al(\fF)}\cup 
	X^\cV_{\be_\al(\fE)}
	=
	\wh X^\cV_{\be_\al(\fF)}\cup 
	\ov X^\cV_{\be_\al(\fE)}\,.
\end{align*}
Hence $\wh X^\cV_{\be_\al(\fE)}$ also pulls back to a Cartier divisor of $f^{-1}(\cV).$
The universality of blowup thus ensures
$f$ factors through $\ti\cV/\cV$,
hence through $\cV^{\tf}/\cV$.

Notice that for every $\al\inn A$ with $\tau_\al\eq\tau_\bullet$,
$\cV^{\tf}/\cV$ is but the identity morphism.
Therefore,
we obtain
\begin{align*}
	\ov f_\cV:f^{-1}(\cV)\lra\cV^{\tf}
	\quad\tn{satisfying}\quad
	\ov\varpi\circ\ov f_\cV=f|_{f^{-1}(\cV)}\qquad
	\forall~\al\inn A,~
	\cV\inn\fV_\al.
\end{align*}
It remains to show such $\ov f_\cV$'s together give rise to a morphism $\fN\!\lra\!\fM^{\tf}$ that $f$ factors through.
In fact,
for every $\al\inn A$ and $\cV\inn\fV_\al$,
\begin{align*}
	f^{-1}\Big(\cV\big\bsl \big(\bigcup_{\fE\in\Xi(\tau_\al)}\!\!\!\!
	X^\cV_{\be_\al(\fE)}\big)\Big)=
	\begin{cases}
		f^{-1}(\cV)\bsl f^{-1}\big(\ov X^{\cV}_{\be_\al(\min(\tau_\al))}\big)\ 
		&\tn{if}~\tau_\al\!\ne\!\tau_\bullet,\\
		f^{-1}(\cV)\ 
		&\tn{if}~\tau_\al\!=\!\tau_\bullet,
	\end{cases}
\end{align*}
is always a dense open subset of $f^{-1}(\cV)$ by assumption.
On this dense open subset,
we have $\ov f_\cV\eq f$, because  $\ov\varpi$ restricts to the identity morphism on the preimage of $\cV\big\bsl \big(\bigcup_{\fE\in\Xi(\tau_\al)}
X^\cV_{\be_\al(\fE)}\big)$.
Hence for every $\al,\al'\inn A$, $\cV\inn\fV_\al$, and $\cV'\inn\fV_{\al'}$,
we have
$\ov f_\cV\eq\ov f_{\cV'}\eq f$ on the dense open subset
\begin{align*}
	f^{-1}\Big((\cV\cap\cV')\big\bsl
	\big(\bigcup_{\fE\in\Xi(\tau_\al)}\!\!\!\!
	X^\cV_{\be_\al(\fE)}
	\cup 
	\bigcup_{\fE'\in\Xi(\tau_{\al'})}\!\!\!\!
	X^{\cV'}_{\be_{\al'}(\fE')}\big)\Big)\,.
\end{align*}
of $f^{-1}(\cV\!\cap\!\cV')$,
so 
$\ov f_\cV\eq\ov f_{\cV'}$ on $f^{-1}(\cV\!\cap\!\cV')$.
This completes the proof of the RL case.

The LR case is parallel, hence is omitted.
\end{proof}

Proposition~\ref{Prp:univ} leads to the following conclusion,
which indicates  $\fM^{\tf}/\fM$ and $\fM^{\rtf}/\fM$ are not isomorphic in general.

\begin{crl}
	\label{Crl:tf_rtf_non_isom}
There exists a morphism $f_1:\fM^{\tf}\!\lra\!\fM^{\rtf}$ (resp.~$f_2:\fM^{\rtf}\!\lra\!\fM^{\tf}$) satisfying $\ov\varpi\eq \ud\varpi\!\circ\!f_1$  (resp.~$\ud\varpi\eq \ov\varpi\!\circ\!f_2$),
if and only if for every $\al\inn A$,
the rooted tree $\tau_\al$ is trivially ordered or linearly ordered.
\end{crl}

\begin{proof}
	$\tau_\al$ being trivially ordered implies $|\Xi(\tau_\al)|\eq 0$ or 1,
	while $\tau_\al$ being linearly ordered implies $|\fE|\eq 1$ for all $\fE\inn\Xi(\tau_\al)$.
	Therefore, $\fM^{\tf}/\fM$ and $\fM^{\rtf}/\fM$ are naturally isomorphic.
	
	Conversely,
	assume there exist $\al\inn A$ such that $\tau_\al$ is neither trivially nor linearly ordered.
	There then exist distinct $e',e'',e'''\inn\tau_\al$ such that $e'\!\prec\! e'''$, whereas $e''$ is incomparable with $e'$ and~$e'''$.
	We choose $e_1,e_2\inn\min(\tau_\al)$ such that
	\begin{align*}
		e_1\prec e'\prec e''',\qquad
		e_2\prec e''.
	\end{align*}
	Then by~(\ref{Eqn:tree_order}),
	we have
	\begin{align*}
		e''\not\in e_1^{\,\succeq},\qquad
		e'\not\in e_2^{\,\succeq},\qquad
		e'''\not\in e_2^{\,\succeq}.
	\end{align*}
	
	For every $\fe\inn e_2^{\,\succeq}$, consider the subset of $\tau_\al$ given by
	\begin{align*}
		\fF_\fe:=
		\max\big(\bigcup_{\fE\in\Xi(\tau_\al)~\tn{s.t.}\;\fE\ni\fe}\hspace{-.2in}
		\fE\ \ \big)\qquad(\ni\fe)\,.
	\end{align*}
	By Corollary~\ref{Crl:Transverse_sections} and Lemma~\ref{Lm:Transverse_sections},
	we have
	$\fF_\fe\inn\Xi(\tau_\al)$.
	Moreover, for every pair $\fe,\fe^\flat\inn e_2^{\,\succeq}$ with $\fe\!\succ\!\fe^\flat$,
	we claim that 
	\begin{align}\label{Eqn:F_e_F_eflat}
		\fF_\fe\succ\fF_{\fe^\flat},\qquad
		\fF_\fe\bsl\fF_{\fe^\flat}=\{\fe\}.
	\end{align}
	Indeed, $\fF_\fe$ cannot be equal to $\fF_{\fe^\flat}$.
	Suppose $\fF_\fe\!\not\succ\!\fF_{\fe^\flat}$. 
	There would then exist $\fe'\inn\fF_{\fe^\flat}$ such that $\fe'\inn\fF_\fe^{\,\succ}$.
	Obviously, $\fe'\!\not\succ\!\fe$, for otherwise $\fe'\!\succ\!\fe^\flat$,
	which would make $\fF_{\fe^\flat}$ no longer a transverse section.
	Consequently, Lemma~\ref{Lm:Transverse_sections} would imply
	\begin{align*}
		\fe\in \max(\fF_\fe\!\cup\!\{\fe'\})
		\in\Xi(\tau_\al)\qquad\tn{and}\qquad
		\fF_{\fe}\prec  \max(\fF_\fe\!\cup\!\{\fe'\}),
	\end{align*}
	which would then contradict the definition of $\fF_{\fe}$.
	Therefore, the inequality in (\ref{Eqn:F_e_F_eflat}) holds.
	To see the equality in (\ref{Eqn:F_e_F_eflat}),
	first observe that $\fe\!\not\in\!\fF_{\fe^\flat}$,
	for otherwise $\fF_{\fe^\flat}$ would contain distinct comparable edges.
	Next, suppose there exists $\fe''\inn\fF_\fe\big\bsl \big(\{\fe\}\!\cup\!\fF_{\fe^\flat}\big)$.
	Then, 
	\begin{align*}
		\fe''\in\fF_{\fe^\flat}^{\,\succ},\qquad
		\fe^\flat\in \max(\fF_{\fe^\flat}\!\cup\!\{\fe''\})
		\in\Xi(\tau_\al),\qquad\tn{and}\qquad
		\fF_{\fe^\flat}\prec  \max(\fF_{\fe^\flat}\!\cup\!\{\fe''\}),
	\end{align*}
	which would contradict the definition of $\fF_{\fe^\flat}$.
	In sum, the equality in (\ref{Eqn:F_e_F_eflat}) is also established.
	
	Notice that $\fF_{\max(e_2^{\succ})}\eq\max(\tau_\al)$. 
	Hence by (\ref{Eqn:F_e_F_eflat}), we obtain an RLS with dominant edges
	\begin{align*}
		\ov\lt=
		\big(\;
		\tau_\al,~
		\ov\bE\eq\{\fF_\fe:\fe\inn e_2^{\,\succ}\},~
		\Dm(\ov\lt)\eq 
		e_2^{\,\succ}
		\,\big)\,.
	\end{align*}
	In this way,  $\fE_\cht\eq\fF_{e_2}$. 
	Moreover, we have 
	\begin{align}
		\label{Eqn:e'_in_cht_prec}
		(e')^\preceq\subset\fE_\cht^\prec
		\qquad\big(\subset\ND(\ov\lt)\,\big)\,,
	\end{align}
	for otherwise we would once again have
	\begin{align*}
		e'''\in\fF_{e_2}^{\,\succ},\qquad
		e_2\in \max(\fF_{e_2}\!\cup\!\{e'''\})
		\in\Xi(\tau_\al),\qquad\tn{and}\qquad
		\fF_{e_2}\prec  \max(\fF_{e_2}\!\cup\!\{e'''\}),
	\end{align*}
	which would contradict the definition of $\fF_{e_2}$.
	
	We are ready to show the nonexistence of $f_1:\fM^{\tf}\!\lra\!\fM^{\rtf}$ with $\ov\varpi\eq\ud\varpi\!\circ\!f_1$.
	Let
	\begin{align*}
		\fE_{e'}:=
		\max\big(\min(\tau_\al)\cup\{e'\}\big),
	\end{align*}
	which is a transverse section by  Lemma~\ref{Lm:Transverse_sections}.
	For every $\fE'\inn\Xi(\tau_\al)$ with $\fE'\!\preceq\!\fE$,
	(\ref{Eqn:e'_in_cht_prec}) and the fact $e'\!\not\in\!e_2^{\,\succeq}$ together imply 
	\begin{align*}
		e_2\in\fE'\subset\fE_\cht^\preceq\qquad\tn{and}\qquad
		\fE'\cap (e')^\succeq\ne\emptyset.
	\end{align*}
	Fix a modular chart $\cV\inn\fV_\al$ and the twisted chart $\cU\eq\cU_{\cV,\ov\lt}\inn\fV_{(\al,\ov\lt)}$.
	By Theorem~\ref{Thm:tf_smooth}~\ref{Cond:smooth_parameters} and (\ref{Eqn:e'_in_cht_prec}),
	the pullback in $\cU$ of the locus
	$
		\ud X_{\be_\al(\fE_{e'})}^\cV
	$ as in (\ref{Eqn:univ})
	is contained in
	\begin{align*}
		\big\{\,
		\xi^{\cU}_{\fF_{e_2}}\!=0\,,\ \ 
		\prod_{e\in (e')^\prec}\!\!\xi^{\cU}_{e}= 0\,
		\big\}\qquad\subset\cU\,.
	\end{align*}
	However, the codimension of the above locus is already 2, so
	$
	\ud X_{\be_\al(\fE_{e'})}^\cV
	$ cannot be pulled back to a Cartier divisor of $\cV^{\tf}$.
	Applying Proposition~\ref{Prp:univ} to $f\eq\ov\varpi$, we thus conclude
	there does not exist any $f_1:\fM^{\tf}\!\lra\!\fM^{\rtf}$ satisfying $\ov\varpi\eq\ud\varpi\!\circ\!f_1$.
	
	The argument for the nonexistence of $f_2:\fM^{\rtf}\!\lra\!\fM^{\tf}$ with $\ud\varpi\eq\ov\varpi\!\circ\!f_2$ is parallel, hence is omitted.
\end{proof}

The following conclusions of Proposition~\ref{Prp:univ} are needed in the construction of $\ti M_2^{\rm tf}(\P^n,d)$,
in that they allow local constructions to glue globally.

Let $\fM$ and $\La$ be as in Proposition~\ref{Prp:univ},
$I_U$ be an index set,
and $(f_{U,t}:U\!\lra\!\fM)_{t\in I_U}$ be a family of morphisms of a stack $U$.
\begin{dfn}\label{Dfn:weak_and_strong_adm}
We call $(f_{U,t}:U\!\lra\!\fM)_{t\in I_U}$ an  \ts{admissible family of morphisms} if
	for every $x\inn U$,
	there exists $\al(x)\inn A$ such that 
	\begin{itemize}[leftmargin=*]
		\item $f_{U,t}(x)\inn\fM_{\al(x)}$ for all $t\inn I_U$, and
		\item 
		for every $t_i\inn I_U$ and $\cV_i\inn\fV_{\al(x)}$ containing $f_{U,t_i}(x)$, $i\eq 1,2$,
		there exists a neighborhood $U'\!\subset\! f_{U,t_1}^{-1}(\cV_1)\!\cap\!f_{U,t_2}^{-1}(\cV_2)$ of $x$ such that
		\begin{align}\label{Eqn:f_U_assump''}
			\big\{f_{U,t_1}^*\!\big(\ze_e^{\cV_{1}}\big)\eq 0\ 
			\forall\,e\inn \be_{\al(x)}\!(\fE)\big\}\cap U'
			\!=
			\big\{f_{U,t_2}^*\!\big(\ze_e^{\cV_{2}}\big)\eq 0\ 
			\forall\,e\inn \be_{\al(x)}\!(\fE)\big\}\cap U'
			\quad 
			\forall~\fE\inn\Xi(\tau_{\al(x)})\,.
		\end{align}
	\end{itemize}

We call $(f_{U,t}:U\!\lra\!\fM)_{t\in I_U}$ a  \ts{strongly admissible family of morphisms} if it  satisfies all the conditions of an admissible family of morphisms, except (\ref{Eqn:f_U_assump''}) is replaced by
\begin{align}\label{Eqn:f_U_assump}
	f_{U,t_1}^*\big(\ze_e^{\cV_{1}}\big)
	\big/ f_{U,t_2}^*\big(\ze_e^{\cV_{2}}\big)
	\in 
	\Ga(\sO^*_{U'})
	\qquad
	\forall~e\in S_{\al(x)}\,.	
\end{align}
\end{dfn}

We remark that (\ref{Eqn:f_U_assump}) implies (\ref{Eqn:f_U_assump''}),
hence a strongly admissible family of morphisms is always admissible.

Given an admissible family of morphisms $(f_{U,t}:U\!\lra\!\fM)_{t\in I_U}$, we write 
\begin{align*}
	U_{t}^{\tf}:=
	U\times_{f_{U,t};\,\fM}
	\fM^{\tf}
	\qquad\tn{and}\qquad
	U_{t}^{\rtf}:=
	U\times_{f_{U,t};\,\fM}
	\fM^{\rtf}\,,
	\qquad
	t\inn I_U\,,
\end{align*} 
and denote by  $U_t^{\tf}/U$ 
and $U_t^{\rtf}/U$ the corresponding natural morphisms.

\begin{lmm}
	\label{Lm:local_tf_isom'}
Let $\fM$ and $\La$ be as in Proposition~\ref{Prp:univ},
$U$ be a stack, and $(f_{U,t}:U\!\lra\!\fM)_{t\in I_U}$ be an admissible family of morphisms.
Then, for every $t_1,t_2\inn I_U$,
there  exist isomorphisms
\begin{align*}
	\ov\phi_{t_2,t_1}:
	U_{t_1}^{\tf}/U
	\;\stackrel{\sim}{\lra}\;
	U_{t_2}^{\tf}/U
	\qquad\tn{and}\qquad
	\ud\phi_{t_2,t_1}:
	U_{t_1}^{\rtf}/U
	\;\stackrel{\sim}{\lra}\;
	U_{t_2}^{\rtf}/U\,.
\end{align*}
\end{lmm}

\begin{proof}
We begin with the RL case.

First, we fix arbitrary $x\inn U$, $t_1,t_2\inn I_U$,
and $\cV_1,\cV_2\inn\fV_{\al(x)}$ with $f_{U,t_i}(x)\inn\cV_i$, $i\eq 1,2$.
We write $\al(x)$ and $f_{U,t_i}$ respectively as $\al$ and $f_i$ for conciseness.
We also fix a neighborhood $U'$ of $x$ satisfying $U'\!\subset\!f_1^{-1}(\cV_1)\!\cap\!f_2^{-1}(\cV_2)$.
The restrictions to $U'$ of $U_{t_i}^{\tf}/U$, $i\eq 1,2$, are respectively denoted by 
\begin{align*}
	\ov{\tn p}_i:(U')^{\tf}_1\lra U',\qquad
	i\eq 1,2.
\end{align*}
For $i\eq 1,2$ let $\tn N_i$ be the union of all the irreducible components of $U'$ whose images under $f_i$ are contained in $	 \bigcup_{\fE\in\Xi(\tau_\al)}	X^{\cV_i}_{\be_\al(\fE)}$.
By (\ref{Eqn:f_U_assump''}), we have
$\tn N_1\eq \tn N_2$, which we denote by $\tn N$.
The union of the remaining irreducible components of $U'$ is denoted by $\tn M$.
Then, $U'\eq \tn M\!\cup\!\tn N$.
We set $${\tn M}_i^{\tf}=\ov{\tn p}_i^{\,-1}(\tn M)\qquad\tn{and}\qquad{\tn N}_i^{\tf}\eq\ov{\tn p}_i^{\,-1}(\tn N).$$

The assumption (\ref{Eqn:f_U_assump''}) implies for every $\fE\inn\Xi(\tau_\al)$,
\begin{align*}
	f_{1}^{-1}\big(\ov X^{\cV_{1}}_{\be_\al(\fE)}\big)
	&=
	\bigcup_{\fE'\in\Xi(\tau_\al)~\tn{s.t.}\;\fE'\succeq\fE}\hspace{-.35in}
	\big\{\,
	f_{1}^*\big(\ze_e^{\cV_{1}}\big)\eq 0\ \ 
	\forall\ 
	e\inn\be_\al(\fE')\,
	\big\}\\
	&=
	\bigcup_{\fE'\in\Xi(\tau_\al)~\tn{s.t.}\;\fE'\succeq\fE}\hspace{-.35in}
	\big\{\,
	f_{2}^*\big(\ze_e^{\cV_{2}}\big)\eq  0\ \ 
	\forall\ 
	e\inn\be_\al(\fE')\,
	\big\}
	=
	f_{2}^{-1}\big(\ov X^{\cV_{2}}_{\be_\al(\fE)}\big)\qquad(\subset U')\,.
\end{align*}
Therefore, we have
\begin{align*}
	(f_{1}\!\circ\!\ov{\tn p}_2)^{-1}\big(\ov X^{\cV_{1}}_{\be_\al(\fE)}\big)
	=
	(f_{2}\!\circ\!\ov{\tn p}_2)^{-1}\big(\ov X^{\cV_{2}}_{\be_\al(\fE)}\big)\qquad
	\big(\,\subset (U')^{\tf}_2\,\big)\,.
\end{align*}

Since $f_{2}\!\circ\!\ov{\tn p}_2$ factors through  $\cV_{2}^{\tf}\!\lra\!\cV_{2}$,
Proposition~\ref{Prp:univ} implies 
the restriction to ${\tn M}_{2}^{\tf}$  of the right-hand side of the above identity is a Cartier divisor of ${\tn M}_{2}^{\tf}$,
and so is the the restriction to ${\tn M}_{2}^{\tf}$  of left-hand side.
Once again, by Proposition~\ref{Prp:univ},
we obtain a morphism 
\begin{align*}
	\ov f_{1,2}:\,
	{\tn M}_{2}^{\tf}
	\lra \cV_{1}^{\tf}
\end{align*}
lifting $f_{1}\!\circ\!\ov{\tn p}_2$.
Moreover,
mimicking the proof of \cite[Lemma 5.2]{HLN},
we see (\ref{Eqn:f_U_assump''}) implies that $\ov f_{1,2}$ extends to a morphism
\begin{align*}
	\ov f_{1,2}:\,
	(U')_{2}^{\tf}
	\lra \cV_{1}^{\tf}\,,
\end{align*}
whose restriction to ${\tn N}_{2}^{\tf}$ factors through an isomorphism $\ov\phi_{1,2;\tn N}:{\tn N}_{2}^{\tf}\!\lra\!{\tn N}_{1}^{\tf}$.

The universality of pullbacks then induces a morphism 
\begin{align*}
	\ov \phi_{1,2}:\,
	(U')_{2}^{\tf}
	\lra 
	(U')_{1}^{\tf}
\end{align*}
that $\ov f_{1,2}$ and $\ov{\tn p}_2$ respectively factor through.
By symmetry,
we obtain 
\begin{align*}
	\ov \phi_{2,1}:\,
	(U')_{1}^{\tf}
	\lra 
	(U')_{2}^{\tf}
\end{align*}
analogously.
The above constructions are summarized in the commutative diagram below.
\begin{center}
	\begin{tikzpicture}{h}
		\draw 
		(2.5,1.2) node {$(U')_{1}^{\tf}\eq U'\!\times_{f_{1},\cV_{1}}\!\cV_{1}^{\tf}$}
		(5.5,1.2) node {$\cV_{1}^{\tf}$}
		(-2.5,1.2) node {$(U')_{2}^{\tf}\eq U'\!\times_{f_{2},\cV_{2}}\!\cV_{2}^{\tf}$}
		(-5.5,1.2) node {$\cV_{2}^{\tf}$}
		(0,0) node {$U'$}
		(5.5,0) node {$\cV_{1}$}
		(-5.5,0) node {$\cV_{2}$}
		(2.7,.2) node {\scriptsize{$f_{1}$}}
		(-2.7,.2) node {\scriptsize{$f_{2}$}}
		(-.9,.55) node {\scriptsize{$\ov{\tn p}_2$}}
		(.96,.52) node {\scriptsize{$\ov{\tn p}_1$}}
		(0,1.45) node {\tiny{$\ov\phi_{1,2}$}}
		(0,.9) node {\tiny{$\ov\phi_{2,1}$}}
		(2.25,1.93) node {\tiny{$\ov f_{1,2}$}};
		\draw[->,>=stealth'] 
		(4.28,1.2)--(5.15,1.2);
		\draw[->,>=stealth'] 
		(-4.28,1.2)--(-5.15,1.2);
		\draw[->,>=stealth'] 
		(.33,0)--(5.15,0);
		\draw[->,>=stealth']
		(-.35,0)--(-5.15,0);
		\draw[->,>=stealth'] 
		(5.4,.9)--(5.4,.3);
		\draw[->,>=stealth'] 
		(-5.6,.9)--(-5.6,.3);
		\draw[->,>=stealth'] 
		(.93,.9)--(.23,.23);
		\draw[->,>=stealth'] 
		(-.9,.9)--(-.23,.23);
		\draw[dashed,->,>=stealth'] 
		(-.71,1.25)--(.71,1.25);
		\draw[dashed,->,>=stealth'] 
		(.71,1.15)--(-.71,1.15);
		\draw[dashed,->,>=stealth'] 
		(-.85,1.55) .. controls (.15,2.4) and (4,2.4) .. (5.2,1.5);
	\end{tikzpicture}		
\end{center}

Notice that the composition $\ov \phi_{2,1}\!\circ\!\ov \phi_{1,2}$ restricts to the identity morphism on 
$${\tn M}_2^{\tf}\cap\Big((f_2\!\circ\!\ov{\tn p}_2)^{-1}\big(\cV_2\big\bsl\!\bigcup_{\fE\in\Xi(\tau_\al)}\!\!\!\!X^{\cV_2}_{\be_\al(\fE)}\big)\Big),
$$
which is dense and open in ${\tn M}_2^{\tf}$,
hence $\ov \phi_{1,2}|_{{\tn M}_2^{\tf}}$ and $\ov \phi_{2,1}|_{{\tn M}_1^{\tf}}$ are inverse to each other.
As $\ov \phi_{1,2}|_{{\tn N}_2^{\tf}}$ and $\ov \phi_{2,1}|_{{\tn N}_1^{\tf}}$ are respectively $\ov\phi_{1,2;\tn N}$ and its inverse,
we conclude that $\ov \phi_{1,2}$ and $\ov \phi_{2,1}$ are inverse to each other.
Since $\{U'\}$ covers $U$,
the RL case of the corollary is established.

The LR case is parallel,
hence is omitted.
\end{proof}

Next, with $\fM$ and $\La$ still as in Proposition~\ref{Prp:univ}, we assume 
$\fN$ is a stack equipped with an open covering $\fU$ such that for every $U\inn\fU$ there exists an admissible family of morphisms $(f_{U,t}:U\!\lra\!\fM)_{t\in I_U}$.
In addition,
assume for every $U,\wh U\inn\fU$ with $U\!\cap\!\wh U\!\ne\!\emptyset$,
we have
\begin{align}\label{Eqn:f_U_assump2}
	U\!\cap\!\wh U\in \fU\,,\qquad
	I_U\!\cap\! I_{U\cap \wh U}\!\ne\!\emptyset\,,\qquad\tn{and}
	\qquad
	f_{U,t}|_{U\cap\wh U}=
	f_{U\cap\wh U,t}\quad \forall~t\inn I_U\!\cap\! I_{U\cap \wh U}\,.
\end{align} 

\begin{crl}
	\label{Crl:local_tf_isom}
Under the above setting,
the natural morphisms
\begin{align*}
	U^{\tf}_t/U\quad
	\big(\,\tn{resp.}\ \,U^{\rtf}_t/U\,\big),\qquad
	U\in\fU,~
	t\in I_U,
\end{align*}
glue to form a  morphism
$\ov\varpi_\fN\!:\fN^{\tf}\!\lra\!\fN$  (resp.~$\ud\varpi_\fN\!:\fN^{\rtf}\!\lra\!\fN$) of a stack $\fN^{\tf}$ (resp.~$\fN^{\rtf}$) such that for every
$U\inn\fU$ and $t\inn I_U$,
$U^{\tf}_t/U$ (resp.~$U^{\rtf}_t/U$) is 
 isomorphic to the natural morphism $U\!\times_{\fN}\!\fN^{\tf}\!\lra\!U$ 
(resp.~$U\!\times_{\fN}\!\fN^{\rtf}\!\lra\!U$),
where $U\!\lra\!\fN$ is the inclusion.

Moreover,
the open covering of $\fN^{\tf}$ (resp.~$\fN^{\rtf}$) consisting of $U^{\tf}_t$ (resp.~$U^{\rtf}_t$) for all $U\inn\fU$ and $t\inn I_U$,
along with the index set
$I_{U_t^{\tf}}\!:=\!I_U$ (resp.~$I_{U_t^{\rtf}}\!:=\!I_U$)
and the natural morphisms 
$$
f^{\tf}_{U,t}:
U_t^{\tf}\lra\fM
^{\tf}\qquad
\big(\,\tn{resp}.~f^{\rtf}_{U,t}:U_t^{\rtf}\lra\fM^{\rtf}\,
\big),$$
satisfies (\ref{Eqn:f_U_assump2}).
\end{crl}

\begin{proof}
Below, we show the  statement of Corollary~\ref{Crl:local_tf_isom} in the RL case. 
The proof of the RL case is parallel, hence is omitted.

By (\ref{Eqn:f_U_assump2}),
for every $U,\wh U\inn\fU$ with $U\!\cap\!\wh U\!\ne\!\emptyset$, 
there exist $t\inn I_U\!\cap\!I_{U\cap\wh U}$ and $\wh t\inn I_U\!\cap\!I_{U\cap\wh U}$.
Then, (\ref{Eqn:f_U_assump2}) and Lemma~\ref{Lm:local_tf_isom'} together imply
\begin{align*}
	U_t^{\tf}\big/(U\cap \wh U)=(U\cap \wh U)_t^{\tf}\big/(U\cap \wh U)
	\simeq 
	(U\cap \wh U)_{\wh t}^{\tf}\big/(U\cap \wh U)
	=\wh U_t^{\tf}\big/(U\cap \wh U).
\end{align*}
Therefore, all these $U_t^{\tf}/U$ can be patched up to form  the desired stack $\fN^{\tf}$ and the morphism $\fN^{\tf}\!\lra\!{\fN}$.
The verification of the open covering consisting of all $U^{\tf}_t$ satisfies (\ref{Eqn:f_U_assump2}) is straightforward.
\end{proof}

In \S\ref{Sec:genus_2_twisted_fields},
we shall construct (\ref{base blowups}) for $g\eq 2$ by adding twisted fields to $\fD_2$ successively in nine steps.
Notice  for every small open  $U\!\subset\!\ov M_2(\P^n,d)$,
the morphisms (\ref{Eqn:Mg_tauto}) is indeed strongly admissible.
This means
once the first step $\fD_{2;\fk 1}^{\tf}/\fD_2$ is constructed,
we immediately have a family of new morphisms $U_H^{\tf}\!\lra\!\fD_2^{\tf}$, $H\inn\bH_U$.
However, in order to apply Lemma~\ref{Lm:local_tf_isom'} and Corollary~\ref{Crl:local_tf_isom} recursively,
we expect this new family to be strongly admissible as well.
This leads us to consider the following statement on strongly admissible families of morphisms, which still applies to (\ref{Eqn:Mg_tauto}).

\begin{crl}
	\label{Crl:f_U_induction}
Let $\fM$ and $\La$ be as in Proposition~\ref{Prp:univ},
$U$ be a stack, and $(f_{U,t}:U\!\lra\!\fM)_{t\in I_U}$ be a strongly admissible family of morphisms.
Then, with
$f^{\tf}_{U,t}, f^{\rtf}_{U,t}$ as in Corollary~\ref{Crl:local_tf_isom} and $\ov\phi_{t_2,t_1}$ and $\ud\phi_{t_2,t_1}$ as in Lemma~\ref{Lm:local_tf_isom'},
for every $t_0\inn I_U$, the families
\begin{align*}
	f^{\tf}_{U,t}\circ\ov\phi_{t,t_0}:
	U^{\tf}_{t_0}\lra \fM^{\tf}\qquad
	\tn{and}\qquad
	f^{\rtf}_{U,t}\circ\ud\phi_{t,t_0}:
	U^{\rtf}_{t_0}\lra \fM^{\rtf},\qquad
	t\in I_U,
\end{align*}
are both strongly admissible families of morphisms.
\end{crl}

\begin{proof}
In the RL setting, we justify the analogue of (\ref{Eqn:f_U_assump}) for the family $\big(f^{\tf}_{U,t}\!\circ\!\ov\phi_{t,t_0}\big)$.
To this end, it suffices to show for every $y\inn U^{\tf}_{t_0}$ over every $x\inn U$,
every $\ov\lt\inn\ov\bT_{\al(x)}(\La)$,
every $\cU$ and $\cU_0\inn\cV_{(\al(x),\ov\lt)}$ satisfying $f_{U,t}^{\tf}(y)\inn\cU$ and $f_{U_0,t}^{\tf}(y)\inn\cU_0$, 
and every $s\inn S_{(\al(x),\ov\lt)}$,
there exists a small neighborhood $U'$ of $y$ in $U^{\tf}_{t_0}$ such that
\begin{align}\label{Eqn:f_U_assump'}
	{(f_{U,t_0}^{\tf})^*\big(\xi_s^{\cU_{0}}\big)}
	\big/
	{(f_{U,t}^{\tf}\!\circ\!\ov\phi_{t,t_0})^*\big(\xi_s^{\cU}\big)}
	\in\Ga\big(\sO^*_{U'}\big)\,.
\end{align}
Notice by Theorem~\ref{Thm:tf_smooth}~\ref{Cond:smooth_parameters},
each twisted parameter can be written as the quotient of the products of the pullbacks of certain original modular parameters.
Therefore,
(\ref{Eqn:f_U_assump'}) follows from (\ref{Eqn:f_U_assump}) and the fact that each $\ov\varpi\!\circ\!f_{U,t}^{\tf}$ factors through $f_{U,t}$.

The verification of the remaining conditions as per Definition~\ref{Dfn:weak_and_strong_adm}, as well as the proof of the LR case, are straightforward and hence are omitted.
\end{proof}

\section{Application to the moduli of genus 1 stable maps: proof of Theorem~\ref{Thm:g1}}\label{Sec:Proof_g1}
Some arguments and calculations are parallel to the corresponding ones in \cite{HL10} and \cite{g1modular}. We provide a sketch.

\ref{Cond:g1_tf_wt} and \ref{Cond:g1_tf_stable_map} are respectively the restatements of Theorem~1.1 and Corollary~1.2 of~\cite{g1modular}.

\ref{Cond:g1_rtf_wt}. Let $\fM$ be the connected component of
$\fM_1^\wt$ with total weight $d$. 
Since the $d\eq 1$ case is trivial,
we assume $d\!>\!1$ and denote by $$\widetilde \fM \lra \fM$$ the iterated blowup of $\fM$ along the proper transforms of $\Theta_d, \cdots, \Theta_2$.
In the proof of Proposition \ref{Prp:RL/LR-cptb},
on any modular chart $\cV$ of $\fM$,
with the modular parameters $\ze_e$ corresponding to the smoothing of the nodes (c.f.~Example~\ref{Eg:G-adim_fixture}),
and with $Z_k$ as in (\ref{Eqn:Z_k}),
we observe that
$$
Z_k=\Th_{d-k+1}\cap\cV,\qquad
1\le k\le d\!-\!1.
$$
As shown in the proof of Proposition \ref{Prp:RL/LR-cptb},
the restriction $\ti\fM/\cV$ is thus isomorphic to $\fM^{\rtf}/\cV$.
Moreover, 
for any other modular chart $\cV'$ of $\fM$,
the isomorphisms
$$
\ti\fM/\cV\lra \fM^{\rtf}/\cV
\qquad\tn{and}\qquad
\ti\fM/\cV'\lra \fM^{\rtf}/\cV'
$$
agree on $\cV\!\cap\!\cV'$,
for they agree on the dense open subset $\cV\!\cap\!\cV'\!\cap\!\fM^\mn$.
Therefore,
the local isomorphisms glue and give rise to a global isomorphism of $\ti\fM/\fM$ to $\fM^{\rtf}/\fM$.
This proves the first statement of~\ref{Cond:g1_rtf_wt}.

For the last statement,
we assume the total weight $d$ of the component $\fM$ satisfies $d\!\ge\!3$.
With notation as in Example~\ref{Eg:genus_1},
we choose $\al\eq(\ga,p_g,\bfw)\inn A$ such that
the dual graph $\ga$ is given by the leftmost graph of Figure~\ref{Fig:Treelike},
the geometric genus function $p_g$ assigns  1 to the vertex $o$ and 0 to the rest vertices,
and the weight function $\bfw$ assigns positive weights to the leaves (as vertices) and 0 to the rest.
For such $\al$, the rooted tree $\tau_\al$ is the same as $\ga$.
The last statement then follows directly from Corollary~\ref{Crl:tf_rtf_non_isom}.

\ref{Cond:g1_rtf_stable_map}.  
Once again we fix an arbitrary connected component $\fM$ of $\fM_1^\wt$ with total weight $d$, and arbitrary $\al\inn A$, $\cV\inn\fV_\al$, $\ud\lt\inn\ud\bT_\al$, and $\cU\inn\fV_{(\al,\ud\lt)}$ with $\ud\varpi(\cU)\!\subset\!\cV$.
We list the leaves of $\tau_\al$ as follows:
$$
\min(\tau_\al)=\{\fe_1,\cdots,\fe_\ell\}\qquad
\tn{for~some}~\ell\!\le\!d.
$$
As shown in~\cite[\S5.2]{HL10}, locally $\ov M_1(\P^n,d)$ is defined by the following equations on the total space $\cE_\cV$ of $\sO_\cV^{\oplus n(d+1)}$:
\begin{align*}
\big[\,\prod_{\fe\succeq\fe_1}\!\ze_\fe\,,\,
\cdots,\,\prod_{\fe\succeq\fe_\ell}\!\ze_\fe\,,\,0,\,
\cdots,\,0\,\big]\cdot
[\,w^i_1,\,\cdots,\,w_\ell^i,\,w_{\ell+1}^i,\cdots,\,w_{d+1}^i\,]^T=0,\qquad
1\!\le\!i\!\le\!n,
\end{align*}
where $w_j^i$'s are free variables on $\cE_\cV$ that, along with the modular parameters on $\cV$, form a subset of a system of local parameters on $\cE_\cV$.
On $\cU$,
the pullbacks of $\prod_{\fe\succeq\fe_k}\!\ze_\fe$, $1\!\le\!k\!\le\!\ell$, are explicitly given in Corollary \ref{Crl:prod_LR}.
Hence by the last statement of Corollary \ref{Crl:prod_LR},
under suitable trivialization of $\sO_\cU^{\oplus n(d+1)}$,
we see that locally, $\ov M_1^{\rtf}(\P^n,d)$ is defined by the following equations on the total space $\cE_\cU$ of $\sO_\cU^{\oplus n(d+1)}$:
\begin{align}\label{Eqn:g1_local_eqn_resolved}
 \prod_{\fE\in\ud\bE}\xi_\fE^\cU\,\cdot \,\ti w_1^i=0,\qquad 1\!\le\!i\!\le\!n,
\end{align}
where $\ti w_1^i$'s are free variables on $\cE_\cU$ that, along with the modular parameters on $\cU$, form a subset of a system of local parameters on $\cE_\cU$.
Hence on  $\cE_\cU$,
$\ov M_1^{\rtf}(\P^n,d)^{\mc}$ is given by
\begin{align}
\label{Eqn:g1_local_main}
\ti w_1^1=\cdots=\ti w_1^n=0,
\end{align}
whereas any other irreducible components of $\ov M_1^{\rtf}(\P^n,d)$ is given by
\begin{align}
\label{Eqn:g1_local_boundary}
\xi_\fE^\cU=0\qquad\tn{for~some}\ \fE\in\ud\bE\,.
\end{align}
This establishes Condition~\ref{Cond:MainSmooth} for $\ov M_1^{\rtf}(\P^n,d)$.

By~(\ref{Eqn:g1_local_main}) and~(\ref{Eqn:g1_local_boundary}),
on $\cE_\cU$, the boundary of
$\ov M_1^{\rtf}(\P^n,d)^{\mc}$ is given by
$$
\prod_{\fE\in\ud\bE}\!\xi_\fE^\cU=\ti w_1^1=\cdots=\ti w_1^n=0\,.
$$
This, along with~(\ref{Eqn:g1_local_main}), gives rise to Condition~\ref{Cond:MainMainComponent}.

Condition~\ref{Cond:MainBirational} for $\ov M_1^{\rtf}(\P^n,d)$ and the properness of $\ud\varpi$ both follow from Theorem~\ref{Thm:tf_smooth_revert}~\ref{Cond:smooth_tf_revert}.

For $k\!\ge\!1$,
by~\cite[Theorem 4.16]{HL10},
the direct image sheaf $\pi_*\ff^*\sO_{\P^n}(k)$ locally is the kernel of the homomorphism $\varphi:\sO_\cV^{\oplus kd+1}\!\lra\!\sO_\cV$:
\begin{align*}
\varphi=\big[\,\prod_{\fe\succeq\fe_1}\!\ze_\fe\,,\,
\cdots,\,\prod_{\fe\succeq\fe_\ell}\!\ze_\fe\,,\,0,\,
\cdots,\,0\,\big]\,.
\end{align*}
Once again by the last statement of Corollary \ref{Crl:prod_LR},
under suitable trivialization of $\sO_\cU^{\oplus kd+1}$,
we see  the pullback of $\varphi$ locally takes the form
$$
\big(\prod_{\fE\in\ud\bE}\xi_\fE^\cU\,\big)\cdot \,[1,0,\cdots,0]\,.
$$
As shown in the proof of~\cite[Theorem 2.11]{HL10} (c.f.~the sentence containing \cite[(5.22)]{HL10}),
the direct image sheaf 
$\ti\pi_*\ti\ff^*\sO_{\P^n}(k)$ locally takes the same form,
so it is locally free and of rank $kd$ on the main component $\ov M_1^{\rtf}(\P^n,d)^\mc$, and is locally free and of rank $kd\!+\!1$ on any other irreducible component.
This establishes Condition~\ref{Cond:MainLocallyFree} for $\ov M_1^{\rtf}(\P^n,d)$.

In sum, the first statement of Part~\ref{Cond:g1_rtf_stable_map} is proved.

The last statement of Part~\ref{Cond:g1_rtf_stable_map} follows from a modified version the proof of Corollary~\ref{Crl:tf_rtf_non_isom}.
We assume $d\!\ge\!3$.
Let
$$
p:\,\ov M_1(\P^n,d)\lra\fM_1^\wt,\qquad
\ov p:\,\ov M_1^{\tf}(\P^n,d)\lra\fM_1^\wt,\qquad\tn{and}\qquad
\ud p:\,\ov M_1^{\rtf}(\P^n,d)\lra\fM_1^\wt
$$
be the morphism (\ref{Eqn:M1wt_tauto}), 
the composition of $p$ and the morphism $\ov M_1^{\tf}(\P^n,d)\!\lra\!\ov M_1(\P^n,d)$ in \ref{Cond:g1_tf_stable_map},
and the composition of $p$ and the morphism $\ov M_1^{\rtf}(\P^n,d)\!\lra\!\ov M_1(\P^n,d)$ in \ref{Cond:g1_rtf_stable_map}, respectively.

We observe that the $RL$ counterpart of (\ref{Eqn:g1_local_eqn_resolved}) is also true, i.e.~locally, $\ov M_1^{\tf}(\P^n,d)$ always takes the form:
\begin{align*}
	\prod_{\fE\in\ov\bE}\xi_\fE^\cU\,\cdot \,\ti w_1^i=0,\qquad 1\!\le\!i\!\le\!n.
\end{align*}
So, $\ov M_1^{\tf}(\P^n,d)^\mc$ locally takes the same form as in (\ref{Eqn:g1_local_main}),
and $\ov M_1^{\tf}(\P^n,d)^\mc\!\cap\!\ov p^{-1}(\Th_2)$ is a Cartier divisor of $\ov M_1^{\tf}(\P^n,d)^\mc$ locally given by
$$
\xi_\fE^\cU=\ti w_1^1=\cdots=\ti w_1^n=0,
$$
where $\fE\eq\max(\tau_\al)$ contains exactly two edges.
Hence if there was a morphism of
\begin{align}\label{Eqn:tf->rtf_no_morph}
	\ov M^{\rtf}_1(\P^n,d)^\mc\big/\ov M_1(\P^n,d)^\mc\lra\ov M^{\tf}_1(\P^n,d)^\mc\big/\ov M_1(\P^n,d)^\mc,
\end{align}
then $\ov M_1^{\tf}(\P^n,d)^\mc\!\cap\!\ov p^{-1}(\Th_2)$ would pull back to $\ov M_1^{\rtf}(\P^n,d)^\mc\!\cap\!\ud p^{-1}(\Th_2)$,
which in turn would be a Cartier divisor of $\ov M_1^{\rtf}(\P^n,d)^\mc$.

To see this is impossible,
we continue with the setting in the proof of the last statement of Part~\ref{Cond:g1_rtf_wt}
so that $\al\eq(\ga,p_g,\bfw)\inn A$ remains the same.
On an arbitrary modular chart $\cV\inn\fV_\al$, we have
\begin{align*}
	&\fM_\al\cap\cV=\{\ze_a\eq\ze_b\eq\ze_c\eq\ze_d\eq 0\},\qquad
	\Th_2\cap\cV=\{\ze_a\eq\ze_d\eq 0\},\\
	&\Th_3\cap\cV=\{\ze_b\eq\ze_c\eq\ze_d\eq 0\},\qquad\tn{and}\qquad
	\Th_k\cap\cV\eq\emptyset\quad\forall~k\!>\!3.
\end{align*}
Consider 
$$
\ud\lt=\big(\tau_\al,\ud\bE,\Dm(\ud\lt)\big)\in\ud\bT_\al\qquad\tn{with}\qquad
\ud\bE=\big\{\fE\eq\{b,c,d\}\big\},\quad
\Dm(\ud\lt)=\{d\}.
$$
On a twisted chart $\cU\inn\fV_{(\al,\ud\lt)}$ with $\ud\varpi(\cU)\!\subset\!\cV$,
by (\ref{Eqn:mod_par_pullback_revert}), we have
\begin{align}\label{Eqn:Th2_local}
	\big(\ud\varpi^{-1}(\Th_2\!\cap\!\cV)\big)\cap\cU
	=\big\{\xi^\cU_\fE\eq\xi^\cU_a\eq 0\big\}.
\end{align}
Then,
$$
\wh \Th_2\;:=\ 
p^{-1}(\Th_2\!\cap\!\cV)\,\times_{(\Th_2\cap\cV)}
\,\big(\cU\cap\ud\varpi^{-1}(\Th_2\!\cap\!\cV)\big)
$$
is an open subset of
$$
p^{-1}(\Th_2\!\cap\!\cV)\times_{(\Th_2\cap\cV)}
\ud\varpi^{-1}(\Th_2\!\cap\!\cV)\,=\,
\ud p^{-1}(\Th_2\!\cap\!\cV)\,.
$$
By (\ref{Eqn:g1_local_eqn_resolved}) and (\ref{Eqn:Th2_local}),
on $\cE_\cU$,
$\ov M_1^{\rtf}(\P^n,d)^\mc\!\cap\!\wh\Th_2$ is then given by
$$
\xi_\fE^\cU=\xi_a^\cU=\ti w_1^1=\cdots=\ti w_1^n=0\,.
$$
Taking~(\ref{Eqn:g1_local_main}) into consideration, we see $\ov M_1^{\rtf}(\P^n,d)^\mc\!\cap\!\ud p^{-1}(\Th_2)$ is not Cartier in $\ov M_1^{\rtf}(\P^n,d)^\mc$,
hence the morphism (\ref{Eqn:tf->rtf_no_morph})
does not exist.

A parallel argument also implies the morphism in (\ref{Eqn:tf->rtf_no_morph}) with the arrow reversed does not exist.
To summarize, the last statement of~\ref{Cond:g1_rtf_stable_map} is established.

\section{Operations on treelike structures}
\label{Sec:Induced_and_derived_tf_stacks}

Starting from $\fD_2$,  progressively,  we will apply Theorems~\ref{Thm:tf_smooth} and \ref{Thm:tf_smooth_revert}
to obtain $\widetilde\fD_2^{\rm tf}$ as motivated in \eqref{base blowups} and stated in Theorem \ref{Thm:Main}.
The choices of the treelike structures that will be added to $\fD_2$ successively are guided
by the local equations of $\ov M_2(\P^n,d)$. 
In~\S\ref{Subsec:M2Pnd_loc_eqn},
we provide a brief review of the local equations of $\ov M_2(\P^n,d)$.
Based on these equations,
we develop certain operations on stacks with LESs and treelike structures in 
\S\ref{Subsec:two_treelike}-\ref{Subsec:Grafted}.

\subsection{Local defining equations of the moduli of genus 2 stable maps: a review}
\label{Subsec:M2Pnd_loc_eqn}
The statements of this subsection are summarized from~\cite{HLN}.

Recall  that $\mdv_2$ is the smooth algebraic stack 
of the  pre-stable pairs $(C,D)$,
where~$C$ are connected genus $2$ nodal curves and $D$ are sums of distinct smooth points on $C$
with multiplicity one. 
Here,
a pair $(C,D)$ in~$\mdv_2$ is said to be \ts{pre-stable} if every rational irreducible component without any divisorial marking contains at least three nodal points.
The connected components $\mdv_2(m)$ of $\mdv_2$ are determined by the degrees $m$ of the divisors.
	
Fix an arbitrary $m\inn \Z_{>0}$.
Consider the algebraic stack $\wh\fD_2(m)$ of the marked nodal curves $(C;\de_1,\ldots,\de_m)$ where $C$ are genus $2$ nodal curves and $\de_i$'s are smooth points on $C$,
satisfying {\it the same pre-stability condition as} $\fD_2(m)$.	
The morphism
\begin{align}\label{Eqn:wh_D}
\wh\fD_2(m)\lra\fD_2(m),\qquad
(C;\,\de_1,\cdots,\de_m)\mapsto
(C, \de_1\!+\!\cdots\!+\!\de_m)
\end{align}
is \'etale, 
hence every {\it sufficiently small} chart $\cV\!\lra\!\fD_2(m)$ can be considered as a chart $\cV\!\lra\!\wh\fD_2(m)$,
for  given a point $(C,D)\inn\cV$, 
one can assign an order to the points of $D$ and extend it to the entire~$\cV$.

\begin{cnv}
	\label{Cnv:chart}
Hereafter, every (small) chart $\cV\!\lra\!\fD_2(m)$ is assumed to factor through (\ref{Eqn:wh_D}).
\end{cnv}

To endow treelike structures to $\fD_2$,
we need two types of local parameters:
those for node-smoothing and those for loci of conjugate points. The former are standard,
so we view the latter.

\subsubsection{Loci of conjugate points}
\label{Subsubsec:conj}

Recall that two distinct points $\de_1$ and $\de_2$ on a smooth genus 2 curve $C$ are said to be \ts{conjugate} if $h^0\big(C,\sO_C(\de_1\!+\!\de_2)\big)\eq 2$.
For distinct $1\!\le\!i,j\!\le\!m$,
we denote by $\cK_{m;i,j}\!\subset\!\cM_{2,m}$ the locus of the smooth genus 2 curves with $m$ marked points whose $i$-th and $j$-th marked points are conjugate, and by $\ov\cK_{m;i,j}\!\subset\!\ov\cM_{2,m}$ the closure of $\cK_{m;i,j}$ in $\ov\cM_{2,m}$.
By~\cite[Lemma~2.14]{HLN}, $\ov\cK_{m;i,j}$ is a Cartier divisor of $\ov\cM_{2,m}$.
We take $\wh\cK_{m;i,j}\!\subset\!\wh\fD_2(m)$
to the pullback of $\ov\cK_{m;i,j}$ under the stabilization 
\begin{align*}
	f_{\rm st} :\, 
	\wh\fD_2(m)\lra\ov\cM_{2,m}\,,
\end{align*}
which is in turn a Cartier divisor of $\wh\fD_2(m)$.
For every $x\eq(C;\de_1,\ldots,\de_m)\inn\wh\fD_2(m)$ and every small chart $\cV\!\lra\!\wh\fD_2(m)$ containing $x$,
let $\cV'\!\lra\!\ov\cM_{2,m}$ be a chart containing $f_{\rm st}(\cV)$,
and $\ka'_{ij}\inn\Ga(\sO_{\cV'})$ be such that 
$\ov\cK_{m;i,j}\!\cap\!\cV'\eq\{\ka'_{ij}\eq 0\}$.
Then, by writing 
\begin{align}\label{Eqn:ka_ij}
	\ka_{ij} := f_{\rm st}^*\,\ka'_{ij}\,,
\end{align}
we have 
$\wh\cK_{m;i,j}\!\cap\!\cV\eq\{\ka_{ij}\eq 0\}$.

In Lemma~\ref{Lm:kappa} below, 
we show the same expressions of $\ka'_{ij}$ in~\cite[Lemma~2.14]{HLN} apply to $\ka_{ij}$ as well.
To this end,
we recall certain terminology and notation about $\fD_2$ from~\cite{HLN}.

For $x\eq (C;\de_1,\ldots,\de_m)\inn\wh\fD_2(m)$,
let $N(C)$ be the set of the nodes of $C$.
By the deformation theory of nodal
curves,
$x$ has a (sufficiently small) neighborhood $\cV$ such that
for every $e\inn N(C)$, there exists  $\ze_e\inn\Ga(\sO_\cV)$ satisfying $\{\ze_e\eq 0\}$ is the locus
where the node~$e$ is not smoothed.
The set $\{\ze_e\}_{e\in N(C)}$ is a subset of a system of local parameters on $\cV$, known as the set of the \ts{node-smoothing parameters}.

	Next, we introduce some notation related to $(C,D)\inn \fD_2$, in order to describe certain subsets of $N(C)$ that are needed for  $\ka_{ij}$ and for the local equations of $\ov M_2(\P^n,d)$.
Let $C_{\core}$ be the smallest (w.r.t.~inclusion) connected genus-2 subcurve of $C$, known as the \ts{core} of $C$.
We write the set of the minimal (w.r.t.~inclusion) connected (arithmetic) genus-1 subcurves of $C_\core$ as $\score(C)$,
whose cardinal is up to three.
In addition,
if $\score(C)$ has two elements and they either are disjoint or meet at one node,
we say $C_\core$ is \ts{separable};
otherwise, $C_\core$ is \ts{inseparable}.
For every $1\!\le\!i\!\le\!m$,
we denote by $C_i$ the irreducible component of $C$ that contains $\de_i$.

Let $\ga_x$ be the \ts{dual graph} of $x$ in the usual sense (which encodes the positions of the nodes and irreducible components of $C$, the geometric genera of the irreducible components,
as well as the distribution of $\de_i$'s on the irreducible components of $C$;
c.f.~\cite[\S23.4]{MirSym}).
Given a union of connected subcurves $\Si'_1\!\cup\!\cdots\!\cup\! \Si'_\ell$ of $C$ {\it satisfying for every $1\!\le\!h\!<\!j\!\le\!\ell$, $\Si'_h\!\cap\!\Si'_j$ either is empty or consists of one node,}
for every $1\!\le\!j\!\le\!\ell$,
we contract the subgraph of $\ga_x$ corresponding to $\Si'_j$ into a single vertex $v_{\Si'_j}$
(more precisely, we remove all the edges corresponding to the nodes of $\Si'_j$, identify all the vertices corresponding to the irreducible components of $\Si'_j$ as a single vertex $v_{\Si'_j}$, and assign all the marked points on $\Si'_j$ to $v_{\Si'_j}$).
The resulting dual graph is written as $\ga_{x;(\Si'_1,\cdots, \Si'_\ell)}$.
We emphasize whenever $\Si'_h$ and $\Si'_j$ meets at a node,
the corresponding edge is {\it not} contracted in the construction of $\ga_{x;(\Si'_1,\cdots, \Si'_\ell)}$.

We are ready to describe the following subsets of the power set of $N(C)$.

\begin{itemize}[leftmargin=*]
\item 
$P_i(x)\eq\{\wp_i(x)\}$.
For every $1\!\le\!i\!\le\!m$,
we denote by $P_i(x)$ the set of the paths of $\ga_{x;(C_\core)}$ connecting $v_{C_\core}$ and the vertex on which $\de_i$ lies.
Here, a \ts{path} is in the usual sense of graph theory, i.e.~a finite sequence of edges which joins a sequence of {\it distinct} vertices, which is allowed to be empty.
Since $\ga_{x;(C_\core)}$ is a tree,
the cardinal of $P_i(x)$ is always 1.
The unique element of $P_i(x)$ is denoted by $\wp_i(x)$, which can naturally be considered as a subset of $N(C)$.

\item 
$P_{T,i}(x)$.
For every $T\inn \score(C)$ and $1\!\le\!i\!\le\!m$,
we denote by $P_{T,i}(x)$ the set of the paths of $\ga_{x;(T)}$ connecting $v_{T}$ and the vertex on which $\de_i$ lies.
The cardinal of $P_{T,i}(x)$ is either 1 or 2 (should $T$ exist).
Each element of $P_{T,i}(x)$ is also considered as a subset of $N(C)$.

\item 
$\wh P_{T,i}(x)\eq\{\wp_{T,i}(x)\}$.
{\it Assume $C_\core$ is separable.}
Then, $\score(C)$ has exactly two elements,
written as $T_1$ and $T_2$,
and they either meet at a node or are disjoint.
On $\ga_{x;(T_1,T_2)}$, for every $T\inn \score(C)$ and $1\!\le\!i\!\le\!m$,
we denote by $\wh P_{T,i}(x)$ the set of the paths of $\ga_{x;(T_1,T_2)}$ connecting $v_{T}$ and the vertex on which $\de_i$ lies.
The cardinal of $\wh P_{T,i}(x)$ is always 1.
The unique element of $\wh P_{T,i}(x)$ is denoted by $\wp_{T,i}(x)$, which is also considered as a subset of $N(C)$.

\item 
$Q_{i}(x)$.
For every $1\!\le\!i\!\le\!m$,
if there exist $T^{(i)}\inn \score(C)$ such that $|P_{T^{(i)},i}(x)|\eq 2$,
then such $T^{(i)}$ must be unique (for otherwise the arithmetic genus of $C_\core$ would be greater than two).
In this case, we write $Q_{i}(x)\!:=\!P_{T^{(i)},i}(x)$.
If such $T^{(i)}\inn \score(C)$ does not exist,
we simply set $Q_{i}(x)\!:=\!\emptyset$.
\end{itemize}

Intuitively,
the path in $P_i(x)$ consists of all the nodes between $C_\core$ and $C_i$,
while the path in $\wh P_{T,i}(x)$ consists of all the nodes between $T$ and $C_i$ {\it after smoothing out the internal nodes of $T_1$ and $T_2$}.
If we only smooth out the internal nodes of $T$,
then each element of $P_{T,i}(x)$ is a path from $T$ to $C_i$;
if $C_\core$ is separable, then
$P_{T,i}(x)\eq \wh P_{T,i}(x)$ if and only if $|P_{T,i}(x)|\eq 1$.
The set $Q_{i}(x)$ is nonempty whenever  $\de_i$, or the rational tail of $C$ containing $\de_i$ if $\de_i\!\notin\!C_\core$,
is on a {non-separating bridge} in the sense of \cite[\S2.5]{HLN},
so there are two ways of connecting $T$ and~$C_i$.
An example for these sets of paths is provided in Figure~\ref{Fig:four_paths}.

\begin{figure}
	\begin{center}
		\begin{tikzpicture}{htb}
			\def\g1{
				(-1,0) ellipse (1 and 0.5)
				(-1.4,0)..controls(-1,0.1)..(-0.6,0)
				(-0.6,0)..controls(-1,-0.1)..(-1.4,0)
				(-0.5,0.05)--(-0.6,0)
				(-1.4,0)--(-1.5,0.05)
			}
			\def\crs{(.05,.05)--(-.05,-.05) 
				(.05,-.05)--(-.05,.05)}
			\draw[xshift=.6cm]\g1;
			\draw[xshift=-1cm]
			(-.7,0) circle (.3cm)
			(-.7,.6) circle (.3cm)
			(-1.3,0) circle (.3cm)
			(-1.724,.424) circle (.3cm)
			(-1.724,-.424) circle (.3cm)
			(-2.149,0) circle (.3cm)
			(-2.749,0) circle (.3cm)
			(-2.149,.829) circle (.3cm)
			(-2.149,-.829) circle (.3cm);
			\draw[thick,xshift=-3.9cm,yshift=.14cm]\crs;
			\draw[thick,xshift=-3.7cm,yshift=-.2cm]\crs;
			\draw[thick,xshift=-3.1cm,yshift=1cm]\crs;
			\draw[thick,xshift=-3cm,yshift=-.95cm]\crs;
			\draw[thick,xshift=-3.3cm,yshift=-.9cm]\crs;
			\draw[thick,xshift=3.25cm,yshift=1cm]\crs;	
			\draw[thick,xshift=-1.6cm,yshift=.77cm]\crs;		
			
			\draw
			(-4,.45) node {\tiny{$\de_1$}}
			(-3.88,-.4) node {\tiny{$\de_2$}}
			(-3.1,1.3) node {\tiny{$\de_3$}}
			(-3.6,-.9) node {\tiny{$\de_4$}}
			(-2.8,-1.1) node {\tiny{$\de_5$}}
			(-1.5,.98) node {\tiny{$\de_6$}}
			(-3.3,0) node {\tiny{$q_8$}}
			(-2.8,.3) node {\tiny{$q_5$}}
			(-2.8,.55) node {\tiny{$q_7$}}
			(-2.8,-.3) node {\tiny{$q_6$}}
			(-2.8,-.55) node {\tiny{$q_9$}}
			(-2.35,-.15) node {\tiny{$q_4$}}
			(-2.35,.15) node {\tiny{$q_3$}}
			(-1.84,-.05) node {\tiny{$q_2$}}
			(-1.66,.17) node {\tiny{$q_{10}$}}
			(-1.23,-.05) node {\tiny{$q_1$}}
			(-.4,-.7) node {\tiny{$T_2$}}
			(-1,-1.2) node {$x$}
			(3,.6) node[right] {\tiny{$\wp_1(x)=\{q_8\}$}}
			(3,.2) node[right] {\tiny{$P_{T_2,1}(x)=\big\{\{q_1,q_2,q_3,q_5,q_8\},\,\{q_1,q_2,q_4,q_6,q_8\}\big\}$}}
			(3,-.2) node[right] {\tiny{$\wp_{T_2,1}(x)=\{q_1,q_2,q_8\}$}}
			(3,-.6) node[right] {\tiny{$Q_{1}(x)=P_{T_2,1}(x)$}}
			(3,-1) node[right] {\tiny{$Q_{13}(x)=\big\{\{q_3\},\,\{q_4,q_6\}\big\}$}}
			(3.3,1) node[right] {\tiny{$:\,\de_i$}};
		\end{tikzpicture}
	\end{center}
	\caption{$P_i(x)$, $P_{T,i}(x)$, $\wh P_{T,i}(x)$, and $Q_i(x)$}\label{Fig:four_paths}
\end{figure}

\begin{lmm}
	\label{Lm:kappa}
Given $x\inn\wh\cK_{m;i,j}$ and a small neighborhood $\cV$ of $x$, we have the following.
\begin{itemize}[leftmargin=*]
	\item If $Q_{i}(x)\eq\emptyset$, then $Q_{j}(x)\eq\emptyset$ as well.
	In such case,
	$\wh \cK_{m;i,j}$ is smooth at
	$x$,
	and $\{\ka_{ij}\}\!\sqcup\!\{\ze_e\}_{e\in N(C)}$ is a subset of a system of local parameters on $\cV$.
	
	\item If $Q_{i}(x)\!\ne\!\emptyset$, then $Q_{j}(x)\!\ne\!\emptyset$ as well, and $T^{(i)}\eq T^{(j)}\inn\score(C)$.
	In such case, 
	the set
	\begin{align*}
		&Q_{ij}(x):=\\
		&\quad
		\begin{cases}
		\min\Big(\big\{\,
		\wp\bsl \wp_{T^{(i)},i}(x)
		:\,
		\wp\inn
		Q_{i}(x)\,
		\big\}\cup 
		\big\{\,
		\wp\bsl \wp_{T^{(j)},j}(x)
		:\,
		\wp\inn
		Q_{j}(x)\,
		\big\}\Big)
		&
		\tn{if}~C_\core~\tn{is~separable,}
		\\
		\min\Big(\big\{\,
		\wp\bsl \wp_{i}(x)
		:\,
		\wp\inn
		Q_{i}(x)\,
		\big\}\cup 
		\big\{\,
		\wp\bsl \wp_{j}(x)
		:\,
		\wp\inn
		Q_{j}(x)\,
		\big\}\Big)
		&
		\tn{if}~C_\core~\tn{is~inseparable},
		\end{cases}
	\end{align*}
	where the minima are taken with respect to inclusion,
	has exactly two elements;
	moreover,
	there exist invertible functions $u_{\wp}$, $\wp\inn Q_{ij}(x)$, on $\cV$ such that
	\begin{align*}
		\ka_{ij}=
		\sum_{\wp\in Q_{ij}(x)}\!\!
		\Big(u_{\wp}
		\prod_{e\in \wp}
		\ze_{e}\Big)\ .
	\end{align*}
\end{itemize}
\end{lmm}

In Figure~\ref{Fig:four_paths},
for example, we have $Q_{12}(x)\eq\big\{\{q_3,q_5\},\,\{q_4,q_6\}\big\},$
$Q_{13}(x)\eq\big\{\{q_3\},\,\{q_4,q_6\}\big\},$
$Q_{14}(x)\eq Q_{15}(x)\eq\big\{\{q_3,q_5\},\,\{q_4\}\big\},$
etc.

\begin{proof}[Proof of Lemma~\ref{Lm:kappa}]
	Let $J\!\subset\!\{1,\ldots,m\}$ be comprised of all $1\!\le\!j\!\le\!m$ such that $\de_j$ lies on a rational component $R_j$ of $C$, and $R_j$ is contracted in $f_{\rm st}(x)$.
	In other words,
	$R_j$ contains only one nodal point, denoted by $q_j$, and does not contain any marked point except $\de_j$.
	Let $\ti x\eq (C;\de_1,\ldots,\de_m;\de_{m+1},\ldots,\de_{m+\ell})\inn\ov\fM_{2;m+\ell}$
	be such that
	$\de_{m+j}\inn R_j$ for all $1\!\le\!j\!\le\!\ell$.
	Then, the image of $\ti x$ under the forgetful morphism
	\begin{align*}
		\ti f:\,
		\ov\cM_{2,m+\ell}\lra 
		\wh\fD_2(m)
	\end{align*}
	that forgets the last $\ell$ marked points is $x$.	
	Let $\ti\cV\!\lra\!\ov\cM_{2,m+\ell}$ be  a small chart containing $\ti x$ satisfying $\ti f(\ti\cV)\!\subset\!\cV$.
	Then for every $\ti y\inn\ti\cV$,
	the underlying nodal curve of $\ti f(\ti y)\inn\cV$ is identical to that of $\ti y$.
	
	For each node $q\inn C$, notice that $\ti f^*\ze_q$ is a node-smoothing parameter corresponding to $q$ on $\ti\cV$.
	Moreover,
	observe that 
	$f_{\rm st}\!\circ\!\ti f$ is exactly the same as the forgetful morphism $\ov\cM_{2,m+\ell}\!\lra\!\ov\cM_{2,m}$ that forgets the last $\ell$ marked points.
	Since $\ov\cK_{m;i,j}\!\subset\!\ov\cM_{2;m}$ are all defined via pullbacks under the forgetful morphisms,
	we have
	\begin{align*}
		\ov\cK_{m+\ell;i,j}\cap\ti\cV=
		\{(f_{\rm st}\!\circ\!\ti f)^*\ka'_{ij}\eq 0\}=
		\{\ti f^*\ka_{ij}\eq 0\}\,.
	\end{align*}
	Therefore,
	the statements of Lemma~\ref{Lm:kappa} all follow from applying~\cite[Lemma~2.14]{HLN} to the corresponding pullbacks to $\ti\cV$.
\end{proof}

\subsubsection{Local equations of $\ov M_2(\P^n,d)$ and $\mathbf{R}\pi_*\ff^*\sO_{\P^n}(k)$}

We are ready to present the local equations of $\ov M_2(\P^n,d)$ and the derived object $\mathbf{R}\pi_*\ff^*\sO_{\P^n}(k)$ (for a fixed positive integer $k$).
By~\cite[Theorem 2.5]{HLN},
for every $y\inn\ov M_2(\P^n,d)$,
a small neighborhood $U$ of $y$ can always be identified with the  kernel of the $m\eq d$ case of a homomorphism
\begin{align}
	\label{Eqn:str_hom_0}
(0\oplus \varphi)^{\oplus n}:\,
\big(\sO_\cV\,\oplus\,\sO_\cV^{\oplus m}\big)^{\oplus n}\lra \big(\sO_\cV^{\oplus 2}\big)^{\oplus n},
\end{align}
where $\cV\!\lra\!\mdv_2$ is a smooth chart containing the image $x\eq(C,\de_1\!+\!\cdots\!+\!\de_m)$ of $y$ in $\fD_2$ (c.f.~\cite[(2.8)]{HLN}).
Over $U$,
${\bf R}\pi_*\ff^*\sO_{\P^n}(k)$ can be represented by the pull back of 
the same homomorphism as in (\ref{Eqn:str_hom_0}), but with $m\eq kd$.
So hereafter,
{\it we assume $m\eq d$ when studying $\ov M_2(\P^n,d)$ and $m\eq kd$ when studying $\mathbf{R}\pi_*\ff^*\sO_{\P^n}(k)$,}
and focus on the degree-$m$ component of $\fD_2$.

By~\cite[Propostion 2.17]{HLN},
the homomorphism 
$\varphi:\sO_\cV^{\oplus m}\!\lra\!\sO_\cV^{\oplus 2}$ takes the form
\begin{align}\label{Eqn:str_hom}
	\varphi=
	\left[\begin{matrix}
		c_{11}\,\chi_{11} &
		c_{12}\,\chi_{12} &
		\cdots &
		c_{1d}\,\chi_{1m} \\
		c_{21}\,\chi_{21} &
		c_{22}\,\chi_{22} &
		\cdots &
		c_{2d}\,\chi_{2m}
	\end{matrix}\right]\,,
\end{align}
where $c_{si}\inn\Ga(\sO_\cV^*)$
and $\chi_{si}\inn\Ga(\sO_\cV)$, $1\!\le\!s\!\le\!2$, $1\!\le\!i\!\le\!m$.
Moreover,
\begin{itemize}[leftmargin=*]
	\item if $C_\core$ is inseparable, then
	\begin{align}
		\label{Eqn:local_eqn_chi1}
		\chi_{1i}=\chi_{2i}=
		\prod_{e\in \wp_i(x)}\!\!\ze_e;
	\end{align}
	
	\item if $C_\core$ is separable, then by writing $\score(C)\eq\{T_1,T_2\}$,
	we have
	\begin{align}
		\label{Eqn:local_eqn_chi2}
		\chi_{si}=
		\prod_{e\in \wp_{T_s,i}(x)}\!\!\!\!\ze_e,\qquad s=1,2.
	\end{align}
\end{itemize}

In order to achieve a modular  resolution of $\ov M_2(\P^n,d)$,
we aim to diagonalize $\varpi$ in the sense of \cite[Definition 2.6]{HLN},
i.e.~to show locally the pullback of $\varphi$, after suitable trivialization, takes the form
\begin{align}\label{Eqn:diag}
	\left[\begin{matrix}
		\ti z_1 & 0 \\
		0 &
		\ti z_2 
	\end{matrix}\right]
	\left[\begin{matrix}
		1 & 0 & 0 &
		\cdots &
		0 \\
		0 &
		1 & 0 &
		\cdots &
		0
	\end{matrix}\right]\,,
\end{align}
where $\ti z_1$ and $\ti z_2$ are products of local parameters satisfying $\ti z_1\!\mid\!\ti z_2$; c.f.~\cite[Definition 2.6]{HLN}.
To this end,
we also need to analyze
$$
\th_{ij}:=\det\left[\begin{matrix}
	c_{1i} & c_{1j}\\ c_{2i} & c_{2j}
\end{matrix}\right],\qquad
1\!\le\!i,j\!\le\!m\,.
$$

For distinct $1\!\le\!i,j\!\le\!m$,
the expression of $\th_{ij}$ is presented in~\cite[Proposition 2.17 ($\tn{SH}_2$)]{HLN}:
\begin{align}\label{Eqn:local_eqn_th}
	\th_{ij}=
	u_{ij}\cdot \ka_{ij} \cdot \prod_{e\in \wp_i(x)\cap \wp_j(x)}\!\!\!\!\!\!\ze_e\ ,
\end{align}
where $u_{ij}\inn\Ga(\sO^*_\cV)$,
and $\ka_{ij}$ is given in Lemma~\ref{Lm:kappa}.
Moreover, for pairwise distinct $1\!\le\!i,j,\ell\!\le\!m$,
the relation between $\ka_{ij}$ and $\ka_{i\ell}$ is provided in~\cite[Proposition 2.17 ($\tn{SH}_3$)]{HLN}:
there exist $u_{ij\ell}\inn\Ga(\sO^*_\cV)$ and $w_{ij\ell}\inn\Ga(\sO_\cV)$ such that
\begin{align}
	\label{Eqn:local_eqn_ka_ijk}
	\ka_{i\ell}=u_{ij\ell}\cdot\ka_{ij}+w_{ij\ell}\,,
\end{align}
satisfying
\begin{itemize}[leftmargin=*]
\item 
if $Q_i(x)$, $Q_j(x)$ and $Q_\ell(x)$ are all nonempty and satisfy
$T^{(i)}\eq T^{(j)}\eq T^{(\ell)}\inn\score(C)$,
then with $Q_{j\ell}(x)$ as in Lemma~\ref{Lm:kappa},
there exist invertible functions $v_{\wp}$, $\wp\inn Q_{j\ell}(x)$, on $\cV$, such that
\begin{align}\label{Eqn:local_eqn_remainder1}
	w_{ij\ell}=\sum_{\wp\in Q_{j\ell}(x)}\!\!
	\Big(v_{\wp}\prod_{e\in \wp}
	\ze_{e}\Big)
	\cdot 
	\Big(\prod_{e\in\wp_j(x)\cap \wp_\ell(x)}\!\!\!\!\!\!\ze_e\,\Big);
\end{align}
moreover,
by writing $Q_{ij}\eq\{\wp_1,\wp_2\}$ and $Q_{j\ell}\eq\{\wp'_1,\wp'_2\}$ so that
$\wp_r\!\cap\!\wp'_r\!\ne\!\emptyset$, $r\eq 1,2$,
then with $u_{\wp_r}\inn\Ga(\sO^*_\cV)$ as in Lemma~\ref{Lm:kappa},
we have
\begin{align}
	\label{Eqn:local_eqn_remainder2}
\det
\left[\begin{matrix}
	u_{\wp_1} & v_{\wp'_1}\\
	u_{\wp_2} & v_{\wp'_2}
\end{matrix}\right]
\in\Ga(\sO^*_\cV)\,;
\end{align}

\item 
otherwise,
there exists an invertible function $v_{ij\ell}$ on $\cV$ such that
\begin{align}\label{Eqn:local_eqn_remainder3}
	w_{ij\ell}=
	v_{ij\ell}\cdot\prod_{e\in\wp_j(x)\cap \wp_\ell(x)}\!\!\!\!\!\!\ze_e.
\end{align}
\end{itemize}

\begin{eg}\label{Eg:abcd}  
Let $x\eq(C,D)\inn\fD_2$, where $C$ is comprised of a smooth irreducible component $C_\core$ of genus 2, as well as four smooth rational irreducible components $R_a$, $R_b$, $R_c$, and $R_f$; c.f.~Figure~\ref{Fig:abcd}.
Let $D\eq\de_1\!+\!\cdots\!+\!\de_d$ be a simple divisor of $C$, satisfying 
\begin{align*}
	&
	D\cap R_a=\{\de_1,\cdots,\de_{h-1}\},&&
	D\cap R_c=\{\de_{h},\cdots,\de_{r-1}\},\\
	&
	D\cap R_f=\{\de_{r},\cdots,\de_d\},&&
	D\cap C_\core=D\cap R_b=\emptyset
\end{align*}
(for some $2\!<\!h\!<\!r\!\le\!d$).
In Figure~\ref{Fig:abcd}, the components $R_a$, $R_c$, and $R_f$ (resp.~$C_\core$ and $R_b$) are shaded (resp.~unshaded) to indicate they contain (resp.~do not contain) elements of $D$.

\begin{figure}[htp]
	\begin{center}
		\begin{tikzpicture}
			\def\g2left{
				(0,0.8) arc (90:270:1.6 and 0.8)
				(-1.04,0.08)--(-0.88,0)
				..controls (-0.64,-0.12)..(-0.4,0)
				--(-0.24,0.08)
				(-0.88,0)..controls (-0.64,0.12)..(-0.4,0)
			}
			
			\draw \g2left
			[xscale=-1] \g2left;
			\draw
			(.8,1.13) circle (0.4);
			\draw[fill=black!42]
			(-.8,1.13) circle (0.4)
			(.8,1.93) circle (0.4)
			(1.6,1.13) circle (0.4);
			\draw 
			(-.8,1.13) node {\tiny{$R_a$}}
			(.66,1) node {\tiny{$R_b$}}
			(.8,1.93) node {\tiny{$R_c$}}
			(1.6,1.13) node {\tiny{$R_f$}}
			(0,-.5) node {\tiny{$C_\core$}}
			(-.64,.5) node {\tiny{$q_a$}}
			(.7,.5) node {\tiny{$q_b$}}
			(.8,1.38) node {\tiny{$q_c$}}
			(1.03,1.13) node {\tiny{$q_f$}};
			
			\filldraw[xshift=5.5cm,yshift=1cm]
			(.35,-.66) circle (1pt)
			(1.05,-.66) circle (1pt)
			(0,0) circle (1pt)
			(.35,.66) circle (1pt)
			(.7,0) circle (1pt);
			\draw[xshift=5.5cm,yshift=1cm]
			(1.05,-.66)--(.35,.66)--(0,0)
			(.7,0)--(.35,-.66);
			\draw[xshift=5.5cm,yshift=1cm]
			(.35,.66) node[above] {\tiny{$o_\al$}}
			(.175,.33) node[left] {\tiny{$a$}}
			(.525,.33) node[right] {\tiny{$b$}}
			(.875,-.33) node[right] {\tiny{$f$}}
			(.525,-.33) node[left] {\tiny{$c$}}
			(-.35,-.66) node {\tiny{$\tau_{\al}$}}
			;
			
			\draw[xshift=8.8cm,yshift=1cm]
			(0,.5) node [left] {\tiny{$\{a,b,c,f\}$}}
			(.5,.5) node [right] {\tiny{$S_\al\eq\{a,b,c,f\}$}}
			(.25,.5) node [above] {\tiny{$\be_\al$}}
			(.18,.25) node [below] {\tiny{$\be_\al(e)\eq e\ \tn{for} \ e\eq a,b,c,f$}}
			;
			
			\draw[xshift=8.8cm,yshift=1cm,->,>=stealth']
			(0,.5)--(.5,.5);
		\end{tikzpicture}
	\end{center}
	\caption{The curve $C$ in Example~\ref{Eg:abcd}}\label{Fig:abcd}
\end{figure}

Assume neither $q_a$ nor $q_b$ is a Weierstrass point of $C_\core$.
This implies 
$$
\ka_{1i}(x)\ne 0\quad\forall~2\le i< h,\qquad
\ka_{1i}(x)= 0\quad\forall~h\le i\le d;
$$
c.f.~\cite[Lemma 2.15]{HLN}.

As in~(\ref{Eqn:str_hom}), locally on a chart $\cV$ containing the point $x$, we have
$$\varphi =\left[\;
\begin{matrix}  
	c_{11} \zeta_a & 
	\cdots &
	c_{1,{h-1}} \zeta_a\  &
	c_{1h} \zeta_b\zeta_c & 
	\cdots &
	c_{1,{r-1}} \zeta_b\zeta_c\  &
	c_{1r} \zeta_b\zeta_f &
	\cdots &
	c_{1d} \zeta_b\zeta_f
	\\
	c_{21} \zeta_a & 
	\cdots &
	c_{2,{h-1}} \zeta_a &
	c_{2h} \zeta_b\zeta_c & 
	\cdots &
	c_{2,{r-1}} \zeta_b\zeta_c &
	c_{2r} \zeta_b\zeta_f &
	\cdots &
	c_{2d} \zeta_b\zeta_f
\end{matrix}\; \right]
$$  
where $\zeta_a$ (resp.~$\zeta_b, \zeta_c, \zeta_d$) is the node-smoothing parameter for $q_a$ (resp.~$q_b,$ $q_c,$ $q_d$).

Since each entry of $\varphi$ contains one among $\ze_a$, $\ze_b\ze_c$ and $\ze_b\ze_f$ as a factor,
we anticipate an LES and a treelike structure $\La$ on $\fD_2$ so that
the stratum $\fM_\al$ containing $x$ 
is locally given by
$$
\fM_\al\cap\cV\,=\,
\big(\,\ze_a=\ze_b=\ze_c=\ze_f=0\,\big)\,,
$$
and the rooted tree $\tau_\al$ and the injection $\be_\al$ are as in Figure~\ref{Fig:abcd}.
In \S\ref{Subsec:Step1}, we will construct $\La$ explicitly.

Assume $\La$ exists for now.
Let $\varpi:(\fD_2)^{\tf}_{\La}\to\fD_2$ be the forgetful morphism as in Theorem~\ref{Thm:tf_smooth}~\ref{Cond:smooth_tf}.
By Corollary~\ref{Crl:prod_RL},
near any lift of $x$,
one of
$$\varpi^*(\ze_a),\quad
\varpi^*(\ze_b\ze_c),\quad\tn{and}\quad
\varpi^*(\ze_b\ze_f)
$$
divides the other two.
For instance,
consider $\ti\al\!:=\!(\al,\ov\lt)\inn\ov\La$ with $\ov\lt$ given in Figure~\ref{Fig:abcd'}.
\begin{figure}[htp]
\begin{center}
\begin{tikzpicture}
	\filldraw
	(.6,.5) circle (1.2pt)
	(-.3,0) circle (1.2pt)
	(.3,1) circle (1.2pt)
	(0,-.5) circle (1.2pt)
	(1.2,-.5) circle (1.2pt);
	\draw
	(1.2,-.5)--(.3,1)--(-.3,0)
	(0,-.5)--(.6,.5);
	\draw[dotted]
	(1.5,0)--(-.6,0)
	(1.5,1)--(-.6,1)
	(1.5,.5)--(-.6,.5);
	\draw[thin,decorate,decoration=brace,xshift=0cm]
	(1.5,.95)--(1.5,.55);
	\draw[thin,decorate,decoration=brace,xshift=0cm]
	(1.5,.45)--(1.5,.05);
	\draw
	(.4,.8) node[right] {\tiny{$b$}}
	(-.1,.3) node[left] {\tiny{$a$}}
	(.1,-.25) node[right] {\tiny{$c$}}
	(1.1,-.25) node[left] {\tiny{$f$}}
	(1.5,.75) node[right] {\tiny{$\fE_1$}}
	(1.5,.25) node[right] {\tiny{$\fE_2$}};
	\draw
	(1.8,-.5) node[right] {\scriptsize{$\ov\lt=(\,\tau_\al\,,\,\{\fE_1,\fE_2\}\,,\,
			\{b,a\}\,)$}}
	;
	
	\filldraw[xshift=9cm,yshift=0cm]
	(0,.5) circle (1.2pt)
	(-.5,0) circle (1.2pt)
	(.5,1) circle (1.2pt)
	(.5,0) circle (1.2pt)
	(1.5,0) circle (1.2pt);
	\draw[xshift=9cm,yshift=0cm]
	(1.5,0)--(.5,1)--(-.5,0)
	(.5,1)--(.5,0);
	\draw[xshift=9cm,yshift=0cm]
	(-.2,.3) node[left] {\tiny{$\fE_{2}$}}
	(.3,.8) node[left] {\tiny{$\fE_{1}$}}
	(.45,.45) node[right] {\tiny{$c$}}
	(1.15,.45) node[right] {\tiny{$f$}}
	(.5,-.5) node {\scriptsize{$\vr_{\ov\lt}=\dg{\ov\lt}$}}
	;
\end{tikzpicture}	
\end{center}
\caption{An element of $\ov\bT_\al(\La)$ and its derived tree in Example~\ref{Eg:abcd}}\label{Fig:abcd'}	
\end{figure} 
In this case, we have
\begin{align*}
	\varpi^*(\ze_a)=u_{a}\xi_{\fE_1}\xi_{\fE_2},\quad
	\varpi^*(\ze_b\ze_c)=u_{c}\xi_{\fE_1}\xi_{\fE_2}\xi_{c},\quad
	\varpi^*(\ze_b\ze_f)=u_{f}\xi_{\fE_1}\xi_{\fE_2}\xi_f,
\end{align*}
where $u_a$, $u_c$, and $u_f$ are invertible.
Therefore, after choosing suitable trivialization and applying~(\ref{Eqn:local_eqn_th}),
we have
\begin{align*}
	\varpi^*\varphi
	=\xi_{\fE_1}\xi_{\fE_2}
	\left[
	\begin{matrix}  
		1 & 0 & \cdots & 0 & 0 & \cdots & 0 & 0 & \cdots & 0
		\\
		0 & 
		\ti\ka_{12}\xi_{\fE_1}\xi_{\fE_2} &
		\cdots &
		\ti\ka_{1,h-1}\xi_{\fE_1}\xi_{\fE_2} &
		\ti\ka_{1h}\xi_c &
		\cdots & 
		\ti\ka_{1,r-1}\xi_c &
		\ti\ka_{1r}\xi_f &
		\cdots & 
		\ti\ka_{1d}\xi_f
	\end{matrix}\right]
\end{align*}
where $\ti\ka_{ij}$ are the pullbacks of $\ka_{ij}$ to $(\fD_2)_\La^{\tf}$.
All $\ti\ka_{ij}$ above are smooth functions in this case.

At this stage, we expect a new treelike structure on $(\fD_2)^{\tf}_\La$ so as to compare the entries of the second row of $\varpi^*\varphi$.
\begin{itemize}[leftmargin=*]
\item 
If $q_a$ and $q_b$ are not conjugate,
then all $\ti\ka_{ij}$ above are invertible,
so we only need to compare
$$\xi_{\fE_1}\xi_{\fE_2},\qquad
\xi_{c},\qquad\tn{and}\qquad
\xi_f.$$
As we will see in \S\ref{Subsec:derived},
the \textsf{derived treelike structure}
$\ov{\dt}\La$ assigns the rooted tree $\vr_{\ov\lt}$ to $(\fD_2)^{\tf}_{\ti\al}$,
where $\vr_{\ov\lt}$ is shown in Figure~\ref{Fig:abcd'}.
By
Theorems~\ref{Thm:tf_smooth}~\ref{Cond:smooth_parameters} and~\ref{Thm:tf_smooth_revert}~\ref{Cond:smooth_parameters_revert}, the pullback of $\varphi$ to either $\big((\fD_2)^{\tf}_\La\big)^{\tf}_{\ov\dt\La}$ or $\big((\fD_2)^{\tf}_\La\big)^{\rtf}_{\ov\dt\La}$ is in the desired form (\ref{Eqn:diag}).

\item
If $q_a$ and $q_b$ are  conjugate, then $\ti\ka_{1h}$ above is not invertible.
By~(\ref{Eqn:local_eqn_ka_ijk}) and (\ref{Eqn:local_eqn_remainder3}),
we have
\begin{align*}
	\ti\ka_{1i}=
	\begin{cases}
		u_i\ti\ka_{1h}+v_i
		&\forall~2\le i<h,\\
		u_i\ti\ka_{1h}+v_i\xi_{\fE_1}\xi_{\fE_2}\xi_c
		&\forall~h< i< r,\\
		u_i\ti\ka_{1h}+v_i\xi_{\fE_1}
		&\forall~r\le i\le d,
	\end{cases}
\end{align*}
where all $u_i$ and $v_i$ are invertible.
After choosing suitable trivialization again, 
we can rewrite $\varpi^*\varphi$ as
\begin{align*}
	\xi_{\fE_1}\xi_{\fE_2}
	\left[
	\begin{matrix}  
		1 & 0 &  0 &  0 & \cdots & 0 & 0 & \cdots & 0
		\\
		0 & 
		\xi_{\fE_1}\xi_{\fE_2} &
		\ti\ka_{1h}\xi_c &
		\ti\ka_{1r}\xi_f &
		\cdots & 
		\ti\ka_{1d}\xi_f &
		0 & \cdots & 0
	\end{matrix}\right].
\end{align*}

The following two examples suggest that
adding RL twisted fields to $(\fD_2)^{\tf}_\La$ with respect to $\ov\dt\La$
may not simplify, if not complicate,
the second row of the pullback of $\varphi$; however,
adding LR twisted fields instead does turn
the second row  into the desired form (\ref{Eqn:diag}).

\begin{figure}[htp]
	\begin{center}
		\begin{tikzpicture}
			\filldraw
			(.6,.5) circle (1.2pt)
			(-.18,-.3) circle (1.2pt)
			(.6,1) circle (1.2pt)
			(0,0) circle (1.2pt)
			(1.5,0) circle (1.2pt);
			\draw
			(-.18,-.3)--(.6,1)--(1.5,0)
			(.6,1)--(.6,.5);
			\draw[dotted]
			(1.8,1)--(-.6,1)
			(1.8,.5)--(-.6,.5);
			\draw[thin,decorate,decoration=brace,xshift=0cm]
			(1.8,.95)--(1.8,.55);
			\draw[xshift=.3cm]
			(.2,.65) node[right] {\tiny{$c$}}
			(0,.3) node[left] {\tiny{$\fE_1$}}
			(-.6,-.1) node {\tiny{$\fE_2$}}
			(1.1,.35) node {\tiny{$f$}}
			(1.55,.75) node[right] {\tiny{$\fF$}}
			(1.6,-.3) node[right] {\scriptsize{$\ov\fr=(\,\vr_{\ov\lt},\,\{\fF\},\,\{c\}\,)$}}
			;
			
			\filldraw[xshift=8cm,yshift=0cm]
			(0,.5) circle (1.2pt)
			(-.5,0) circle (1.2pt)
			(.5,1) circle (1.2pt)
			(.5,0) circle (1.2pt)
			(1.5,0) circle (1.2pt);
			\draw[xshift=8cm,yshift=0cm]
			(1.5,0)--(.5,1)--(-.5,0)
			(.5,1)--(.5,0);
			\draw[xshift=8cm,yshift=0cm]
			(-.15,.3) node[left] {\tiny{$\fE_2$}}
			(.35,.8) node[left] {\tiny{$\fE_1$}}
			(.65,.35) node {\tiny{$\fF$}}
			(1.4,.35) node {\tiny{$f$}}
			(1.7,-.25) node[right] {\scriptsize{$\dg{\ov\fr}$}}
			;
		\end{tikzpicture}	
	\end{center}
	\caption{An element of $\ov\bT_{\ti\al}(\ov\dt\La)$ and its derived tree in Example~\ref{Eg:abcd}}\label{Fig:abcd''}	
\end{figure}
\begin{itemize}[leftmargin=*]
\item 
Consider $\wc\al\eq(\al,\ov\lt,\ov\fr)\inn\ov{\ov\dt\La}$,
where $\ov\fr$ is given in Figure~\ref{Fig:abcd''}.
We denote by
$$\varpi':\big((\fD_2)^{\tf}_\La\big)^{\tf}_{\ov\dt\La}\lra
(\fD_2)^{\tf}_\La$$
the forgetful morphism of Theorem~\ref{Thm:tf_smooth}.
Then, locally we have
$$
(\varpi')^*(\xi_{\fE_1}\xi_{\fE_2})=u'_a\wc\xi_{\fE_1}\xi_{\fF}\xi_{\fE_2},\qquad
(\varpi')^*(\xi_{c})=u'_c\xi_{\fF},\qquad
(\varpi')^*(\xi_{\fE_1}\xi_{f})=u'_f\wc\xi_{\fE_1}\wc\xi_{f}(\xi_{\fF})^2,
$$
where $u'_a$, $u'_c$, and $u'_f$ are invertible functions,
and $\xi_{\fF}$, $\wc\xi_{\fE1}$, $\xi_{\fE_2}$, and $\wc\xi_f$ are local parameters.
Therefore, the second row of $(\varpi\!\circ\!\varpi')^*\varphi$,
after extracting the common factor $\xi_{\fF}$ and taking suitable trivialization,
can be written as
$$
\left[
\begin{matrix}
	0 &
	\wc\xi_{\fE_1}\xi_{\fE_2} &
	\ti\ka_{1h} & 
	\wc\xi_{\fE_1}\xi_{\fF}\wc\xi_f & 
	0 &
	\cdots & 
	0
\end{matrix}
\right].
$$
In order to proceed,
one needs to compare $$
 \wc\xi_{\fE_1}\xi_{\fE_2},\qquad\wc\xi_{\fE_1}\xi_{\fF}\wc\xi_f,\qquad
 \tn{and}\qquad
 \ti \ka_{1h}.
$$
While the extra parameter $\ti \ka_{1h}$ can be dealt with using the technique of \ts{grafting} introduced in \S\ref{Subsec:Grafted},
neither $\dg{\ov\lt}$ nor $\dg{\ov\fr}$, respectively illustrated in Figures~\ref{Fig:abcd'} and~\ref{Fig:abcd''},
appear promising for the comparison of $\xi_{\fE_2}$ and $\xi_{\fF}\wc\xi_f$.
In sum,
adding RL-twisted fields to $(\fD_2)^{\tf}_\La$ with respect to $\ov\dt\La$ could make
the second row of the pullback of $\varphi$ even more complicated.

\begin{figure}[htp]
\begin{center}
\begin{tikzpicture}
	\filldraw
	(-.5,1) circle (1.2pt)
	(-.5,0) circle (1.2pt)
	(.5,1) circle (1.2pt)
	(.5,0) circle (1.2pt)
	(1.5,0) circle (1.2pt);
	\draw
	(-.5,0)--(-.5,1)--(.5,1)--(1.5,0)
	;
	\draw[ultra thick]
	(.5,1)--(.5,0);
	\draw[dotted]
	(1.8,1)--(-.6,1)
	(1.8,0)--(-.6,0);
	\draw[thin,decorate,decoration=brace,xshift=0cm]
	(1.85,.95)--(1.85,.05);
	\draw
	(.62,.35) node {\tiny{$c$}}
	(0,.85) node {\tiny{$\fE_1$}}
	(-.3,.35) node {\tiny{$\fE_2$}}
	(1.4,.35) node {\tiny{$f$}}
	(1.85,.5) node[right] {\tiny{$\fF$}}
	(.7,-.4) node[below] {\scriptsize{$\ud\ls=(\,\vr_{\ov\lt},\,\{\fF\},\,\{c\}\,)$}}
	;
	
	\filldraw[xshift=4cm,yshift=0cm]
	(0,.5) circle (1.2pt)
	(-.5,0) circle (1.2pt)
	(.5,1) circle (1.2pt)
	(.5,0) circle (1.2pt)
	(1.5,0) circle (1.2pt);
	\draw[xshift=4cm,yshift=0cm]
	(1.5,0)--(.5,1)--(-.5,0)
	(.5,1)--(.5,0);
	\draw[xshift=4cm,yshift=0cm]
	(-.15,.3) node[left] {\tiny{$\fE_2$}}
	(.35,.8) node[left] {\tiny{$\fE_1$}}
	(.65,.35) node {\tiny{$\fF$}}
	(1.4,.35) node {\tiny{$f$}}
	(.5,-.4) node[below] {\scriptsize{$\dg{\ud\ls}$}}
	;
	
	\filldraw[xshift=8cm,yshift=0cm]
	(.6,.5) circle (1.2pt)
	(-.3,0) circle (1.2pt)
	(.3,1) circle (1.2pt)
	(.3,0) circle (1.2pt)
	(1.2,-.5) circle (1.2pt);
	\draw[xshift=8cm,yshift=0cm]
	(1.2,-.5)--(.3,1)--(-.3,0)
	(.3,0)--(.6,.5);
	\draw[dotted,xshift=8cm,yshift=0cm]
	(1.5,0)--(-.6,0)
	(1.5,1)--(-.6,1)
	(1.5,.5)--(-.6,.5);
	\draw[thin,decorate,decoration=brace,xshift=8cm]
	(1.5,.95)--(1.5,.55);
	\draw[thin,decorate,decoration=brace,xshift=8cm]
	(1.5,.45)--(1.5,.05);
	\draw[xshift=8cm,yshift=0cm]
	(.4,.8) node[right] {\tiny{$b$}}
	(-.1,.3) node[left] {\tiny{$a$}}
	(.1,.25) node[right] {\tiny{$c$}}
	(1.1,-.25) node[left] {\tiny{$f$}}
	(1.55,.75) node[right] {\tiny{$\fE_1$}}
	(1.55,.25) node[right] {\tiny{$\fE_2$}};
	\draw[xshift=8cm,yshift=0cm]
	(.2,-.4) node[below] {\scriptsize{$\ov\lt_{\Dm(\ud\ls)}$}}
	;
	
	\filldraw[xshift=12cm,yshift=0cm]
	(0,.5) circle (1.2pt)
	(-.5,0) circle (1.2pt)
	(0,1) circle (1.2pt)
	(.5,0) circle (1.2pt);
	\draw[xshift=12cm,yshift=0cm]
	(.5,0)--(0,.5)--(-.5,0)
	(0,1)--(0,.5);
	\draw[xshift=12cm,yshift=0cm]
	(-.2,.3) node[left] {\tiny{$\fE_{2}$}}
	(.1,.75) node[left] {\tiny{$\fE_{1}$}}
	(.2,.3) node[right]{\tiny{$f$}}
	(0,-.5) node {\scriptsize{$\dg{(\ov\lt_{\Dm(\ud\ls)})}$}}
	;
\end{tikzpicture}			
\end{center}
\caption{An element of $\ud\bT_{\ti\al}(\ov\dt\La)$ and its derived tree in Example~\ref{Eg:abcd}}\label{Fig:abcd'''}	
\end{figure}
\item 
If instead, the LR-twisted fields are  added to $(\fD_2)^{\tf}_\La$ with respect to $\ov\dt\La$,
then the aforementioned unpleasant situation can be avoided.
For instance,
consider $\wh\al\eq(\al,\ov\lt,\ud\ls)\inn\ud{\ov\dt\La}$,
where $\ud\ls$ is given in Figure~\ref{Fig:abcd'''}.
We denote by
$$\wh\varpi:\big((\fD_2)^{\tf}_\La\big)^{\rtf}_{\ov\dt\La}\lra
(\fD_2)^{\tf}_\La$$
the forgetful morphism of Theorem~\ref{Thm:tf_smooth_revert}.
Then, locally we have
$$
\wh\varpi^*(\xi_{\fE_1}\xi_{\fE_2})=u'_a\xi_{\fE_1}\wc\xi_{\fE_2}\xi_{\fF},\qquad
\wh\varpi^*(\xi_{c})=u'_c\xi_{\fF},\qquad
\wh\varpi^*(\xi_{f})=u'_f\wc\xi_{f}\xi_{\fF},
$$
where $u'_a$, $u'_c$, and $u'_f$ are invertible functions,
and $\xi_{\fF}$, $\xi_{\fE1}$, $\wc\xi_{\fE_2}$, and $\wc\xi_f$ are local parameters.
Therefore, the second row of $(\varpi\!\circ\!\wh\varpi)^*\varphi$,
after extracting common factors and taking suitable trivialization,
can be written as
$$
\left[
\begin{matrix}
	0 &
	\xi_{\fE_1}\wc\xi_{\fE_2} &
	\ti\ka_{1h} & 
	\xi_{\fE_1}\wc\xi_f & 
	0 &
	\cdots & 
	0
\end{matrix}
\right].
$$
To proceed,
one needs to compare 
$$
\xi_{\fE_1}\wc\xi_{\fE_2},\qquad\xi_{\fE_1}\wc\xi_f,\qquad
\tn{and}\qquad
\ti \ka_{1h}.
$$
The extra parameter $\ti \ka_{1h}$ can be dealt with via grafting (c.f.~\S\ref{Subsec:Grafted}),
and
the \ts{doubly derived tree} of $\ov\lt_{(\Dm(\ud\ls))}$,
illustrated in Figure~\ref{Fig:abcd'''} and formally introduced in \S\ref{Subsec:2nd_derived_treelike}, 
provides the model for comparing $\xi_{\fE_1}\wc\xi_{\fE_2}$ and $\xi_{\fE_1}\wc\xi_f$.
\end{itemize}
\end{itemize}
\end{eg}

\subsection{Proper transforms of treelike structures}\label{Subsec:two_treelike}

In this subsection,
we establish a statement about two treelike structures $\La$ and $\La^\dag$ on a {\set} stratified stack $\fM$:
if we add RL (resp.~LR) twisted fields to $\fM$ with respect to $\La$,
there is a standard way to convert $\La^\dag$ to a new treelike structure on $\fM^{\tf}_\La$ (resp.~$\fM^{\rtf}_\La$).
This will play a crucial role in the construction of new treelike structures after \S\ref{Subsec:Step4}.

We start with the RL case.
Let $\fM$ be a stack with an LES as in Definition~\ref{Dfn:G-adim_fixture} 
and 
\begin{align*}
\La=(\tau_\al,\be_\al)_{\al\in A}
\qquad\tn{and}\qquad
\La^\dag=(\tau^\dag_\al,\be^\dag_\al)_{\al\in A}
\end{align*}
be two treelike structures on $\fM$ as in Definition~\ref{Dfn:Treelike_structure}.
Notice that for every $\al\inn A$, the targets of $\be_\al$ and $\be^\dag_\al$ are both $S_\al$.
After adding RL twisted fields to $\fM$ with respect to $\La$,
for every 
\begin{align}\label{Eqn:ti_al}
	\ti\al:=\big(\al,\ov\lt\big)\in \ov\La,\qquad
\tn{where}\quad \ov\lt=\big(\tau_\al,\ov\bE,\Dm(\ov\lt)\big)\in\ov\bT_\al(\La), 
\end{align}
we set
\begin{equation}\begin{split}\label{Eqn:treelike_new}
	&\tau^\dag_{\ti\al} :=
	\tau^\dag_{\al_{(\fB_{\ti\al})}},
	\qquad\tn{where}\qquad
	\fB_{\ti\al}:=\be_\al\big(\Dm(\ov\lt)\big).
\end{split}\end{equation}
The map $\be^\dag_{\ti\al}:\tau^\dag_{\ti\al}\hookrightarrow S_{\ti\al}$ is taken to be the composition:
\begin{align}\label{Eqn:be'_ti_al}
	\be^\dag_{\ti\al}:\,
	\tau^\dag_{\ti\al}
	\stackrel{\be^\dag_{\al_{(\fB_{\ti\al})}}}{\longrightarrow}
	S_{\al_{(\fB_{\ti\al})}}
	\stackrel{\iota_{\al;\fB_{\ti\al}}}{\lra}
	S_\al\bsl\fB_{\ti\al}
	\,\hookrightarrow\,
	\ov\bE\sqcup\big(S_\al\bsl\fB_{\ti\al}\big)
	= S_{\ti\al}.
\end{align} 
Here, $\iota_{\al;\fB_{\ti\al}}$ is as in Corollary~\ref{Crl:M_strata_local}, and
the arrow $\hookrightarrow$ is the inclusion (as a subset).
The injectivity of such defined $\be^\dag_{\ti\al}$ follows from 
that of $\be^\dag_{\al_{(\fB_{\ti\al})}}$ and $\iota_{\al;\fB_{\ti\al}}$.

\begin{prp}\label{Prp:two_treelike_strs}
	The indexed family 
	$$
	\ov\PT_{\La}(\La^\dag):=\big(\,\tau^\dag_{\ti\al},\,\be^\dag_{\ti\al}\,\big)_{\ti\al\in \ov\La}
	$$ gives a treelike structure on $\fM^{\tf}_\La$,
	called the \ts{RL proper transform} of $\La^\dag$ \ts{with respect to} $\La$.
\end{prp}

\begin{proof}
Fix $\ti\al\eq(\al,\ov\lt)\inn\ov\La$ as in (\ref{Eqn:ti_al}) such that $\ti\al\!\ne\!(0,\ov\lt_\bullet)$.
Also fix
$$
S'\,\subset\, S_{\ti\al}=\ov\bE\sqcup\!\big(S_\al\bsl\fB_{\ti\al}\big).
$$
By (\ref{Eqn:be'_ti_al}), we have
\begin{align*}
	(\be^\dag_{\ti\al})^{-1}(S')=
	\big(\be^\dag_{\al_{(\fB_{\ti\al})}}\big)^{-1}\circ 
	\iota_{\al;\fB_{\ti\al}}^{-1}(R'),\qquad
	\tn{where}\quad
	R':= S'\!\cap(S_\al\bsl\fB_{\ti\al}).
\end{align*}
By Theorem~\ref{Thm:tf_smooth}~\ref{Cond:smooth_Z},
there uniquely exist $J\!\subset\!\bbI(\ov\lt)$ and $E\!\subset\!S_\al$ such that
\begin{align}\label{Eqn:ti_al_(s)}
\ti\al_{(S')}=\big(\,\al_{(E)},\,\ov\lt'\,\big)\,,\qquad\tn{where}\quad
\ov\lt':=\phi^{-1}_{\al;E}(\ov\lt_{(J)});
\end{align}
c.f.~(\ref{Eqn:C_J}) and Propositions~\ref{Prp:M_strata_local} and~\ref{Prp:Treelike_contraction} for notation.
We aim to construct an isomorphism
\begin{align}\label{Eqn:phi_ti_al_s}
	&
	\phi_{\ti\al;S'}:\,\tau^\dag_{\ti\al_{(S')}}\lra
	\tau^\dag_{\ti\al}\big\bsl\,
	\big((\be^\dag_{\ti\al})^{-1}(S')\big)^{\wedge}=
	\tau^\dag_{\al_{(\fB_{\ti\al})}}
	\Big\bsl\,
	\Big(\big(\be^\dag_{\al_{(\fB_{\ti\al})}}\big)^{-1}\circ 
	\iota^{-1}_{\al;\fB_{\ti\al}}(R')\Big)^{\!\wedge},\\
	\label{Eqn:phi_ti_al_s'}
	&
	\tn{satisfying}\qquad
	\iota_{\ti\al;S'}\circ\be^\dag_{\ti\al_{(S')}}
	=\be^\dag_{\ti\al}\circ \phi_{\ti\al;S'}.
\end{align}
Here, $\iota_{\ti\al;S'}:S_{\ti\al_{(S')}}\!\lra\!S_{\ti\al}\bsl S'$ is described in Theorem~\ref{Thm:tf_smooth}~\ref{Cond:smooth_Z}.
Since $\ti\al$ and $S'$ are arbitrary,
(\ref{Eqn:phi_ti_al_s}) and (\ref{Eqn:phi_ti_al_s'}) will guarantee $\ov\PT_\La(\La^\dag)$ is a treelike structure.

We start with analyzing $\tau^\dag_{\ti\al_{(S')}}$, which by (\ref{Eqn:treelike_new}) is given by 
\begin{align*}
	\tau^\dag_{\ti\al_{(S')}}=
	\tau^\dag_{(\al_{(E)})_{\big(\be_{\al_{(E)}}(\Dm(\ov\lt'))\big)}}.
\end{align*}
By Corollary~\ref{Crl:M_strata_local},
we have
\begin{align*}
	(\al_{(E)})_{\big(\be_{\al_{(E)}}(\Dm(\ov\lt'))\big)}
	=\al_{\big(E\;\sqcup\; \iota_{\al;E}\circ\be_{\al_{(E)}}(\Dm(\ov\lt'))\big)}.
\end{align*}
Notice that Lemma~\ref{Lm:RLS_LRS_isom} and the expression of $\ov\lt'$ in (\ref{Eqn:ti_al_(s)}) imply that 
\begin{align*}
	\Dm(\ov\lt')=
	\phi^{-1}_{\al;E}\big(\Dm(\ov\lt_{(J)})\big).
\end{align*}
The fact that $\La$ is a treelike structure and the expression of $E$ in Theorem~\ref{Thm:tf_smooth}~\ref{Cond:smooth_Z} then yield
\begin{align}\label{Eqn:E_cup_dom_t'}
	E\,\sqcup\, \iota_{\al;E}\circ\be_{\al_{(E)}}\big(\Dm(\ov\lt')\big)
	=
	\be_\al\big(\Dm(\ov\lt_{(J)})\big)\sqcup \be_{\al}(\ov\sfE_J)\sqcup 
	\Big(S'\cap\big(S_\al\bsl\be_\al(\tau_\al)\big)\Big).
\end{align}

To further analyze $\be_\al\big(\Dm(\ov\lt_{(J)})\big)$,
observe the expression of
$J$ in  Theorem~\ref{Thm:tf_smooth}~\ref{Cond:smooth_Z} implies
\begin{align*}
	\be_\al\big(J\cap \ND(\ov\lt)\big)=
	S'\cap\be_\al\big(\ND(\ov\lt)\big).
\end{align*}
So by
Lemma~\ref{Lm:I(t_J)} and the injectivity of $\be_\al$, we have
\begin{align*}
	\be_\al\big(\Dm(\ov\lt_{(J)})\big)\sqcup
	\be_\al(\ov\sfE_J)=
	\fB_{\ti\al}\sqcup 
	\Big(S'\cap\be_\al\big(\ND(\ov\lt)\big)\Big).
\end{align*}
Therefore, (\ref{Eqn:E_cup_dom_t'}) can be rewritten as
\begin{align}\label{Eqn:E_cup_dom_t''}
	E\,\sqcup\, \iota_{\al;E}\circ\be_{\al_{(E)}}\big(\Dm(\ov\lt')\big)
	=
	\fB_{\ti\al}\sqcup 
	\big(S'\cap(S_\al\bsl\fB_{\ti\al})\big)
	=\fB_{\ti\al}\sqcup R'.
\end{align}
Hence we have
\begin{align}\label{Eqn:al_(E)_(dom_t')}
	(\al_{(E)})_{\big(\be_{\al_{(E)}}(\Dm(\ov\lt'))\big)}
	=\al_{(\fB_{\ti\al}\sqcup R')}
	=
	(\al_{(\fB_{\ti\al})})_{\big(\iota^{-1}_{\al;\fB_{\ti\al}}(R')\big)},
\end{align}
where the last equality once again follows from Corollary~\ref{Crl:M_strata_local}.

By (\ref{Eqn:al_(E)_(dom_t')}), we can rewrite $\tau^\dag_{\ti\al_{(S')}}$ as
\begin{align*}
	\tau^\dag_{\ti\al_{(S')}}=
	\tau^\dag_{(\al_{(\fB_{\ti\al})})_{\big(\iota^{-1}_{\al;\fB_{\ti\al}}(R')\big)}},
\end{align*}
and take the proposed isomorphism $\phi_{\ti\al;S'}$ of (\ref{Eqn:phi_ti_al_s}) to be
\begin{align}\label{Eqn:phi_ti_al}
	\phi_{\ti\al;S'}:=
	\phi_{\al_{(\fB_{\ti\al})};\;\iota^{-1}_{\al;\fB_{\ti\al}}(R')}:\,
	\tau^\dag_{(\al_{(\fB_{\ti\al})})_{\big(\iota^{-1}_{\al;\fB_{\ti\al}}\!(R')\big)}}\,\lra\,
	\tau^\dag_{\al_{(\fB_{\ti\al})}}\Big\bsl\,
	\Big(\big(\be^\dag_{\al_{(\fB_{\ti\al})}}\big)^{-1}\circ 
	\iota^{-1}_{\al;\fB_{\ti\al}}(R')\Big)^{\!\wedge},
\end{align}
where $\phi_{\al_{(\fB_{\ti\al})};\;\iota^{-1}_{\al;\fB_{\ti\al}}(R')}$ is the isomorphism given by the treelike structure $\La^\dag$  as per Definition~\ref{Dfn:Treelike_structure}.

It remains to verify such defined $\phi_{\ti\al;S'}$ satisfies (\ref{Eqn:phi_ti_al_s'}).
On the one hand,
from the expression of $\iota_{\ti\al;S'}$ in Theorem~\ref{Thm:tf_smooth}~\ref{Cond:smooth_Z}, (\ref{Eqn:be'_ti_al}) for $\ti\al_{(S')}$, (\ref{Eqn:al_(E)_(dom_t')}), Corollary~\ref{Crl:M_strata_local}, and (\ref{Eqn:E_cup_dom_t''}),
we conclude that
\begin{align*}
	\iota_{\ti\al;S'}\circ 
	\be^\dag_{\ti\al_{(S')}}
	&=
	\iota_{\al;E}\circ 
	\iota_{\al_{(E)};\be_{\al_{(E)}}(\Dm(\ov\lt'))}\circ 
	\be^\dag_{\al_{(\fB_{\ti\al}\sqcup R')}}
	\\
	&=
	\iota_{\al\,;\,(E\;\sqcup\; \iota_{\al;E}\circ\be_{\al_{(E)}}(\Dm(\ov\lt')))}\circ 
	\be^\dag_{\al_{(\fB_{\ti\al}\sqcup R')}}
	=\iota_{\al\,;\,(\fB_{\ti\al}\sqcup R')}\circ 
	\be^\dag_{\al_{(\fB_{\ti\al}\sqcup R')}}.
\end{align*}
On the other hand,
by (\ref{Eqn:be'_ti_al}),
the fact $\La^\dag$ is a treelike structure,
(\ref{Eqn:al_(E)_(dom_t')}),
and Corollary~\ref{Crl:M_strata_local},
we have
\begin{align*}
	\be^\dag_{\ti\al}\circ \phi_{\ti\al;S'}
	&=
	\iota_{\al;\fB_{\ti\al}}\circ
	\be^\dag_{\al_{(\fB_{\ti\al})}}\circ 
	\phi_{\al_{(\fB_{\ti\al})};\;\iota^{-1}_{\al;\fB_{\ti\al}}(R')}
	\\
	&=
	\iota_{\al;\fB_{\ti\al}}\circ
	\iota_{\al_{(\fB_{\ti\al})};\;\iota^{-1}_{\al;\fB_{\ti\al}}(R')}\circ
	\be^\dag_{\al_{(\fB_{\ti\al}\sqcup R')}}
	=\iota_{\al\,;\,(\fB_{\ti\al}\sqcup R')}\circ 
	\be^\dag_{\al_{(\fB_{\ti\al}\sqcup R')}}.
\end{align*}
Therefore, $\phi_{\ti\al;S'}$ satisfies (\ref{Eqn:phi_ti_al_s'}).
This completes the proof of Proposition~\ref{Prp:two_treelike_strs}.
\end{proof}

The result above has its analogue in the LR case.
Let $\fM$, $\La$ and $\La^\dag$ be as in Proposition~\ref{Prp:two_treelike_strs}.
After adding LR twisted fields to $\fM$ with respect to $\La$, we set
\begin{equation*}\begin{split}
&
\fB_{\ut\al}:=\be_\al\big(\Dm(\ud\lt)\big),\quad
\tau^\dag_{\ut\al} :=
\tau^\dag_{\al_{(\fB_{\ut\al})}},
\quad
\be^\dag_{\ut\al} :=
\iota_{\al;\fB_{\ut\al}}\!\circ
\be^\dag_{\al_{(\fB_{\ut\al})}}
\!:\,
\tau^\dag_{\ut\al}\hookrightarrow S_{(\al,\ud\lt)},\qquad
\forall~\ut\al\!:=\!(\al,\ud\lt)\inn \ud\La.		
\end{split}\end{equation*}
Here, $\iota_{\al;\fB_{\ut\al}}$ is as in Corollary~\ref{Crl:M_strata_local}.

\begin{prp}\label{Prp:two_treelike_strs_revert}
	The indexed family 
	$$
	\ud\PT_{\La}(\La^\dag):=\big(\,\tau^\dag_{\ut\al},\,\be^\dag_{\ut\al}\,\big)_{\ut\al\in \ud\La}
	$$ gives a treelike structure on $\fM^{\rtf}_\La$,
	called the \ts{LR proper transform} of $\La^\dag$ \ts{with respect to} $\La$.
\end{prp}

The proof of Proposition~\ref{Prp:two_treelike_strs_revert} is largely parallel to that of Proposition~\ref{Prp:two_treelike_strs},
so we omit the details.

The reason $\ov\PT_\La(\La^\dag)$ and $\ud\PT_\La(\La^\dag)$ are called ``proper transforms'' is as follows.
In the RL case, according to Theorem~\ref{Thm:tf_smooth}~\ref{Cond:smooth_local_blowup},
locally on any modular chart $\cV$ of $\fM$, $\fM^{\tf}_\La/\cV$ and $\fM^{\tf}_{\La^\dag}/\cV$ are isomorphic to sequences of blowups of $\cV$, respectively.
By direct computation,
one can show that locally on any twisted chart $\cU$ of $\fM^{\tf}_\La$, $(\fM^{\tf}_{\La})^{\tf}_{\ov\PT_\La(\La^\dag)}\big/\cU$ is isomorphic to the sequence of blowups of $\cU$ along the {\it proper transforms} of the blowup centers of $\fM^{\tf}_{\La^\dag}\big/\cU$.
This observation, along with its counterpart on the LR scenario, justifies the terminology.

\subsection{Secondary proper transforms}\label{Subsec:2nd_PT}

We continue with the notation from \S\ref{Subsec:two_treelike}.
In \S\ref{Sec:genus_2_twisted_fields},
we will encounter situations when 
there are treelike structures 
\begin{alignat*}{2}
	&\La=(\tau_\al,\be_\al)_{\al\in A}\qquad\tn{and}\qquad
	\La^\dag=(\tau^\dag_\al,\be^\dag_\al)_{\al\in A}\qquad\quad
	&&\tn{on}\quad\fM\,,\\
	&\La^\ddag=
	\big(\,\tau^\ddag_{(\al,\ud\ls)}\,,\,
	\be^\ddag_{(\al,\ud\ls)}\,\big)_{(\al,\ud\ls)\in \ud{\La^\dag}}
	\qquad&&\tn{on}\quad\fM^{\rtf}_{\La^\dag}\,.
\end{alignat*}
Here, recall $\ud{\La^\dag}$ consists of all the pairs $(\al,\ud\ls)$ with $\al\inn A$ and $\ud\ls\inn \ud\bT_\al(\La^\dag)$;
c.f.~Theorem~\ref{Thm:tf_smooth_revert}~\ref{Cond:smooth_tf_revert}.
Proposition~\ref{Prp:two_treelike_strs} guarantees
the proper transform
$\ov{\PT}_{\La}(\La^\dag)$ is a treelike structure on $\fM^{\tf}_\La$.
In this subsection,
we show $\La^\ddag$ can be proper transformed in a similar way into a treelike structure
\begin{align*}
	\ov{\PT}_{\La}(\La^\ddag)\qquad\tn{on}\quad
	(\fM^{\tf}_\La)^{\rtf}_{\ov{\PT}_{\La}(\La^\dag)}\,.
\end{align*} 
Such operation will provide models for the operations of \S\ref{Subsec:Step3}-\S\ref{Subsec:Step6-9}.

\begin{rmk}
We only consider $(\fM^{\tf}_\La)^{\rtf}_{\ov{\PT}_{\La}(\La^\dag)}$
because this is the scenario that occurs in \S\ref{Sec:genus_2_twisted_fields}.
Statements parallel to Proposition~\ref{Prp:2nd_PT} below should be expected for $(\fM^{\tf}_\La)^{\tf}_{\ov{\PT}_{\La}(\La^\dag)}$, $(\fM^{\rtf}_\La)^{\tf}_{\ud{\PT}_{\La}(\La^\dag)}$, and $(\fM^{\rtf}_\La)^{\rtf}_{\ud{\PT}_{\La}(\La^\dag)}$ as well.
\end{rmk}

By Proposition~\ref{Prp:two_treelike_strs},
$\ov\PT_\La(\La^\dag)$ consists of 
\begin{align*}
	&
	\big(\,\tau^\dag_{\ti\al}=
	\tau_{\al_{(\fB_{\ti\al})}}^\dag\,,\;\be^\dag_{\ti\al}
	=
	\iota_{\al;\fB_{\ti\al}}\!\circ\!
	\be^\dag_{\al_{(\fB_{\ti\al})}}\,\big)\,,\qquad
	\ti\al=\big(\,\al,\,
	\ov\lt\eq(\tau_\al,\ov\bE,\Dm(\ov\lt))\,\big)\in \ov\La\,.
\end{align*}
Here, recall  $\fB_{\ti\al}:=\be_\al\big(\Dm(\ov\lt)\big).$
Notice that
\begin{align*}
	\big(\,\al_{(\fB_{\ti\al})}\,,\,
	\ud\fr\,\big)\,\in\, 
	\ud{\La^\dag}\qquad
	&\forall\ 
	\big(\ti\al\eq(\al,\ov\lt),\,\ud\fr\big)\in \ud{\ov\PT_\La(\La^\dag)}\quad
	\Big(\tn{i.e.}\ 
	\forall~
	\ti\al\inn\ov\La,\ 
	\ud\fr\inn\ud\bT_{\ti\al}\big(\ov\PT_\La(\La^\dag)\big)\Big).
\end{align*}
Consequently, the rooted trees
\begin{align}\label{Eqn:tau_dagg}
	\tau^\ddag_{\wh\al}:=
	\tau^\ddag_{(\,\al_{(\fB_{\ti\al})}\,,\,\ud\fr\,)},
	\qquad
	\wh\al=(\al,\ov\lt,\ud\fr)\in \ud{\ov\PT_\La(\La^\dag)}\,,
\end{align}
are well defined.

In addition, 
if we write $\ud\fr\eq\big(\tau^\dag_{\ti\al},\ud\bF,\Dm(\ud\fr)\big)$,
then
\begin{align*}
	S_{\wh\al}=
	\ud\bF\sqcup\ov\bE\sqcup 
	\Big(S_\al\big\bsl 
	\big(\fB_{\ti\al}\sqcup\be^\dag_{\ti\al}(\Dm(\ud\fr))\big)\Big),
\end{align*}
because by (\ref{Eqn:be'_ti_al}),
the image of $\be^\dag_{\ti\al}$ in $S_{\ti\al}$ is disjoint from $\ov\bE$.
Let 
\begin{equation}\label{Eqn:j_wh_al}
\begin{split}
	&\mfi_{\wh\al}:\,
	S_{(\al_{(\fB_{\ti\al})},\,\ud\fr)}
	=\ud\bF\sqcup \big(S_{\al_{(\fB_{\ti\al})}}
	\bsl \,\be^\dag_{\al_{(\fB_{\ti\al})}}\!
	(\Dm(\ud\fr))\big)
	\,\hookrightarrow\,
	S_{\wh\al}\qquad
	\tn{satifying}\\
	&
	\mfi_{\wh\al}(\fF)=\fF\quad\forall\ 
	\fF\in\ud\bF\,,\qquad
	\mfi_{\wh\al}(e)=\iota_{\al;\fB_{\ti\al}}(e)\quad\forall\ 
	e\in S_{\al_{(\fB_{\ti\al})}}
	\bsl\, \be^\dag_{\al_{(\fB_{\ti\al})}}\!(\Dm(\ud\fr)).
\end{split}\end{equation}
Here, $\iota_{\al;\fB_{\ti\al}}\!:S_{\al_{(\fB_{\ti\al})}}\!\lra\!S_\al\bsl\fB_{\ti\al}$ is the bijection in Corollary~\ref{Crl:M_strata_local}.
The expression of $\be^\dag_{\ti\al}$ in  (\ref{Eqn:be'_ti_al}) implies
\begin{align*}
	\mfi_{\wh\al}
	\big(S_{\al_{(\fB_{\ti\al})}}
	\bsl \,\be^\dag_{\al_{(\fB_{\ti\al})}}\!
	(\Dm(\ud\fr))\big)
	=\iota_{\al;\fB_{\ti\al}}
	\big(S_{\al_{(\fB_{\ti\al})}}
	\bsl \,\be^\dag_{\al_{(\fB_{\ti\al})}}\!
	(\Dm(\ud\fr))\big)
	=
	S_\al\big\bsl 
	\big(\fB_{\ti\al}\sqcup\be^\dag_{\ti\al}(\Dm(\ud\fr))\big).
\end{align*}
We thus set
\begin{align}\label{Eqn:be_dagg}
	\be^\ddag_{\wh\al}:=
	\mfi_{\wh\al}\circ\be^\ddag_{(\al_{(\fB_{\ti\al})},\,\ud\fr)}\,
	:\,
	\tau^\ddag_{\wh\al}\hookrightarrow 
	S_{\wh\al},
\end{align}
whose injectivity follows from that of the two maps on the right-hand side.

\begin{prp}
\label{Prp:2nd_PT}
With $\tau^\ddag_{\wh\al}$ and $\be^\ddag_{\wh\al}$ introduced in (\ref{Eqn:tau_dagg}) and (\ref{Eqn:be_dagg}), respectively,
\begin{align*}
	\ov\PT_\La(\La^\ddag):=
	\big(\,
	\tau^\ddag_{\wh\al}\,,\;
	\be^\ddag_{\wh\al}\,\big)_{\wh\al\,\in\, \ud{\ov\PT_\La(\La^\dag)}}
\end{align*}
is a treelike structure on $(\fM^{\tf}_\La)^{\rtf}_{\ov{\PT}_{\La}(\La^\dag)}$,
called the \ts{secondary RL proper transform} of $\La^\ddag$ \ts{with respect to} $\La$.
\end{prp}

The general structure of the proof of Proposition~\ref{Prp:2nd_PT} is parallel to that of Proposition~\ref{Prp:two_treelike_strs},
but the details are  more complicated.

\begin{proof}[Proof of Proposition~\ref{Prp:2nd_PT}]
Fix 
\begin{alignat*}{2}
&\ti\al=\big(\al,\ov\lt\big)\in \ov\La\,\bsl\{(0,\ov\lt_\bullet)\},\qquad
&&
\tn{where}\quad
\ov\lt=\big(\tau_\al,\ov\bE,\Dm(\ov\lt)\big)\in\ov\bT_\al(\La),
\\
&\wh\al=\big(\ti\al,\ud\fr\big)\in \ud{\ov\PT_\La(\La^\dag)}\,,\qquad
&&
\tn{where}\quad
\ud\fr=\big(\tau^\dag_{\ti\al},\ud\bF,\Dm(\ud\fr)\big)\in\ud\bT_{\ti\al}\big(\ov\PT_\La(\La^\dag)\big),
\\
&
S'\,\subset\, S_{\wh\al}=\ud\bF\sqcup\ov\bE\sqcup 
\Big(S_\al\big\bsl 
\big(\fB_{\ti\al}\sqcup\be^\dag_{\ti\al}(\Dm(\ud\fr))\big)\Big).\hspace{-1.2in}
&&
\end{alignat*}
For conciseness,
we write
\begin{align*}
	\ut\al':=\big(\,\al_{(\fB_{\ti\al})},\,\ud\fr\,\big).
\end{align*}
Then by (\ref{Eqn:be_dagg}), we have
\begin{align*}
	&
	(\be^\ddag_{\wh\al})^{-1}(S')=
	\big(\be^\ddag_{\ut\al'}\big)^{-1}
	\Bigg(\big(S'\!\cap\!\ud\bF\big)\,\sqcup\, 
	\iota_{\al;\fB_{\ti\al}}^{-1}\bigg(S'\!\cap
	\Big(S_\al\big\bsl
	\big(\fB_{\ti\al}\sqcup\be^\dag_{\ti\al}(\Dm(\ud\fr))\big)\Big)\bigg)\Bigg).
\end{align*} 
Let
\begin{alignat*}{2}
	J&:=(S'\cap\ud\bF)\,\sqcup\,(\be^\dag_{\ti\al})^{-1}
	\big(S'\cap\be^\dag_{\ti\al}(\ND(\ud\fr))\big)\qquad&&\subset\bbI(\ud\fr),\\
	E&:=
	\be^\dag_{\ti\al}(\ud\sfE_J)\,\sqcup\, 
	\big(S'\cap (S_{\ti\al}\bsl\be^\dag_{\ti\al}(\tau^\dag_{\ti\al}))\big)&&\subset S_{\ti\al},\\
	J_0&:=
	(E\cap\ov\bE)\,\sqcup\,\be_\al^{-1}
	\big(E\cap\be_\al(\ND(\ov\lt))\big)&&
	\subset\bbI(\ov\lt),\\
	E_0&:=
	\be_{\al}(\ov\sfE_{J_0})\,\sqcup\, 
	\big(E\cap (S_\al\bsl\be_{\al}(\tau_\al))\big)&&\subset S_{\al}.
\end{alignat*}
By Theorem~\ref{Thm:tf_smooth}~\ref{Cond:smooth_Z} and Theorem~\ref{Thm:tf_smooth_revert}~\ref{Cond:smooth_Z_revert}, we then have
\begin{align}\label{Eqn:ti_al_(s)_ddag}
	\wh\al_{(S')}=\big(\,\ti\al_{(E)},\,\ud\fr'\,\big)
	=\big(\,\al_{(E_0)},\,\phi^{-1}_{\al;E_0}(\ov\lt_{(J_0)}),\,\ud\fr'\,\big)\,,\qquad\tn{where}\quad
	\ud\fr':=\phi^{-1}_{\ti\al;E}(\ud\fr_{(J)});
\end{align}
c.f.~(\ref{Eqn:C_J}), (\ref{Eqn:C_J_revert}), and Propositions~\ref{Prp:M_strata_local} and~\ref{Prp:Treelike_contraction} for notation.
We aim to construct an isomorphism
\begin{align}\label{Eqn:phi_ti_al_s_ddag}
	&\phi_{\wh\al;S'}:\,\tau^\ddag_{\wh\al_{(S')}}\lra\,
	\tau^\ddag_{\wh\al}\big\bsl\,
	\big((\be^\ddag_{\wh\al})^{-1}(S')\big)^{\wedge}
	\qquad\tn{
satisfying}\\
	\label{Eqn:phi_ti_al_s'_ddag}
	&\iota_{\wh\al;S'}\circ\be^\ddag_{\wh\al_{(S')}}
	=\be^\ddag_{\wh\al}\circ \phi_{\wh\al;S'}\,.
\end{align}
Here, $\iota_{\wh\al;S'}:S_{\wh\al_{(S')}}\!\lra\!S_{\wh\al}\bsl S'$ is described in Theorem~\ref{Thm:tf_smooth_revert}~\ref{Cond:smooth_Z_revert}.
Since $\wh\al$ and $S'$ are arbitrary,
(\ref{Eqn:phi_ti_al_s_ddag}) and (\ref{Eqn:phi_ti_al_s'_ddag}) will guarantee $\ov\PT_\La(\La^\ddag)$ is a treelike structure.

We start with analyzing $\tau^\ddag_{\wh\al_{(S')}}$. 

On the one hand,
by (\ref{Eqn:ti_al_(s)_ddag}), (\ref{Eqn:tau_dagg}), Corollary~\ref{Crl:M_strata_local},
Lemma~\ref{Lm:RLS_LRS_isom},
and the fact that $\La$ is a treelike structure,
we have 
\begin{align*}
	\tau^\ddag_{\wh\al_{(S')}}=
	\tau^\ddag_{(\al_{(E_0\sqcup\be_\al(\Dm(\ov\lt_{(J_0)}))},\,\ud\fr')},\qquad
	\be^\ddag_{\wh\al_{(S')}}=
	\mfi_{\wh\al_{(S')}}\circ
	\be^\ddag_{(\al_{(E_0\sqcup\be_\al(\Dm(\ov\lt_{(J_0)}))},\,\ud\fr')}\,.
\end{align*}
The expressions of $E_0$ and $J_0$ in the previous paragraph,
along with
Lemma~\ref{Lm:I(t_J)} and the injectivity of $\be_\al$, then implies
\begin{align}\label{Eqn:tau_wh_al}
\tau^\ddag_{\wh\al_{(S')}}=
\tau^\ddag_{\big(\al_{\big(\fB_{\ti\al}\sqcup\,(E\cap(S_\al\bsl\fB_{\ti\al}))\big)},\,
	\phi^{-1}_{\ti\al;E}(\ud\fr_{(J)}\big)},
\quad
\be^\ddag_{\wh\al_{(S')}}=
\mfi_{\wh\al_{(S')}}\circ
\be^\ddag_{\big(\al_{\big(\fB_{\ti\al}\sqcup\,(E\cap(S_\al\bsl\fB_{\ti\al}))\big)},\,
\phi^{-1}_{\ti\al;E}(\ud\fr_{(J)}\big)}.	
\end{align}

On the other hand,
by (\ref{Eqn:j_wh_al}),
we have
\begin{align*}
	\mfi^{-1}_{\wh\al}(S')
	\,=\,(S'\cap \ud\bF)\;\sqcup\; \iota_{\al;\fB_{\ti\al}}^{-1}
	\bigg(S'\cap 
	\Big(S_\al\big\bsl 
	\big(	\fB_{\ti\al}\sqcup\be^\dag_{\ti\al}(\Dm(\ud\fr))\big)\Big)\bigg)
	\quad\big(\subset 
	S_{\ut\al'}\big)\,.
\end{align*}
By
Theorem~\ref{Thm:tf_smooth_revert}~\ref{Cond:smooth_Z_revert},
$\mfi^{-1}_{\wh\al}(S')$
determines $J'\!\subset\!\bbI(\ud\fr)$ and $E'\!\subset\!S_{\al_{(\fB_{\ti\al})}}$,
satisfying
\begin{align*}
	J' 
	&= 
	\big(S'\cap 
	\ud\bF\big)\;\sqcup\;
	(\be^\dag_{\al_{(\fB_{\ti\al})}})^{-1}
	\big(\iota_{\al;\fB_{\ti\al}}^{-1}(S')\cap\be^\dag_{\al_{(\fB_{\ti\al})}}\!(\ND(\ud\fr))\big)
	\\
	&=
	\big(S'\cap 
	\ud\bF\big)\;\sqcup\;
	\big((\be^\dag_{\ti\al})^{-1}(S')\cap\ND(\ud\fr)\big)=J\,,\\
	E'&=
	\be^\dag_{\al_{(\fB_{\ti\al})}}
	(\ud\sfE_J)
	\;\sqcup\;
	\iota_{\al;\fB_{\ti\al}}^{-1}
	\Big(S'\cap
	S_{\ti\al}\bsl
	\big(\fB_{\ti\al}\sqcup\be^\dag_{\ti\al}(\tau^\dag_{\ti\al})\big)\Big).
\end{align*}
Comparing the expressions of $E$ and $E'$,
we obtain
\begin{align*}
	\iota_{\al;\fB_{\ti\al}}(E')= 
	E\cap\big(S_\al\bsl\fB_{\ti\al}\big)
	=: R'.
\end{align*}
Therefore, Theorem~\ref{Thm:tf_smooth_revert}~\ref{Cond:smooth_Z_revert}, Corollary~\ref{Crl:M_strata_local},
and (\ref{Eqn:phi_ti_al}) imply
\begin{equation}\begin{split}\label{Eqn:ut_al'}
	\ut\al'_{\big(\mfi^{-1}_{\wh\al}(S')\big)}=
	\big(\,
	\al_{(\fB_{\ti\al}\sqcup R')}\,,\,
	\phi^{-1}_{\al_{(\fB_{\ti\al})};\,\iota_{\al;\fB_{\ti\al}}^{-1}\!(R')}(\ud\fr_{(J)})\,
	\big) =
	\big(\,
	\al_{(\fB_{\ti\al}\sqcup\,R')}\,,\,
	\phi^{-1}_{\ti\al;E}(\ud\fr_{(J)})\,
	\big)
	\quad
	\big(\subset A^{\ud\der}_{\La^\dag}\big)\,.
\end{split}\end{equation}
Combining the above equality and (\ref{Eqn:tau_wh_al}),
we have
\begin{align}\label{Eqn:tau_ddag}
	\tau^\ddag_{\wh\al_{(S')}}=
	\tau^\ddag_{\;\ut\al'_{(\mfi^{-1}_{\wh\al}(S'))}}\,.
\end{align}

As $\La^\ddag$ is a treelike structure,
there exists an isomorphism
\begin{align*}
	\phi_{\ut\al';\,\mfi^{-1}_{\wh\al}(S')}:\,
	\tau^\ddag_{\;\ut\al'_{(\mfi^{-1}_{\wh\al}(S'))}}\,
	\lra\ 
	&\tau^\ddag_{\ut\al'}\big\bsl\,
	\big((\be^\ddag_{\;\ut\al'})^{-1}\circ \mfi^{-1}_{\wh\al}(S')
	\big)^\wedge =\,\tau^\ddag_{\wh\al}\big\bsl\,
	\big((\be^\ddag_{\wh\al})^{-1}(S')\big)^{\wedge},
\end{align*}
satisfying
\begin{align}\label{Eqn:phi_ut_al'}
	\iota_{\ut\al';\,\mfi^{-1}_{\wh\al}(S')}\circ
	\be^\ddag_{\;\ut\al'_{(\mfi^{-1}_{\wh\al}(S'))}}
	=
	\be^\ddag_{\ut\al'}
	\circ 
	\phi_{\ut\al';\,\mfi^{-1}_{\wh\al}(S')}.
\end{align}
By (\ref{Eqn:tau_ddag}),
we simply take the proposed isomorphism $\phi_{\ti\al;S'}$ of (\ref{Eqn:phi_ti_al_s_ddag}) to be
\begin{align*}
	\phi_{\wh\al;S'}:=\phi_{\ut\al';\,\mfi^{-1}_{\wh\al}(S')}.
\end{align*}

It remains to verify (\ref{Eqn:phi_ti_al_s'_ddag}) for such $\phi_{\ti\al;S'}$.
To this end,
by (\ref{Eqn:tau_wh_al}), (\ref{Eqn:be_dagg}), (\ref{Eqn:phi_ut_al'}), (\ref{Eqn:ut_al'}) and (\ref{Eqn:phi_ti_al_s}),
it suffices to show
\begin{align}\label{Eqn:iota_jmath}
	\iota_{\wh\al;S'}\circ\mfi_{\wh\al_{(S')}}&=
	\mfi_{\wh\al}\circ\iota_{\ut\al';\,\mfi^{-1}_{\wh\al}(S')}:
	\\
	S_{\ut\al'_{(\mfi^{-1}_{\wh\al}(S'))}}
	&=\phi^{-1}_{\ti\al;E}(\ud\bF\bsl S')\sqcup 
	\big(S_{\al_{(\fB_{\ti\al}\sqcup R')}}
	\!\big\bsl 
	\be^\dag_{\al_{(\fB_{\ti\al}\sqcup R')}}\!
	(\tau^\dag_{\al_{(\fB_{\ti\al}\sqcup R')}})\big)\ 
	\hookrightarrow\ 
	S_{\wh\al}\bsl S'.\nonumber
\end{align}
Indeed,
from (\ref{Eqn:j_wh_al}), (\ref{Eqn:ti_al_(s)_ddag}),
Theorem~\ref{Thm:tf_smooth_revert}~\ref{Cond:smooth_Z_revert},
and Corollary~\ref{Crl:M_strata_local},
we conclude that
\begin{alignat*}{1}
	&\iota_{\wh\al;S'}\circ\mfi_{\wh\al_{(S')}}(\fF')
	=\phi_{\ti\al;E}(\fF')
	=\phi_{\al_{(\fB_{\ti\al})};\,E'}(\fF')
	=\mfi_{\wh\al}\circ\iota_{\ut\al';\,\mfi^{-1}_{\wh\al}(S')}(\fF')
	\qquad
	\forall~\fF'\in \phi^{-1}_{\ti\al;E}(\ud\bF\bsl S');\\
	&\iota_{\wh\al;S'}\circ\mfi_{\wh\al_{(S')}}(e)=
	\iota_{\al;\,\fB_{\ti\al}\sqcup R'}(e)
	=
	\mfi_{\wh\al}\circ\iota_{\ut\al';\,\mfi^{-1}_{\wh\al}(S')}(e)\qquad
	\forall~e\in S_{\al_{(\fB_{\ti\al}\sqcup R')}}\!\!
	\big\bsl 
	\be^\dag_{\al_{(\fB_{\ti\al}\sqcup R')}}\!
	(\tau^\dag_{\al_{(\fB_{\ti\al}\sqcup R')}}).
\end{alignat*}
This establishes (\ref{Eqn:iota_jmath}) and completes the proof of Proposition~\ref{Prp:2nd_PT}.
\end{proof}

\subsection{Derived treelike structures}
\label{Subsec:derived}

When we add twisted fields to $\fD_2$,
certain entries of the pullback of $\varphi$ of (\ref{Eqn:str_hom})  (after suitable trivialization of the source and target of $\varphi$) contain some twisted parameters $\xi_\fE$, $\fE\inn\ov\bE$ or $\ud\bE$, of orders higher than one.
Inevitably,
one need to acquire suitable treelike structures whose modular monomials contain such modular parameters as well as those labeled by $\ND(\ov\lt)$ or $\ND(\ud\lt)$.
This leads to the constructions in \S\ref{Subsec:derived} and \S\ref{Subsec:2nd_derived_treelike}.

For every $\ov\lt\eq\big(\tau,\ov\bE,\Dm(\ov\lt)\big)\inn\ov\bT$,
recall  the index set $$\bbI(\ov\lt)= \ov\bE\sqcup\ND(\ov\lt)$$
given in (\ref{Eqn:bbI}) that labels part of the modular parameters in Theorem~\ref{Thm:tf_smooth}~\ref{Cond:smooth_parameters}.

The tree order on $\tau$ further induces a tree order on $\bbI(\ov\lt)$ so that it becomes a rooted tree (i.e.~a poset in the sense of Definition~\ref{Dfn:Rooted_tree}) $\dg{\ov\lt}$ as follows:
$$
 \dg{\ov\lt}:=
 \big(\,\bbI(\ov\lt)\,,\,\preceq\big)\,,
$$
where $\preceq$ is the relation on $\dg{\ov\lt}$ given by the restriction of the partial order (\ref{Eqn:transverse_sections_order}) on $\ov\bE$, the restriction of the tree order (\ref{Eqn:tree_order}) on $\ND(\ov\lt)$,
as well as the following:
for every $e\inn\ND(\ov\lt)$ and $\fE\inn\ov\bE$,
\begin{equation}\begin{split}\label{Eqn:derived_order}
e\prec\fE
\qquad\Longleftrightarrow\qquad
&
e\in\big(\fE\cap\Dm^*(\ov\lt)\big)^{\!\prec},\\
&
\tn{where}\quad \Dm^*(\ov\lt):=\big(\dot\fE_\cht\cap\min(\tau)\big)^\succeq\ \ \subset\Dm(\ov\lt)\,;
\end{split}\end{equation}
c.f.~(\ref{Eqn:E_R}) for notation.
Particularly,
$\fE\!\not\preceq\!e$ for any $e\inn\ND(\ov\lt)$ and $\fE\inn\ov\bE$.

It is a direct check that $\preceq$ is a partial order on $\bbI(\ov\lt)$. 
Moreover, it satisfies (\ref{Eqn:tree_order}).
To see this,
notice that for any $e,e'\inn\ND(\ov\lt)$ and $\fE\inn\ov\bE$ such that 
$
e\!\prec\! e'$ and $e\!\prec\! \fE,$
by (\ref{Eqn:derived_order}),
we can find $e''\inn\fE\!\cap\big(\dot\fE_\cht\!\cap\!\min(\tau)\big)^\succeq$ such that $e\!\prec\!e''$.
Then, (\ref{Eqn:tree_order}) implies $e'$ and $e''$ are comparable.
Moreover, $e''\inn\Dm(\ov\lt)$ implies  $e'\!\not\succeq\!e''$, for otherwise $e'$ would be in $\Dm(\ov\lt)$ due to Lemma~\ref{Lm:dominant}.
Hence $e'\!\prec\!e''$, which in turn implies $e'\!\prec\!\fE$.
(The remaining cases for checking $\preceq$ on $\bbI(\ov\lt)$ satisfies (\ref{Eqn:tree_order}) are straightforward.)
To summarize, we observe that $\preceq$ is a tree order on $\bbI(\ov\lt)$,
i.e.~$\dg{\ov\lt}$ is a rooted tree.

Similarly,
for every $\ud\lt\eq\big(\tau,\ud\bE,\ND(\ud\lt)\big)\inn\ud\bT$,
recall the index set 
$$\bbI(\ud\lt)= \ud\bE\sqcup\ND(\ud\lt)$$
as in (\ref{Eqn:bbI}) that labels part of the modular parameters in Theorem~\ref{Thm:tf_smooth_revert}~\ref{Cond:smooth_parameters_revert}.
On $\bbI(\ud\lt)$,
the relation~$\preceq$ is likewise given by the restriction of the partial order (\ref{Eqn:transverse_sections_order}) on $\ud\bE$, the restriction of the tree order (\ref{Eqn:tree_order}) on $\ND(\ud\lt)$,
as well as the following:
for every $e\inn\ND(\ud\lt)$ and $\fE\inn\ud\bE$,
\begin{equation}\begin{split}\label{Eqn:derived_order_revert}
	e\prec\fE
	\qquad\Longleftrightarrow\qquad
	&e\in \big(\,
	\fE\cap\Dm^*(\ud\lt)\,\big)^{\!\prec},\\
	&\tn{where}\quad
	\Dm^*(\ud\lt):=
	\bigcup_{\fe\in\udt\fE_\cht\cap\max(\tau)}\!\!\!\!\!\!\!\!\!\!\Om(\fe)^\succeq\ \ 
	\subset\Dm(\ud\lt);
\end{split}\end{equation}
c.f.~(\ref{Eqn:Omega_e}) for $\Om(\fe)$.
Particularly,
$\fE\!\not\preceq\!e$ for any $e\inn\ND(\ud\lt)$ and $\fE\inn\ud\bE$.
Mimicking the argument in the RL case,
we observe that $\prec$ is a partial order on $\bbI(\ud\lt)$ satisfying (\ref{Eqn:tree_order}),
hence 
$$
\dg{\ud\lt}:=\big(\,\bbI(\ud\lt)\,,\,\preceq\big)
$$
is a rooted tree. 

Intuitively, 
by calling each root-to-leaf path of $\tau$ whose edges are all dominant a \ts{dominant root-to-leaf path},
we see
$\Dm^*(\ov\lt)$ (resp.~$\Dm^*(\ud\lt)$) is simply the union of the dominant root-to-leaf paths.
For every $e\inn\ND(\ov\lt)$ (resp.~$e\inn\ND(\ud\lt)$),
$e'\inn e^\succ$ (in $\tau$) that belongs to a dominant root-to-leaf path,
and $\fE\inn\ov\bE$ (resp.~$\fE\inn\ud\bE$) containing $e'$,
we have $e\!\prec\!\fE$ in $\dg{\ov\lt}$ (resp.~$\dg{\ud\lt}$),
and vice versa.

\begin{dfn}\label{Dfn:Derived_tree}
We call the rooted tree $\dg{\ov\lt}\inn\bT$ (resp.~$\dg{\ud\lt}\inn \bT$) the \ts{derived tree of} $\ov\lt$ (resp.~$\ud\lt$).
The elements of $\ov\bE$ (resp.~$\ud\bE$) are called the \ts{exceptional edges of} $\dg{\ov\lt}$ (resp.~$\dg{\ud\lt}$).
\end{dfn}

\begin{figure}[htp]
	\begin{center}
		\begin{tikzpicture}
			\filldraw[xshift=6cm]
			(0,.1) circle (1.2pt)
			(0,.4) circle (1.2pt)
			(0,1) circle (1.2pt)
			(-.3,.1) circle (1.2pt)
			(.5,-.1) circle (1.2pt)
			(1.1,-.1) circle (1.2pt)
			(.3,.7) circle (1.2pt);
			
			\draw[xshift=6cm]
			(0,.1)--(0,1)--(1.1,-.1)
			(-.3,.1)--(0,.4)--(.5,-.1);
			\draw[xshift=6cm, dotted]
			(-.5,.1)--(1.2,.1)
			(-.5,.4)--(1.2,.4)
			(-.5,.7)--(1.2,.7)
			(-.5,1)--(1.2,1);
			\draw[thin,decorate,decoration=brace,xshift=6cm]
			(1.2,.95)--(1.2,.75);
			\draw[thin,decorate,decoration=brace,xshift=6cm]
			(1.2,.65)--(1.2,.45);
			\draw[thin,decorate,decoration=brace,xshift=6cm]
			(1.2,.35)--(1.2,.15);
			
			\draw[xshift=6cm]
			(-.1,.55) node {\tiny{$a$}}
			(-.1,.2) node {\tiny{$c$}}
			(-.35,.25) node {\tiny{$b$}}
			(.34,.25) node {\tiny{$d$}}
			(.33,.89) node {\tiny{$f$}}
			(.92,.25) node {\tiny{$g$}}
			(1.2,.85) node[right] {\tiny{$\fE_1$}}
			(1.2,.25) node[right] {\tiny{$\fE_3$}}
			(1.2,.55) node[right] {\tiny{$\fE_2$}};
			
			\filldraw[xshift=11.5cm] 
			(0,.1) circle (1.2pt)
			(0,.4) circle (1.2pt)
			(0,1) circle (1.2pt)
			(-.3,.1) circle (1.2pt)
			;
			
			\draw[xshift=11.5cm]
			(0,.1)--(0,1)
			(-.3,.1)--(0,.4)
			(-.1,.55) node {\tiny{$a$}}
			(-.1,.2) node {\tiny{$c$}}
			(-.35,.25) node {\tiny{$b$}}
			;
			
			\draw[xshift=11.5cm, dotted]
			(-.5,.1)--(.2,.1)
			(-.5,.4)--(.2,.4)
			(-.5,.7)--(.2,.7)
			(-.5,1)--(.2,1)
			;
			
			\draw[loosely dashed, yshift=.3cm]
			(10.7,.9) rectangle (12,-.4);

			\filldraw[xshift=14.5cm] 
			(0,.1) circle (1.2pt)
			(0,.4) circle (1.2pt)
			(0,.7) circle (1.2pt)
			(0,1) circle (1.2pt) 
			(.5,-.1) circle (1.2pt)
			(1.1,-.1) circle (1.2pt)
			;
			
			\draw[xshift=14.5cm]
			(0,.1)--(0,1)
			(0,.4)--(.5,-.1)
			(0,1)--(1.1,-.1)
			(0,.85) node[left] {\tiny{$\fE_1$}}
			(0,.25) node[left] {\tiny{$\fE_3$}}
			(0,.55) node[left] {\tiny{$\fE_2$}}
			(.34,.25) node {\tiny{$d$}}
			(.92,.25) node {\tiny{$g$}}
			; 
			
			\filldraw[xshift=17.5cm] 
			(0,.4) circle (1.2pt)
			(0,.7) circle (1.2pt)
			(0,1) circle (1.2pt) 
			(1.1,-.1) circle (1.2pt)
			;
			
			\draw[xshift=17.5cm]
			(0,.4)--(0,1)--(1.1,-.1)
			(0,.85) node[left] {\tiny{$\fE_1$}}
			(0,.55) node[left] {\tiny{$\fE_2$}}
			(.92,.25) node {\tiny{$g$}}
			; 
			
			\draw[yshift=-.5cm]
			(6.75,0) node {\tiny{$\ov\lt=\big(\tau,\big\{\{a,f\},\{a,g\},\{b,c,d,g\}\big\},\,\{f,a,b,c\}\big)$}}
			(11.35,0) node {\tiny{$\Dm^*(\ov\lt)$}}
			(14.9,0) node {\tiny{$\dg{\ov\lt}$}}
			(17.9,0) node {\tiny{$\vr_{\ov\lt}$}}
			;
		\end{tikzpicture}
	\end{center}
	\caption{A derived tree}\label{Fig:derived}
\end{figure}  

\begin{eg}\label{Eg:derived}
Consider 
$$
\ov\lt=
\big(\,
\tau\,;\ 
\ov\bE\eq \big\{\fE_1\eq\{a,f\},\,\fE_2\eq\{a,g\},\,\fE_3\eq\{b,c,d,g\}\big\}\,;\ 
\Dm(\ov\lt)\eq\{f,a,b,c\}\,\big)\,\in\,\ov\bT\,;
$$
c.f.~the leftmost graph of Figure~\ref{Fig:derived}.
In this case, $\bbI(\ov\lt)\eq\{\fE_1,\fE_2,\fE_3,d,g\}$,
and $\Dm^*(\ov\lt)\eq\{a,b,c\}$.
By~(\ref{Eqn:derived_order}), 
the tree order on $\dg{\ov\lt}$ is thus given by
$
d\!\prec\!\fE_2\prec\!\fE_1$ and $\fE_3\!\prec\!\fE_2\!\prec\!\fE_1.$
An illustration of $\dg{\ov\lt}$ is given in Figure~\ref{Fig:derived}.
\end{eg}

A more delicate example of derived trees can be found in Example~\ref{Eg:2nd_der_treelike}.

\begin{lmm}
	\label{Lm:Derived_min}
	For every $\ov\lt\inn\ov\bT$, we have
	$\fE_{\cht}\inn\min(\dg{\ov\lt})$.
	For every $\ud\lt\inn\ud\bT$, we have
	$\fE_{1}\inn\min(\dg{\ud\lt})$.
\end{lmm}

\begin{proof}
	For $\ov\lt\inn\ov\bT$,
	\begin{align*}
		\big(\dot\fE_\cht\!\cap\!\min(\tau)\big)^\succeq\cap\fE_\cht\,=\,
		\dot\fE_\cht\cap\min(\tau)\,\subset\,
		\min(\tau),
	\end{align*}
	hence the first statement of Lemma~\ref{Lm:Derived_min} follows from (\ref{Eqn:derived_order}).
	
	The argument for $\ud\lt\inn\ud\bT$ is simpler:
	$\fE_1\eq\min(\tau)$ implies
	$\fE_1^\prec\eq\emptyset$,
	hence the last statement of Lemma~\ref{Lm:Derived_min} follows from (\ref{Eqn:derived_order_revert}).
\end{proof}

Let $\phi:\tau\!\lra\!\tau'$ be an isomorphism of rooted trees $\tau$ and $\tau'$.
For every $\ov\lt\eq \big(\tau,\ov\bE,\Dm(\ov\lt)\big)\inn\ov\bT$ (resp.~$\ud\lt\eq \big(\tau,\ud\bE,\Dm(\ud\lt)\big)\inn\ud\bT$), 
let $\phi(\ov\lt)\inn\ov\bT$ (resp.~$\phi(\ud\lt)\inn\ud\bT$) be as in Lemma~\ref{Lm:RLS_LRS_isom}.

\begin{lmm}
	\label{Lm:Derived_isom}
With notation as above,
the bijections
\begin{align*}
&\dg{\ov\lt}\lra\dg{\phi(\ov\lt)},\qquad
\fE\mapsto \phi(\fE)\quad\forall~\fE\inn\ov\bE,\qquad
e\mapsto \phi(e)\quad\forall~e\inn\ND(\ov\lt),\\
&\dg{\ud\lt}\lra\dg{\phi(\ud\lt)},\qquad
\fE\mapsto \phi(\fE)\quad\forall~\fE\inn\ud\bE,\qquad
e\mapsto \phi(e)\quad\forall~e\inn\ND(\ud\lt),
\end{align*}
are respectively isomorphisms between the derived trees, which are still denoted by $\phi$.
\end{lmm}

\begin{proof}
It is a direct check that the maps above preserve the tree orders on $\dg{\ov\lt}$ and $\dg{\phi(\ov\lt)}$, and on $\dg{\ud\lt}$ and $\dg{\phi(\ud\lt)}$, respectively.
We omit the details.
\end{proof}

The idea behind the notion of the derived trees comes from the second and third rounds of the blowups in~\cite{HLN},
where the blowup centers lie within the exceptional divisors of the previous rounds.
A simple example is provided below;
c.f.~Example~\ref{Eg:abcd} for a more thorough example.

\begin{eg}\label{Eg:derived_eqn}
Let $\fM$, $\La$ and $\fM^{\tf}$ be as in Theorem~\ref{Thm:tf_smooth},
and $(\al,\ov\lt)\inn \ov\La$ be such that $\ov\lt$ is the same as in Example~\ref{Eg:derived};
c.f.~the leftmost graph of Figure~\ref{Fig:derived}.
We further assume $S_\al\eq \tau_\al$ and $\be_\al\eq\tn{Id}$.
Let $\cV\inn\fV_\al$ 
and $\cU\inn\fV_{(\al,\ov\lt)}$ be as in Theorem~\ref{Thm:tf_smooth}~\ref{Cond:smooth_parameters}.
Consider the following equation on $\cV\!\times\!\A^4$:
$$
 \ze_a^2\ze_b^2w_1+\ze_a^2\ze_cw_2+\ze^2_a\ze_dw_3+\ze_f\ze_gw_4=0,
$$
where $w_1,\ldots,w_4$ are coordinates on $\A^4$.
Such equations appear as part of the local equations of $\ov M_2(\P^n,d)$; c.f.~\S\ref{Subsec:M2Pnd_loc_eqn}.

By Theorem~\ref{Thm:tf_smooth}~\ref{Cond:smooth_parameters}, 
we obtain the pullback of the equation above to $\cU$:
$$
 (\xi_{\fE_1}\xi_{\fE_2}\xi_{\fE_3})\big(
 \xi_{\fE_1}\xi_{\fE_2}(u_cw_2+\xi_{\fE_3}w_1+\xi_dw_3)+\xi_gw_4\big)=0,
$$
where $u_c$ is a  unit on $\cU$.
After changing variables we rewrite the above equation as
$$
 (\xi_{\fE_1}\xi_{\fE_2}\xi_{\fE_3})\big(
 \xi_{\fE_1}\xi_{\fE_2}w_2'+\xi_gw'_4\big)=0.
$$
Albeit less singular than the original equation, the new equation is still singular.
To further desingularize it,
we need to compare $\xi_{\fE_1}\xi_{\fE_2}$ and $\xi_g$,
hence need the rightmost rooted tree $\vr_{\ov\lt}$ of Figure~\ref{Fig:derived}.
Notice all the data of $\vr_{\ov\lt}$ is encoded in the derived tree $\dg{\ov\lt}$.
\end{eg}

In the remainder of this subsection,
we use the derived trees $\dg{\ov\lt}$ to construct a new treelike structure on $\fM^{\tf}_\La$ that lays the foundation of the construction of $\fM^{\fk 5}/\fM^{\fk 4}$ in~\S\ref{Sec:genus_2_twisted_fields}.

For every $\ti\al\eq(\al,\ov\lt)\inn \ov\La$ as in (\ref{Eqn:ti_al}), let
\begin{equation}\begin{split}\label{Eqn:beta_ti_al}
&
\be_{\ti\al}:\,\dg{\ov\lt}=\ov\bE\sqcup\ND(\ov\lt)\ \hookrightarrow\ 
S_{\ti\al}=
\ov\bE\sqcup\be_\al\big(\ND(\ov\lt)\big)\sqcup \big(S_\al\big\bsl\be_\al(\tau_\al)\big),
\\
&
\be_{\ti\al}|_{\ov\bE} = \tn{Id}_{\ov\bE},\qquad
\be_{\ti\al}|_{\ND(\ov\lt)}=\be_\al|_{\ND(\ov\lt)}\,.
\end{split}\end{equation}
Although it is natural to anticipate 
\begin{align}\label{Eqn:false_derived_tree_str}
 \big(\,\dg{\ov\lt},\,\be_{\ti\al}\,\big)_{\ti\al=(\al,\ov\lt)\in\ov\La}
\end{align}
gives a treelike structure on $\fM^{\tf}_\La$,
Example~\ref{Eg:Derived_CE} in \S\ref{Sec:eg} shows the contrary.
Nevertheless, the derived trees can be slightly modified to yield a treelike structure on $\fM^{\tf}_\La$.

For every $\ov\lt\eq\big(\tau,\ov\bE,\Dm(\ov\lt)\big)\inn\ov\bT$ and $k\inn\Z_{>0}$,
let
\begin{align}\label{Eqn:mu_k}
\ny_{k}(\ov\lt):=
\big\{\,\fE\in\ov\bE\,:\; \big|\,\Dm^*\!(\ov\lt)\cap\fE\,\big|\ge k\,\big\}\,.
\end{align}
The tree order (\ref{Eqn:tree_order}) ensures that
for every $\fE\inn\ny_k(\ov\lt)$ and $\fE'\inn\ov\bE$ with $\fE'\!\prec\!\fE$,
we have $\fE'\inn\ny_k(\ov\lt)$.

\begin{prp}\label{Prp:Derived_treelike_1}
Let $\fM$ and $\La$ be as in Theorem~\ref{Thm:tf_smooth}.
For each $\ti\al\eq(\al,\ov\lt)\inn \ov\La$,
let 
$$
\vr_{\ov\lt}:=
\dg{\ov\lt}\big\bsl\,\big(\ny_{2}(\ov\lt)\big)^{\wedge},\qquad
\be_{1;\ti\al}=\be_{\ti\al}\big|_{\vr_{\ov\lt}}:\vr_{\ov\lt}\hookrightarrow S_{\ti\al},
$$
where $\be_{\ti\al}$ is as in~(\ref{Eqn:beta_ti_al}).
Then,
$$
\ov\dt \La:=\big(\,\varrho_{\ov\lt},\,\be_{1;\ti\al}\,\big)_{\ti\al=(\al,\ov\lt)\in \ov\La}
$$
is a treelike structure on $\fM^{\tf}_\La$,
called the \ts{derived treelike structure with respect to} $\La$.
\end{prp}

\begin{proof}
Fix 
\begin{align*}
	&\ti\al=\big(\al,\ov\lt\big)\in \ov\La\big\bsl\big\{(0,\ov\lt_\bullet)\big\},\qquad\tn{where}\quad
	\ov\lt=(\tau_\al,\ov\bE,\Dm(\ov\lt))
	\in\ov\bT_\al(\La);\\
	&s\,\in\, S_{\ti\al}=\ov\bE\sqcup\be_\al\big(\ND(\ov\lt)\big)\sqcup \big(S_\al\big\bsl\be_\al(\tau_\al)\big).
\end{align*}
Taking $S'\eq\{s\}$ in Theorem~\ref{Thm:tf_smooth}~\ref{Cond:smooth_Z},
we obtain
$$
J=\big(\{s\}\!\cap\!\ov\bE\big)\sqcup\be_\al^{-1}\big(\{s\}\!\cap\!\be_\al(\ND(\ov\lt))\big)\ \subset\bbI(\ov\lt),
\qquad
E=\be_\al(\ov\sfE_J)\sqcup\big(\{s\}\!\cap\!(S_\al\bsl\be_\al(\tau_\al))\big)\ \subset S_\al,
$$
so that
$$
 \ti\al_{(s)}=\big(\,\al_{(E)},\,\phi_{\al;E}^{-1}(\ov\lt_{(J)})\,\big)\,,
$$
where $\phi_{\al;E}\!:\tau_{\al_{(E)}}\!\lra\!\tau_\al\bsl\,(\ov\sfE_J)^\wedge$ is the isomorphism in Proposition~\ref{Prp:M_strata_local}.

In order to show Proposition~\ref{Prp:Derived_treelike_1}, we aim to find an isomorphism 
\begin{align}\label{Eqn:Derived_1}
&\phi_{1}:\,\vr_{\phi_{\al;E}^{-1}(\ov\lt_{(J)})}\lra\vr_{\ov\lt}\big\bsl\,\big(\be_{1;\ti\al}^{-1}(s)\big)^{\wedge}
\qquad
\tn{s.t.}\qquad
\iota_{\ti\al;\,s}\circ\be_{\ti\al_{(s)}}=\be_{\ti\al}\circ\phi_{1}\,.
\end{align}
Obviously, $\phi_1$ depends on $\ti\al$ and $s$, which are fixed in this proof, 
so we do not write them explicitly as subscripts of $\phi_1$.
In (\ref{Eqn:Derived_1}),
as well as in every edge contraction below,
the adjoint is always with respect to the rooted tree from which the adjoint set is contracted.

By Lemma~\ref{Lm:Derived_isom},
we have
$$
\vr_{\phi_{\al;E}^{-1}(\ov\lt_{(J)})}=
\phi_{\al;E}^{-1}\big(\vr_{(\ov\lt_{(J)})}\big)=
\phi_{\al;E}^{-1}\Big(\dg{\ov\lt_{(J)}}\!\big\bsl\,
\big(\ny_2(\ov\lt_{(J)})\big)^{\wedge}\Big)\,.
$$
In addition, by Lemma~\ref{Lm:E^} and (\ref{Eqn:beta_ti_al}),
we have
$$
\vr_{\ov\lt}\big\bsl\,\big(\be_{1;\ti\al}^{-1}(s)\big)^{\wedge}=
\dg{\ov\lt}\big\bsl\,\big(\ny_2(\ov\lt)\sqcup\be_{1;\ti\al}^{-1}(s)\big)^\wedge
=
\dg{\ov\lt}\big\bsl\,\big(\ny_2(\ov\lt)\cup J\big)^\wedge\,.
$$
Therefore, to establish~(\ref{Eqn:Derived_1}), it suffices to show
\begin{align}\label{Eqn:Derived_1'}
\dg{\ov\lt_{(J)}}\!\big\bsl\,
\big(\ny_2(\ov\lt_{(J)})\big)^{\wedge}=
\dg{\ov\lt}\big\bsl\,\big(\ny_2(\ov\lt)\cup J\big)^\wedge,
\end{align}
then by setting 
\begin{align}\label{Eqn:phi_1}
	\phi_{1}:=\phi_{\al;E} 
\end{align}
and applying the last statement of Theorem~\ref{Thm:tf_smooth}~\ref{Cond:smooth_Z} as well as Proposition~\ref{Prp:Treelike_contraction},
we obtain the equality in (\ref{Eqn:Derived_1}).

Below, we prove (\ref{Eqn:Derived_1'}) by considering all the possibilities for $s\in S_{\ti\al}$.

\textbf{Case A}: $s\!\in\!S_\al\bsl\be_\al(\tau_\al)$. 
In this case,
$J\eq\emptyset$ and $E\eq \{s\}$,
so $\ov\lt_{(J)}\eq\ov\lt$,
and (\ref{Eqn:Derived_1'}) holds.

\textbf{Case B}:
$s\inn \be_\al\big(\ND(\ov\lt)\big)$.
In this case, we write $e\!:=\!\be_\al^{-1}(s)$.
Then,
$$J= \{e\}\qquad\tn{and}\qquad E=\be_\al(\ov\sfE_J).$$
By Lemma~\ref{Lm:dominant}~\ref{Cond:dom_complement},
we further divide Case B into two sub-cases:

\begin{itemize}
[leftmargin=*]
\item If $e\inn\fE^\prec_\cht$,
then
\begin{align*}
\ov\lt_{(J)}=\big(\,\tau_\al\bsl\, e^\wedge,\,\ov\bE,\,\Dm(\ov\lt)\,\big)
\qquad\Longrightarrow\qquad
&\dg{\ov\lt_{(J)}}=\dg{\ov\lt}\bsl\,e^{\wedge}\quad\tn{and}\quad
\ny_2\big(\ov\lt_{(J)}\big)=\ny_2\big(\ov\lt\big),
\end{align*}
hence by Lemma~\ref{Lm:E^}, either side of (\ref{Eqn:Derived_1'}) is equal to $\dg{\ov\lt}\big\bsl\,\big(\ny_2(\ov\lt)\!\cup\! \{e\}\big)^\wedge$.

\item If $e\inn\fE_\cht\bsl\dot\fE_\cht$, 
we write
\begin{align*}
\fE':=
\max\{\,\fE\in\ov\bE:\,\fE\not\succ_{\dg{\ov\lt}}\!e\,\},
\end{align*}
where the maximum is with respect to the order (\ref{Eqn:transverse_sections_order}) of $\Xi(\tau_\al)$. 
(Lemma~\ref{Lm:Derived_min} guarantees the existence of $\fE'$.)
For such $e$,
\begin{itemize}[leftmargin=*]
\item 
if $e\!\not\in\!\min(\tau_\al)$ or $\fE'\inn\ny_2(\ov\lt)$,
then 
$$ \dg{\ov\lt_{(J)}}=\dg{\ov\lt}\bsl\{e\},\qquad
 \ny_2(\ov\lt_{(J)})=\ny_2(\ov\lt),\qquad
 \tn{and}\qquad
 J^\wedge\subset \big(\ny_2(\ov\lt)\big)^\wedge\ \tn{in}\ \dg{\ov\lt},
$$
so either side of (\ref{Eqn:Derived_1'}) is equal to $\dg{\ov\lt}\big\bsl\,\big(\ny_2(\ov\lt)\big)^\wedge$;

\item 
if $e\!\in\!\min(\tau_\al)$ and $\fE'\!\not\in\!\ny_2(\ov\lt)$,
then
\begin{align*} &\dg{\ov\lt_{(J)}}=\dg{\ov\lt}\bsl\{e\},\qquad
\ny_2(\ov\lt_{(J)})=\{\,\fE\inn\ov\bE:\,\fE\!\succeq\!\fE'\,\}\supset\ny_2(\ov\lt),\qquad
\tn{and}\\
&J^\wedge= \{\,\fE\inn\ov\bE:\,\fE\!\succeq\!\fE'\,\}^\wedge\supset \big(\ny_2(\ov\lt)\big)^\wedge
\quad\tn{in}\quad\dg{\ov\lt},
\end{align*}
so either side of (\ref{Eqn:Derived_1'}) is equal to $\dg{\ov\lt}\big\bsl\,e^\wedge$.
\end{itemize}
\end{itemize}

\textbf{Case C}:
$s\inn \ov\bE$.
In this case, we write $\fE\!:=\!s$.
Then,
$$J= \{\fE\}\qquad\tn{and}\qquad E=\be_\al(\ov\sfE_J).$$
We further divide Case C into two sub-cases:

\begin{itemize}
	[leftmargin=*]
\item
if $\fE\inn\ny_2(\ov\lt)$,
then 
$$ \dg{\ov\lt_{(J)}}=\dg{\ov\lt}\bsl\,\fE^\wedge,\qquad
\ny_2(\ov\lt_{(J)})=\ny_2(\ov\lt)\bsl\{\fE\},\qquad
\tn{and}\qquad
J\subset\ny_2(\ov\lt),
$$
so by Lemma~\ref{Lm:E^}, either side of (\ref{Eqn:Derived_1'}) is equal to $\dg{\ov\lt}\big\bsl\,\big(\ny_2(\ov\lt)\big)^\wedge$;

\item 
if $\fE\!\not\in\!\ny_2(\ov\lt)$,
then
\begin{align*} &\dg{\ov\lt_{(J)}}=\dg{\ov\lt}\bsl\,\fE^\wedge,\qquad
\ny_2(\ov\lt_{(J)})=\ny_2(\ov\lt),\qquad
\tn{and}\qquad
J^\wedge\cap \big(\ny_2(\ov\lt)\big)^\wedge=\emptyset
\quad\tn{in}\quad\dg{\ov\lt},	
\end{align*}
so by Lemma~\ref{Lm:E^}, either side of (\ref{Eqn:Derived_1'}) is equal to $\dg{\ov\lt}\big\bsl\,\big(\ny_2(\ov\lt)\!\sqcup\!\{\fE\}\big)^\wedge$.
\qedhere
\end{itemize}
\end{proof}

\begin{eg}
In Example~\ref{Eg:derived},
the set $\ny_2(\ov\lt)\eq\{\fE_3\}$ because
\begin{align*}
&
\Big(\big(\dot\fE_3\!\cap\!\min(\tau)\big)^\succeq\Big)\cap\fE_1=
\Big(\big(\dot\fE_3\!\cap\!\min(\tau)\big)^\succeq\Big)\cap\fE_2=
\{a\}\,,\quad
\Big(\big(\dot\fE_3\!\cap\!\min(\tau)\big)^\succeq\Big)\cap\fE_3=\{b,c\}\,.
\end{align*}
Therefore, $\vr_{\ov\lt}\eq\dg{\ov\lt}\bsl\,\{\fE_3\}^\wedge\eq \dg{\ov\lt}\bsl\,\{\fE_3,d\}$;
c.f.~Figure~\ref{Fig:derived} and Example~\ref{Eg:derived_eqn}.
\end{eg}

\begin{rmk}
In the LR case, the derived trees $\dg{\ud\lt}$ can also be modified to give rise to
the $LR$ counterpart $\ud\dt\La$ of Proposition~\ref{Prp:Derived_treelike_1}.
Such operation, however, is not needed for the construction of $\ti\fD_2^{\tn{tf}}$ in \S\ref{Sec:genus_2_twisted_fields},
so we shall not dig further into this topic in this paper.
\end{rmk}

\subsection{Doubly and second-order derived treelike structures}
\label{Subsec:2nd_derived_treelike}

In this subsection,
we construct two treelike structures on 
$(\fM^{\tf}_{\La})^{\rtf}_{\ov\dt \La}$
that provides models for $\fM^{\fk 8}/\fM^{\fk 7}$ and $\fM^{\fk 7}/\fM^{\fk 6}$ in \S\ref{Sec:genus_2_twisted_fields}, respectively.
(The reason for choosing the leaf-to-root direction when adding twisted fields w.r.t.~$\ov\dt \La$ has been explained in \S\ref{Subsec:M2Pnd_loc_eqn}.) 
An example for these two treelike structures is provided in Example~\ref{Eg:2nd_der_treelike} of~\S\ref{Sec:eg}.

We begin with some notation.
For each 
\begin{equation}\label{Eqn:wh_al}
	\begin{split}
	\wh\al=\big(\ti\al,\ud\ls\big)
	=\big(\al,\ov\lt,\ud\ls\big)\in \ud{\ov\dt\La},\qquad
	\tn{where}\quad
	&\ov\lt=\big(\tau,\ov\bE,\Dm(\ov\lt)\big)\in\ov\bT_\al(\La),\\
	&\ud\ls=\big(\vr_{\ov\lt},\ud\bF,\Dm(\ud\ls)\big)\in\ud\bT_{\ti\al}(\ov\dt\La),
	\end{split}
\end{equation}
notice that
$$
J:=\Dm(\ud\ls)\,\subset\,\vr_{\ov\lt}\,\subset\,\bbI(\ov\lt)\,=\,\ov\bE\sqcup\ND(\ov\lt).
$$
Here, we treat $\vr_{\ov\lt}$ as a set instead of a poset.
By Theorem~\ref{Thm:tf_smooth_revert}~\ref{Cond:smooth_parameters_revert}, we have
\begin{align*}
	S_{\wh\al}
	&=
	\ud\bF\sqcup\be_{1;\ti\al}\big(\ND(\ud\ls)\big)\sqcup \big(S_{\ti\al}\bsl\be_{1;\ti\al}(\vr_{\ov\lt})\big)
	\\
	&=
	\ud\bF\sqcup 
	\big(\ov\bE\!\cap\!\ND(\ud\ls)\big)\sqcup 
	\be_\al\big(\ND(\ov\lt)\!\cap\!\ND(\ud\ls)\big)\sqcup
	\ny_2(\ov\lt)\sqcup 
	\be_\al\big(\ny_2(\ov\lt)^\wedge\!\cap\!\ND(\ov\lt)\big)\sqcup \big(S_\al\bsl\be_\al(\tau_\al)\big)\,.
\end{align*}
Here the adjoint in $\ny_2(\ov\lt)^\wedge$ is with respect to $\dg{\ov\lt}$,
yet its intersection with $\ND(\ov\lt)$ is also a subset of $\tau_\al$, hence it is legitimate to consider its image in $S_\al$ under $\be_\al$.

Notice that
\begin{align*}
	\big(\ov\bE\!\cap\!\ND(\ud\ls)\big)\sqcup\ny_2(\ov\lt)=\ov\bE\bsl J\,,\qquad
	\ND(\ov\lt)\!\cap\!\big(\ny_2(\ov\lt)^\wedge\!\sqcup\!\ND(\ud\ls)\big)=\ND(\ov\lt)\bsl J\,.
\end{align*}
The expression of $S_{\wh\al}$ can thus be rewritten as:
\begin{align}\label{Eqn:S_wh_al}
	S_{\wh\al}
	&=
	\underbrace{\ud\bF}_{\subset \Xi(\vr_{\ov\lt})}\sqcup 
	\underbrace{\big(\ov\bE\bsl J\big)}_{\subset \Xi(\tau_\al)}\sqcup 
	\underbrace{\be_\al\big(\ND(\ov\lt)\bsl J\big)}_{\subset\be_\al(\tau_\al)}\sqcup \big(S_\al\bsl\be_\al(\tau_\al)\big)\,.	
\end{align}
Meanwhile, by Lemma~\ref{Lm:I(t_J)}, the underlying set of $\dg{\ov\lt_{(J)}}$, which is $\bbI\big(\ov\lt_{(J)}\big)$, satisfies
\begin{align*}
	\dg{\ov\lt_{(J)}}
	=(\ov\bE\bsl J)\sqcup 
	\big(\ND(\ov\lt)\big\bsl (J\cup \ov\sfE_J^\wedge)\big).
\end{align*}
Therefore,
there exists a natural injection (of sets)
\begin{align}\label{Eqn:beta_wh_al}
	\be_{\wh\al}=\be_{\ti\al}|_{\dg{\ov\lt_{(J)}}}:\,
	\dg{\ov\lt_{(J)}}\hookrightarrow S_{\wh\al}\,,
\end{align}
where $\be_{\ti\al}$ is as in (\ref{Eqn:beta_ti_al}).
In other words,
\begin{align*}
	\be_{\wh\al}|_{\ov\bE\bsl J}=\tn{Id}_{\ov\bE\bsl J},\qquad
	\be_{\wh\al}|_{\ND(\ov\lt)\bsl (J\cup\ov\sfE_J^\wedge)}=\be_\al|_{\ND(\ov\lt)\bsl (J\cup\ov\sfE_J^\wedge)}.
\end{align*}

\begin{prp}\label{Prp:Derived_treelike_2}
	Let $\fM$ and $\La$ be as in Theorem~\ref{Thm:tf_smooth} and $\ov\dt \La$ be as in Proposition~\ref{Prp:Derived_treelike_1}.
	For each $\wh\al\eq(\al,\ov\lt,\ud\ls)\inn \ud{\ov\dt\La}$ as in (\ref{Eqn:wh_al}), 
	let 
	$$
	J:=\Dm(\ud\ls),\qquad
	\tau_{\wh\al}:=
	\dg{\ov\lt_{(J)}}\big\bsl\,\big(\ny_{3}(\ov\lt_{(J)})\big)^{\wedge},\qquad
	\be_{2;\wh\al}:=\be_{\wh\al}\big|_{\tau_{\wh\al}}=\be_{\ti\al}\big|_{\tau_{\wh\al}}:\tau_{\wh\al}\hookrightarrow S_{\wh\al},
	$$
	where $\be_{\wh\al}$ is as in~(\ref{Eqn:beta_wh_al}).
	Then,
	$$
	\dt_2\La:=\big(\,\tau_{\wh\al},\,\be_{2;\wh\al}\,\big)_{\wh\al=(\al,\ov\lt,\ud\ls)\,\in\, \ud{\ov\dt\La}}
	$$
	is a treelike structure on $(\fM^{\tf}_\La)^{\rtf}_{\ov\dt \La}$,
	known as the \ts{doubly derived treelike structure with respect to} $\La$.
\end{prp}

\begin{eg}
	\label{Eg:derived_treelike_2}
	Let $\fM$, $\La$ and $\fM^{\tf}_\La$ be as in Theorem~\ref{Thm:tf_smooth}.
	Assume there exist $(\al,\ov\lt)$ and $(\al,\ov\lt')$ in $\ov\La$ such that $\ov\lt\eq\big(\tau_\al,\ov\bE,\Dm(\ov\lt)\big)$ and $\ov\lt'\eq\big(\tau_\al,\ov\bE,\Dm(\ov\lt')\big)$
	are respectively
	given by the first two graphs of Figure~\ref{Fig:derived_2}.
	In other words,
	$\tau_\al\eq\{b,c\}$ with the trivial tree order,
	and $\ov\bE\eq\{\fE\}\eq\big\{\{b,c\}\big\}$ for both $\ov\lt$ and $\ov\lt'$.
	The only difference is $\Dm(\ov\lt)\eq\{b\}$ while $\Dm(\ov\lt')\eq\{b,c\}$.
	We further assume $S_\al\eq \tau_\al$ and $\be_\al\eq\tn{Id}$.
\begin{figure}
\begin{center}
\begin{tikzpicture}
	\filldraw
	(.117,.22) circle (1pt)
	(.35,.66) circle (1pt)
	(.7,0) circle (1pt);
	\draw
	(.117,.22)--(.35,.66)--(.7,0);
	\draw[dotted]
	(-.2,.22)--(.9,.22)
	(-.2,.66)--(.9,.66);
	\draw[thin,decorate,decoration=brace,xshift=0cm]
	(.95,.61)--(.95,.27);
	\draw
	(.233,.44) node[left] {\tiny{$b$}}
	(.525,.33) node[right] {\tiny{$c$}}
	(.95,.44) node[right] {\tiny{$\fE$}}
	(.35,-.25) node {\tiny{$\ov\lt$}}
	;
	
	\filldraw[xshift=2.3cm]
	(0,0) circle (1pt)
	(.35,.66) circle (1pt)
	(.7,0) circle (1pt);
	\draw[xshift=2.3cm]
	(0,0)--(.35,.66)--(.7,0);
	\draw[xshift=2.3cm, dotted]
	(-.2,0)--(.9,0)
	(-.2,.66)--(.9,.66);
	\draw[thin,decorate,decoration=brace,xshift=2.3cm]
	(.95,.61)--(.95,.05);
	\draw[xshift=2.3cm]
	(.175,.33) node[left] {\tiny{$b$}}
	(.525,.33) node[right] {\tiny{$c$}}
	(.95,.33) node[right] {\tiny{$\fE$}}
	(.35,-.25) node {\tiny{$\ov\lt'$}}
	;
	
	\filldraw[xshift=4.6cm]
	(0,0) circle (1pt)
	(.35,.66) circle (1pt)
	(.7,0) circle (1pt);
	\draw[xshift=4.6cm]
	(0,0)--(.35,.66)--(.7,0);
	\draw[xshift=4.6cm]
	(.175,.33) node[left] {\tiny{$\fE$}}
	(.525,.33) node[right] {\tiny{$c$}}
	(.35,-.25) node {\tiny{$\ga$}}
	;
	
	\filldraw[xshift=6.9cm]
	(0,0) circle (1pt)
	(0,.66) circle (1pt)
	;
	\draw[xshift=6.9cm]
	(0,0)--(0,.66);
	\draw[xshift=6.9cm]
	(.05,.33) node[left] {\tiny{$\fE$}}
	(0,-.25) node {\tiny{$\ga'$}};
	
	\filldraw[xshift=8.5cm]
	(0,0) circle (1pt)
	(.35,.66) circle (1pt)
	(.7,0) circle (1pt);
	\draw[xshift=8.5cm]
	(0,0)--(.35,.66)--(.7,0);
	\draw[ultra thick,xshift=8.5cm]
	(0,0)--(.35,.66);
	\draw[xshift=8.5cm, dotted]
	(-.2,0)--(.9,0)
	(-.2,.66)--(.9,.66);
	\draw[thin,decorate,decoration=brace,xshift=8.5cm]
	(.95,.61)--(.95,.05);
	\draw[xshift=8.5cm]
	(.175,.33) node[left] {\tiny{$\fE$}}
	(.525,.33) node[right] {\tiny{$c$}}
	(.95,.33) node[right] {\tiny{$\fF$}}
	(.35,-.25) node {\tiny{$\ud\ls^1$}}
	;
	
	\filldraw[xshift=10.8cm]
	(0,0) circle (1pt)
	(.35,.66) circle (1pt)
	(.7,0) circle (1pt);
	\draw[xshift=10.8cm]
	(0,0)--(.35,.66)--(.7,0);
	\draw[ultra thick,xshift=10.8cm]
	(0,0)--(.35,.66)
	(.35,.66)--(.7,0);
	\draw[xshift=10.8cm, dotted]
	(-.2,0)--(.9,0)
	(-.2,.66)--(.9,.66);
	\draw[thin,decorate,decoration=brace,xshift=10.8cm]
	(.95,.61)--(.95,.05);
	\draw[xshift=10.8cm]
	(.175,.33) node[left] {\tiny{$\fE$}}
	(.525,.33) node[right] {\tiny{$c$}}
	(.95,.33) node[right] {\tiny{$\fF$}}
	(.35,-.25) node {\tiny{$\ud\ls^2$}}
	;
	
	\filldraw[xshift=13.1cm]
	(0,0) circle (1pt)
	(.35,.66) circle (1pt)
	(.7,0) circle (1pt);
	\draw[xshift=13.1cm]
	(0,0)--(.35,.66)--(.7,0);
	\draw[ultra thick,xshift=13.1cm]
	(.35,.66)--(.7,0);
	\draw[xshift=13.1cm, dotted]
	(-.2,0)--(.9,0)
	(-.2,.66)--(.9,.66);
	\draw[thin,decorate,decoration=brace,xshift=13.1cm]
	(.95,.61)--(.95,.05);
	\draw[xshift=13.1cm]
	(.175,.33) node[left] {\tiny{$\fE$}}
	(.525,.33) node[right] {\tiny{$c$}}
	(.95,.33) node[right] {\tiny{$\fF$}}
	(.35,-.25) node {\tiny{$\ud\ls^3$}}
	;
\end{tikzpicture}		
\end{center}
\caption{Examples of doubly derived trees}\label{Fig:derived_2}		
\end{figure}

Let $\ga\eq\{\fE,c\}$ and $\ga'\eq\{\fE\}$ be the rooted trees respectively given in Figure~\ref{Fig:derived_2}. 
It is direct to check that $\dg{\ov\lt}\eq\vr_{\ov\lt}\eq\ga$, $\dg{\ov\lt'}\eq\ga'$ and $\vr_{\ov\lt'}\eq\tau_\bullet$.
There are three possible LRS with dominant edges based on $\vr_{\ov\lt}$:
\begin{align*}
	&\ud\ls^i=\big(\ga,\{\fF\},\Dm(\ud\ls^i)\big),\quad
	i\eq 1,2,3,\qquad
	\tn{where}\\
	&
	\fF=\{\fE,c\},\quad
	\Dm(\ud\ls^1)=\{\fE\},\quad
	\Dm(\ud\ls^2)=\{\fE,c\},\quad
	\Dm(\ud\ls^3)=\{c\}.
\end{align*}
Based on $\vr_{\ov\lt'}$,
the only LRS with dominant edges is the trivial one: $\ud\lt_\bullet$.

We are ready to analyze the rooted trees assigned by $\dt_2\La$ to the strata of $(\fM^{\tf}_\La)^{\rtf}_{\ov\dt \La}$ labeled by the following indices:
\begin{align*}
	\wh\al^i:=\big(\al,\ov\lt,\ud\ls^i\big),\ \ 
	i=1,2,3,\quad
	\tn{and}\quad
	\wh\al':=\big(\al,\ov\lt',\ud\lt_\bullet\big)\quad
	\in\ud{\ov\dt\La}.
\end{align*}  
For $1\!\le\!i\!\le\!3$,
let $J^i\!:=\!\Dm(\ud\ls^i)$.
It is then direct to check $
	\ov\lt_{(J^1)}\eq\ov\lt_{(J^2)}\eq\ov\lt_\bullet,$ 
	$\ov\lt_{(J^3)}\eq \ov\lt',$ and
	$\ny_3\big(\ov\lt_{(J^3)}\big)\eq 
	\ny_3\big(\ov\lt'\big)\eq \emptyset.$
Therefore, 
\begin{align}\label{Eqn:tau_wh_al_eg}
	\tau_{\wh\al^1}=\tau_{\wh\al^2}= \tau_\bullet,\qquad
	\tau_{\wh\al^3}= \tau_{\wh\al'}=\ga'\,.
\end{align}

Notice that $\wh\al^3_{(\fF)}\eq\wh\al'$ and $\wh\al^3_{(\fE)}\eq\wh\al^2.$
Since $\dt_2\La$ is claimed to be a treelike structure in Proposition~\ref{Prp:Derived_treelike_2},
one should particularly expect 
\begin{align*}
	\tau_{\wh\al^3}=\tau_{\wh\al'},\qquad
	\tau_{\wh\al^3}\bsl\{\fE\}^\wedge=\tau_{\wh\al^2},
\end{align*}
which are indeed true because of (\ref{Eqn:tau_wh_al_eg}).
\end{eg}

\begin{proof}[Proof of Proposition~\ref{Prp:Derived_treelike_2}]
Throughout the proof,
we fix $\wh\al\inn\ud{\ov\dt\La}$ as in (\ref{Eqn:wh_al}).
According to	Theorem~\ref{Thm:tf_smooth}~\ref{Cond:smooth_Z} and Theorem~\ref{Thm:tf_smooth_revert}~\ref{Cond:smooth_Z_revert},
each $S'\!\subset\! S_{\wh\al}$ uniquely determines
\begin{alignat*}{2}
	J'&:=(S'\cap\ud\bF)\,\sqcup\,\be^{-1}_{\ti\al}
	\big(S'\cap\be_{\ti\al}(\ND(\ud\ls))\big)\qquad&&\subset\bbI(\ud\ls),\\
	E'&:=
	\be_{\ti\al}(\ud\sfE_{J'})\,\sqcup\, 
	\big(S'\cap (S_{\ti\al}\bsl\be_{\ti\al}(\vr_{\ov\lt}))\big)&&\subset S_{\ti\al},\\
	J_0&:=
	(E'\cap\ov\bE)\,\sqcup\,\be_\al^{-1}
	\big(E'\cap\be_\al(\ND(\ov\lt))\big)&&
	\subset\bbI(\ov\lt),\\
	E_0&:=
	\be_{\al}(\ov\sfE_{J_0})\,\sqcup\, 
	\big(E'\cap (S_\al\bsl\be_{\al}(\tau_\al))\big)&&\subset S_{\al},
\end{alignat*}
such that
\begin{align*}
	\wh\al_{(S')}=\big(\,\ti\al_{(E')},\,\ud\ls'\,\big)
	=\big(\,\al_{(E_0)},\,\ov\lt',\,\ud\ls'\,\big)\,,\qquad\tn{where}\quad
	\ov\lt':=\phi^{-1}_{\al;E_0}(\ov\lt_{(J_0)}),\quad
	\ud\ls':=\phi^{-1}_{\ti\al;E'}(\ud\ls_{(J')}).
\end{align*}

In order to show Proposition~\ref{Prp:Derived_treelike_2},
we aim to find an isomorphism 
\begin{equation}\begin{split}\label{Eqn:Derived_2}
	&\phi_{2}:\,\dg{\ov\lt'_{(\Dm(\ud\ls'))}}\!\big\bsl\, \big(\ny_{3}(\ov\lt'_{(\Dm(\ud\ls'))})\big)^{\wedge}
	\lra
	\Big(\dg{\ov\lt_{(J)}}\!\big\bsl\big(\ny_{3}(\ov\lt_{(J)})\big)^{\wedge}\Big)
	\Big\bsl\,\big(\be_{2;\wh\al}^{-1}(S')\big)^{\wedge}
	\\
	&
	\tn{satisfying}\qquad
	\iota_{\wh\al;\,S'}\circ\be_{2;\wh\al_{(S')}}=\be_{2;\wh\al}\circ\phi_{2}	
\end{split}\end{equation}	
for each $S'\!\subset\! S_{\wh\al}$ with $|S'|\eq 1$.
By Lemma~\ref{Lm:Derived_isom} and (\ref{Eqn:phi_1}) (where $\phi_{\al;E_0}$ is denoted by $\phi_1$),
we have
$$
\dg{\ov\lt'_{(\Dm(\ud\ls'))}}\!\big\bsl\, \big(\ny_{3}(\ov\lt'_{(\Dm(\ud\ls'))})\big)^{\wedge}
=
\phi_{\al;E_0}^{-1}\bigg(\Big(\dg{\big((\ov\lt_{(J_0)})_{(\Dm(\ud\ls_{(J')}))}\big)}\Big)\Big\bsl\,
\Big(\ny_3\big((\ov\lt_{(J_0)})_{(\Dm(\ud\ls_{(J')}))}\big)\Big)^{\!\!\wedge}\bigg)\,.
$$
In addition, by Lemma~\ref{Lm:E^},
we have
$$
\Big(\dg{\ov\lt_{(J)}}\!\big\bsl\big(\ny_{3}(\ov\lt_{(J)})\big)^{\!\wedge}\Big)
\Big\bsl\,\big(\be_{2;\wh\al}^{-1}(S')\big)^{\!\wedge}
=
\dg{\ov\lt_{(J)}}\big\bsl\,
\big(\ny_{3}(\ov\lt_{(J)})\sqcup \be_{2;\wh\al}^{-1}(S') \big)^{\!\wedge}
\,.
$$
Therefore, to establish~(\ref{Eqn:Derived_2}), it suffices to show
\begin{align}\label{Eqn:Derived_2'}
	\Big(\dg{\big((\ov\lt_{(J_0)})_{(\Dm(\ud\ls_{(J')}))}\big)}\Big)\Big\bsl\,
	\Big(\ny_3\big((\ov\lt_{(J_0)})_{(\Dm(\ud\ls_{(J')}))}\big)\Big)^{\!\!\wedge}
	=
	\dg{\ov\lt_{(J)}}\big\bsl\,
	\big(\ny_{3}(\ov\lt_{(J)})\sqcup \be_{2;\wh\al}^{-1}(S') \big)^{\!\wedge}
	\,,
\end{align}
then by setting $\phi_{2}\!:=\!\phi_{\al;E_0}$
and applying the last statements of Theorems~\ref{Thm:tf_smooth}~\ref{Cond:smooth_Z} and~\ref{Thm:tf_smooth_revert}~\ref{Cond:smooth_Z_revert}, respectively, as well as Proposition~\ref{Prp:Treelike_contraction},
we obtain the equality in (\ref{Eqn:Derived_2}).

Below, we prove (\ref{Eqn:Derived_2'}) by considering all the possibilities for $S'\!\subset\! S_{\wh\al}$ with $|S'|\eq 1$,
based on the expression (\ref{Eqn:S_wh_al}) of $S_{\wh\al}$.

\textbf{Case A}: $S'\eq\{e\}\inn S_\al\bsl\be_\al(\tau_\al)$.
In this case, $J'\eq J_0\eq \emptyset$ and
$\be_{2;\wh\al}^{-1}(S')\eq\emptyset$, so (\ref{Eqn:Derived_2'}) holds trivially.

\textbf{Case B.i}: $S'\eq\{\be_\al(e)\}$ for some $e\inn\ND(\ov\lt)\!\cap\!\ND(\ud\ls).$
In this case, 
\begin{align*}
	J'=\{e\},\qquad
	J_0=\begin{cases}
		 \{e\}& \tn{if}~e\inn\fF_\cht^\succ,\\
		\emptyset & \tn{otherwise},
	\end{cases}
	\qquad
	\be^{-1}_{2;\wh\al}(S')=
	\begin{cases}
		\{e\}& \tn{if}~e\!\not\in\!\ov\sfE_J^\wedge,\\
		\emptyset & \tn{if}~e\!\in\!\ov\sfE_J^\wedge.
	\end{cases}
\end{align*}
Therefore, 
\begin{align}\label{Eqn:t_J_0_dom}
	(\ov\lt_{(J_0)})_{(\Dm(\ud\ls_{(J')}))}
	=
	\ov\lt_{(J\sqcup\{e\})}
	=
	\begin{cases}
		(\ov\lt_{(J)})_{(e)}& \tn{if}~e\!\not\in\!\ov\sfE_J^\wedge,\\
		\ov\lt_{(J)} & \tn{if}~e\!\in\!\ov\sfE_J^\wedge.
	\end{cases}
\end{align}
\begin{itemize}[leftmargin=*]
\item If $e\!\in\!\ov\sfE_J^\wedge$,
both sides of (\ref{Eqn:Derived_2'}) are equal to $\dg{\ov\lt_{(J)}}\big\bsl
\big(\ny_{3}(\ov\lt_{(J)})\big)^\wedge$.

\item 
If $e\!\not\in\!\ov\sfE_J^\wedge$,
then $e\inn\ND(\ov\lt_{(J)})$.
There are two sub-cases.

\begin{itemize}[leftmargin=*]
\item 
If
$e\inn\big(\fE_{\cht}(\ov\lt_{(J)})\big)^\prec$ or $e\inn \big(\fE_{\cht}(\ov\lt_{(J)})\bsl\min(\tau_{(J)})\big)$, then
\begin{align*}
	\dg{\big((\ov\lt_{(J)})_{(e)}\big)}=
	\dg{\ov\lt_{(J)}}\bsl e^\wedge,\qquad
	\ny_3 \big((\ov\lt_{(J)})_{(e)}\big)=
	\ny_3 (\ov\lt_{(J)}),
\end{align*}
thus (\ref{Eqn:Derived_2'}) holds.

\item 
If $e\inn\fE_{\cht}(\ov\lt_{(J)})\!\cap\!\min(\tau_{(J)})$, observe  the construction of the underlying tree $\vr_{\ov\lt}$ of $\ud\ls$ in Proposition~\ref{Prp:Derived_treelike_1} ensures that 
\begin{align*}
	\ny_2\big(\ov\lt_{(J)}\big)=\ov\bE\bsl J.
\end{align*} 
In other words,
every element of $\ov\bE(\ov\lt_{(J)})$ contains at least two edges of $\tau_{(J)}$ contained in $\big(\dot\fE_\cht(\ov\lt_{(J)})\!\cap\!\min(\tau_{(J)})\big)^\succeq$.
Therefore,
\begin{align*}
	\ny_3 \big((\ov\lt_{(J)})_{(e)}\big)=
	\ny_3 (\ov\lt_{(J)})\cup 
	\ov{\mathbf{H}},\quad
	\tn{where}\quad
	\ov{\mathbf{H}}:=
	\big\{\,
	\fE\inn\ov\bE\bsl J:\,
	\fE\!\cap\!\Dm^*(\ov\lt_{(J)})\!\cap\!\{e\}^\succeq\!=\!\emptyset\,\big\}.
\end{align*}
Meanwhile,
notice that the tree orders on $\dg{\big((\ov\lt_{(J)})_{(e)}\big)}$ and $\dg{\ov\lt_{(J)}}\bsl\{e\}$ can only be different on $\ov{\mathbf{H}}^\wedge$;
moreover,
notice that $\ov{\mathbf{H}}^\wedge\eq\{e\}^\wedge$ in $\ov\lt_{(J)}$. 
Consequently,
\begin{align*}
	\tn{LHS~of}~(\ref{Eqn:Derived_2'})
	&=
	\dg{\big((\ov\lt_{(J)})_{(e)}\big)}\big\bsl\,
	\big(\ny_3 (\ov\lt_{(J)})\!\cup\! 
	\ov{\mathbf{H}}\big)^\wedge
	=
	\dg{\ov\lt_{(J)}}\Big\bsl\,
	\Big(\big(\ny_3 (\ov\lt_{(J)})\!\cup\! 
	\ov{\mathbf{H}}\big)^{\!\wedge}\cup \{e\}\Big)
	\\
	&=
	\dg{\ov\lt_{(J)}}\big\bsl\,
	\big(\ny_3 (\ov\lt_{(J)})\!\cup\! \{e\}\big)^{\!\wedge}
	=
	\tn{RHS~of}~(\ref{Eqn:Derived_2'}).
\end{align*}
\end{itemize} 
\end{itemize}

\textbf{Case B.ii}: $S'\eq\{\be_\al(e)\}$ for some $e\in\ND(\ov\lt)\cap\ny_2(\ov\lt)^{\!\wedge}.$
In this case, 
\begin{align*}
	J'=\emptyset,\qquad
	J_0=\{e\},
	\qquad
	\be^{-1}_{2;\wh\al}(S')=
	\begin{cases}
		\{e\}& \tn{if}~e\!\not\in\!\ov\sfE_J^\wedge,\\
		\emptyset & \tn{if}~e\!\in\!\ov\sfE_J^\wedge,
	\end{cases}
\end{align*}
so we still have (\ref{Eqn:t_J_0_dom}).
Therefore,
the argument of Case B.i applies to the current situation verbatim.

\textbf{Case C.i}: $S'\eq\{\fE\}\!\subset\!\ny_2(\ov\lt)$.
In this case, 
\begin{align*}
	J'=\emptyset,\qquad
	J_0=\{\fE\}\qquad\Longrightarrow\qquad
	(\ov\lt_{(J_0)})_{(\Dm(\ud\ls_{(J')}))}
	=
	\ov\lt_{(J\sqcup\{\fE\})}\,.
\end{align*}
\begin{itemize}[leftmargin=*]
	\item 
	If $\fE\inn\ny_2(\ov\lt)\!\cap\!\ny_3(\ov\lt_{(J)})$,
	then 
	\begin{align*}
	\be^{-1}_{2;\wh\al}(S')=\emptyset,\qquad
	\dg{\ov\lt_{(J\sqcup\{\fE\})}}\!=\dg{\ov\lt_{(J)}}\!\big\bsl\{\fE\}^\wedge,\qquad
	\ny_3(\ov\lt_{(J\sqcup\{\fE\})})=
	\ny_3(\ov\lt_{(J)})\bsl\{\fE\}=
	\ny_3(\ov\lt_{(J)})\big\bsl\{\fE\}^\wedge,
	\end{align*}
	where the last adjoint of $\{\fE\}$ is taken in $\dg{\ov\lt_{(J)}}$.
	Therefore,
	by Lemma~\ref{Lm:E^} we have
	\begin{align*}
		\dg{\ov\lt_{(J\sqcup\{\fE\})}}\big\bsl\,
		\big(\ny_3(\ov\lt_{(J\sqcup\{\fE\})})\big)^{\!\wedge}
		=
		\big(\dg{\ov\lt_{(J)}}\big\bsl\{\fE\}^\wedge\big)
		\Big\bsl
		\big(\ny_3(\ov\lt_{(J)})\big\bsl\{\fE\}^\wedge\big)^{\!\wedge}
		=\dg{\ov\lt_{(J)}}\big\bsl
		\big(\ny_3(\ov\lt_{(J)})\big)^\wedge.
	\end{align*}
	\item If $\fE\inn\ny_2(\ov\lt)\bsl\ny_3(\ov\lt_{(J)})$,
	then 
	\begin{align*}
		\be^{-1}_{2;\wh\al}(S')=\{\fE\},\qquad
		\dg{\ov\lt_{(J\sqcup\{\fE\})}}\!=\dg{\ov\lt_{(J)}}\!\big\bsl\{\fE\}^\wedge,\qquad
		\ny_3(\ov\lt_{(J\sqcup\{\fE\})})=
		\ny_3(\ov\lt_{(J)}),
	\end{align*}
	hence
	\begin{align*}
		\dg{\ov\lt_{(J\sqcup\{\fE\})}}\big\bsl\,
		\big(\ny_3(\ov\lt_{(J\sqcup\{\fE\})})\big)^{\!\wedge}
		=
		\big(\dg{\ov\lt_{(J)}}\bsl\{\fE\}^\wedge\big)
		\big\bsl
		\big(\ny_3(\ov\lt_{(J)})\big)^{\!\wedge}
		=\dg{\ov\lt_{(J)}}\!\big\bsl\,
		\big(\ny_3(\ov\lt_{(J)})\sqcup\{\fE\}\big)^\wedge.
	\end{align*}
\end{itemize}

\textbf{Case C.ii}:
$S'\eq\{\fE\}\!\subset\!\ov\bE\!\cap\!\ND(\ud\ls)$.
In this case, 
\begin{align*}
	J'=\{\fE\},\qquad
	J_0=\begin{cases}
		\{\fE\} & \tn{if}~\fE\inn\fF_\cht^\succ,\\
		\emptyset & \tn{otherwise},
	\end{cases}
	\qquad\Longrightarrow\qquad
	(\ov\lt_{(J_0)})_{(\Dm(\ud\ls_{(J')}))}
	=
	\ov\lt_{(J\sqcup\{\fE\})}\,.
\end{align*}
Thus, we have $\dg{\ov\lt_{(J\sqcup\{\fE\})}}\!=\dg{\ov\lt_{(J)}}\!\big\bsl\{\fE\}^\wedge$ and
\begin{align*}
	\ny_3(\ov\lt_{(J\sqcup\{\fE\})})=
	\begin{cases}
		\ny_3(\ov\lt_{(J)})\bsl\{\fE\} & 
		\tn{if}~\fE\inn \ny_3(\ov\lt_{(J)}),\\
		\ny_3(\ov\lt_{(J)}) & 
		\tn{if}~\fE\!\not\in\! \ny_3(\ov\lt_{(J)});
	\end{cases}\qquad
	\be^{-1}_{2;\wh\al}(S')=
	\begin{cases}
		\emptyset &\tn{if}~\fE\inn \ny_3(\ov\lt_{(J)}),\\
		\{\fE\} & 
		\tn{if}~\fE\!\not\in\! \ny_3(\ov\lt_{(J)}).
	\end{cases}
\end{align*}
Therefore,
\begin{align*}
	\dg{\ov\lt_{(J\sqcup\{\fE\})}}\big\bsl\,
	\big(\ny_3(\ov\lt_{(J\sqcup\{\fE\})})\big)^{\!\wedge}
	&=
	\begin{cases}
	\dg{\ov\lt_{(J)}}\big\bsl
	\big(\ny_3(\ov\lt_{(J)})\big)^\wedge
	&\tn{if}~\fE\inn \ny_3(\ov\lt_{(J)}),\\
	\dg{\ov\lt_{(J)}}\!\big\bsl\,
	\big(\ny_3(\ov\lt_{(J)})\sqcup\{\fE\}\big)^\wedge
	& 
	\tn{if}~\fE\!\not\in\! \ny_3(\ov\lt_{(J)}),
	\end{cases}
	\\
	&=
	\dg{\ov\lt_{(J)}}\!\big\bsl\,
	\big(\ny_3(\ov\lt_{(J)})\sqcup\be^{-1}_{2;\wh\al}(S')\big)^{\!\wedge}\,.
\end{align*}

\textbf{Case D}:
$S'\eq\{\fF\}\!\subset\!\ud\bF$.
In this case, 
\begin{align*}
	J'=\{\fF\},\quad
	J_0=\ud\sfE_{\{\fF\}}
	\qquad\Longrightarrow\qquad
	(\ov\lt_{(J_0)})_{(\Dm(\ud\ls_{(J')}))}
	=
	\ov\lt_{(J)}\,,\quad 
	\be^{-1}_{2;\wh\al}(S')=\emptyset\,.
\end{align*}
Consequently, both sides of (\ref{Eqn:Derived_2'}) are equal to $\dg{\ov\lt_{(J)}}\big\bsl
\big(\ny_{3}(\ov\lt_{(J)})\big)^\wedge$.
\end{proof}

Next, we introduce another treelike structure on $(\fM^{\tf}_{\La})^{\rtf}_{\ov\dt \La}$ as follows.
For each $\wh\al$ as in (\ref{Eqn:wh_al}),
let
\begin{alignat*}{2}
	&\ny_{2}^{\ov\lt}(\ud\ls):=
	\big\{\,\fF\in\ud\bF:\,
	\big(\,\fF\cap\Dm^*\!(\ud\ls)\,\big)\big\bsl\ov\bE\ne\emptyset\,
	\big\}&&\subset\ud\bF\subset\dg{\ud\ls}\,,\\
	&
	\nu^{\ov\lt}(\ud\ls):= 
	\big\{\,
	e\in\ND(\ud\ls)\cap\fF_\cht^\prec:\,
	\{e\}^{\succ}\bsl\ov\bE\ne\emptyset\,\big\}\qquad&&\subset\ND(\ud\ls)\subset\dg{\ud\ls}\,.
\end{alignat*}
Here, we emphasize $\{e\}^{\succ}$ is taken within $\vr_{\ov\lt}$, the underlying rooted tree of $\ud\ls$.
We then set
\begin{align*}
\vs_{\wh\al}:= 
\dg{\ud\ls}\big\bsl\,
\big(\,\ny_{2}^{\ov\lt}(\ud\ls)\cup 
\nu^{\ov\lt}(\ud\ls)\,\big)^{\!\wedge}\,.
\end{align*}
Notice that as a set,
\begin{align*}
	\vs_{\wh\al}
	\,\subset\,
	\bbI(\ud\ls)\bsl \ov\bE 
	\,=\,
	\ud\bF\sqcup \big(\ND(\ud\ls)\bsl\ov\bE\big)
	\,\subset\, 
	\ud\bF\sqcup \big(\vr_{\ov\lt}\bsl\ov\bE\big)
	\,\subset\, 
	\ud\bF\sqcup\ND(\ov\lt)\,.
\end{align*}
Therefore,
there exists a unique injection
\begin{align*}
	&\be'_{\wh\al}:\,\vs_{\wh\al}\hookrightarrow S_{\wh\al}\qquad
	\tn{satisfying}\qquad
	\be'_{\wh\al}|_{\vs_{\wh\al}\cap \ud\bF}= 
	\tn{Id}|_{\vs_{\wh\al}\cap \ud\bF}\,,\quad
	\be'_{\wh\al}|_{\vs_{\wh\al}\cap\ND(\ov\lt)}=\be_\al|_{\vs_{\wh\al}\cap\ND(\ov\lt)}
	\,.
\end{align*}

\begin{prp}\label{Prp:Derived_treelike_second_order_1}
With notation as above,
\begin{align*}
	\ud\partial\ov\dt \La:=\big(\,\vs_{\wh\al}\,,\,\be'_{\wh\al}\,\big)_{\wh\al=(\al,\ov\lt,\ud\ls)\,\in\,\ud{\ov\dt\La}}
\end{align*}
is a treelike structure on $(\fM^{\tf}_{\La})^{\rtf}_{\ov\dt\La}$,		known as the \ts{second-order derived treelike structure with respect to} $\La$.
\end{prp}

\begin{rmk}
	The $\ud\partial$ notation in Proposition~\ref{Prp:Derived_treelike_second_order_1} suggests the treelike structure is ``partially'' derived from $\ov\dt\La$, in that $\ny_2(\ud\ls)$,
	if defined mimicking (\ref{Eqn:mu_k}) with $\ov\lt$ replaced by $\ud\ls$, is a subset of $\ny_2^{\ov\lt}(\ud\ls)$.
	Therefore, 
	$\vs_{\wh\al}$ is an edge contraction from $\dg{\ud\ls}\big\bsl\big(\ny_2(\ud\ls)\big)^{\!\wedge}$. 
\end{rmk}

\begin{eg}
	\label{Eg:2nd_derived_treelike}
	Let $\fM$, $\La$, $\al$, $\ov\lt$, $\ov\lt'$, $\ud\ls^i$ and $\wh\al^i$, $1\!\le\!i\!\le\!3$, and $\wh\al'$ be the same as in Example~\ref{Eg:derived_treelike_2}.
	Set 
	\begin{align*}
		\ti\ga:=\{\fF,c\},\quad \ti\ga':=\{\fF\}, \quad 
		\ti\ga'':=\{\fF,\fE\}\quad\in\,\bT
	\end{align*}
	so that the tree orders on them are all trivial.
	In other words,
	$\ti\ga$ and $\ti\ga''$ are the rooted tree $\ga$ in Figure~\ref{Fig:derived_2} with $\fE$ and $c$  respectively replaced by $\fF$,
	and $\ti\ga'$ is $\ga'$ in Figure~\ref{Fig:derived_2} with $\fE$ replaced by $\fF$.
	
	It is a direct check that
	\begin{align*}
		&
		\dg{\ud\ls^1}=\ti\ga,\quad
		\dg{\ud\ls^2}=\ti\ga',\quad
		\dg{\ud\ls^3}=\ti\ga'';\\
		&
		\ny^{\ov\lt}_2(\ud\ls^2)=\ny^{\ov\lt}_2(\ud\ls^3)=\{\fF\},\quad
		\ny^{\ov\lt}_2(\ud\ls^1)=
		\nu^{\ov\lt}_2(\ud\ls^1)=
		\nu^{\ov\lt}_2(\ud\ls^2)=
		\nu^{\ov\lt}_2(\ud\ls^3)=\emptyset.
	\end{align*}
	Therefore, $\vs_{\wh\al^1}\eq\ti\ga$ and $\vs_{\wh\al^2}\eq\vs_{\wh\al^3}\eq\vs_{\wh\al'}\eq \tau_\bullet$.
	
	Notice that $\wh\al^1_{(c)}\eq\wh\al^2$ and $\wh\al^1_{(\fF)}\eq(\al_{(b)},\ov\lt_\bullet,\ud\lt_\bullet)$.
	Since $\ud\partial\ov\dt\La$ is claimed to be a treelike structure in Proposition~\ref{Prp:Derived_treelike_second_order_1},
	one should particularly expect 
	\begin{align*}
		\vs_{\wh\al^1}\bsl\, c^\wedge=\vs_{\wh\al^2},\qquad
		\vs_{\wh\al^1}\bsl\{\fF\}^\wedge=\tau_{\bullet},
	\end{align*}
	which are indeed true according to the preceding paragraph.
\end{eg}

\begin{proof}[Proof of Proposition~\ref{Prp:Derived_treelike_second_order_1}]
Throughout the proof,
we once again fix $\wh\al\inn\ud{\ov\dt\La}$ as in (\ref{Eqn:wh_al}).

Given $S'\!\subset\! S_{\wh\al}$,
the induced subsets
$J'$, $S'$, $J_0$ and $S_0$,
as well as the index $\wh\al_{(S')}\inn\ud{\ov\dt\La}$, are identical with those in the proof of Proposition~\ref{Prp:Derived_treelike_2}.
Mimicking the argument in the paragraph containing (\ref{Eqn:Derived_2}),
we conclude that in order to establish the statement of Proposition~\ref{Prp:Derived_treelike_second_order_1}, it suffices to justify the following analogue of (\ref{Eqn:Derived_2'}) for all $S'\!\subset\!S_{\wh\al}$ with $|S'|\eq 1$:
\begin{align}
	\label{Eqn:Derived_2nd}
	\dg{(\ud\ls_{(J')})}\big\bsl\,
	\big(\ny_2^{\ov\lt_{(J_0)}}\!(\ud\ls_{(J')})\cup
	\nu^{\ov\lt_{(J_0)}}\!(\ud\ls_{(J')})\big)^{\!\wedge}
	=
	\dg{\ud\ls}\big\bsl\,
	\big(\ny_2^{\ov\lt}(\ud\ls)\cup
	\nu^{\ov\lt}(\ud\ls)
	\sqcup(\be'_{\wh\al})^{-1}(S')\big)^{\!\wedge}.
\end{align}

Below, we prove (\ref{Eqn:Derived_2nd}) by considering all the possibilities for $S'$,
based on the expression (\ref{Eqn:S_wh_al}) of $S_{\wh\al}$.

\textbf{Case A}: $S'\eq\{e\}\inn S_\al\bsl\be_\al(\tau_\al)$.
In this case, $J'\eq J_0\eq \emptyset$ and
$(\be'_{\wh\al})^{-1}(S')\eq\emptyset$, so (\ref{Eqn:Derived_2nd}) holds trivially.

\textbf{Case B.i}: $S'\eq\{\be_\al(e)\}$ for some $e\inn\ND(\ov\lt)\!\cap\!\ND(\ud\ls).$
In this case, 
\begin{align*}
	J'=\{e\},\quad
	J_0=\begin{cases}
		\{e\}& \tn{if}~e\inn\fF_\cht^\succ,\\
		\emptyset & \tn{otherwise},
	\end{cases}
	\qquad
	\Longrightarrow\qquad
	\dg{\ud\ls_{(J')}}=
	\dg{\ud\ls}\bsl 
	e\,.
\end{align*}
\begin{itemize}[leftmargin=*]
	\item 
	If $e\in\fF_\cht^\succ\cup\big(\fF_\cht\bsl\max(\vr_{\ov\lt})\big)$,
	then
	\begin{align*}
		\dg{\ud\ls}\bsl 
		e=
		\dg{\ud\ls}\bsl 
		e^\wedge,\qquad
		\ny_2^{\ov\lt_{(J_0)}}\!(\ud\ls_{(J')})
		=
		\ny_2^{\ov\lt}(\ud\ls),\qquad
		\nu^{\ov\lt_{(J_0)}}\!(\ud\ls_{(J')})
		=
		\nu^{\ov\lt}(\ud\ls),
		\qquad
		(\be'_{\wh\al})^{-1}(S')
		=
		\{e\},
	\end{align*}
	hence both sides of (\ref{Eqn:Derived_2nd}) are equal to
	$\dg{\ud\ls}\big\bsl\,
	\big(\ny_2^{\ov\lt}(\ud\ls)\cup
	\nu^{\ov\lt}(\ud\ls)
	\sqcup\{e\}\big)^{\!\wedge}.$
	
	\item 
	If $e\inn\fF_\cht\!\cap\!\max(\vr_{\ov\lt})$,
	then $(\be'_{\wh\al})^{-1}(S')\eq\{e\}$.
	On the one hand,
	we have $\ny_2^{\ov\lt_{(J_0)}}\!(\ud\ls_{(J')})\eq \ud\bF$,
	which implies the left hand side of (\ref{Eqn:Derived_2nd}) is $\tau_\bullet$.
	On the other hand,
	notice that for every $e'\!\prec\!e$,
	we have $e'\inn\nu^{\ov\lt}(\ud\ls)$,
	hence we can always choose $e''\inn\min(\vr_{\ov\lt})\!\cap\!\{e\}^\prec$ so that
	\begin{align*}
		\nu^{\ov\lt}(\ud\ls)\sqcup\{e\}
		\;\supset\; \{e''\}^\succeq.
	\end{align*}
	(Here, notice that $\{e''\}^\succeq$ is taken in $\vr_{\ov\lt}$.)
	Therefore,
	the right hand side of (\ref{Eqn:Derived_2nd}) is also $\tau_\bullet$.
	
	\item 
	If $e\inn\fF_\cht^\prec\big\bsl\nu^{\ov\lt}(\ud\ls)$,
	then $\{e\}^\succ\!\subset\!\ov\bE$.
	Hence
	\begin{align*}
		\ny_2^{\ov\lt_{(J_0)}}\!(\ud\ls_{(J')})
		=
		\ny_2^{\ov\lt}(\ud\ls)\cup
		\big\{\fF\inn\ud\bF:\,\fF\!\cap\!\{e\}^\preceq\!\ne\!\emptyset\big\},\quad
		\nu^{\ov\lt_{(J_0)}}\!(\ud\ls_{(J')})
		\eq
		\nu^{\ov\lt}(\ud\ls),\quad
		(\be'_{\wh\al})^{-1}(S')
		\eq
		\{e\}.
	\end{align*}
	Notice that 
	\begin{align*}
		e\in\big\{\fF\inn\ud\bF:\,\fF\!\cap\!\{e\}^\preceq\!\ne\!\emptyset\big\}^\wedge
		\quad\tn{in}\ \ 
		\dg{\ud\ls}\,;\qquad
		\big\{\fF\inn\ud\bF:\,\fF\!\cap\!\{e\}^\preceq\!\ne\!\emptyset\big\}^\wedge\!= e^\wedge\quad
		\tn{in}\ \ 
		\dg{\ud\ls}\big\bsl
		\big(\nu^{\ov\lt}(\ud\ls)\big)^\wedge.
	\end{align*}
	Therefore,
	\begin{align*}
		\tn{LHS~of}~(\ref{Eqn:Derived_2nd})
		&=
		\dg{\ud\ls}\Big\bsl
		\Big(\ny_2^{\ov\lt}(\ud\ls)\cup
		\nu^{\ov\lt}(\ud\ls)\cup
		\big\{\fF\inn\ud\bF:\,\fF\!\cap\!\{e\}^\preceq\big\}\Big)^{\!\wedge}\\
		&=
		\dg{\ud\ls}\big\bsl
		\big(\ny_2^{\ov\lt}(\ud\ls)\cup
		\nu^{\ov\lt}(\ud\ls)\cup
		\{e\}\big)^{\!\wedge}
		=
		\tn{RHS~of}~(\ref{Eqn:Derived_2nd})\,.
	\end{align*}
	
	\item 
	If $e\inn\nu^{\ov\lt}(\ud\ls)$,
	then 
	\begin{align*}
		\dg{\ud\ls}\bsl 
		e=
		\dg{\ud\ls}\bsl 
		e^\wedge,\qquad
		\ny_2^{\ov\lt_{(J_0)}}\!(\ud\ls_{(J')})
		=
		\ny_2^{\ov\lt}(\ud\ls),\qquad
		\nu^{\ov\lt_{(J_0)}}\!(\ud\ls_{(J')})
		=
		\nu^{\ov\lt}(\ud\ls)\bsl\{e\},
		\qquad
		(\be'_{\wh\al})^{-1}(S')
		=
		\emptyset,
	\end{align*}
	hence both sides of (\ref{Eqn:Derived_2nd}) are equal to
	$\dg{\ud\ls}\big\bsl\,
	\big(\ny_2^{\ov\lt}(\ud\ls)\cup
	\nu^{\ov\lt}(\ud\ls)\big)^{\!\wedge}.$
\end{itemize}

\textbf{Case B.ii}: $S'\eq\{\be_\al(e)\}$ for some $e\inn\ND(\ov\lt)\!\cap\!\ny_2(\ov\lt)^\wedge.$
In this case, 
$J'\eq\emptyset,$ $J_0\eq\{e\}$, and $(\be'_{\wh\al})^{-1}(S')
\eq\emptyset$,
thus both sides of (\ref{Eqn:Derived_2nd}) are equal to
$\dg{\ud\ls}\big\bsl\,
\big(\ny_2^{\ov\lt}(\ud\ls)\cup
\nu^{\ov\lt}(\ud\ls)\big)^{\!\wedge}.$

\textbf{Case C.i}:
$S'\eq\{\fE\}\!\subset\!\ny_2(\ov\lt)$.
This case is analogous to Case B.ii:
we still have $J'\eq\emptyset$  and $(\be'_{\wh\al})^{-1}(S')
\eq\emptyset$,
so both sides of (\ref{Eqn:Derived_2nd}) are still equal to
$\dg{\ud\ls}\big\bsl\,
\big(\ny_2^{\ov\lt}(\ud\ls)\cup
\nu^{\ov\lt}(\ud\ls)\big)^{\!\wedge}.$

\textbf{Case C.ii}:
$S'\eq\{\fE\}\!\subset\!\ov\bE\!\cap\!\ND(\ud\ls)$.
In this case,  we have
\begin{align*}
	J'=\{\fE\},\qquad
	J_0=\begin{cases}
		\{\fE\} & \tn{if}~\fE\inn\fF_\cht^\succ,\\
		\emptyset & \tn{otherwise},
	\end{cases}
	\qquad
	(\be'_{\wh\al})^{-1}(S')
	=\emptyset\,.
\end{align*}
\begin{itemize}[leftmargin=*]
	\item 
	If $\fE\inn\fF_\cht^\succeq$,
	then $\ny_2^{\ov\lt}(\ud\ls_{(\fE)})\eq\ny_2^{\ov\lt}(\ud\ls)\eq \ud\bF$,
	thus  both sides of (\ref{Eqn:Derived_2nd}) are $\tau_\bullet$.
	
	\item 
	If $\fE\inn\fF_\cht^\prec$,
	then we have
	\begin{align*}
		\dg{\ud\ls_{(J')}}=
		\dg{\ud\ls}\bsl 
		\{\fE\}^\wedge,\qquad
		\ny_2^{\ov\lt_{(J_0)}}\!(\ud\ls_{(J')})
		=
		\ny_2^{\ov\lt}(\ud\ls),\qquad
		\nu^{\ov\lt_{(J_0)}}\!(\ud\ls_{(J')})
		=
		\nu^{\ov\lt}(\ud\ls).
	\end{align*}
	Observe that $\fE\inn \big(\ny_2^{\ov\lt}(\ud\ls)\big)^{\!\wedge}$ in $\dg{\ud\ls}$, so both sides of (\ref{Eqn:Derived_2nd}) are equal to
	$\dg{\ud\ls}\big\bsl\,
	\big(\ny_2^{\ov\lt}(\ud\ls)\cup
	\nu^{\ov\lt}(\ud\ls)\big)^{\!\wedge}.$
\end{itemize}

\textbf{Case D}:
$S'\eq\{\fF\}\!\subset\!\ud\bF$.
In this case, 
\begin{align*}
	J'=\{\fF\},\quad
	J_0=\ud\sfE_{\{\fF\}}
	\qquad
	(\be'_{\wh\al})^{-1}(S')
	=
	\begin{cases}
		\{\fF\} & \tn{if}~\fF\!\not\in\!\ny_2^{\ov\lt}(\ud\ls),\\
		\emptyset & \tn{if}~\fF\!\in\!\ny_2^{\ov\lt}(\ud\ls).
	\end{cases}.
\end{align*}
\begin{itemize}[leftmargin=*]
	\item 
	If $\fF\!\not\in\!\ny_2^{\ov\lt}(\ud\ls)$,
	then 
	\begin{align*}
		\dg{\ud\ls_{(J')}}=
		\dg{\ud\ls}\bsl 
		\fF^\wedge,\qquad
		\ny_2^{\ov\lt_{(J_0)}}\!(\ud\ls_{(J')})
		=
		\ny_2^{\ov\lt}(\ud\ls),\qquad
		\nu^{\ov\lt_{(J_0)}}\!(\ud\ls_{(J')})
		=
		\nu^{\ov\lt}(\ud\ls).
	\end{align*}
	so both sides of (\ref{Eqn:Derived_2nd}) are equal to $\dg{\ud\ls}\big\bsl\,
	\big(\ny_2^{\ov\lt}(\ud\ls)\cup
	\nu^{\ov\lt}(\ud\ls)\cup\{\fF\}\big)^{\!\wedge}.$
	
	\item 
	If $\fF\!\in\!\ny_2^{\ov\lt}(\ud\ls)$,
	then
	\begin{align*}
		\dg{\ud\ls_{(J')}}=
		\dg{\ud\ls}\bsl 
		\fF^\wedge,\qquad
		\ny_2^{\ov\lt_{(J_0)}}\!(\ud\ls_{(J')})
		=
		\ny_2^{\ov\lt}(\ud\ls)\bsl\{\fF\},\qquad
		\nu^{\ov\lt_{(J_0)}}\!(\ud\ls_{(J')})
		=
		\nu^{\ov\lt}(\ud\ls),
	\end{align*}
	hence both sides of (\ref{Eqn:Derived_2nd}) are equal to
	$\dg{\ud\ls}\big\bsl\,
	\big(\ny_2^{\ov\lt}(\ud\ls)\cup
	\nu^{\ov\lt}(\ud\ls)\big)^{\!\wedge}.$
	\qedhere
\end{itemize}
\end{proof}

\subsection{Grafted stratification}
\label{Subsec:Grafted}

For every rooted tree $\tau\eq(E,\preceq)$, 
let $\ex(\tau)$ be the rooted tree
\begin{gather*}
 \ex(\tau)=\big(\,E\sqcup\{e_\ex\}\,,\,\preceq\,\big),
\end{gather*}
where $\preceq$ is the relation on $\ex(\tau)$ given by the tree order on $E$ as well as the following:
every $e\inn E$ is {\it not} comparable with $e_\ex$.
It is a direct check that $\preceq$ is a tree order as per Definition~\ref{Dfn:Rooted_tree},
and $e_\ex$ is both maximal and minimal.

Intuitively,
$\ex(\tau)$ is obtained from $\tau$ by {\it grafting} a new vertex onto the root via a new  edge $e_\ex$.
An example is illustrated in Figure~\ref{Fig:grafted}.

\begin{figure}[htp]
\begin{center}
\begin{tikzpicture}
\filldraw[xshift=2cm]
 (0,0) circle (1pt)
 (.35,.5) circle (1pt)
 (.7,0) circle (1pt);
\filldraw[xshift=6cm]
 (0,0) circle (1pt)
 (.35,.5) circle (1pt)
 (.7,0) circle (1pt)
 (1.4,0) circle (1pt);
 
\draw[xshift=2cm]
 (0,0)--(.35,.5)--(.7,0);
\draw[xshift=6cm]
 (0,0)--(.35,.5)--(.7,0)
 (.35,.5)--(1.4,0)
 ;
 
\draw[xshift=2cm]
 (.05,.25) node {\tiny{$a$}}
 (.65,.25) node {\tiny{$b$}}
 (.35,.5) node[above] {\tiny{$o$}}
 ;
\draw[xshift=6cm]
 (.05,.25) node {\tiny{$a$}}
 (.45,.2) node {\tiny{$b$}}
 (1.1,.3) node {\tiny{$e_\ex$}}
 (.35,.5) node[above] {\tiny{$o$}}
 ;
\draw[yshift=-.5cm]
 (2.35,0) node {\tiny{$\tau$}}
 (6.7,0) node {\tiny{$\ex(\tau)$}};
\end{tikzpicture}
\end{center}
\caption{An example of grafting}\label{Fig:grafted}
\end{figure}  

The idea of grafting comes from the third round of the blowups in~\cite{HLN},
where the conjugate or Weierstrass loci of the underlying marked curves give another 1-dimensional condition.
The grafted edge $e_\ex$ 
reflects this extra condition.
In the remainder of this subsection,
we apply this idea to LESs and treelike structures,
which provides models for \S\ref{Subsec:Step6-9}.

Let $\fM$ be a stack with an LES as in Definition~\ref{Dfn:G-adim_fixture}, and $\De\!\subset\!\fM$ be the boundary of $\fM$ as in~(\ref{Eqn:boundary}) whose strata are labeled by the elements of $B\eq A\bsl \{0\}$.
Assume that 
there exists a closed substack 
$$K \subset \De$$
satisfying that for every $\al\inn B$
and $\cV\inn\fV_{\al}$,
there exist $\ka^\cV\inn\Ga(\sO_\cV)$ and
a {\it possibly empty} set $P(K)$ of subsets of $S_\al$ (i.e.~$P(K)$ is a subset of the power set of $S_\al$), satisfying
\begin{align}\label{Eqn:K_local}
	K\cap\cV=
	\bigcup_{E\in P(K)}\!\!\!
	\big\{\,\ze_e^\cV\eq 0\ \,\forall\,e\inn E\ ;\ \ka^\cV\eq 0\,\big\}\,;
\end{align}
moreover, 
if $K\!\cap\!\cV\!\ne\!\emptyset$,
then $\ka^\cV$ and $\ze_e^\cV,$ $e\inn S_\al$, form a subset of local parameters on $\cV$.

We set
\begin{align}\label{Eqn:grafted_boundary}
	B_\ex:=
	\{\al\inn B:\,\fM_\al\!\cap\!K\!\ne\!\emptyset\},\qquad
	A_\ex:=B_\ex\sqcup\{0\}\ \ (\subset A).
\end{align}
For each $\al\inn B_\ex$,
let
\begin{align}\label{Eqn:grafted_boundary'}
	\fM_{\ex;\al}:= \fM_{\al} \cap K,\qquad
	S_{\ex;\al} := S_{\al}\sqcup\{e_\ex\}\,.
\end{align}
The choice of $B_\ex$ implies 
\begin{align*}
	\bigsqcup_{\al\in B_\ex}\!\!\fM_{\ex;\al} = K,
\end{align*}
so we take
$$
\fM_{\ex;0}:=
\fM\,\big\bsl\big(\!
\bigsqcup_{\al\in B_\ex}\!\!\fM_{\ex;\al}\big)
=\fM\bsl K
\,,\qquad
S_{\ex;0}:=\emptyset\,.
$$

\begin{lmm}\label{Lm:grafted_strat}
With notation as above,
$$
\fM=\bigsqcup_{\al\in A_\ex}\!\fM_{\ex;\al}
$$
is an LES of $\fM$,
known as the \ts{grafted stratification of $\fM$ (with respect to the LES $\fM\eq\bigsqcup_{\al\in A}\fM_\al$ and $K$)}.
The modular charts are given by
\begin{align*}
	\fV_{\ex;\al}:=
	\big\{\,\cV\inn\fV_\al:\,
	K\!\cap\!\cV\!\neq\!\emptyset\,
	\big\},\qquad\al\in B_{\ex},
\end{align*}
while the modular parameters $\{\ze_{\ex;e}^\cV\}_{e\in S_{\ex;\al}}$,
where $\al\inn B_{\ex}$ and $\cV\inn\fV_{\ex;\al}$, satisfy
\begin{align}\label{Eqn:grafted_parameter}
	\ze_{\ex;e}^\cV=\ze_e^\cV\quad\forall\ 
	e\inn S_\al,\qquad
	\ze_{\ex;e_{\ex}}^\cV=\ka^\cV\,.
\end{align}
For every $\al,\al'\inn B_\ex$,
$\cV\inn\fV_{\ex;\al}$, and $\cV'\inn\fV_{\ex;\al'}$ satisfying
$\fM_{\ex;\al'}\!\cap\!\cV\!\cap\!\cV'\!\ne\!\emptyset$,
the injection  $\iota_{\ex;\cV,\cV'}:S_{\ex;\al'}\hookrightarrow S_{\ex;\al}$ is given by
\begin{align}\label{Eqn:grafted_injection}
	\iota_{\ex;\cV,\cV'}(e)=\iota_{\cV,\cV'}(e)\quad\forall\ e\inn S_\al,\qquad
	\iota_{\ex;\cV,\cV'}(e_\ex)=e_\ex.
\end{align}
\end{lmm}

\begin{proof}
The verification of (\ref{Eqn:loc_Euc_transition}) under the current setting is straightforward.
 
To show (\ref{Eqn:loc_Euc_equation}) for the grafted stratification of $\fM$,
observe that
\begin{align*}
	\fM_{\ex;\al'}\cap\cV\cap\cV'=
	K\cap\fM_{\al'}\cap\cV\cap\cV'
	\ne\emptyset.
\end{align*}
Applying (\ref{Eqn:loc_Euc_equation}) to the original stratification $\fM\eq\bigsqcup_{\al\in A}\fM_\al$,
and taking (\ref{Eqn:K_local}) into consideration,
we have
\begin{align*}
	K\cap\fM_{\al'}&\cap\cV\cap\cV'
	\\
	=~&
	\big\{~\ze_e^{\cV}\!\!\eq 0\ \;
	\forall\,e\in \iota_{\cV,\cV'}(S_{\al'})\,;\ \;
	\ze_e^{\cV}\!\!\!\neq\! 0\ \;
	\forall\,e\in S_\al\big\bsl\iota_{\cV,\cV'}(S_{\al'})~\big\}\,\cap \\
	&\Big(\!\bigcup_{E\in P(K)\ \tn{s.t.}\; E\subset \iota_{\cV,\cV'}(S_{\al'})}
	\hspace{-.52in} 
	\big\{\,\ze_e^\cV\eq 0\ \,\forall\,e\inn E\ ;\ \ka^\cV\eq 0\,\big\}\Big) \cap\,\cV'\\
	=~&
	\big\{~\ze_e^{\cV}\!\!\eq 0\ \;
	\forall\,e\in \iota_{\cV,\cV'}(S_{\al'})\,;\ \;
	\ka^\cV\eq 0\,;\ \;
	\ze_e^{\cV}\!\!\!\neq\! 0\ \;
	\forall\,e\in S_{\ex;\al}\big\bsl\iota_{\ex;\cV,\cV'}(S_{\ex;\al'})~\big\}\,\cap\,\cV'\,.
\end{align*}
This, along with (\ref{Eqn:grafted_parameter}) and (\ref{Eqn:grafted_injection}),
establishes (\ref{Eqn:loc_Euc_equation}) for the grafted stratification of $\fM$.
\end{proof}

Next, let $\La\eq(\tau_\al,\be_\al)_{\al\in A}$ be a treelike structure on $\fM$, satisfying
\begin{equation}\label{Eqn:grafted_treelike}
	\tau_{\al}\ne\tau_\bullet\qquad\tn{if~and~only~if}\qquad
	\al\in B_{\ex}\quad\big(\,\tn{i.e.}\ \al\inn B\ \ \tn{and}\ \ \fM_\al\!\cap\!K\!\ne\!\emptyset\,\big)\,.
\end{equation}

\begin{prp}\label{Prp:grafted_treelike}
Under the assumption (\ref{Eqn:K_local}), $\La$ induces a treelike structure 
\begin{align*}
	\La_\ex=\big(\,\tau_{\ex;\al}\,,\,\be_{\ex;\al}\,\big)_{\al\in A_\ex}
\end{align*}
on $\fM\eq\bigsqcup_{\al\in A_{\ex}}\!\fM_{\ex;\al}$ such that 
$\tau_{\ex;0}\eq \tau_\bullet$, 
$\be_{\ex;0}:\emptyset\!\lra\!\emptyset$, and
for every $\al\inn B_\ex$, 
\begin{alignat*}{2}
	&\tau_{\ex;\al}=\ex(\tau_\al)\qquad& 
	&\be_{\ex;\al}:\tau_{\ex;\al}\hookrightarrow S_{\ex;\al}\ \ \tn{s.t.}\ \ 
	\be_{\ex;\al}|_{\tau_\al}\eq\be_\al,\ \ 
	\be_{\ex;\al}(e_\ex)\eq e_\ex\,.
\end{alignat*}
Such defined $\La_\ex$ is called the \ts{grafted treelike structure on $\fM$ (with respect to $\La$ and $K$).}
\end{prp}

\begin{proof}
Fix $\al\inn B_{\ex}$ and $e\inn S_{\ex;\al}\eq S_\al\!\sqcup\!\{e_\ex\}$.
They together determine $\al_{\ex;(e)}\inn A_{\ex}$ as per Proposition~\ref{Prp:M_strata_local}; i.e.
\begin{align*}
	\fM_{\ex;\al_{\ex;(e)}}\!\cap\cV\,\supset \,
	\big\{\;\ze_{\ex;e'}^{\cV}\eq 0\ \;
	\forall\,e'\inn S_{\ex;\al}\bsl\{e\}\,;\ \;
	\ze_{\ex;e}^{\cV}\!\neq\! 0\;\big\}\qquad\forall\ 
	\cV\inn\fV_{\al}\ \tn{with}\ K\!\cap\!\cV\!\ne\!\emptyset\,.
\end{align*}
In addition, the choice of $B_\ex$ implies $\al\inn B$, so when $e\inn S_\al$,
we also have $\al_{(e)}\inn A$ as per Proposition~\ref{Prp:M_strata_local}.

Below, we divide the proof into three cases.

\textbf{Case 1}: $e\inn S_\al$ and $\al_{(e)}\inn B_\ex$. 
Then,
$\fM_{\al_{(e)}}\!\cap\!K\!\ne\!\emptyset$,
so there exist $\cV\inn\fV_{\ex;\al}$ such that
\begin{align*}
	\big\{\;\ze_{e'}^{\cV}\eq 0\ \;
	\forall\,e'\in S_\al\bsl\{e\}\,;\ \;
	\ka^\cV\eq 0\,;\ \;
	\ze_e^{\cV}\!\neq\! 0\;\big\}\,\subset \,
	\fM_{\al_{(e)}}\!\cap K\cap\cV\,=\,
	\fM_{\ex;\al_{(e)}}\!\cap \cV.
\end{align*}
Hence $\al_{\ex;(e)}\eq\al_{(e)}$.

Since $\La$ is a treelike structure,
there exists an isomorphism
\begin{align*}
	\phi_{\al;e}:\,\tau_{\al_{(e)}}\lra 
	\tau_\al\big\bsl\,\big(\be_\al^{-1}(e)\big)^\wedge.
\end{align*}
The assumptions $\al_{(e)}\inn B_\ex$ and (\ref{Eqn:grafted_treelike}) imply $\tau_{\al_{(e)}}\!\ne\!\tau_\bullet$,
so $\big(\be_\al^{-1}(e)\big)^\wedge\!\ne\!\tau_\al$.
We thus construct the isomorphism
\begin{align*}
	&
	\phi^\ex_{\al;e}:\,
	\tau_{\ex;\al_{(e)}}=\ex\big(\tau_{\al_{(e)}}\big)\,\lra\,
	\ex\Big(\tau_\al\big\bsl\,\big(\be_\al^{-1}(e)\big)^\wedge\Big)
	=\ex(\tau_\al)\big\bsl\,\big(\be_\al^{-1}(e)\big)^\wedge
	=
	\tau_{\ex;\al}\big\bsl\,\big(\be_\al^{-1}(e)\big)^\wedge\\
	&\tn{satisfying}\qquad\phi^\ex_{\al;e}|_{\tau_{\al_{(e)}}}\!=\phi_{\al;e},\qquad
	\phi^\ex_{\al;e}(e_\ex)=e_\ex\,.
\end{align*}
It is then a direct check that $\phi^\ex_{\al;e}$ satisfies (\ref{Eqn:Treelike}).

\textbf{Case 2}: $e\inn S_\al$ and $\al_{(e)}\!\not\in\! B_\ex$.
In this case, 
$\fM_{\al_{(e)}}\!\cap\!K\!=\!\emptyset$,
so for every $\cV\inn\fV_{\ex;\al}$,
the locus
\begin{align*}
	\big\{\;\ze_{e'}^{\cV}\eq 0\ \;
	\forall\,e'\in S_\al\bsl\{e\}\,;\ \;
	\ka^\cV\eq 0\,;\ \;
	\ze_e^{\cV}\!\neq\! 0\;\big\}\,\subset \,
	\fM_{\al_{(e)}}\!\cap\cV
\end{align*}
is disjoint from $K$,
hence is in $\fM_{\ex;0}$ ($=\!\fM\bsl K$).
Therefore, $\al_{\ex;(e)}\eq 0$.

The assumptions $\al_{(e)}\!\not\in\! B_\ex$ and (\ref{Eqn:grafted_treelike}) imply $\tau_{\al_{(e)}}\!=\!\tau_\bullet$,
so $\big(\be_\al^{-1}(e)\big)^\wedge\!=\!\tau_\al$, 
and consequently the adjoint set of $\be_\al^{-1}(e)$ in $\tau_{\ex;\al}\eq\ex(\tau_\al)$ is $\tau_{\ex;\al}$ itself.
We thus take
\begin{align*}
	\phi^\ex_{\al;e}:\,
	\tau_{\ex;\al_{(e)}}=\tau_{\ex;0}=\tau_\bullet\,\lra\,
	\tau_\bullet
	=
	\tau_{\ex;\al}\big\bsl\,\big(\be_\al^{-1}(e)\big)^\wedge
\end{align*}
to be the trivial isomorphism,
which obviously satisfies (\ref{Eqn:Treelike}).

\textbf{Case 3}: $e\eq e_\ex$.
In this case, notice that for every $\cV\inn\fV_{\ex;\al}$,
$\ze^\cV_{\ex;e}\eq\ka^\cV$,
hence 
$\fM_{\ex;\al_{\ex;(e)}}\!\cap\!\cV$ is disjoint from $K$,
so $\al_{\ex;(e)}\eq 0$ and $\tau_{\ex;\al_{(e)}}\eq\tau_\bullet$.
In addition,
observe that $e_\ex$ is both maximal and minimal in $\ex(\tau_{\al})$,
hence $\tau_{\ex;\al}\big\bsl\,\big(\be_\al^{-1}(e)\big)^\wedge\eq\ex(\tau_{\al})\bsl (e_\ex)^\wedge\eq\tau_\bullet$.
As in Case 2,
the trivial isomorphism $\phi^\ex_{\al;e}$ satisfies (\ref{Eqn:Treelike}).
\end{proof}

\section{Application to the moduli of genus 2 stable maps}
\label{Sec:genus_2_twisted_fields}
 In~\S\ref{Subsec:Step1}-\S\ref{Subsec:Step6-9},
we apply the STF theory to $\fD_2$ in nine times and obtain $\ti{\fD}_2^{\rm tf}\!:=\!\fM^{\fk{10}}$.
The construction is guided by the local equations of $\ov M_2(\P^n,d)$ and $\mathbf{R}\pi_*\ff^*\sO_{\P^n}(k)$ reviewed in \S\ref{Subsec:M2Pnd_loc_eqn}.
In each step,
the explicit expression of the pullback of the structural homomorphism~(\ref{Eqn:str_hom_0}) is provided,
whose entries suggest the treelike structure of the succeeding step.

The outline of the construction  is illustrated in Figure~\ref{Fig:blowup}.
Here,
each arrow represents a forgetful morphism determined by the corresponding treelike structure.
Only the treelike structures whose arrows are  dashed need to be explicitly constructed in~\S\ref{Sec:genus_2_twisted_fields}.
The solid arrows follow directly from the statements of \S\ref{Sec:Induced_and_derived_tf_stacks}, indicated beside each arrow.
For example,
once the treelike structure $\La^{\fk 1}$ is established on $\fM^{\fk 1}$,
which gives rise to the dashed arrow $\fM^{\fk 2}\!\lra\!\fM^{\fk 1}$ at the bottom of Figure~\ref{Fig:blowup},
its RL derived treelike structure $\ov{\tn d}\La^{\fk 1}$ on $\fM^{\fk 2}$ is fully described in Proposition~\ref{Prp:Derived_treelike_1},
illustrated by the red arrow toward $\fM^{\fk 2}$.

\begin{figure}[htp]
	\begin{center}	
		\begin{tikzpicture}
			[
			base/.style={
				rectangle,minimum size=8mm,
				ultra thick,draw=black},
			rle/.style={
				rectangle,minimum size=6mm,
				very thick,draw=black,
				top color=white,bottom color=black!10},
			rln/.style={
				rectangle,minimum size=6mm,
				very thick,draw=black,dashed,
				top color=white,bottom color=black!10},
			lre/.style={
				rectangle,minimum size=6mm,
				very thick,draw=black,
				top color=black!10,bottom color=white},
			lrn/.style={
				rectangle,minimum size=6mm,
				very thick,draw=black,dashed,
				top color=black!10,bottom color=white},
			Prln/.style={
				rectangle,minimum size=6mm,rounded corners=3mm,
				very thick,dashed,draw=black!20,
				top color=white,bottom color=black!10},
			Prle/.style={
				rectangle,minimum size=6mm,rounded corners=3mm,
				very thick,draw=black!20,
				top color=white,bottom color=black!10},
			Plrn/.style={
				rectangle,minimum size=6mm,rounded corners=3mm,
				very thick,dashed,draw=black!20,
				top color=black!10,bottom color=white},
			Plre/.style={
				rectangle,minimum size=6mm,rounded corners=3mm,
				very thick,draw=black!20,
				top color=black!10,bottom color=white}
			]
			\node (M1) [node distance=4mm and 2mm, base,] {\scriptsize{$~\fM^{\fk 1}=\fD_2~$}};
			\node (r1p1) [node distance=4mm and 2mm, above=of M1] {};
			\node (M2) [node distance=4mm and 2mm, rln,above=of r1p1] {\scriptsize{$\fM^{\fk 2}=(\fM^{\fk 1})^{\tf}_{\La^{\fk 1}}$}};
			\node (r1p2) [node distance=8mm and 4mm,above=of M2] {};
			\node (M3) [node distance=8mm and 4mm,rln,above=of r1p2] {\scriptsize{$\fM^{\fk 3}=(\fM^{\fk 2})^{\tf}_{\La^{\fk 2}}$}};
			\node (r1p34) [node distance=8mm and 2mm,above=of M3] {};
			\node (M4) [node distance=8mm and 2mm,rln,above=of r1p34] {\scriptsize{$\fM^{\fk 4}=(\fM^{\fk 3})^{\tf}_{\La^{\fk 3}}$}};
			\node (r1p5) [node distance=4mm and 2mm,above=of M4] {};
			\node (M5) [node distance=4mm and 2mm,rle,above=of r1p5] {\scriptsize{$\fM^{\fk 5}=(\fM^{\fk 4})^{\tf}_{{\PT}(\Up^{\fk 4})}$}};
			\node (r2) [node distance=4mm and 2mm,above=of M5] {};
			\node (M6) [node distance=4mm and 2mm,lre,above=of r2] {\scriptsize{$\fM^{\fk 6}=(\fM^{\fk 5})^{\rtf}_{{\PT}(\ov{\tn d}\La^{\fk 1})}$}};
			\node (r3p1) [node distance=4mm and 2mm,above=of M6] {};
			\node (M7) [node distance=4mm and 2mm,rln,above=of r3p1] {\scriptsize{$\fM^{\fk 7}=(\fM^{\fk 6})^{\tf}_{({\PT}({\tn d}_2\La^{\fk 1}))_{\ex}}$}};
			\node (r3p2) [node distance=4mm and 2mm,above=of M7] {};
			\node (M8) [node distance=4mm and 2mm,rle,above=of r3p2] {\scriptsize{$\fM^{\fk 8}=(\fM^{\fk 7})^{\tf}_{({\PT}(\ud{\partial}\ov{\tn d}\La^{\fk 1}))_{\ex}}$}};
			\node (r3p3) [node distance=4mm and 2mm,above=of M8] {};
			\node (M9) [node distance=4mm and 2mm,rle,above=of r3p3] {\scriptsize{$\fM^{\fk 9}=(\fM^{\fk 8})^{\tf}_{({\PT}(\Up^{\fk 8}))_{\ex}}$}};
			\node (r3p4) [node distance=4mm and 2mm,above=of M9] {};
			\node (M10) [node distance=4mm and 2mm,rle,above=of r3p4] {\scriptsize{$\fM^{\fk{10}}=(\fM^{\fk 9})^{\tf}_{({\PT}({\tn d}\Up^{\fk 4}))_{\ex}}$}};
			
			\node (Pr2) [node distance=2mm and 6mm,right=of M2] {};
			\node (PM6) [node distance=2mm and 6mm,Plre,right=of Pr2] {\scriptsize{$\fN^{\fk 6}_{\fk 2}=(\fM^{\fk 2})^{\rtf}_{\ov{\tn d}\La^{\fk 1}}$}};
			\node (Pr3p1) [node distance=-2mm and 4mm,below right=of PM6] {};
			\node (PM7) [node distance=-2mm and 4mm,Prle,below right=of Pr3p1] {\scriptsize{$(\fN^{\fk 6}_{\fk 2})^{\tf}_{{\tn d}_2\La^{\fk 1}}$}};
			\node (Pr3p2) [node distance=-2mm and 4mm,above right=of PM6] {};
			\node (PM8) [node distance=-2mm and 4mm,Prle,above right=of Pr3p2] {\scriptsize{$(\fN^{\fk 6}_{\fk 2})^{\tf}_{\ud{\partial}\ov{\tn d}\La^{\fk 1}}$}};
			\node (Pr1p5) [node distance=2mm and 6mm,left=of M2] {};
			\node (PM5) [node distance=2mm and 6mm,Prln,left=of Pr1p5] {\scriptsize{$\fN^{\fk 5}_{\fk 2}=(\fM^{\fk 2})^{\tf}_{\Up^{\fk 4}}$}};
			\node (Pr3p3) [node distance=-2mm and 4mm,above left=of PM5] {};
			\node (PM9) [node distance=-2mm and 4mm,Prln,above left=of Pr3p3] {\scriptsize{$(\fN^{\fk 5}_{\fk 2})^{\tf}_{\Up^{\fk 8}}$}};
			\node (Pr3p4) [node distance=-2mm and 4mm,below left=of PM5] {};
			\node (PM10) [node distance=-2mm and 4mm,Prle,below left=of Pr3p4] {\scriptsize{$(\fN^{\fk 5}_{\fk 2})^{\tf}_{\ov{\tn d}\Up^{\fk 4}}$}};
			
			\node (Pr2') [node distance=2mm and 4.5mm,right=of M3] {};
			\node (PM6') [node distance=2mm and 4.5mm,Plre,right=of Pr2'] {\scriptsize{$\fN^{\fk 6}_{\fk 3}=(\fM^{\fk 3})^{\rtf}_{{\PT}(\ov{\tn d}\La^{\fk 1})}$}};
			\node (Pr3p1') [node distance=-2mm and 4mm,below right=of PM6'] {};
			\node (PM7') [node distance=-2mm and 4mm,Prle,below right=of Pr3p1'] {\scriptsize{$(\fN^{\fk 6}_{\fk 3})^{\tf}_{{\PT}({\tn d}_2\La^{\fk 1})}$}};
			\node (Pr3p2') [node distance=-2mm and 4mm,above right=of PM6'] {};
			\node (PM8') [node distance=-2mm and 4mm,Prle,above right=of Pr3p2'] {\scriptsize{$(\fN^{\fk 6}_{\fk 3})^{\tf}_{{\PT}(\ud{\partial}\ov{\tn d}\La^{\fk 1})}$}};
			\node (Pr1p5') [node distance=2mm and 4.5mm,left=of M3] {};
			\node (PM5') [node distance=2mm and 4.5mm,Prle,left=of Pr1p5'] {\scriptsize{$\fN^{\fk 5}_{\fk 3}=(\fM^{\fk 3})^{\tf}_{{\PT}(\Up^{\fk 4})}$}};
			\node (Pr3p3') [node distance=-2mm and 4mm,above left=of PM5'] {};
			\node (PM9') [node distance=-2mm and 4mm,Prle,above left=of Pr3p3'] {\scriptsize{$(\fN^{\fk 5}_{\fk 3})^{\tf}_{{\PT}(\Up^{\fk 8})}$}};
			\node (Pr3p4') [node distance=-2mm and 4mm,below left=of PM5'] {};
			\node (PM10') [node distance=-2mm and 4mm,Prle,below left=of Pr3p4'] {\scriptsize{$(\fN^{\fk 5}_{\fk 3})^{\tf}_{{\PT}(\ov{\tn d}\Up^{\fk 4})}$}};
			
			\node (Pr2'') [node distance=2mm and 4.5mm,right=of M4] {};
			\node (PM6'') [node distance=2mm and 4.5mm,Plre,right=of Pr2''] {\scriptsize{$\fN^{\fk 6}_{\fk 4}=(\fM^{\fk 4})^{\rtf}_{{\PT}(\ov{\tn d}\La^{\fk 1})}$}};
			\node (Pr3p1'') [node distance=-2mm and 4mm,below right=of PM6''] {};
			\node (PM7'') [node distance=-2mm and 4mm,Prle,below right=of Pr3p1''] {\scriptsize{$(\fN^{\fk 6}_{\fk 4})^{\tf}_{{\PT}({\tn d}_2\La^{\fk 1})}$}};
			\node (Pr3p2'') [node distance=-2mm and 4mm,above right=of PM6''] {};
			\node (PM8'') [node distance=-2mm and 4mm,Prle,above right=of Pr3p2''] {\scriptsize{$(\fN^{\fk 6}_{\fk 4})^{\tf}_{{\PT}(\ud{\partial}\ov{\tn d}\La^{\fk 1})}$}};
			
			\node (Pr3p3'') [node distance=-2mm and 6mm,above left=of M5] {};
			\node (PM9'') [node distance=-2mm and 6mm,Prle,above left=of Pr3p3''] {\scriptsize{$(\fM^{\fk 5})^{\tf}_{{\PT}(\Up^{\fk 8})}$}};
			\node (Pr3p4'') [node distance=-2mm and 6mm,below left=of M5] {};
			\node (PM10'') [node distance=-2mm and 6mm,Prle,below left=of Pr3p4''] {\scriptsize{$(\fM^{\fk 5})^{\tf}_{{\PT}(\ov{\tn d}\Up^{\fk 4})}$}};
			
			\node (Pr3p2''') [node distance=-2mm and 15mm,left=of M6] {};
			\node (PM8''') [node distance=-2mm and 15mm,Prln,left=of Pr3p2'''] {\scriptsize{$(\fM^{\fk 6})^{\tf}_{({\PT}(\ud{\partial}\ov{\tn d}\La^{\fk 1}))_{\ex}}$}};
			\node (Pr3p3''') [node distance=-2mm and 15mm,above right=of M6] {};
			\node (PM9''') [node distance=-2mm and 15mm,Prln,above right=of Pr3p3'''] {\scriptsize{$(\fM^{\fk 6})^{\tf}_{({\PT}(\Up^{\fk 8}))_{\ex}}$}};
			\node (Pr3p4''') [node distance=-2mm and 15mm,below right=of M6] {};
			\node (PM10''') [node distance=-2mm and 15mm,Prln,below right=of Pr3p4'''] {\scriptsize{$(\fM^{\fk 6})^{\tf}_{({\PT}(\ov{\tn d}\Up^{\fk 4}))_{\ex}}$}};
			
			\node (Pr3p3'''') [node distance=-2mm and 6mm,right=of M7] {};
			\node (PM9'''') [node distance=-2mm and 6mm,Prle,right=of Pr3p3''''] {\scriptsize{$(\fM^{\fk 7})^{\tf}_{({\PT}(\Up^{\fk 8}))_{\ex}}$}};
			\node (Pr3p4'''') [node distance=-2mm and 6mm,left=of M7] {};
			\node (PM10'''') [node distance=-2mm and 6mm,Prle,left=of Pr3p4''''] {\scriptsize{$(\fM^{\fk 7})^{\tf}_{({\PT}(\ov{\tn d}\Up^{\fk 4}))_{\ex}}$}};
			
			\node (Pr3p4''''') [node distance=-2mm and 6mm,left=of M8] {};
			\node (PM10''''') [node distance=-2mm and 6mm,Prle,left=of Pr3p4'''''] {\scriptsize{$(\fM^{\fk 8})^{\tf}_{({\PT}(\ov{\tn d}\Up^{\fk 4}))_{\ex}}$}};
			
			\node [node distance=-1mm and -1.5mm, left=of r1p1, black!50] {\tiny{$(\mathsf{r}_1\mathsf{p}_1)$}};
			\node [node distance=-1mm and -1.5mm, left=of r1p2, black!50] {\tiny{$(\mathsf{r}_1\mathsf{p}_2)$}};
			\node [node distance=-1mm and -1.5mm, left=of r1p34, black!50] {\tiny{$(\mathsf{r}_1\mathsf{p}_4)$ then $(\mathsf{r}_1\mathsf{p}_3)$}};
			\node [node distance=-1mm and -1.5mm, left=of r1p5, black!50] {\tiny{$(\mathsf{r}_1\mathsf{p}_5)$}};
			\node [node distance=-1mm and -1.5mm, left=of r2, black!50] {\tiny{$(\mathsf{r}_2)$}};
			\node [node distance=-1mm and -1.5mm, left=of r3p1, black!50] {\tiny{$(\mathsf{r}_3\mathsf{p}_1)$}};
			\node [node distance=-1mm and -1.5mm, left=of r3p2, black!50] {\tiny{$(\mathsf{r}_3\mathsf{p}_2)$}};
			\node [node distance=-1mm and -1.5mm, left=of r3p3, black!50] {\tiny{$(\mathsf{r}_3\mathsf{p}_3)$}};
			\node [node distance=-1mm and -1.5mm, left=of r3p4, black!50] {\tiny{$(\mathsf{r}_3\mathsf{p}_4)$}};
			
			\node [node distance=-1mm and -1.5mm, above=of Pr2, black!50] {\tiny{Prp.~\ref{Prp:Derived_treelike_1}}};
			\node [node distance=-1mm and -6mm, below left=of Pr3p1, black!50] {\tiny{Prp.~\ref{Prp:Derived_treelike_2}}};
			\node [node distance=-1mm and -6mm, above left=of Pr3p2, black!50] {\tiny{Prp.~\ref{Prp:Derived_treelike_second_order_1}}};
			\node [node distance=-1mm and -6mm, below right=of Pr3p4, black!50] {\tiny{Prp.~\ref{Prp:Derived_treelike_1}}};

			\node [node distance=-1mm and -1.5mm, above=of Pr2', black!50] {\tiny{Prp.~\ref{Prp:two_treelike_strs}}};
			\node [node distance=-1mm and -6mm, below left=of Pr3p1', black!50] {\tiny{Prp.~\ref{Prp:2nd_PT}}};
			\node [node distance=-1mm and -6mm, above left=of Pr3p2', black!50] {\tiny{Prp.~\ref{Prp:2nd_PT}}};
			\node [node distance=-1mm and -2.5mm, above=of Pr1p5', black!50] {\tiny{Prp.~\ref{Prp:two_treelike_strs}}};
			\node [node distance=-1mm and -6mm, below right=of Pr3p4', black!50] {\tiny{Prp.~\ref{Prp:2nd_PT}}};
			\node [node distance=-1mm and -6mm, above right=of Pr3p3', black!50] {\tiny{Prp.~\ref{Prp:2nd_PT}}};
			
			\node [node distance=-1mm and -1.5mm, above=of Pr2'', black!50] {\tiny{Prp.~\ref{Prp:two_treelike_strs}}};
			\node [node distance=-1mm and -6mm, below left=of Pr3p1'', black!50] {\tiny{Prp.~\ref{Prp:2nd_PT}}};
			\node [node distance=-1mm and -6mm, above left=of Pr3p2'', black!50] {\tiny{Prp.~\ref{Prp:2nd_PT}}};
			
			\node [node distance=-1mm and -1.5mm, right=of r1p5, black!50] {\tiny{Prp.~\ref{Prp:two_treelike_strs}}};
			\node [node distance=-1mm and -8mm, below right=of Pr3p4'', black!50] {\tiny{Prp.~\ref{Prp:2nd_PT}}};
			\node [node distance=-1mm and -8mm, above right=of Pr3p3'', black!50] {\tiny{Prp.~\ref{Prp:2nd_PT}}};
			
			\node [node distance=-1mm and -1.5mm, right=of r2, black!50] {\tiny{Prp.~\ref{Prp:two_treelike_strs}}};
			
			\node [node distance=-1mm and -8mm, above=of Pr3p2''', black!50] {\tiny{Prp.~\ref{Prp:2nd_PT}~\&~\ref{Prp:grafted_treelike}}};
			\node [node distance=-1mm and -8mm, below left=of Pr3p4''', black!50] {\tiny{Prp.~\ref{Prp:two_treelike_strs}~\&~\ref{Prp:grafted_treelike}}};
			\node [node distance=-1mm and -8mm, above left=of Pr3p3''', black!50] {\tiny{Prp.~\ref{Prp:two_treelike_strs}~\&~\ref{Prp:grafted_treelike}}};
			
			\node [node distance=-1mm and -1.5mm, right=of r3p1, black!50] {\tiny{Prp.~\ref{Prp:2nd_PT}~\&~\ref{Prp:grafted_treelike}}};
			
			\node [node distance=-1mm and -8mm, above=of Pr3p4'''', black!50] {\tiny{Prp.~\ref{Prp:two_treelike_strs}}};
			\node [node distance=-1mm and -8mm, above=of Pr3p3'''', black!50] {\tiny{Prp.~\ref{Prp:two_treelike_strs}}};
			
			\node [node distance=-1mm and -1.5mm, right=of r3p2, black!50] {\tiny{Prp.~\ref{Prp:two_treelike_strs}}};
			
			\node [node distance=-1mm and -8mm, above=of Pr3p4''''', black!50] {\tiny{Prp.~\ref{Prp:two_treelike_strs}}};
			
			\node [node distance=-1mm and -1.5mm, right=of r3p3, black!50] {\tiny{Prp.~\ref{Prp:two_treelike_strs}}};
			
			\node [node distance=-1mm and -1.5mm, right=of r3p4, black!50] {\tiny{Prp.~\ref{Prp:two_treelike_strs}}};
			
			\path (M2) edge[->,>=stealth', very thick, dashed, style={shorten >=1pt}] (M1);
			
			\path (M3) edge[->,>=stealth', very  thick, dashed, style={shorten >=1pt}] (M2);
			
			\path (M4) edge[->,>=stealth', very  thick, dashed, style={shorten >=1pt}] (M3);
			
			\path (M5) edge[->,>=stealth', very  thick, orange, style={shorten >=1pt}] (M4);
			
			\path (M6) edge[->,>=stealth', very  thick, red, style={shorten >=1pt}] (M5);
			
			\path (PM6) edge[->,>=stealth', thick, red, style={shorten >=1pt}] (M2);
			
			\path (PM7) edge[->,>=stealth', thick, blue, style={shorten >=1pt}] (PM6);
			
			\path (PM8) edge[->,>=stealth', thick, green, style={shorten >=1pt}] (PM6);
			
			\path (PM5) edge[->,>=stealth', thick, dashed, orange, style={shorten >=1pt}] (M2);
			
			\path (PM9) edge[->,>=stealth', thick, dashed, cyan, style={shorten >=1pt}] (PM5);
			
			\path (PM10) edge[->,>=stealth', thick, olive, style={shorten >=1pt}] (PM5);
			
			\path (PM6') edge[->,>=stealth', thick, red, style={shorten >=1pt}] (M3);
			
			\path (PM7') edge[->,>=stealth', thick, blue, style={shorten >=1pt}] (PM6');
			
			\path (PM8') edge[->,>=stealth', thick, green, style={shorten >=1pt}] (PM6');
			
			\path (PM5') edge[->,>=stealth', thick, orange, style={shorten >=1pt}] (M3);
			
			\path (PM9') edge[->,>=stealth', thick, cyan, style={shorten >=1pt}] (PM5');
			
			\path (PM10') edge[->,>=stealth', thick, olive, style={shorten >=1pt}] (PM5');
			
			\path (PM6'') edge[->,>=stealth', thick, red, style={shorten >=1pt}] (M4);
			
			\path (PM7'') edge[->,>=stealth', thick, blue, style={shorten >=1pt}] (PM6'');
			
			\path (PM8'') edge[->,>=stealth', thick, green, style={shorten >=1pt}] (PM6'');
			
			\path (PM9'') edge[->,>=stealth', thick, cyan, style={shorten >=1pt}] (M5);
			
			\path (PM10'') edge[->,>=stealth', thick, olive, style={shorten >=1pt}] (M5);
			
			\path (M7) edge[->,>=stealth', very thick, dashed, blue, style={shorten >=1pt}] (M6);
			
			\path (PM8''') edge[->,>=stealth', thick, dashed, green, style={shorten >=1pt}] (M6);
			
			\path (PM9''') edge[->,>=stealth', thick, dashed, cyan, style={shorten >=1pt}] (M6);
			
			\path (PM10''') edge[->,>=stealth', thick, dashed, olive, style={shorten >=1pt}] (M6);
			
			\path (M8) edge[->,>=stealth', very thick, green, style={shorten >=1pt}] (M7);
			
			\path (PM9'''') edge[->,>=stealth', thick, cyan, style={shorten >=1pt}] (M7);
			
			\path (PM10'''') edge[->,>=stealth', thick, olive, style={shorten >=1pt}] (M7);
			
			\path (M9) edge[->,>=stealth', very thick, cyan, style={shorten >=1pt}] (M8);
			
			\path (PM10''''') edge[->,>=stealth', thick, olive, style={shorten >=1pt}] (M8);
			
			\path (M10) edge[->,>=stealth', very thick, olive, style={shorten >=1pt}] (M9);
		\end{tikzpicture}		
	\end{center}
	\caption{Outline of the construction of $\ti\fD_2^{\rm tf}$}\label{Fig:blowup}	
\end{figure} 

On the STF of an intermediate step,
there may be several treelike structures.
One can add the twisted fields with respect to each of them and obtain a new STF, 
on which there may exist even more treelike structures.
For example, in Figure~\ref{Fig:blowup},
consider $\fM^{\fk 2}$ 
and the treelike structures $\La^{\fk 2}$, $\Up^{\fk 4}$ and $\ov{\tn d}\La^{\fk 1}$ on it.
One can add the RL twisted fields with respect to each of the three treelike structures
and obtain $\fM^{\fk 3}$, $\fN_{\fk 2}^{\fk 5}$ and $\fN_{\fk 2}^{\fk 6}$, respectively.
On each of them,
locally, the kernel of the pullback of the structural homomorphism (\ref{Eqn:str_hom_0}) becomes less singular.
One can further consider the treelike structures $\La^{\fk 3}$ on $\fM^{\fk 3}$,
$\Up^{\fk 8}$ and $\ov{\tn d}\Up^{\fk 4}$ on $\fN_{\fk 2}^{\fk 5}$,
and ${\tn d}_2\La^{\fk 1}$ and $\ud\partial\ov{\tn d}\La^{\fk 1}$ on $\fN_{\fk 2}^{\fk 6}$,
and add RL twisted fields respectively so that the kernel of the pullback of (\ref{Eqn:str_hom_0}) becomes even less singular.
On the set of these eight treelike structures,
there is a natural partial order whose minimal elements
with $\La^{\fk 2}$, $\Up^{\fk 4}$ and $\ov{\tn d}\La^{\fk 1}$.
In order to achieve $\ti\fD_2^{\rm tf}$,
we need to extend this partial order into a linear order (such an extension is often not unique)
and add the proper transforms (including the secondary proper transforms) of these treelike structures successively.
In Figure~\ref{Fig:blowup},
the morphisms given by the proper transforms 
of the same treelike structure are shown in the same color.
For each proper transform in Figure~\ref{Fig:blowup}, 
the treelike structure that it is taken with respect to,
as well as whether it is RL or LR,
are implicitly encoded in the arrows and not shown explicitly.

The vertical sequence of arrows in the middle of Figure~\ref{Fig:blowup} is the nine-step construction of $\ti\fD_2^{\rm tf}$.
Beside each vertical arrow, $(\mathsf{r}_i\mathsf{p}_j)$ refers to the round and phase of the blowups introduced in~\cite{HLN} that the arrow corresponds to.
As explained in \S\ref{Sec:Intro},
the order of the blowups performed in~\cite{HLN} is not identical to the current construction,
which is a typical phenomenon of a resolution problem.

The final stack $\ti{\fD}_2^{\rm tf}$ leads to the modular resolution~$\ti M_2^{\rm tf}(\P^n,d)$ of $\ov M_2(\P^n,d)$ as in Theorem~\ref{Thm:Main}.
The proof of Theorem~\ref{Thm:Main} is provided in~\S\ref{Subsec:Proof_Main}.

\subsection{Step~$\fk 1$}\label{Subsec:Step1}
We take the initial base moduli to be $\fD_2$,
renamed as $\fM^{\fk 1}$,
which has a natural LES given by the dual graphs of the curves and the distribution of the divisorial marking.

More precisely,
for every $x\eq(C,D)\inn\fD_2$,
let $\ga_x$ be its \ts{dual graph} in the usual sense that encodes the positions of the nodes and irreducible components of $C$, the geometric genera of the irreducible components,
as well as the number of the points of $D$ on each irreducible component;
c.f.~\cite[\S23.4]{MirSym}.
For another $x'\eq (C',D')\inn\fD_2$,
$\ga_{x'}$ and $\ga_{x}$ are said to be \ts{equivalent},
written as $\ga_{x'}\!\sim\! \ga_{x}$, 
if there exists an isomorphism of dual graphs between them (i.e.~there is an isomorphism of the underlying graphs that preserves the geometric genera and the distribution of the divisorial marked points).
If a dual graph is edgeless, then it contains exactly one vertex. 

Let
\begin{alignat*}{2}
	&
	B^{\fk 1}:=
	\big\{\,\ga_{x}\,:\ 
	x\in \fD_2,~\ga_x~\tn{is~not~edgeless}\,\big\}\big/\sim\,,\qquad&&
	A^{\fk 1}:=
	\{0\}\sqcup B^{\fk 1}\,,\\
	&
	\fM^{\fk 1}_{\al^{\fk 1}}:=
	\big\{\,
	x\in\fD_2\,:\,
	\ga_x\in\al
	\,\big\}\quad
	\forall~\al\in B^{\fk 1},\qquad&&
	\fM_0=\fM^{\fk 1}\Big\bsl
	\Big(\bigsqcup_{\al\in B^{\fk 1}}
	\!\fM^{\fk 1}_\al\Big)\,.
\end{alignat*}
For every $(C,D)\inn\fM^{\fk 1}$, recall $N(C)$ consists of all the nodes of $C$;
c.f.~\S\ref{Subsec:M2Pnd_loc_eqn}.
For every $\al\inn B^{\fk 1}$,
consider the following subset of the sections of the restriction of the universal curve $\cC/\fD_2$ to $\fM^{\fk 1}_{\al}$:
\begin{align}\label{Eqn:S_al}
	S_{\al}\subset \Ga(\fM^{\fk 1}_{\al};\cC)
	\qquad\tn{satisfying}\qquad
	\big\{\,e(C,D):\,e\inn S_{\al}\,\big\}
	=N(C)\quad
	\forall~(C,D)\inn\fM^{\fk 1}\,.
\end{align}
For every $e\inn S_{\al^{\fk 1}}$ and $(C,D)\inn\fM^{\fk 1}$,
we may simply write the node $e(C,D)$  as $e$ should no ambiguity occur (because each node of $C$ uniquely determines an element of $S_{\al}$).

For every $\al\inn B^{\fk 1}$ and $x\eq (C,D)\inn \fM^{\fk 1}_{\al}$,
we take
an affine smooth chart $\cV_x\!\lra\!\fD_2$ containing~$x$ and satisfying the degeneracy of any point in $\cV_x$ is no worse than $x$.
Recall on $\cV_x$, there are the node-smoothing parameters $\ze_e^{\cV_x}$, $e\inn S_{\al}$, as in \S\ref{Subsec:M2Pnd_loc_eqn}.
The topology of $\fD_2$ then implies $\fV_{\al}:=\{\cV_x\}_{x\in\fM^{\fk 1}_{\al}}$,
$\al\inn B^{\fk 1}$,
along with the the node-smoothing parameters,
satisfy (\ref{Eqn:loc_Euc_transition}) and (\ref{Eqn:loc_Euc_equation}).
In other words,
$\fM\eq\bigcup_{\al\in A^{\fk 1}}\fM^{\fk 1}_{\al}$ is an LES as per Definition~\ref{Dfn:G-adim_fixture},
and the node-smoothing parameters on $\fD_2$ serve as the modular parameters.

Recall the subsets $\wp_i(x)$,
$\wp_{T_s,i}$, $Q_i(x)$,
of $N(C)$, as well as the subset $P_{T,i}(x)$ of the power set of $N(C)$, introduced in \S\ref{Subsec:M2Pnd_loc_eqn};
c.f.~Figure~\ref{Fig:four_paths}.
Given $\al\inn A^{\fk 1}$ and $1\!\le\!i\!\le\!m$,
there exists  $\wp_{\al;i}\!\subset\!S_\al$ satisfying
\begin{align*}
	\big\{\,e(x) :\, e\inn \wp_{\al;i}\,\big\}
	= \wp_i(x)\qquad
	\forall~x\in\fM^{\fk 1}_\al\,.
\end{align*}
Similarly,
if $\score(C)$ has an element $T$ for some $x\eq(C,D)\inn\fM^{\fk 1}_{\al}$,
then there exists a subset $P_{\al;T,i}$ of the power set of $S_{\al}$, satisfying
\begin{align*}
	\big\{\,\{e(x):e\inn\wp\}
	 :\, \wp\inn P_{\al^{\fk 1};T,i}\,\big\}
	= P_{T,i}(x)\,.
\end{align*}
Likewise,
if $C_\core$ is separable for some $x\eq(C,D)\inn\fM^{\fk 1}_\al$,
then by writing $\score(C)\eq\{T_1, T_2\}$,
there exists subsets $\wp_{\al;s,i}\!\subset\!S_\al$, $s\eq 1,2$, satisfying
\begin{align*}
	\big\{\,e(x) :\, e\inn \wp_{\al;s,i}\,\big\}
	= \wp_{T_s,i}(x)\qquad
	s\eq 1,2.
\end{align*}
Finally, if $Q_i(x)\!\ne\!\emptyset$ (resp.~$Q_{ij}(x)\!\ne\!\emptyset$) for some $x\inn\fM^{\fk 1}_\al$,
then there exists a subset $Q_{\al;i}$ (resp.~$Q_{\al;ij}$) of the power set of $S_\al$ satisfying
\begin{align*}
	\big\{\,\{\,e(x):e\inn \wp\,\}\,:\,\wp\in Q_{\al;i}\,\big\}
	=Q_i(x)\qquad
	\Big(\tn{resp.}~
	\big\{\,\{\,e(x):e\inn \wp\,\}\,:\,\wp\in Q_{\al;ij}\,\big\}
	=Q_{ij}(x)
	\Big)\,.
\end{align*}

From (\ref{Eqn:local_eqn_chi1}), (\ref{Eqn:local_eqn_chi2}),
and the fact that $\wp_{\al;s,i}\!\supset\! \wp_{\al;i}$,
we conclude that either term of the $i$-th column of (\ref{Eqn:str_hom}) has a factor
$\prod_{e\in \wp_{\al;i}}\!\ze_e$,
so we aim to construct a treelike structure~$(\tau_{\al},\be_{\al})$, $\al\inn A^{\fk 1}$, on $\fM^{\fk 1}$, so that the root-to-leaf paths are minimal $\wp_{\al;i}$'s.
More precisely,
for every $\al\inn A^{\fk 1}$,
we say $\wp_{\al;j}$ is \ts{minimal} if it is minimal among all $\wp_{\al;i}$, $1\!\le\!i\!\le\!m$, with respect to the inclusion of subsets.
For instance,
in Example~\ref{Eg:abcd},
all $\wp_{\al;i}$'s are minimal.
Another simple example is when some
$\wp_{\al;i_0}$ is empty, hence is minimal;
in such situation, there does not exist any nonempty minimal $\wp_{\al;i}$.

We then construct a rooted tree $\tau_\al$ whose root-to-leaf paths are exactly given by $\wp_{\al;i}$, $1\!\le\!i\!\le\!m$.
Concretely, we take
\begin{align*}
	I_{\al,\min}:=
	\{\,1\le i\le m\,:\,
	\wp_{\al;i}~\tn{is~minimal}\,
	\},\qquad
	\tau_\al:=
	\bigcup_{i\in I_{\al,\min}}\!\!\!
	\wp_{\al;i}\ \ (\subset S_\al).
\end{align*}
On the set $\tau_\al$,
consider the relation $\preceq$ given by
$e\!\preceq\! e'$ if and only if for any $x\eq (C,D)\inn\fM^{\fk 1}_\al$, any connected subcurve of $C$ containing $e(x)$ and $C_\core$ must contain $e'(x)$.
The fact that $C_\core$ is of (arithmetic) genus 2 ensures that $\preceq$ satisfies~(\ref{Eqn:tree_order}),
which turns $\tau_\al$ into a rooted tree in the sense of Definition~\ref{Dfn:Rooted_tree}.
Such $\tau_\al$ is equivalent to the \ts{terminally weighted tree} of~\cite{HLN}.
Finally, since $\tau_\al\!\subset\!S_\al$ as sets,
we simply take $\be_\al:\tau_{\al}\!\hookrightarrow\!S_\al$ to be the inclusion,
and then set
\begin{equation*}\begin{split}
 \La^{\fk 1}
 &:=
 \big(\,\tau_{\al},\,
  \be_{\al}\,
 \big)_{
  {\al}\in A^{\fk 1}
  }.
\end{split}\end{equation*}

\begin{prp}
\label{Prp:G1-level}
$\La^{\fk 1}$ is a treelike structure on 
$\fM^{\fk 1}$ as per Definition~\ref{Dfn:Treelike_structure}.
\end{prp}

\begin{proof}
The proof is parallel to its genus~1 counterpart in Example~\ref{Eg:genus_1},
with $\fM^\wt_1$ replaced by $\fD_2$,
hence we omit the details.
\end{proof}

Proposition~\ref{Prp:G1-level} and Theorem~\ref{Thm:tf_smooth} together imply the following:

\begin{crl}
\label{Crl:G1-tf} 
Let $$\fM^{\fk 2}:=(\fM^{\fk 1})_{\La^{\fk 1}}^{\tf}$$ be constructed as in Theorem~\ref{Thm:tf_smooth}~\ref{Cond:smooth_tf}.
Then, $\fM^{\fk 2}$ is a smooth algebraic stack that is endowed with an LES, whose index set of the strata is 
\begin{align*}
	A^{\fk 2}:=
	\ov{\La^{\fk 1}}=
	\big\{\,(\al,\ov\lt):\,
	\al\inn A^{\fk 1},\,
	\ov\lt\inn\ov\bT_{\al}(\La^{\fk 1})\,\big\}\,,
\end{align*}
as well as a proper and birational morphism $(\fM^{\fk 1})_{\La^{\fk 1}}^{\tf}\!\lra\!\fM^{\fk 1}\eq\fD_2$ (known as the forgetful morphism) that is isomorphic on $(\fM^{\fk 1})^{\mn}$.
\end{crl} 

To study the pullback of the structural homomorphism $\varphi$ of (\ref{Eqn:str_hom}),
we fix arbitrary $\al\inn B^{\fk 1}$, $x\eq(C,D)\inn\fM^{\fk 1}_\al$,
and $\ov\lt\inn\ov\bT_\al(\La^{\fk 1})$,
and write 
\begin{align}\label{Eqn:al_2}
	\al^{\fk 2}:=(\al,\ov\lt)
	=\big(\,\al,\,
	(\tau_\al,\ov\bE,\Dm(\ov\lt))\,\big)\in A^{\fk 2}\,.
\end{align}
We then fix a sufficiently small chart $\cU^{(\fk 2)}\inn\fV_{\al^{\fk 2}}$ containing a fixed lift $\ti x^{(\fk 2)}$ of $x$ in $\fM^{\fk 2}_{\al^{\fk 2}}$.
On $\cU^{(\fk 2)}$,
we denote by $\xi^{(\fk 2)}_s$, $s\inn S_{\al^{\fk 2}}$, the twisted parameters,
and by  $\ti\varphi^{(\fk 2)}$ the restriction to $\cU^{(\fk 2)}$ of the pullback of $\varphi$.
For every $1\!\le\!i,j\!\le\!m$,
let 
\begin{align*}
	\ov \bE_{\al;ij}:=
	\big\{\,\fE\in\ov\bE\,:\,
	\fE\cap \wp_{\al;i}\cap \wp_{\al;j}\ne\emptyset
	\,\big\}\,.
\end{align*}

\begin{crl}
	\label{Crl:Step2_str_hom}
With notation as above, we have the following.
\begin{enumerate}[leftmargin=*,label*=($\fk 2$\alph*)]
	\item \label{Cond:Step2_sep}
	Assume $C_\core$ is separable, and $\wp_{\al;s,i}\bsl \Dm(\ov\lt)\!\ne\!\emptyset$, $s\eq 1,2$, whenever $\wp_{\al;i}\!\subset\!\Dm(\ov\lt)$.
	Then, 
	\begin{align*}
		\ti\varphi^{(\fk 2)}=
		\Big(\prod_{\fE\in\ov\bE}\xi^{(\fk 2)}_\fE	\Big)\cdot
		\left[
		\begin{matrix}
			\prod_{e\in \wp_{\al;1,1}\bsl\Dm(\ov\lt)}\xi^{(\fk 2)}_e  & \cdots & 
			\prod_{e\in \wp_{\al;1,m}\bsl\Dm(\ov\lt)}\xi^{(\fk 2)}_e
			\\
			\prod_{e\in \wp_{\al;2,1}\bsl\Dm(\ov\lt)}\xi^{(\fk 2)}_e & \cdots &
			\prod_{e\in \wp_{\al;2,m}\bsl\Dm(\ov\lt)}\xi^{(\fk 2)}_e
		\end{matrix}
		\right]\,.
	\end{align*}
	
	\item \label{Cond:Step2_other}
	In all other scenarios, under suitable trivialization of (\ref{Eqn:str_hom_0}) and after re-ordering $\de_i$'s and $T_s$'s if necessary, we have
	\begin{align*}
		\ti\varphi^{(\fk 2)}=
		\Big(\prod_{\fE\in\ov\bE}\xi^{(\fk 2)}_\fE	\Big)\cdot
		\left[
		\begin{matrix}\,
			1 & 0 & \cdots & 0\\
			0 & \tw^{(\fk 2)}_{12} & \cdots &
			\tw^{(\fk 2)}_{1m}\,
		\end{matrix}
		\right]\,,
	\end{align*}
	where each $\tw_{1i}^{(\fk 2)}$ takes the form
	\begin{align*}
		\frac{\tw_{1i}^{(\fk 2)}}{
			\ti\ka_{1i}^{(\fk 2)}
			\cdot
			\Big(\prod_{\fE\in \ov\bE_{\al;1i}}
			\xi^{(\fk 2)}_{\fE}\Big)}
		=
		\begin{cases}
			\prod_{e\in \wp_{\al;2,i}\bsl\Dm(\ov\lt)}\xi^{(\fk 2)}_e
			&
			\tn{if}~C_\core~\tn{is~separable},
			\\
			\prod_{e\in \wp_{\al;i}\bsl\Dm(\ov\lt)}\xi^{(\fk 2)}_e
			&
			\tn{if}~C_\core~\tn{is~inseparable},
		\end{cases}
	\end{align*}
	in which $\ti\ka_{1i}^{(\fk 2)}$ denotes the pullback of $\ka_{1i}$ of (\ref{Eqn:ka_ij}).
\end{enumerate}
\end{crl}
\begin{proof}
	\ref{Cond:Step2_sep} follows directly from (\ref{Eqn:local_eqn_chi2}) and Corollary~\ref{Crl:prod_RL}.
	
	Under the assumption of \ref{Cond:Step2_other},
	each entry of $\ti\varphi^{(\fk 2)}$ still contains $\prod_{\fE\in\ov\bE}\xi^{(\fk 2)}_\fE$ as a factor.
	Moreover, there exists $1\!\le\!j\!\le\!m$ (and $1\!\le\!s\!\le\!2$ if $C_\core$ is separable) such that 
	$\wp_{\al;j}\bsl\Dm(\ov\lt)\eq\emptyset$
	(and $\wp_{\al;s,j}\bsl\Dm(\ov\lt)\eq \emptyset$, too, if $C_\core$ is separable).
	Therefore, we can re-order the $\de_i$'s if necessary so that $\de_j$ becomes $\de_1$, and obtain the desired form of $\ti\varphi^{(\fk 2)}$ after elementary row and column operations (particularly, switching the rows is equivalent to switching $T_1$ and $T_2$).
	The expressions of $\tw_{1i}^{(\fk 2)}$,
	$i\!>\!1$, follow from (\ref{Eqn:local_eqn_th}) and Theorem~\ref{Thm:tf_smooth}~\ref{Cond:smooth_parameters}.
\end{proof}

\subsection{Step~$\fk 2$ and prototypes for Steps~$\fk 4$ and $\fk 5$}\label{Subsec:Step2}

At this stage,
three treelike structures on $\fM^{\fk 2}$ should be taken into consideration:
\begin{itemize}[leftmargin=*]
	\item 
	a new treelike structure $\La^{\fk 2}$,
	which aims to handle the terms in \ref{Cond:Step2_sep} of Corollary~\ref{Crl:Step2_str_hom};
	\item 
	a new treelike structure $\Up^{\fk 4}$,
	which aims to handle the factors $\prod_{e\in \wp_{\al;i}\bsl\Dm(\ov\lt)}\!\xi_e^{(\fk 2)}$ the entries of $\ti\varphi^{(\fk 2)}$ in Corollary~\ref{Crl:Step2_str_hom},
	when $\wp_{\al;i}\bsl\Dm(\ov\lt)$, $1\!<\!i\!\le\!m$, are all nonempty; and
	\item 
	the derived treelike structure $\ov{\tn d}\La^{\fk 1}$ as in Proposition~\ref{Prp:Derived_treelike_1}, which aims to handle the  terms 
	\begin{align*}
		\Big(\prod_{\fE\in \ov\bE_{\al;1i}}\!\!
		\xi^{(\fk 2)}_{\fE}\ \Big)
		\cdot
		\Big(\prod_{e\in \wp_{\al;i}\bsl\Dm(\ov\lt)}\!\!\!\!\xi^{(\fk 2)}_e\ \Big)\
	\end{align*}
	of $\ti\th^{(\fk 2)}_{1i}$ in \ref{Cond:Step2_other} of Corollary~\ref{Crl:Step2_str_hom}.
\end{itemize} 
In Figure~\ref{Fig:blowup}, they are respectively illustrated  as the dashed black, dashed orange, and solid red arrows toward $\fM^{\fk 2}$.
We shall use $\La^{\fk 2}$ for $\fM^{\fk 3}/\fM^{\fk 2}$,
and the rest two as prototypes for $\fM^{\fk 5}/\fM^{\fk 4}$ and $\fM^{\fk 6}/\fM^{\fk 5}$, respectively.

First of all, observe that the derived treelike structure $\ov{\tn d}\La^{\fk 1}$ is completely determined by Proposition~\ref{Prp:Derived_treelike_1},
so only $\La^{\fk 2}$ and $\Up^{\fk 4}$ need to be constructed.

First, we consider $\La^{\fk 2}$.
Let $\al^{\fk 2}\eq(\al,\ov\lt)\inn A^{\fk 2}$ as in (\ref{Eqn:al_2}) and $x\eq(C,D)\inn\fM^{\fk 1}_\al$ be as in  Corollary~\ref{Crl:Step2_str_hom}.
We still begin with describing the root-to-leaf paths of $\tau_{\al^{\fk 2}}$.
\begin{itemize}[leftmargin=*,label=-]
	\item If $C_\core$ is inseparable,
	so is any other point of $\fM^{\fk 1}_\al$;
	we simply take $\tau_{\al^{\fk 2}}\eq\emptyset$ as a set.
	
	\item 
	If $C_\core$ is separable,
	so is any other point of $\fM^{\fk 1}_\al$;
	we write
	\begin{align*}
		\wp_{\al^{\fk 2};s,i}:=
		\wp_{\al;s,i}\bsl \Dm(\ov\lt),\qquad
		1\le s\le 2,~1\le i\le m.
	\end{align*}
	Once again, we say $\wp_{\al^{\fk 2};s,i}$ is \ts{minimal} if it is minimal with respect to the inclusion of subsets.
	We then set
	\begin{align*}
		I_{\al^{\fk 2},\min}:=
		\big\{\,
		(s,i)\in\{1,2\}\!\times\!\{1,\cdots,m\}\,:\,
		\wp_{\al^{\fk 2};s,i}~\tn{is~minimal}\,\big\},\qquad
		\tau_{\al^{\fk 2}}
		:=\!
		\bigcup_{(s,i)\in I_{\al^{\fk 2},\min}}
		\!\!\!\!\!\!\!  \wp_{\al^{\fk 2};s,i}\,.
	\end{align*}
\end{itemize}

To see $\tau_{\al^{\fk 2}}$ has a natural tree order, we need the following lemma.
\begin{lmm}
	\label{Lm:r1p2_tree}
	If $C_\core$ is separable,
	then for every $1\!\le\!i,j\!\le\!m$ such that both $(1,i)$ and $(2,j)$ belong to $I_{\al^{\fk 2},\min}$,
	the paths $\wp_{\al^{\fk 2};1,i}$ and $ \wp_{\al^{\fk 2};2,j}$ are disjoint.
\end{lmm}

\begin{proof}
	W.l.o.g.~we assume both $\wp_{\al^{\fk 2};1,i}$ and $ \wp_{\al^{\fk 2};2,j}$ are contained in the set of the nodes between $T_1$ and $T_2$,
	for otherwise the statement is obviously true.
	Let $N'\!:=\!\wp_{\al^{\fk 2};1,i}\!\cap\! \wp_{\al^{\fk 2};2,j}$.
	Suppose $N'\!\ne\!\emptyset$,
	then we have $\wp_{\al^{\fk 2};2,i}\eq \wp_{\al^{\fk 2};2,j}\bsl N'\!\subsetneq\! \wp_{\al^{\fk 2};2,j}$,
	which contradicts the minimality of $\wp_{\al^{\fk 2};2,j}$.
	Therefore, $N'\!=\!\emptyset$.
\end{proof}

By Lemma~\ref{Lm:r1p2_tree}, we have
\begin{align*}
	\tau_{\al^{\fk 2}}
	=
	\tau_{\al^{\fk 2},1}
	\sqcup
	\tau_{\al^{\fk 2},2}\,,
	\qquad
	\tn{where}\quad
	\tau_{\al^{\fk 2},1}:=
	\bigcup_{(1,i)\in I_{\al^{\fk 2},\min}}
	\!\!\!\!\!\! \wp_{\al^{\fk 2};1,i},\,\quad
	\tau_{\al^{\fk 2},2}:=
	\bigcup_{(2,j)\in I_{\al^{\fk 2},\min}}
	\!\!\!\!\!\!\wp_{\al^{\fk 2};2,j}\,.
\end{align*}
Each $\tau_{\al^{\fk 2},s}$ is a rooted tree,
whose tree order $\prec_s$ is given as follows: $e\!\preceq_s\!e'$ if and only if for every connected subcurve $C'$ of $C$ containing $T_s$ and $e(x)$,
we have $e'(x)\inn C'$.
The tree orders $\prec_1$ and $\prec_2$ together give rise to a tree order $\prec$ on $\tau_{\al^{\fk 2}}$,
satisfying $e\!\preceq_s\!e'$ if and only if there exists $s\inn\{1,2\}$ such that $e,e'\inn\tau_{\al^{\fk 2},s}$ and $e\!\preceq_s\!e'$.
In other words,
\begin{align*}
	\tau_{\al^{\fk 2}}\in\bT
	\qquad\forall~\al^{\fk 2}\inn A^{\fk 2}\,.
\end{align*}
Moreover, 
we observe that as a set (instead of a poset),
\begin{align}\label{Eqn:tau_r1p2}
	\tau_{\al^{\fk 2}}
	\subset S_\al\,\bsl\,\Dm(\ov\lt)
	= S_\al\bsl\be_\al(\Dm(\ov\lt))
	\subset \ov\bE\sqcup 
	\big(S_\al\bsl\be_\al(\Dm(\ov\lt))\big)
	=S_{\al^{\fk 2}}\,;
\end{align}
c.f.~(\ref{Eqn:al_2}) and Theorem~\ref{Thm:tf_smooth}~\ref{Cond:smooth_parameters}.
We thus take
$
\be_{\al^{\fk 2}}:\tau_{\al^{\fk 2}}
\!\hookrightarrow \!
S_{\al^{\fk 2}}
$
to be the natural inclusion,
and then set
\begin{align*}
	\La^{\fk 2}:=
	\big(\tau_{\al^{\fk 2}},\,\be_{\al^{\fk 2}}\big)_{\al^{\fk 2}\in A^{\fk 2}}\,.
\end{align*}

\begin{prp}
	\label{Prp:r1p2_treelike}
	$\La^{\fk 2}$ is a treelike structure on $\fM^{\fk 2}$.
\end{prp}

\begin{proof}
	Let $\al^{\fk 2}$ and $x$ be as above. 
	W.l.o.g.~we assume $\tau_{\al^{\fk 2}}\!\ne\!\tau_\bullet$,
	which implies $C_\core$ is separable.
	Given $e\inn S_{\al^{\fk 2}}\eq \ov\bE\!\sqcup\!\big(S_\al\bsl\be_\al(\Dm(\ov\lt))\big)$,
	let $(\al',\ov\lt')\!=\!\al^{\fk 2}_{(e)}$ and $\tau'\!=\!\tau_{\al^{\fk 2}_{(e)}}$.
	
	If $e\inn\fE_\cht\bsl\dot\fE_\cht$,
	then $\al'\eq\al$ and $\ov\lt'
	\eq\ov\lt_{(e)}$;
	in this case, it is straightforward to check $\tau'\eq\tau_{\al^{\fk 2}}\bsl e^\wedge$.
	Otherwise, we have
	$\al'\eq\al_{(e)}$.
	Recall $\phi_{\al;e}:\tau_{\al'}\!\lra\!\tau_\al\bsl e^\wedge$ denotes the bijection in Definition~\ref{Dfn:Treelike_structure}.
	Then, $\Dm(\ov\lt')\eq\phi_{\al;e}^{-1}\big(\Dm(\ov\lt)\big)$.
	It is then a direct  check that
	$\tau'$ is isomorphic to $\tau_{\al^{\fk 2}}\big\bsl \big(\be_\al^{-1}(e)\big)^\wedge$, and the isomorphism satisfies (\ref{Eqn:Treelike}).
\end{proof}

Next, we construct $\Up^{\fk 4}$.
Fix $\al^{\fk 2}\inn A^{\fk 2}$ as in (\ref{Eqn:al_2}).
Once again, we begin with describing the root-to-leaf paths of the rooted tree that $\Up^{\fk 4}$ assigns to $\al^{\fk 2}$.
Let
\begin{align*}
	\wp_{\al^{\fk 2};i}:=
	\wp_{\al;i}\,\bsl\, \Dm(\ov\lt),\qquad
	1\le i\le m.
\end{align*}
The construction of $\tau_\al$ in Step~${\fk 1}$, along with Definition~\ref{Dfn:Enhanced_rtl_seq}, implies at least one $\wp_{\al^{\fk 2};i}$ is empty.
We choose $1\!\le\! j_1\!\le\! m$ satisfying $\wp_{\al^{\fk 2};j_1}\eq \emptyset$, and set
\begin{align*}
	J_{\al^{\fk 2},\min}:=
	\{\,
	1\le i\le m:\,
	i\ne j_1,\,
	\wp_{\al^{\fk 2};i}~\tn{is~minimal}\,\},\qquad
	\u_{\al^{\fk 2}}
	:=
	\bigcup_{i\in J_{\al^{\fk 2},\min}}
	\!\!\!\! \wp_{\al^{\fk 2};i}\,.
\end{align*}
Here, the minimality is among all  $\wp_{\al^{\fk 2};i}$ except $i\eq j_1$ (w.r.t.~the inclusion of subsets).
Notice that if $\u_{\al^{\fk 2}}$ is nonempty,
then $\wp_{\al^{\fk 2};i}\!\neq\!\emptyset$ for all $i$ distinct from $j_1$,
hence the definition of $\u_{\al^{\fk 2}}$ does not depend on the choice of $j_1$ (even though $J_{\al^{\fk 2},\min}$ does).

The reason that $\u_{\al^{\fk 2}}$ is naturally a rooted tree is parallel to that of $\tau_\al$ in Step~$\fk 1$.
Moreover,
we observe that as a set (instead of a poset),
\begin{align}\label{Eqn:tau_r1p5}
	\u_{\al^{\fk 2}}
	\subset S_\al\,\bsl\,\Dm(\ov\lt)
	= S_\al\bsl\be_\al(\Dm(\ov\lt))
	\subset \ov\bE\sqcup 
	\big(S_\al\bsl\be_\al(\Dm(\ov\lt))\big)
	=S_{\al^{\fk 2}}\,.
\end{align}
We thus take
$
	\mathfrak j_{\al^{\fk 2}}:\u_{\al^{\fk 2}}
	\!\hookrightarrow \!
	S_{\al^{\fk 2}}
$
to be the natural inclusion,
and then set
\begin{align*}
	\Up^{\fk 4}:=
	\big(\u_{\al^{\fk 2}},\,\mathfrak j_{\al^{\fk 2}}\big)_{\al^{\fk 2}\in A^{\fk 2}}\,.
\end{align*}

\begin{prp}
	\label{Prp:r1p5_treelike_proto}
	$\Up^{\fk 4}$ is a treelike structure on $\fM^{\fk 2}$.
\end{prp}

\begin{proof}
	The  proof is analogous to that of Proposition~\ref{Prp:G1-level},
	hence is omitted.
\end{proof}

At this stage,
we can proceed with adding twisted fields to $\fM^{\fk 2}$ with respect to any among $\La^{\fk 2}$, $\Up^{\fk 4}$, or $\ov{\tn d}\La^{\fk 1}$,
and then taking the proper transform of the other two.
Mimicking the order of the blowups of~\cite{HLN},
we choose $\La^{\fk 2}$ in Step~$\fk 2$,
and reserve the proper transforms of $\Up^{\fk 4}$ and $\ov{\tn d}\La^{\fk 1}$ as the treelike structures in Steps~$\fk 4$ and $\fk 5$.

Theorems~\ref{Thm:tf_smooth}~\ref{Cond:smooth_tf} and~\ref{Thm:tf_smooth_revert}~\ref{Cond:smooth_tf_revert}
and Proposition~\ref{Prp:two_treelike_strs} together imply the following.
\begin{crl}
\label{Crl:G2-tf} 
Let 
$$
\fM^{\fk 3}:=(\fM^{\fk 2})_{\La^{\fk 2}}^{\tf},\qquad
\fN^{\fk 5}_{\fk 2}:=(\fM^{\fk 2})_{\Up^{\fk 4}}^{\tf},\qquad\tn{and}\qquad
\fN^{\fk 6}_{\fk 2}:=(\fM^{\fk 2})_{\ov{\tn d}\La^{\fk 1}}^{\rtf}
$$ be constructed as in Theorems~\ref{Thm:tf_smooth}~\ref{Cond:smooth_tf} and~\ref{Thm:tf_smooth_revert}~\ref{Cond:smooth_tf_revert}, respectively.
Then, $\fM^{\fk 3}$, $\fN^{\fk 5}_{\fk 2}$, and $\fN^{\fk 6}_{\fk 2}$ are smooth algebraic stacks that are respectively endowed with LESs as well as proper and birational  morphisms  (known as the forgetful morphism)
\begin{align*}
	\fM^{\fk 3}\lra\fM^{\fk 2},\qquad
	\fN^{\fk 5}_{\fk 2}
	\lra\fM^{\fk 2},\qquad\tn{and}\qquad
	\fN^{\fk 6}_{\fk 2}\lra\fM^{\fk 2}
\end{align*}
that are all isomorphic on $(\fM^{\fk 2})^{\mn}$.
The index set of the strata of $\fM^{\fk 3}$ 
is
\begin{align*}
	A^{\fk 3}:=
	\ov{\La^{\fk 2}}=
	\big\{\,(\al,\ov\lt,\ov\lt^{\fk 2}):\,
	(\al,
	\ov\lt)\inn A^{\fk 2},\,
	\ov\lt^{\fk 2}\inn\ov\bT_{(\al,\ov\lt)}(\La^{\fk 2})
	\,\big\}\,.
\end{align*}
Moreover, the proper transforms
\begin{align*}
	\ov\PT_{\La^{\fk 2}}\big(\Up^{\fk 4}\big)\qquad
	\tn{and}
	\qquad
	\ov\PT_{\La^{\fk 2}}\big(\ov{\tn d}\La^{\fk 1}\big)
\end{align*}
as in Proposition~\ref{Prp:two_treelike_strs}
are both treelike structures on $\fM^{\fk 3}$.
\end{crl} 

Let  $\al^{\fk 2}\inn A^{\fk 2}$ be as in (\ref{Eqn:al_2}),
and $x\eq(C,D)\inn\fM^{\fk 1}_\al$ and its lift $\ti x^{(\fk 2)}$ be fixed as in the setting of Corollary~\ref{Crl:Step2_str_hom}.
We then fix arbitrary
\begin{alignat}{2}\label{Eqn:al_3}
	\al^{\fk 3}&:=(\al^{\fk 2},\ov\lt^{\fk 2})
	=\big(\,\al^{\fk 2},\,
	(\tau_{\al^{\fk 2}},\ov\bE^{\fk 2},\Dm(\ov\lt^{\fk 2}))\,\big)&&\in A^{\fk 3}\,,
\end{alignat} 
as well as a sufficiently small chart $\cU^{(\fk 3)}\inn\fV_{\al^{\fk 3}}$ containing a fixed lift $\ti x^{(\fk 3)}$ of $\ti x^{(\fk 2)}$ in $\fM^{{\fk 3}}_{\al^{\fk 3}}$.
We denote by $\xi^{(\fk 3)}_s$, $s\inn S_{\al^{\fk 3}}$, the twisted parameters on $\cU^{(\fk 3)}$,
and by  $\ti\varphi^{(\fk 3)}$ the restriction to $\cU^{(\fk 3)}$ of the pullback of the structural homomorphism~$\varphi$.

\begin{crl}
	\label{Crl:Step3_str_hom}
	Under suitable trivialization of (\ref{Eqn:str_hom_0}) and after re-ordering $\de_i$'s and $T_s$'s if necessary, we have
	\begin{align*}
		\ti\varphi^{(\fk 3)}=
		\Big(\!\prod_{\fE\in\ov\bE\sqcup \ov\bE^{\fk 2}}\!\!\xi^{(\fk 3)}_\fE\Big)\cdot
		\left[
		\begin{matrix}\,
			1 & 0 & \cdots & 0\\
			0 & \tw_{12}^{(\fk 3)} & \cdots &
			\tw_{1m}^{(\fk 3)}\,
		\end{matrix}
		\right]\,,
	\end{align*}
	where each $\tw_{1i}^{(\fk 3)}$ takes the form
	\begin{align*}
		\frac{\tw_{1i}^{(\fk 3)}}{
		\ti\ka_{1i}^{(\fk 3)}
		\cdot
		\Big(\prod_{\fE\in \ov\bE_{\al;1i}}
		\xi^{(\fk 3)}_{\fE}\Big)}
		=
		\begin{cases}
		\prod_{e\in \wp_{\al^{\fk 2};2,i}\bsl\Dm(\ov\lt^{\fk 2})}\xi^{(\fk 3)}_e
		&
		\tn{if}~C_\core~\tn{is~separable},
		\\
		\prod_{e\in \wp_{\al^{\fk 2};i}\bsl\Dm(\ov\lt^{\fk 2})}\xi^{(\fk 3)}_e
		&
		\tn{if}~C_\core~\tn{is~inseparable},
		\end{cases}
	\end{align*}
	in which $\ti\ka_{1i}^{(\fk 3)}$ denotes the pullback of $\ka_{1i}$ of (\ref{Eqn:ka_ij}).
\end{crl}

\begin{proof}
	Consider the underlying stratum $\fM^{\fk 2}_{\al^{\fk 2}}$.
	If the assumption in \ref{Cond:Step2_sep} of Corollary~\ref{Crl:Step2_str_hom} is satisfied,
	comparing the root-to-leaf paths of $\tau_{\al^{\fk 2}}$ and the terms of $\ti\varphi^{(\fk 2)}$ after factoring out $\prod_{\fE\in\ov\bE}\xi^{(\fk 2)}_\fE$,
	we can mimic the proof of Corollary~\ref{Crl:Step2_str_hom}~\ref{Cond:Step2_other} and obtain the desired form of $\ti\varphi^{(\fk 3)}$.
	If the assumption in \ref{Cond:Step2_other} of Corollary~\ref{Crl:Step2_str_hom} is satisfied,
	then $\tau_{\al^{\fk 2}}\eq\tau_\bullet$,
	hence $\ov\bE^{\fk 2}\eq\emptyset$,
	so the expression of $\ti\varphi^{(\fk 2)}$ follows directly from Corollary~\ref{Crl:Step2_str_hom}~\ref{Cond:Step2_other}.
\end{proof}

\subsection{Step~$\fk 3$ and prototypes for Steps~$\fk 6$-$\fk 9$}\label{Subsec:Step3}

In Corollary~\ref{Crl:G2-tf},
we obtain three STFs: $\fM^{\fk 3}$, $\fN^{\fk 5}_{\fk 2}$, and $\fN^{\fk 6}_{\fk 2}$.
On these STFs, the following treelike structures should be taken into consideration:

\begin{itemize}[leftmargin=*]
	\item 
	a new treelike structure $\La^{\fk 3}$ on $\fM^{\fk 3}$,
	which aims to handle the factors $\ti\ka_{1i}^{(\fk 3)}\prod_{e\in \wp_{\al^{\fk 2};2,i}\bsl\Dm(\ov\lt^{\fk 2})}\xi^{(\fk 3)}_e$ 
	and $\ti\ka_{1i}^{(\fk 3)}\prod_{e\in \wp_{\al^{\fk 2};i}\bsl\Dm(\ov\lt^{\fk 2})}\xi^{(\fk 3)}_e$ respectively
	in $\tw^{(\fk 3)}_{1i}$'s of Corollary~\ref{Crl:Step3_str_hom};
	\item 
	a new treelike structure $\Up^{\fk 8}$ on $\fN^{\fk 5}_{\fk 2}$,
	which aims to handle the factors $\prod_{e\in \wp_{\al^{\fk 2};i}\bsl\Dm(\ov\lt^{\fk 2})}\!\xi_e^{(\fk 3)}$ of $\tw^{(\fk 3)}_{1i}$'s of $\ti\varphi^{(\fk 3)}$ in Corollary~\ref{Crl:Step3_str_hom},
	when there is exactly one $1\!<\!i\!\le\!m$ with $\wp_{\al^{\fk 2};i}\bsl\Dm(\ov\lt^{\fk 2})\eq\emptyset$;
	\item 
	the derived treelike structure $\ov{\tn d}\Up^{\fk 4}$ on $\fN^{\fk 5}_{\fk 2}$ as in Proposition~\ref{Prp:Derived_treelike_1};
	\item 
	the doubly derived treelike structure ${\tn d}_2\La^{\fk 1}$ on $\fN^{\fk 6}_{\fk 2}$ as in Proposition~\ref{Prp:Derived_treelike_2}; and
	\item 
	the second-order derived treelike structure $\ud{\partial}\ov{\tn d}\La^{\fk 1}$  on $\fN^{\fk 6}_{\fk 2}$ as in Proposition~\ref{Prp:Derived_treelike_second_order_1}.
\end{itemize} 
The last three treelike structures are all for handling the last terms in~(\ref{Eqn:local_eqn_remainder1}) and~(\ref{Eqn:local_eqn_remainder3}),
yet in different situations.
In Figure~\ref{Fig:blowup}, the treelike structures above are respectively illustrated  as the dashed black arrow toward $\fM^{\fk 3}$, the dashed light blue and solid dark green arrows toward $\fN^{\fk 5}_{\fk 2}$,
and the solid dark blue and light green arrows toward $\fN^{\fk 6}_{\fk 2}$.
We shall use $\La^{\fk 3}$ for $\fM^{\fk 4}/\fM^{\fk 3}$,
and the rest four as prototypes for $\fM^{\fk 9}/\fM^{\fk 8}$, $\fM^{\fk{10}}/\fM^{\fk 9}$, $\fM^{\fk 7}/\fM^{\fk 6}$, and $\fM^{\fk 8}/\fM^{\fk 7}$, respectively.

Once again, notice that $\ov{\tn d}\Up^{\fk 4}$, ${\tn d}_2\La^{\fk 1}$ and $\ud{\partial}\ov{\tn d}\La^{\fk 1}$ are completely determined by Propositions~\ref{Prp:Derived_treelike_1}, \ref{Prp:Derived_treelike_2} and~\ref{Prp:Derived_treelike_second_order_1}, respectively,
so only $\La^{\fk 3}$ and $\Up^{\fk 8}$ need to be constructed.

We begin with $\La^{\fk 3}$ on $\fM^{\fk 3}$.
Let $\al^{\fk 3}\inn A^{\fk 3}$ be as in (\ref{Eqn:al_3}), and
$x\eq(C,D)\inn\fM^{\fk 1}_\al$ along with its lifts $\ti x^{(\fk 2)}$ and $\ti x^{(\fk 3)}$ be as in the setting of Corollary~\ref{Crl:Step3_str_hom}.
We first describe the root-to-leaf paths of the rooted tree that $\La^{\fk 3}$ assigns to $\al^{\fk 3}$.

If $C_\core$ is inseparable, then $\score(C)$ contains at most one element,
denoted by $T_2$, satisfying
\begin{align}\label{Eqn:p3}
	\wp^{\fk 3}:=
	\wp\big\bsl\big(\Dm(\ov\lt)\!\sqcup\!\Dm(\ov\lt^{\fk 2})\big)
	\ne\emptyset\qquad
	\forall~\wp\in \bigcup_{1\le i\le m}\!\!\!P_{\al;T_2,i}.
\end{align}
If $C_\core$ is separable, then $\score(C)$ contains at most one element,
also denoted by $T_2$, satisfying
\begin{align*}
	\wp_{T_2;i}\big\bsl
	\big(\Dm(\ov\lt)\sqcup\Dm(\ov\lt^{\fk 2})\big)\ne\emptyset
	\qquad
	\forall~1\!\le\!i\!\le\!m\,;
\end{align*}
this condition implies (\ref{Eqn:p3}),
so we define $\wp^{\fk 3}$ in the same way.
In either situation,
we construct $\tau_{\al^{\fk 3}}$ as follows.

\begin{itemize}[leftmargin=*]
	\item 
If such $T_2$ does not exist,
we simply take $\tau_{\al^{\fk 3}}\eq\emptyset$ as a set.

\item
If such $T_2$ exists, we say $\wp^{\fk 3}$  is \ts{minimal} if it is minimal (w.r.t.~the inclusion) among all the contracted paths in (\ref{Eqn:p3}).
We then set
\begin{equation*}
	\begin{split}
	&I_{\al^{\fk 3},\min}:=
	\big\{\,\wp\in \bigcup_{1\le i\le m}\!\!\!P_{\al;T_2,i}\,:\,
	\wp^{\fk 3}~\tn{is~minimal}\,\big\},\qquad
	\tau_{\al^{\fk 3}}:=
	\bigcup_{\wp\,\in I_{\al^{\fk 3},\min}}
	\!\!\!\!\!\!
	\wp^{\fk 3}	\,.
	\end{split}
\end{equation*}
\end{itemize}

Notice that the set $\tau_{\al^{\fk 3}}$ in each case above is determined by $\al^{\fk 3}$ and independent of the choice of $x\inn\fM^{\fk 1}_{\al}$ or its lifts.

\begin{lmm}\label{Lm:r1p3_tree}
	Assume that $\tau_{\al^{\fk 3}}\!\ne\!\emptyset$, and there exists $1\!\le\!i\!\le\!m$ such that $P_{\al;T_2,i}\!\cap\!I_{\al^{\fk 3},\min}$ contains two distinct paths $\wp$ and $\wp'$.
	Then,
	$\wp^{\fk 3}\!\cap\!\wp_{\al;i}\eq(\wp')^{\fk 3}\!\cap\!\wp_{\al;i}\!=\!\emptyset$.
\end{lmm}

\begin{proof}
	Since $\tau_{\al^{\fk 3}}\!\ne\!\emptyset$, the aforementioned $T_2$ exists.
	We choose $1\!\le\!j\!\le\!m$ such that $\wp_{\al;j}\!\subset\!\Dm(\ov\lt)\!\sqcup\!\Dm(\ov\lt^{\fk 2})$.
	For each $\wp''\inn P_{\al;T_2,j}$,
	the topology of the genus 2 nodal curves ensures $\wp''\bsl \wp_{\al;j}$ is a subset of either $\wp$ or $wp'$,
	say $\wp$.
	Therefore,
	\begin{align*}
		(\wp'')^{\fk 3}=
		\big(\wp''\bsl \wp_{\al;j}\big)\big\bsl
		\big(\Dm(\ov\lt)\sqcup\Dm(\ov\lt^{\fk 2})\big)
		\subset 
		\wp^{\fk 3}.
	\end{align*}
	Notice that $(\wp'')^{\fk 3}$ is disjoint from $\wp_{\al;i},$
	hence if $\wp^{\fk 3}\!\cap\!\wp_{\al;i}\!\neq\!\emptyset$,
	then $(\wp'')^{\fk 3}\!\subsetneq\!\wp^{\fk 3}$,
	which would contradict the minimality of $\wp^{\fk 3}$.
	Therefore,
	$\wp^{\fk 3}\!\cap\!\wp_{\al;i}\!=\!\emptyset$,
	and thus $(\wp')^{\fk 3}\!\cap\!\wp_{\al;i}\eq \wp^{\fk 3}\!\cap\!\wp_{\al;i}\!=\!\emptyset$.
\end{proof}

Lemma~\ref{Lm:r1p3_tree} determines a partial order $\preceq$ on $\tau_{\al^{\fk 3}}$ so that $e\!\preceq\!e'$ if and only if every path $\wp$ of $\ga_{x,(T_2)}$ starting from $v_{T_2}$ and satisfying
$
	e\inn \wp\!\subset\! \tau_{\al^{\fk 3}}
$ must contain $e'$.
It is straightforward that such $\preceq$ is indeed a tree order satisfying (\ref{Eqn:tree_order}),
i.e.
$$\tau_{\al^{\fk 3}}\in\bT\qquad
\forall~\al^{\fk 3}\inn A^{\fk 3}.$$
Moreover,
we observe that as a set (instead of a poset),
\begin{equation}\label{Eqn:tau_r1p3}
	\begin{split}
	\tau_{\al^{\fk 3}}
	\subset S_\al\,\big\bsl\,\big(\Dm(\ov\lt)\!\sqcup\!\Dm(\ov\lt^{\fk 2})\big)
	&= S_\al\,\big\bsl\,\big(\be_\al(\Dm(\ov\lt))\!\sqcup\!\be_{\al^{\fk 2}}(\Dm(\ov\lt^{\fk 2}))\big)\\
	&\subset
	\ov\bE\sqcup 
	\ov\bE^{\fk 2}\sqcup 
	\Big(S_{\al}\big\bsl\big(\be_\al(\Dm(\ov\lt))\sqcup \be_{\al^{\fk 2}}(\Dm(\ov\lt^{\fk 2}))\big)\Big)
	= S_{\al^{\fk 3}}\,;
	\end{split}
\end{equation}
c.f.~(\ref{Eqn:al_3}).
We thus take
$
\be_{\al^{\fk 3}}:\tau_{\al^{\fk 3}}
\!\hookrightarrow \!
S_{\al^{\fk 3}}
$
to be the natural inclusion,
and then set
\begin{align*}
	\La^{\fk 3}:=
	\big(\tau_{\al^{\fk 3}},\,\be_{\al^{\fk 3}}\big)_{\al^{\fk 3}\in A^{\fk 3}}\,.
\end{align*}

\begin{prp}
	\label{Prp:r1p3_treelike}
	$\La^{\fk 3}$ is a treelike structure on $\fM^{\fk 3}$.
\end{prp}

\begin{proof}
	The proof is parallel to that of Proposition~\ref{Prp:r1p2_treelike},
	hence is omitted.
\end{proof}

Next, we construct $\Up^{\fk 8}$ on $\fN^{\fk 5}_{\fk 2}$,
mimicking the construction of $\Up^{\fk 4}$ on $\fM^{\fk 2}$.
Recall the index set of the strata of $\fN^{\fk 5}_{\fk 2}$,
as in Theorem~\ref{Thm:tf_smooth}~\ref{Cond:smooth_tf},
is 
\begin{align*}
	\wh A^{\fk 5}:=
	\ov{\Up^{\fk 4}}=
	\big\{\,(\al,\ov\lt,\ov\ls):\,
	(\al,
	\ov\lt)\inn A^{\fk 2},\,
	\ov\ls\inn\ov\bT_{(\al,\ov\lt)}(\Up^{\fk 4})
	\,\big\}\,.
\end{align*}
We fix arbitrary
\begin{alignat}{2}\label{Eqn:al_5}
	\al^{\fk 2}&:=(\al,\ov\lt)
	=\big(\,\al,\,
	(\tau_\al,\ov\bE,\Dm(\ov\lt))\,\big)&&\in A^{\fk 2}\,,\nonumber\\
	\wh\al^{\fk 5}&:=(\al^{\fk 2},\ov\ls)
	=\big(\,\al^{\fk 2},\,
	(\u_{\al^{\fk 2}},\ov\bF,\Dm(\ov\ls))\,\big)&&\in\wh A^{\fk 5}\,.
\end{alignat}
By Theorem~\ref{Thm:tf_smooth}~\ref{Cond:smooth_parameters} and (\ref{Eqn:tau_r1p5}), we have
\begin{align*}
	S_{\wh\al^{\fk 5}}=
	\ov\bF\sqcup 
	\big(S_{\al^{\fk 2}}\bsl\mathfrak j_{\al^{\fk 2}}(\Dm(\ov\ls))\big)=
	\ov\bE\sqcup 
	\ov\bF\sqcup 
	\Big(S_{\al}\big\bsl\big(\be_\al(\Dm(\ov\lt))\sqcup \mathfrak j_{\al^{\fk 2}}(\Dm(\ov\ls))\big)\Big)
	\,.
\end{align*}
Let
\begin{align*}
	\wp_{\wh\al^{\fk 5};i}:=
	\wp_{\al;i}\,\big\bsl\, \big(\Dm(\ov\lt)\!\sqcup\!\Dm(\ov\ls)\big),\qquad
	1\le i\le m.
\end{align*}
The construction of $\u_{\al^{\fk 2}}$, along with Definition~\ref{Dfn:Enhanced_rtl_seq}, implies at least two $\wp_{\wh\al^{\fk 5};i}$ are empty.
We choose $1\!\le\! j_1\!<\!j_2\!\le\! m$ satisfying $\wp_{\wh\al^{\fk 5};j_1}\eq \wp_{\wh\al^{\fk 5};j_2}\eq \emptyset$, and set
\begin{align*}
	J_{\wh\al^{\fk 5},\min}:=
	\{\,
	i\inn\{1,\cdots,m\}\bsl\{j_1,j_2\}:\,
	\wp_{\wh\al^{\fk 5};i}~\tn{is~minimal}\,\},\qquad
	\u_{\wh\al^{\fk 5}}
	:=
	\bigcup_{i\in J_{\wh\al^{\fk 5},\min}}
	\!\!\!\! \wp_{\wh\al^{\fk 5};i}\,.
\end{align*}
As always, the minimality is  among all $\wp_{\wh\al^{\fk 5};i}$ except $i\eq j_1, j_2$ (w.r.t.~inclusion).
Notice that if $\u_{\al^{\fk 2}}$ is nonempty,
then $\wp_{\wh\al^{\fk 5};i}\!\neq\!\emptyset$ for all $i$ distinct from $j_1$ and $j_2$,
hence the definition of $\u_{\al^{\fk 5}}$ does not depend on the choice of $j_1$ and $j_2$ (even though $J_{\wh \al^{\fk 5},\min}$ does).

The reason that $\u_{\wh\al^{\fk 5}}$ is naturally a rooted tree is parallel to that of $\tau_\al$ and $\u_{\al^{\fk 2}}$.
Moreover,
we observe that as a set (instead of a poset),
\begin{align*}
	\u_{\wh\al^{\fk 5}}
	\subset S_\al\,\big\bsl\, \big(\Dm(\ov\lt)\!\sqcup\!\Dm(\ov\ls)\big)
	= S_\al\,\big\bsl\, \big(\be_\al(\Dm(\ov\lt))\!\sqcup\!\mathfrak j_{\al^{\fk 2}}(\Dm(\ov\ls))\big)
	\subset S_{\wh\al^{\fk 5}}\,.
\end{align*}
We thus take
$
\mathfrak j_{\wh \al^{\fk 5}}:\u_{\wh \al^{\fk 5}}
\!\hookrightarrow \!
S_{\wh\al^{\fk 5}}
$
to be the natural inclusion,
and then set
\begin{align*}
	\Up^{\fk 8}:=
	\big(\u_{\wh\al^{\fk 5}},\,\mathfrak j_{\wh\al^{\fk 5}}\big)_{\wh\al^{\fk 5}\in \wh A^{\fk 5}}\,.
\end{align*}
The following statement is the analogue of Propositions~\ref{Prp:G1-level} and~\ref{Prp:r1p5_treelike_proto};
its proof is parallel and hence is omitted.

\begin{prp}
	\label{Prp:r3p3_treelike_proto}
	$\Up^{\fk 8}$ is a treelike structure on $\fN^{\fk 5}_{\fk 2}$.
\end{prp}

At this stage,
we can proceed with adding twisted fields to $\fM^{\fk 3}$ with respect to any among $\La^{\fk 3}$, $\ov\PT_{\La^{\fk 2}}(\Up^{\fk 4})$, or $\ov\PT_{\La^{\fk 2}}(\ov{\tn d}\La^{\fk 1})$,
and then taking the proper transform of the other two.
Based on the order of the blowups in~\cite{HLN},
we choose $\La^{\fk 3}$ in Step~$\fk 3$,
and reserve the proper transforms of $\Up^{\fk 4}$ and $\ov{\tn d}\La^{\fk 1}$ for Steps~$\fk 4$ and $\fk 5$.

Theorems~\ref{Thm:tf_smooth}~\ref{Cond:smooth_tf} and~\ref{Thm:tf_smooth_revert}~\ref{Cond:smooth_tf_revert}
and Propositions~\ref{Prp:two_treelike_strs} and~\ref{Prp:2nd_PT} together imply the following.
\begin{crl}
	\label{Crl:G3-tf} 
	Let 
	$$
	\fM^{\fk 4}:=(\fM^{\fk 3})_{\La^{\fk 3}}^{\tf},\qquad
	\fN^{\fk 5}_{\fk 3}:=(\fM^{\fk 3})_{\ov\PT_{\La^{\fk 2}}(\Up^{\fk 4})}^{\tf},\qquad\tn{and}\qquad
	\fN^{\fk 6}_{\fk 3}:=(\fM^{\fk 3})_{\ov\PT_{\La^{\fk 2}}(\ov{\tn d}\La^{\fk 1})}^{\rtf}
	$$ be constructed as in Theorems~\ref{Thm:tf_smooth}~\ref{Cond:smooth_tf} and~\ref{Thm:tf_smooth_revert}~\ref{Cond:smooth_tf_revert}, respectively.
	Then, $\fM^{\fk 4}$, $\fN^{\fk 5}_{\fk 3}$, and $\fN^{\fk 6}_{\fk 3}$ are smooth algebraic stacks that are respectively endowed with LESs as well as proper and birational  morphisms  (known as the forgetful morphism)
	\begin{align*}
		\fM^{\fk 4}\lra\fM^{\fk 3},\qquad
		\fN^{\fk 5}_{\fk 3}
		\lra\fM^{\fk 3},\qquad\tn{and}\qquad
		\fN^{\fk 6}_{\fk 3}\lra\fM^{\fk 3}
	\end{align*}
	that are all isomorphic on $(\fM^{\fk 3})^{\mn}$.
	The index set of the LES of $\fM^{\fk 4}$ is
	\begin{align*}
		A^{\fk 4}:=
		\ov{\La^{\fk 3}}=
		\big\{\,(\al,\ov\lt,\ov\lt^{\fk 2},\ov\lt^{\fk 3}):\,
		(\al,
		\ov\lt,\ov\lt^{\fk 2})\inn A^{\fk 3},\,
		\ov\lt^{\fk 3}\inn\ov\bT_{(\al,\ov\lt,\ov\lt^{\fk 2})}(\La^{\fk 3})
		\,\big\}\,.
	\end{align*}
	Moreover, 
	the secondary RL proper transforms
	\begin{align*}
		\ov\PT_{\La^{\fk 2}}(\Up^{\fk 8}),
		\quad
		\ov\PT_{\La^{\fk 2}}(\ov{\tn d}\Up^{\fk 4})\quad\tn{on}\quad
		\fN^{\fk 5}_{\fk 2}
		\qquad\tn{and}\qquad
		\ov\PT_{\La^{\fk 2}}({\tn d}_2\La^{\fk 1}),
		\quad
		\ov\PT_{\La^{\fk 2}}(\ud{\partial}\ov{\tn d}\La^{\fk 1})
		\quad\tn{on}\quad
		\fN^{\fk 6}_{\fk 2}
	\end{align*}
	as in Proposition~\ref{Prp:2nd_PT}
	are respectively treelike structures.
	Furthermore, the RL proper transforms
	\begin{align*}
		\ov\PT_{\La^{\fk 3}}\big(\ov\PT_{\La^{\fk 2}}\big(\Up^{\fk 4}\big)\big)\qquad
		\tn{and}
		\qquad
		\ov\PT_{\La^{\fk 3}}\big(\ov\PT_{\La^{\fk 2}}\big(\ov{\tn d}\La^{\fk 1}\big)\big)
	\end{align*}
	as in Proposition~\ref{Prp:two_treelike_strs}
	are both treelike structures on $\fM^{\fk 4}$.
\end{crl}

Let  $\al^{\fk 3}\inn A^{\fk 3}$ be as in (\ref{Eqn:al_3}),
and $x\eq(C,D)\inn\fM^{\fk 1}_\al$ and its lifts $\ti x^{(\fk 2)}$ and $\ti x^{(\fk 3)}$ be fixed as in the setting of Corollary~\ref{Crl:Step3_str_hom}.
We then fix arbitrary
\begin{alignat}{2}\label{Eqn:al_4}
	\al^{\fk 4}&:=(\al^{\fk 3},\ov\lt^{\fk 3})
	=\big(\,\al^{\fk 3},\,
	(\tau_{\al^{\fk 3}},\ov\bE^{\fk 3},\Dm(\ov\lt^{\fk 3}))\,\big)&&\in A^{\fk 4}\,,
\end{alignat}  
as well as a sufficiently small chart $\cU^{(\fk 4)}\inn\fV_{\al^{\fk 4}}$  containing a fixed lift $\ti x^{(\fk 4)}$ of $\ti x^{(\fk 3)}$ in $\fM^{{\fk 4}}_{\al^{\fk 4}}$.
We denote by $\xi^{(\fk 4)}_s$, $s\inn S_{\al^{\fk 4}}$, the twisted parameters on $\cU^{(\fk 4)}$,
and by  $\ti\varphi^{(\fk 4)}$ the restriction to $\cU^{(\fk 4)}$ of the pullback of the structural homomorphism~$\varphi$.
For each $1\!<\!i\!\le\!m,$
we denote by $\wc\ka_{1i}^{(\fk 4)}$ the proper transform of $\ka_{1i}$ of (\ref{Eqn:ka_ij});
i.e.~on $\cU^{(\fk 4)}$, the proper transform of the conjugate locus $\wh\cK_{m;1i}$ is given by $\{\wc\ka_{1i}^{(\fk 4)}\eq 0\}$.

\begin{crl}
	\label{Crl:Step4_str_hom}
	Under suitable trivialization of (\ref{Eqn:str_hom_0}) and after re-ordering $\de_i$'s and $T_s$'s if necessary, we have
	\begin{align*}
		&	\ti\varphi^{(\fk 4)}=
		\left[
		\begin{matrix}\,
			\prod_{\fE\in\ov\bE\sqcup \ov\bE^{\fk 2}}\xi^{(\fk 4)}_\fE & 0\\
			0 & 
			\prod_{\fE\in\ov\bE\sqcup \ov\bE^{\fk 2}\sqcup \ov\bE^{\fk 3}}\xi^{(\fk 4)}_\fE\,
		\end{matrix}
		\right]
		\cdot
		\left[
		\begin{matrix}\,
			1 & 0 & \cdots & 0\\
			0 & \wc\tw_{12}^{(\fk 4)} & \cdots &
			\wc\tw_{1m}^{(\fk 4)}\,
		\end{matrix}
		\right]\,,
	\end{align*}
	where each $\wc\tw_{1i}^{(\fk 4)}$, $1\!<\!i\!\le\!m,$ satisfies
	\begin{align*}
		\frac
		{\wc\tw_{1i}^{(\fk 4)}}
		{\wc\ka_{1i}^{(\fk 4)}
		\cdot
		\Big(\prod_{\fE\in \ov\bE_{\al;1i}}
		\!\xi^{(\fk 4)}_{\fE}\,\Big)}
		=
		\begin{cases}
			\prod_{e\in\wp_{\al;2,i}\bsl
				(\Dm(\ov\lt)\sqcup\Dm(\ov\lt^{\fk 2})\sqcup\Dm(\ov\lt^{\fk 3}))}\xi^{(\fk 4)}_e
			&
			\tn{if}~C_\core~ \tn{is~separable},
			\\
			\prod_{e\in\wp_{\al;i}\bsl
				(\Dm(\ov\lt)\sqcup\Dm(\ov\lt^{\fk 2})\sqcup\Dm(\ov\lt^{\fk 3}))}\xi^{(\fk 4)}_e
			&
			\tn{if}~C_\core~ \tn{is~inseparable}.
		\end{cases}	
	\end{align*}
	Moreover,
	\begin{enumerate}[leftmargin=*,label*=($\fk 4$\alph*)]
		\item \label{Cond:Step4_nonep_bridge}
		If $Q_{\al;1}\!\ne\!\emptyset$, $Q_{\al;i}\!\ne\!\emptyset$, and $T^{(i)}\eq T_2$,
		then 
		\begin{align}\label{Eqn:Step4_kappa}
			\wc\ka_{1i}^{(\fk 4)}=
			\sum_{\wp\in Q_{\al;1i}}
			\!
			\Big(\ti u_\wp^{(\fk 4)}
			\hspace{-.1in} \prod_{e\in\wp\bsl(\Dm(\ov\lt)\sqcup\Dm(\ov\lt^{\fk 2})\sqcup\Dm(\ov\lt^{\fk 3}))}
			\hspace{-.3in}\xi^{(\fk 4)}_e\ \  \Big),
		\end{align}
		where each $\ti u_{\wp}^{(\fk 4)}$ denotes the pullback of the coefficient $u_\wp$ of Lemma~\ref{Lm:kappa}.
		
		\item \label{Cond:Step4_other}
		Otherwise, 
		we have $\wc\ka_{1i}^{(\fk 4)}\eq \ti\ka_{1i}^{(\fk 4)}$, the pullback of $\ka_{1i}$ of (\ref{Eqn:ka_ij}),
		which is either a unit, 
		or a local parameter satisfying that
		$
			\{\ti\ka_{1i}^{(\fk 4)}\}\!\sqcup\!
			\{\xi^{(\fk 4)}_s\}_{s\in S_{\al^{\fk 4}}}
		$
		forms a subset of a system of local parameters.	
	\end{enumerate}
\end{crl}

\begin{proof}
Under the assumption of~\ref{Cond:Step4_nonep_bridge},
the expression of $\wc\ka^{(\fk 4)}_{1i}$
follows immediately from  Lemma~\ref{Lm:kappa}.
To obtain the desired form of $\wc\tw^{(\fk 4)}_{1i}$,
we first consider the case when $C_\core$ is inseparable.
For each $\wp\inn Q_{\al;1i}$,
there exists $\wp'\inn P_{\al;T_2,1}\!\sqcup\!P_{\al;T_2,i}$ such that $\wp\!\sqcup\!\wp_{\al;i}\!\subset\!\wp'$.
\begin{itemize}[leftmargin=*]
	\item If $\wp\!\sqcup\!\wp_{\al;i}\!=\!\wp'$,
	then $\big(\wp\sqcup\wp_{\al;i}\big)\big\bsl\big(\Dm(\ov\lt)\sqcup\Dm(\ov\lt^{\fk 2})\sqcup\Dm(\ov\lt^{\fk 3})\big)\eq(\wp')^{\fk 3}.$
	
	\item  If
	$\wp\!\sqcup\!\wp_{\al;i}\!\subsetneq\!\wp'$,
	then there exists $\wp_1\inn Q_{\al;1}$ such that $\wp\!\sqcup\!\wp_{\al;1}\eq\wp_1$.
	Since $\wp_{\al;1}\!\subset\!\Dm(\ov\lt)$,
	we have
	$\wp\big\bsl \big(\Dm(\ov\lt)\sqcup\Dm(\ov\lt^{\fk 2})\sqcup\Dm(\ov\lt^{\fk 3})\big) \eq(\wp_1)^{\fk 3}$.
\end{itemize}
Either way, the construction of $\tau_{\al^{\fk 3}}$ implies $\big(\wp\!\sqcup\!\wp_{\al;i}\big)\big\bsl\big(\Dm(\ov\lt)\!\sqcup\!\Dm(\ov\lt^{\fk 2})\!\sqcup\!\Dm(\ov\lt^{\fk 3})\big)\eq(\wp')^{\fk 3}$,
which contains a root-to-leaf path $\wp''$ of $\tau_{\al^{\fk 3}}$ as a subset.
Hence the expression of $\wc\tw^{(\fk 4)}_{1i}$ follows from Corollary~\ref{Crl:Step3_str_hom} and the expression of  $\wc\ka^{(\fk 4)}_{1i}$.

The same argument applies to the case when $C_\core$ is separable in~\ref{Cond:Step4_nonep_bridge} verbatim, expect $\wp_{\al;i}$ is replaced by $\wp_{\al;2,i}$.
Under the assumption of~\ref{Cond:Step4_other},
the expressions of  $\wc\ka^{(\fk 4)}_{1i}$ and  $\wc\tw^{(\fk 4)}_{1i}$ follow directly from Corollary~\ref{Crl:Step3_str_hom}.
\end{proof}

\subsection{Steps $\fk 4$ and $\fk 5$}\label{Subsec:Step4}

In these two steps,
no new treelike structures need to be introduced.
We simply start with $\fM^{\fk 4}$ and add RL twisted fields with respect to the RL proper transform of $\Up^{\fk 4}$,
then add LR twisted fields with respect to the RL proper transform of $\ov{\tn d}\La^{\fk 1}$.

First, we observe in Corollary~\ref{Crl:Step4_str_hom},
each term $\wc\tw^{(\fk 4)}_{1i}$ contains a factor 
$
	\prod_{e\in\wp_{\al^{\fk 4};i}}
	\xi^{(\fk 4)}_e,
$
where
\begin{align*}
	\wp_{\al^{\fk 4};i}:=
	\wp_{\al;i}\big\bsl
	\big(\Dm(\ov\lt)\sqcup\Dm(\ov\lt^{\fk 2})\sqcup\Dm(\ov\lt^{\fk 3})\big)\,,
	\quad
	1\!<\!i\!\le\!m.
\end{align*}

\begin{lmm}\label{Lm:Step_4_lemma}
	If $C_\core$ is separable,
	i.e.~$\score(C)\eq\{T_1,T_2\}$,
	we have the following.
	\begin{itemize}[leftmargin=*]
		\item If there exists $1\!<\!i\!\le\!m$ so that
		$\wp_{\al^{\fk 4};i}\eq\emptyset$, then $\wp_{\al;2,i}\big\bsl
		\big(\Dm(\ov\lt)\!\sqcup\!\Dm(\ov\lt^{\fk 2})\!\sqcup\!\Dm(\ov\lt^{\fk 3})\big)
		\eq\emptyset$.
				
		\item If 
		$\wp_{\al^{\fk 4};i}\!\neq\!\emptyset$ for all $1\!<\!i\!\le\!m$,
		then 
		$\wp_{\al;2,i}\big\bsl
			\big(\Dm(\ov\lt)\!\sqcup\!\Dm(\ov\lt^{\fk 2})\!\sqcup\!\Dm(\ov\lt^{\fk 3})\big)
			\eq
			\wp_{\al^{\fk 4};i}$
		for all $1\!\le\!i\!\le\!m$.
	\end{itemize}
\end{lmm}

\begin{proof}
	In the former case, notice that $\wp_{\al;2,i}\bsl \wp_{\al;i}$,
	if nonempty,
	is comprised of certain nodes lying between $T_1$ and $T_2$,
	each of which is closer to the root than any element of $\wp_{\al;i}$ in $\tau_{\al^{\fk 2}}$ and $\tau_{\al^{\fk 3}}$.
	Therefore,
	$(\wp_{\al;2,i}\bsl \wp_{\al;i})\!\subset\!\Dm(\ov\lt)\!\sqcup\!\Dm(\ov\lt^{\fk 2})\!\sqcup\!\Dm(\ov\lt^{\fk 3})$ as well,
	thus $\wp_{\al;2,i}\big\bsl
	\big(\Dm(\ov\lt)\!\sqcup\!\Dm(\ov\lt^{\fk 2})\!\sqcup\!\Dm(\ov\lt^{\fk 3})\big)
	\eq\emptyset$.
	
	In the latter case,
	regardless of the position of $\de_1$ on $C$,
	we always have
	$
	(\wp_{\al;1,1}\!\sqcup\!\wp_{\al;2,1})\!\subset\!
	\Dm(\ov\lt)\!\sqcup\!\Dm(\ov\lt^{\fk 2})\!\sqcup\!\Dm(\ov\lt^{\fk 3}).
	$
	Particularly,
	every node between $T_1$ and $T_2$ is contained in $\Dm(\ov\lt)\!\sqcup\!\Dm(\ov\lt^{\fk 2})\!\sqcup\!\Dm(\ov\lt^{\fk 3})$. 
	Once again, this implies 
	$(\wp_{\al;2,i}\bsl \wp_{\al;i})\!\subset\!\Dm(\ov\lt)\!\sqcup\!\Dm(\ov\lt^{\fk 2})\!\sqcup\!\Dm(\ov\lt^{\fk 3})$,
	which leads to the proposed equalities.
\end{proof}

From Lemma~\ref{Lm:Step_4_lemma}, 
we see
in order to simplify the second row of the matrix of Corollary~\ref{Crl:Step4_str_hom},
it suffices to compare the factors $\prod_{e\in\wp_{\al^{\fk 4};i}}
\xi_e^{(\fk 4)}$, $1\!<\!i\!\le\!m.$
For this purpose, 
the proper transform of the treelike structure $\Up^{\fk 4}$ of Proposition~\ref{Prp:r1p5_treelike_proto} is exactly what we need. 
More concretely, in Step~$\fk 4$,
we set
\begin{align*}
	\La^{\fk 4}:=\ov\PT_{\La^{\fk 3}}(\ov\PT_{\La^{\fk 2}}(\Up^{\fk 4}))\,.
\end{align*}
As shown in Corollary~\ref{Crl:G3-tf},
$\La^{\fk 4}$ is a treelike structure on $\fM^{\fk 4}$.

Theorems~\ref{Thm:tf_smooth}~\ref{Cond:smooth_tf} and~\ref{Thm:tf_smooth_revert}~\ref{Cond:smooth_tf_revert}
and Propositions~\ref{Prp:two_treelike_strs} and~\ref{Prp:2nd_PT} together imply the following.

\begin{crl}
	\label{Crl:G4-tf} 
	Let 
	$$	
	\fM^{\fk 5}:=(\fM^{\fk 4})_{\La^{\fk 4}}^{\tf}\qquad\tn{and}\qquad
	\fN^{\fk 6}_{\fk 4}:=(\fM^{\fk 4})_{\ov\PT_{\La^{\fk 3}}(\ov\PT_{\La^{\fk 2}}(\ov{\tn d}\La^{\fk 1}))}^{\rtf}
	$$ be constructed as in Theorems~\ref{Thm:tf_smooth}~\ref{Cond:smooth_tf} and~\ref{Thm:tf_smooth_revert}~\ref{Cond:smooth_tf_revert}, respectively.
	Then, $\fM^{\fk 5}$ and $\fN^{\fk 6}_{\fk 4}$ are smooth algebraic stacks that are respectively endowed with LESs as well as proper and birational  morphisms  (known as the forgetful morphism)  
	\begin{align*}
		\fM^{\fk 5}\lra\fM^{\fk 4}\qquad\tn{and}\qquad
		\fN^{\fk 6}_{\fk 4}\lra\fM^{\fk 4}
	\end{align*}
	that are both isomorphic on $(\fM^{\fk 4})^{\mn}$.
	The index set of the LES of $\fM^{\fk 5}$ is
	\begin{align*}
		A^{\fk 5}:=
		\ov{\La^{\fk 4}}=
		\big\{\,(\al,\ov\lt,\ov\lt^{\fk 2},\ov\lt^{\fk 3},\ov\lt^{\fk 4}):\,
		(\al,
		\ov\lt,\ov\lt^{\fk 2},\ov\lt^{\fk 3})\inn A^{\fk 4},\,
		\ov\lt^{\fk 4}\inn\ov\bT_{(\al,\ov\lt,\ov\lt^{\fk 2},\ov\lt^{\fk 3})}(\La^{\fk 4})
		\,\big\}\,.
	\end{align*}
	Moreover, 
	the secondary RL proper transforms
	\begin{alignat*}{2}
		&\ov\PT_{\La^{\fk 3}}\big(\ov\PT_{\La^{\fk 2}}(\Up^{\fk 8})\big),
		\quad
		\ov\PT_{\La^{\fk 3}}
		\big(\ov\PT_{\La^{\fk 2}}(\ov{\tn d}\Up^{\fk 4})\big)\quad
		&&\tn{on}\quad
		\fM^{\fk 5}
		\qquad\tn{and}\\
		&\ov\PT_{\La^{\fk 3}}\big(\ov\PT_{\La^{\fk 2}}({\tn d}_2\La^{\fk 1})\big),
		\quad
		\ov\PT_{\La^{\fk 3}}\big(\ov\PT_{\La^{\fk 2}}(\ud{\partial}\ov{\tn d}\La^{\fk 1})\big)
		\quad
		&&\tn{on}\quad
		\fN^{\fk 6}_{\fk 4}
	\end{alignat*}
	as in Proposition~\ref{Prp:2nd_PT}
	are respectively treelike structures.
	Furthermore, the RL proper transform
	\begin{align*}
		\ov\PT_{\La^{\fk 4}}
		\big(\ov\PT_{\La^{\fk 3}}\big(\ov\PT_{\La^{\fk 2}}\big(\ov{\tn d}\La^{\fk 1}\big)\big)\big)
	\end{align*}
	as in Proposition~\ref{Prp:two_treelike_strs}
	is a treelike structure on $\fM^{\fk 5}$.
\end{crl}

Let  $\al^{\fk 4}\inn A^{\fk 4}$ be as in (\ref{Eqn:al_4}),
and $x\eq(C,D)\inn\fM^{\fk 1}_\al$ and its lifts $\ti x^{(\fk 2)}$, $\ti x^{(\fk 3)}$ and $\ti x^{(\fk 4)}$ be fixed as in the setting of Corollary~\ref{Crl:Step4_str_hom}.
We then fix arbitrary
\begin{alignat}{2}\label{Eqn:al_5'}
	\al^{\fk 5}&:=(\al^{\fk 4},\ov\lt^{\fk 4})
	=\big(\,\al^{\fk 4},\,
	(\tau_{\al^{\fk 4}},\ov\bE^{\fk 4},\Dm(\ov\lt^{\fk 4}))\,\big)&&\in A^{\fk 5}\,,
\end{alignat}  
as well as a sufficiently small chart $\cU^{(\fk 5)}\inn\fV_{\al^{\fk 5}}$  containing a fixed lift $\ti x^{(\fk 5)}$ of $\ti x^{(\fk 4)}$ in $\fM^{{\fk 5}}_{\al^{\fk 5}}$.
We denote by $\xi^{(\fk 5)}_s$, $s\inn S_{\al^{\fk 5}}$, the twisted parameters on $\cU^{(\fk 5)}$,
and by  $\ti\varphi^{(\fk 5)}$ the restriction to $\cU^{(\fk 5)}$ of the pullback of the structural homomorphism~$\varphi$.
For every $1\!\le\!i\!\le\!m$,
we write 
\begin{align*}
	\wp_{\al^{\fk 5};i}:=
	\wp_{\al;i}\big\bsl
	\big(\Dm(\ov\lt)\sqcup\Dm(\ov\lt^{\fk 2})\sqcup\Dm(\ov\lt^{\fk 3})\sqcup\Dm(\ov\lt^{\fk 4})\big)\,.
\end{align*}
For each $1\!<\!i\!\le\!m,$
we denote by $\wc\ka_{1i}^{(\fk 5)}$ the proper transform of $\ka_{1i}$ of (\ref{Eqn:ka_ij});
i.e.~on $\cU^{(\fk 5)}$, the proper transform of the conjugate locus $\wh\cK_{m;1i}$ is given by $\{\wc\ka_{1i}^{(\fk 5)}\eq 0\}$.

\begin{crl}
	\label{Crl:Step5_str_hom}
	Under suitable trivialization of (\ref{Eqn:str_hom_0}) and after re-ordering $\de_i$'s if necessary, we have
	\begin{align*}
		&	\ti\varphi^{(\fk 5)}=
		\left[
		\begin{matrix}\,
			\prod_{\fE\in\ov\bE\sqcup \ov\bE^{\fk 2}}\xi^{(\fk 5)}_\fE & 0\\
			0 & 
			\prod_{\fE\in\ov\bE\sqcup \ov\bE^{\fk 2}\sqcup\ov\bE^{\fk 3}\sqcup \ov\bE^{\fk 4}}\xi^{(\fk 5)}_\fE\,
		\end{matrix}
		\right]
		\cdot
		\left[
		\begin{matrix}\,
			1 & 0 & \cdots & 0\\
			0 & \wc\tw_{12}^{(\fk 5)} & \cdots &
			\wc\tw_{1m}^{(\fk 5)}\,
		\end{matrix}
		\right]\,,\\
		&	\tn{where}\qquad
		\wc\tw_{1i}^{(\fk 5)}=
		\wc\ka_{1i}^{(\fk 5)}\cdot
		\Big(\prod_{e\in\wp_{\al^{\fk 5};i}}
		\!\!\xi^{(\fk 5)}_{e}\ \Big)
		\cdot
		\Big(\prod_{\fE\in \ov\bE_{\al;1i}}
		\!\!\!\xi^{(\fk 5)}_{\fE}\ \Big)\,,
		\qquad
		1\!<\!i\!\le\!m.
	\end{align*}	
	Moreover, $\wc\ka_{1i}^{(\fk 5)}$ equals the pullback of $\wc\ka_{1i}^{(\fk 4)}$ up to a unit.
	Furthermore,
	there exists at least one $1\!<\!i\!\le\!m$ such that on $\cU^{(\fk 5)}$, $\prod_{e\in\wp_{\al^{\fk 5};i}}\!
	\xi^{(\fk 5)}_{e}$ is a unit.
\end{crl}

\begin{proof}
	Under the assumption of Corollary~\ref{Crl:Step4_str_hom}~\ref{Cond:Step4_nonep_bridge},
	notice the rooted tree $\tau_{\al^{\fk 4}}$ assigned by $\La^{\fk 4}$ is a subset of $\bigsqcup_{1<i\le m}\wp_{\al^{\fk 4};i}$,
	which is disjoint from $Q_{\al;1i}$,
	hence the pullback and the proper transform of $\{\wc\ka_{1i}^{(\fk 4)}\eq 0\}$ are identical.
	The same conclusion holds under the assumption of Corollary~\ref{Crl:Step4_str_hom}~\ref{Cond:Step4_other} trivially.
	
	The expression of $\ti\varphi^{(\fk 5)}$ is obtained by pulling back $\ti\varphi^{(\fk 4)}$ in Corollary~\ref{Crl:Step4_str_hom}.
	We only need to verify the case when $C_\core$ is separable and there exists some $1\!<\! i\!\le\!m$ such that 
	\begin{align}\label{Eqn:Step4_sep}
		\wp_{\al^{\fk 4};2,i}:=
		\wp_{\al;2,i}\big\bsl
		\big(\Dm(\ov\lt)\sqcup\Dm(\ov\lt^{\fk 2})\sqcup\Dm(\ov\lt^{\fk 3})\big)
		\ne \wp_{\al^{\fk 4};i}\,;
	\end{align}
	in this case,
	the factor $\prod_{e\in\wp_{\al^{\fk 5};i}}\!
	\xi^{(\fk 5)}_{e}$ should a priori be 
	$\prod_{e\in\wp_{\al^{\fk 4};2,i}\bsl\Dm(\ov\lt^{\fk 4})}\!
	\xi^{(\fk 5)}_{e}$.
	Nonetheless, by writing $\score(C)\eq\{T_1,T_2\}$,
	we see
	(\ref{Eqn:Step4_sep}) implies the existence of a node $e$ between $T_1$ and $T_2$ satisfying $e\!\notin\!
	\Dm(\ov\lt)\!\sqcup\!\Dm(\ov\lt^{\fk 2})\!\sqcup\!\Dm(\ov\lt^{\fk 3})$,
	so particularly, $\wp_{\al^{\fk 4};2,1}\!\ne\!\emptyset$. 
	Hence by Lemma~\ref{Lm:Step_4_lemma},
	there exists $j\!>\!1$ such that $\wp_{\al^{\fk 4};j}\!=\!\wp_{\al^{\fk 4};2,j}\eq\emptyset$.
	Consequently,
	we see $\de_1$ and $\de_j$ are on different connected components of $C\bsl\{e\}$,
	hence $\ka_{ij}$ is invertible and 
	$\ov\bE_{\al;1i}\eq\emptyset$.
	In sum,
	$\wc\tw^{(\fk 4)}_{1j}$ is invertible, and so is~$\wc\tw^{(\fk 5)}_{1j}$,
	which implies the proposed expressions of the remaining $\wc\tw^{(\fk 5)}_{1\ell}$'s can be achieved after suitable elementary column operations.
	
	The last statement of Corollary~\ref{Crl:Step5_str_hom} follows directly from the fact that every root-to-leaf path of $\tau_{\al^{\fk 4}}$ is equal to some $\wp_{\al^{\fk 4};i}$.
\end{proof}

Next,
in Step~$\fk 5$,
we aim to compare the factors $\big(\prod_{e\in\wp_{\al^{\fk 5};i}}
\xi^{(\fk 5)}_{e}\big)
\cdot
\big(\prod_{\fE\in \ov\bE_{\al;1i}}
\xi^{(\fk 5)}_{\fE}\big)$
by considering the minimal (w.r.t.~inclusion) elements of $\big\{\wp_{\al^{\fk 5};i}\!\sqcup\!\ov\bE_{\al;1i}: 1\!\le \!i\!\le\!m\big\}$.
By Corollary~\ref{Crl:Step2_str_hom}, $\wp_{\al;1}$ is a dominant root-to-leaf path,
i.e.~$\wp_{\al;1}\!\subset\!\Dm^*(\ov\lt)$.
Hence for each $i\!>\!1$,
$\wp_{\al^{\fk 5};i}\!\sqcup\!\ov\bE_{\al;1i}$
contains a root-to-leaf path of the proper transform of $\varrho_{\ov\lt}$ as a subset.
Conversely,
for each root-to-leaf path~$\wp$ of the proper transform of $\varrho_{\ov\lt}$,
there exists $i\!>\!1$ such that $\wp_{\al^{\fk 5};i}\!\sqcup\!\ov\bE_{\al;1i}\eq\wp$.
Therefore,
the proper transform of the derived treelike structure $\ov{\tn d}\La^{\fk 1}$ is exactly what we need.
More concretely,
we set
\begin{align*}
	\La^{\fk 5}:=\ov\PT_{\La^{\fk 4}}
	\big(\ov\PT_{\La^{\fk 3}}\big(\ov\PT_{\La^{\fk 2}}\big(\ov{\tn d}\La^{\fk 1}\big)\big)\big)\,,
\end{align*}
which is a treelike structure on $\fM^{\fk 5}$ by Corollary~\ref{Crl:G4-tf}.

Theorems~\ref{Thm:tf_smooth}~\ref{Cond:smooth_tf} and~\ref{Thm:tf_smooth_revert}~\ref{Cond:smooth_tf_revert}
and Propositions~\ref{Prp:two_treelike_strs} and~\ref{Prp:2nd_PT} together imply the following.

\begin{crl}
	\label{Crl:G5-tf} 
	Let 
	$$	
	\fM^{\fk 6}:=(\fM^{\fk 5})_{\La^{\fk 5}}^{\rtf}
	$$ be constructed as in Theorem~\ref{Thm:tf_smooth_revert}~\ref{Cond:smooth_tf_revert}.
	Then, $\fM^{\fk 6}$ is a  smooth algebraic stack endowed with an LES as well as a proper and birational  morphism (known as the forgetful morphism)   $\fM^{\fk 6}\!\lra\!\fM^{\fk 5}$,
	which is isomorphic on $(\fM^{\fk 5})^{\mn}$.
	The index set of the LES of $\fM^{\fk 6}$ is
	\begin{align*}
		A^{\fk 6}:=
		\ud{\La^{\fk 5}}=
		\big\{\,(\al,\ov\lt,\ov\lt^{\fk 2},\ov\lt^{\fk 3},\ov\lt^{\fk 4},\ud\ls):\ 
		(\al,
		\ov\lt,\ov\lt^{\fk 2},\ov\lt^{\fk 3},\ov\lt^{\fk 4})\inn A^{\fk 5},\ 
		\ud\ls\inn\ud\bT_{(\al,\ov\lt,\ov\lt^{\fk 2},\ov\lt^{\fk 3},\ov\lt^{\fk 4})}(\La^{\fk 5})
		\,\big\}\,.
	\end{align*}
	Moreover, 
	the LR proper transforms 
	\begin{align*}
		\wc\Up^{\fk 8}:=
		\ud\PT_{\La^{\fk 5}}\big(\ov\PT_{\La^{\fk 3}}\big(\ov\PT_{\La^{\fk 2}}(\Up^{\fk 8})\big)\big)\qquad
		\tn{and}\qquad
		\wc\Up^{\fk 9}:=
		\ud\PT_{\La^{\fk 5}}\big(
		\ov\PT_{\La^{\fk 3}}
		\big(\ov\PT_{\La^{\fk 2}}(\ov{\tn d}\Up^{\fk 4})\big)\big)
	\end{align*}
	as in Proposition~\ref{Prp:two_treelike_strs_revert}, as well as
	the secondary RL proper transforms
	\begin{align*}
		\wc\Up^{\fk 6}:=
		\ov\PT_{\La^{\fk 4}}\big(
		\ov\PT_{\La^{\fk 3}}\big(\ov\PT_{\La^{\fk 2}}({\tn d}_2\La^{\fk 1})\big)\big)
		\qquad
		\tn{and}\qquad
		\wc\Up^{\fk 7}:=
		\ov\PT_{\La^{\fk 4}}
		\big(\ov\PT_{\La^{\fk 3}}\big(\ov\PT_{\La^{\fk 2}}(\ud{\partial}\ov{\tn d}\La^{\fk 1})\big)\big)
	\end{align*}
	as in Proposition~\ref{Prp:2nd_PT},
	are respectively treelike structures on $\fM^{\fk 6}$.
\end{crl} 

Let  $\al^{\fk 5}\inn A^{\fk 5}$ be as in (\ref{Eqn:al_5}),
and $x\eq(C,D)\inn\fM^{\fk 1}_\al$ and its lifts $\ti x^{(\fk 2)}$, $\ti x^{(\fk 3)}$, $\ti x^{(\fk 4)}$ and $\ti x^{(\fk 5)}$ be fixed as in the setting of Corollary~\ref{Crl:Step5_str_hom}.
We then fix arbitrary
\begin{alignat}{2}\label{Eqn:al_6}
	\al^{\fk 6}&:=(\al^{\fk 5},\ud\ls)
	=\big(\,\al^{\fk 5},\,
	(\tau_{\al^{\fk 5}},\ud\bF,\Dm(\ud\ls))\,\big)&&\in A^{\fk 6}\,,
\end{alignat}  
as well as a sufficiently small chart $\cU^{(\fk 6)}\inn\fV_{\al^{\fk 6}}$  containing a fixed lift $\ti x^{(\fk 6)}$ of $\ti x^{(\fk 5)}$ in $\fM^{{\fk 6}}_{\al^{\fk 6}}$.
We denote by $\xi^{(\fk 6)}_s$, $s\inn S_{\al^{\fk 6}}$, the twisted parameters on $\cU^{(\fk 6)}$,
and by  $\ti\varphi^{(\fk 6)}$ the restriction to $\cU^{(\fk 6)}$ of the pullback of the structural homomorphism~$\varphi$.

For every $1\!\le\!i\!\le\!m$,
we write 
\begin{align*}
	\wp_{\al^{\fk 6};i}:=
	\wp_{\al;i}\big\bsl
	\big(\Dm(\ov\lt)\sqcup\Dm(\ov\lt^{\fk 2})\sqcup\Dm(\ov\lt^{\fk 3})\sqcup\Dm(\ov\lt^{\fk 4})\sqcup\Dm(\ud\ls)\big)\,.
\end{align*}
For every $1\!<\!i\!\le\!m$,
we denote by $\wc\ka_{1i}^{(\fk 6)}$ the proper transform of $\ka_{1i}$ of (\ref{Eqn:ka_ij});
i.e.~on $\cU^{(\fk 6)}$, the proper transform of the conjugate locus $\wh\cK_{m;1i}$ is given by $\{\wc\ka_{1i}^{(\fk 6)}\eq 0\}$.

\begin{crl}
	\label{Crl:Step6_str_hom}
	Under suitable trivialization of (\ref{Eqn:str_hom_0}) and after re-ordering $\de_i$'s if necessary, we have
	\begin{align*}
		&
		\ti\varphi^{(\fk 6)}=
		\left[
		\begin{matrix}\,
			b_{11} & 0\\
			0 & b_{22}
			\,
		\end{matrix}
		\right]
		\!\cdot\!
		\left[
		\begin{matrix}\,
			1 & 0 & 0 & \cdots & 0\\
			0 & \wc\ka_{12}^{(\fk 6)} & \wc\tw_{13}^{(\fk 6)} & \cdots &
			\wc\tw_{1m}^{(\fk 6)}\,
		\end{matrix}
		\right]\,,\qquad\tn{where}
		\\
		&
		b_{11}=\prod_{\fE\in(\ov\bE\bsl\Dm(\ud\ls))\sqcup \ov\bE^{\fk 2}\sqcup \ud\bF}\!\!\!\!\!\!\!\!\xi^{(\fk 6)}_\fE\ \ ,\qquad
		b_{22}=\Big(\prod_{\fE\in(\ov\bE\bsl\Dm(\ud\ls))\sqcup \ov\bE^{\fk 2}\sqcup\ov\bE^{\fk 3}\sqcup \ov\bE^{\fk 4}}
		\!\!\!\!\!\!\!\!\xi^{(\fk 6)}_\fE\ \Big)\cdot
		\Big(\prod_{\fF\in\ud\bF}
		\xi^{(\fk 6)}_\fF 
		\Big)^{\!2},
		\\
		&
		\wc\tw_{1i}^{(\fk 6)}=
		\wc\ka_{1i}^{(\fk 6)}
		\Big(\!\prod_{e\in\wp_{\al^{\fk 6};i}}
		\!\!\xi^{(\fk 6)}_{e}\Big)
		\Big(\!\prod_{\fE\in \ov\bE_{\al;1i}\bsl\Dm(\ud\ls)}
		\!\!\!\!\!\!\!\!\!\!\!\xi^{(\fk 6)}_{\fE}\Big)
		,
		\qquad
		2\!<\!i\!\le\!m.
	\end{align*}	
	Moreover,
	$\wc\ka_{12}^{(\fk 6)}$  is either a unit, 
	or a local parameter satisfying that
	$
	\{\wc\ka_{12}^{(\fk 6)}\}\!\sqcup\!
	\{\xi^{(\fk 6)}_s\}_{s\in S_{\al^{\fk 6}}}
	$
	forms a subset of a system of local parameters.	
\end{crl}

\begin{proof}
	By the choice of $\La^{\fk 5}$,
	there exists at least one $i\!>\!1$,
	say $i\eq 2$, so that 
	$\wp_{\al^{\fk 5};2}\!\sqcup\!\ov\bE_{\al;12}\!\subset\!\Dm(\ud\ls)$.
	In addition,
	Corollaries~\ref{Crl:Step4_str_hom} and~\ref{Crl:Step5_str_hom} imply $\wc\ka_{12}^{(\fk 6)}$ equals the pullback of $\wc\ka_{12}^{(\fk 4)}$ up to a unit.
	The expression of $\ti\varphi^{(\fk 6)}$ is thus obtained by pulling back $\ti\varphi^{(\fk 5)}$ of Corollary~\ref{Crl:Step5_str_hom}.
	
	To show the last statement of Corollary~\ref{Crl:Step6_str_hom},
	it suffices to assume $\de_1$ and $\de_2$ satisfy the assumption~\ref{Cond:Step4_nonep_bridge} of Corollary~\ref{Crl:Step4_str_hom},
	and in addition, 
	$\wp_{\al^{\fk 6};i}\!\ne\!\emptyset$ whenever $Q_{\al;i}\eq\emptyset$ or $T^{(i)}\!\ne\!T_2$.
	Under this assumption,
	after re-ordering $\de_i$'s if necessary,
	we can find at least one path $\wp'\inn P_{\al;T_2,1}$ such that $\wp'\!\subset\!\Dm(\ov\lt)\!\sqcup\!\Dm(\ov\lt^{\fk 2})\!\sqcup\!\Dm(\ov\lt^{\fk 3})$,
	hence there exists at least one path $\wp\inn  Q_{\al;12}$ such that $\wp\!\subset\!\wp'\!\subset\!\Dm(\ov\lt)\!\sqcup\!\Dm(\ov\lt^{\fk 2})\!\sqcup\!\Dm(\ov\lt^{\fk 3})$.
	By~(\ref{Eqn:Step4_kappa}),
	$\wc\ka_{1i}^{(\fk 4)}\eq \ti\ka_{1i}^{(\fk 4)}$ is either a unit, 
	or a local parameter satisfying that
	$
	\{\ti\ka_{1i}^{(\fk 4)}\}\!\sqcup\!
	\{\xi^{(\fk 4)}_s\}_{s\in S_{\al^{\fk 4}}}
	$
	forms a subset of a system of local parameters.	
	Along with the last statement of Theorem~\ref{Thm:tf_smooth}~\ref{Cond:smooth_parameters},
	this leads to the last statement of Corollary~\ref{Crl:Step6_str_hom}.
\end{proof}

\subsection{Steps $\fk 6$-$\fk 9$}\label{Subsec:Step6-9}

In this subsection,
we aim to compare the terms 
$\wc\tw_{1i}^{(\fk 6)}$, $i\!>\!2$,
in the matrix of $\ti\varphi^{(\fk 6)}$ of Corollary~\ref{Crl:Step6_str_hom},
utilizing the four treelike structures on $\fM^{\fk 6}$ introduced in Corollary~\ref{Crl:G5-tf}.
The term $\wc\tw_{12}^{(\fk 6)}$ is then handled by considering the grafted stratification and treelike structures with respect to the aforementioned four treelike structures.

First of all, we observe that by
Corollary~\ref{Crl:Step6_str_hom}, $\ti\varphi^{(\fk 6)}$ is already diagonalized if $\wc\ka_{12}^{(\fk 6)}$ is invertible;
so hereafter,
{\it we assume $\wc\ka_{12}^{(\fk 6)}$ is not invertible,}
i.e.~it is a local parameter,
which implies $\de_1$ and $\de_2$ are conjugate.
Particularly, we have
\begin{align*}
	\fE\cap\wp_{\al;2}=\emptyset\qquad
	\forall~\fE\inn\ov\bE^{\fk 2}\sqcup\ov\bE^{\fk 3}.
\end{align*}
Along with (\ref{Eqn:local_eqn_remainder1}), (\ref{Eqn:local_eqn_remainder2}),  (\ref{Eqn:local_eqn_remainder3}), and Corollary~\ref{Crl:Step4_str_hom}~\ref{Cond:Step4_nonep_bridge},
this implies for every $2\!<\!i\!\le\!m$,
\begin{align}\label{Eqn:Step6_ka_diff}
	\wc\ka_{1i}^{({\fk 6})} =
	u_{12i}^{({\fk 6})}\wc\ka_{12}^{({\fk 6})}+
	w_{12i}^{({\fk 6})}\!\!\!\!\!\!\prod_{e\in\wp_{\al;2}\cap\wp_{\al;i}}\!\!\!\!\!\!\!\ti\ze_e^{({\fk 6})}
	=\,u_{12i}^{({\fk 6})}\wc\ka_{12}^{({\fk 6})}+
	w_{12i}^{({\fk 6})}
	\Big(\!\prod_{\fF\in \ud\bF_{\al^{\fk 6};2i}}
	\!\!\xi^{(\fk 6)}_{\fF}\Big)
	\Big(\!\prod_{\fE\in (\ov\bE_{\al;2i}\bsl\Dm(\ud\ls))\sqcup \ov\bE^{\fk 4}_{\al;2i}}
	\!\!\!\!\!
	\!\xi^{(\fk 6)}_{\fE}\Big)
	,
\end{align}
where $u_{12i}^{({\fk 6})}$ and $v_{12i}^{({\fk 6})}$ are units on $\cU^{(\fk 6)}$,
$\ti\ze_e^{({\fk 6})}$ denotes the pullback of the node-smoothing parameter $\ze_e$, and
\begin{align*}
	\ud\bF_{\al^{\fk 6};2i}:=
	\big\{\,\fF\inn\ud \bF:\,
	\fF \!\cap\!\big(\ov\bE_{\al;2i}\!\sqcup\!(\wp_{\al^{\fk 5};2}\!\cap\!\wp_{\al^{\fk 5};i})\big)\!\ne\!\emptyset
	\,\big\},\quad
	\ov\bE^{\fk 4}_{\al;2i}:=
	\big\{\,
	\fE\inn\ov\bE^{\fk 4}:\,
	\fE\!\cap\!\wp_{\al;2}\!\cap\!\wp_{\al;i}\!\ne\!\emptyset
	\,\big\}.
\end{align*}
Therefore, in Corollary~\ref{Crl:Step6_str_hom},
we can apply elementary column operations to $\ti\varphi^{(\fk 6)}$ (i.e.~we take another trivialization of (\ref{Eqn:str_hom_0}))
and obtain
\begin{align}\label{Eqn:phi^6}
	\ti\varphi^{(\fk 6)}=
	\left[
	\begin{matrix}\,
		b_{11} & 0\\
		0 & b_{22}
		\,
	\end{matrix}
	\right]
	\!\cdot\!
	\left[
	\begin{matrix}\,
		1 & 0 & 0 & \cdots & 0\\
		0 & \wc\ka_{12}^{(\fk 6)} & \wh\tw_{13}^{(\fk 6)} & \cdots &
		\wh\tw_{1m}^{(\fk 6)}\,
	\end{matrix}
	\right]\,,
\end{align}
where $b_{11}$, $b_{22}$, and $\wc\ka_{12}^{(\fk 6)}$ are the same as in Corollary~\ref{Crl:Step6_str_hom},
and $\wh\tw_{1i}^{(\fk 6)}$'s are given by
\begin{align*}
	\wh\tw_{1i}^{(\fk 6)}=
	\Big(\!\prod_{e\in\wp_{\al^{\fk 6};i}}
	\!\!\!\xi^{(\fk 6)}_{e}\Big)
	\Big(\!\prod_{\fE\in \ov\bE_{\al;1i}\bsl\Dm(\ud\ls)}
	\!\!\!\!\!\!\!\!\!\!\!\xi^{(\fk 6)}_{\fE}\Big)
	\Big(\!\prod_{\fE\in \ov\bE_{\al;2i}\bsl\Dm(\ud\ls)}
	\!\!\!\!\!\!\!\!\!\!\!\xi^{(\fk 6)}_{\fE}\Big)
	\Big(\!\prod_{\fE\in \ov\bE^{\fk 4}_{\al;2i}}
	\!\!\xi^{(\fk 6)}_{\fE}\ \,\Big)
	\Big(\!\prod_{\fF\in \ud\bF_{\al^{\fk 6};2i}}
	\!\!\xi^{(\fk 6)}_{\fF}\Big)
	\,,
	\qquad
	2\!<\!i\!\le\!m.
\end{align*}
Notice that $\ov\bE_{\al;1i}\!\cap\!\ov\bE_{\al;2i}\!\subset\!\Dm(\ud\ls)$,
so
\begin{align}\label{Eqn:Step6_wh_tw}
	\wh\tw_{1i}^{(\fk 6)}=
	\Big(\!\prod_{e\in\wp_{\al^{\fk 6};i}}
	\!\!\!\xi^{(\fk 6)}_{e}\Big)
	\Big(\!\prod_{\fE\in (\ov\bE_{\al;1i}\cup\ov\bE_{\al;2i})\bsl\Dm(\ud\ls)}
	\!\!\!\!\!\!\!\!\!\!\!\!\!\!\xi^{(\fk 6)}_{\fE}\ \ \Big)
	\Big(\!\prod_{\fE\in \ov\bE^{\fk 4}_{\al;2i}}
	\!\!\!\xi^{(\fk 6)}_{\fE}\,\Big)
	\Big(\!\prod_{\fF\in \ud\bF_{\al^{\fk 6};2i}}
	\!\!\!\!\xi^{(\fk 6)}_{\fF}\,\Big)
	\,,
	\qquad
	2\!<\!i\!\le\!m.
\end{align}

Recall the treelike structures $\wc\Up^{\fk h}$, $\fk 6\!\le\!\fk h\!\le\!\fk 9$, on $\fM^{\fk 6}$ introduced in Corollary~\ref{Crl:G5-tf}.
Comparing with their prototypes
$\tn d_2\La^{\fk 1}$ in Proposition~\ref{Prp:Derived_treelike_2} and
$\ud{\partial}\ov{\tn d}\La^{\fk 1}$ in Proposition~\ref{Prp:Derived_treelike_second_order_1}, as well as 
$\Up^{\fk 8}$ and $\ov{\tn d}\Up^{\fk 4}$ in \S\ref{Subsec:Step3},
we observe
\begin{itemize}[leftmargin=*]
	\item 
	the minimal elements (w.r.t.~inclusion) of $\big\{\wp_{\al^{\fk 6};i}\!\sqcup\!\big((\ov\bE_{\al;1i}\!\cup\!\ov\bE_{\al;2i})\bsl\Dm(\ud\ls)\big):1\!\le\!i\!\le \!m\big\}$  are exactly the root-to-leaf paths of the rooted tree assigned by $\wc\Up^{\fk 6}$;
	
	\item 
	the minimal elements (w.r.t.~inclusion) of $\{\wp_{\al^{\fk 6};i}\!\sqcup\!\ud\bF_{\al^{\fk 6};2i}:1\!\le\!i\!\le \!m\}$ are exactly the root-to-leaf paths of the rooted tree assigned by $\wc\Up^{\fk 7}$;
	
\item 
the minimal elements (w.r.t.~inclusion) of $\{\wp_{\al^{\fk 6};i}:2\!<\!i\!\le \!m\}$ are exactly the root-to-leaf paths of the rooted tree assigned by $\wc\Up^{\fk 8}$;

\item 
the minimal elements (w.r.t.~inclusion) of $\{\wp_{\al^{\fk 6};i}\!\sqcup\!\ov\bE^{\fk 4}_{\al;2i}:1\!<\!i\!\le \!m\}$ are exactly the root-to-leaf paths of the rooted tree assigned by $\wc\Up^{\fk 9}$.
\end{itemize}
The above observation gives rise to the following.

\begin{lmm}
	\label{Lm:Step6-9}
Let $(\fM')^{\fk{7}}/\fM^{\fk 6}$,
$(\fM')^{\fk{8}}/(\fM')^{\fk 7}$,
$(\fM')^{\fk{9}}/(\fM')^{\fk 8}$,
and  $(\fM')^{\fk{10}}/(\fM')^{\fk 9}$ be respectively obtained by successively adding the RL twisted fields as per Theorem~\ref{Thm:tf_smooth}  with respect to the proper transforms of 
$\wc\Up^{\fk 6}$, $\wc\Up^{\fk 7}$, $\wc\Up^{\fk 8}$, and $\wc\Up^{\fk 9}$.
Then, locally on $\fM^{\fk{10}}$,
the set of the pullbacks of the regular functions $\wh\tw_{1i}^{(\fk 6)}$, $2\!<\!i\!\le\!m$,
always has a divisibly minimal element in the sense of Definition~\ref{Dfn:Tauto_monom}.
\end{lmm}

\begin{proof}
	The statement follows directly from (\ref{Eqn:Step6_wh_tw}) and Corollary~\ref{Crl:prod_RL}.
	Since the calculation is largely parallel to those in the preceding subsections,
	we omit further details.
\end{proof}

\begin{rmk}
	The above argument suggests that
one can start with $\fM^{\fk 6}$ and successively add RL twisted fields with respect to the proper transforms of $\wc\Up^{\fk h}$, $\fk 6\!\le\!\fk h\!\le\!\fk 9$, in an {\it arbitrary order}.
The pullback of $\{\wh\tw_{1i}^{(\fk 6)}:2\!<\!i\!\le\!m\}$ always has a divisibly minimal element.
\end{rmk}

It remains to handle the term $\wc\ka_{12}^{(\fk 6)}$ of (\ref{Eqn:phi^6}),
using the technique of \S\ref{Subsec:Grafted}.
For every $\al^{\fk 6}\inn A^{\fk 6}$ as in (\ref{Eqn:al_6}),
let
\begin{align*}
	I^{\cj}_{\al^\fk 6}\!:=
	\Big\{
	(i,j)\,:\ 
	1\!\le\!i,j\!\le\!m,\ \ 
	i\!\ne\!j,\ \ \,
	&\wp_{\al,i}\subset
	\Dm(\ov\lt)\!\sqcup\!\Dm(\ov\lt^{\fk 2}),\\
	& 
	\ov\bE_{\al;ij}
	\!\sqcup\! \wp_{\al,j}
	\subset 
	\Dm(\ov\lt)\!\sqcup\!\!\bigsqcup_{\fk 2\le i\le\fk 4}\!\!\!\Dm(\ov\lt^{i})\!\sqcup\!
	\Dm(\ud\ls) 
	\Big\}.
\end{align*}
The choices of $\La^{\fk 4}$ and $\La^{\fk 5}$ imply $I^{\cj}_{\al^\fk 6}\!\ne\!\emptyset$ for all $\al^{\fk 6}\inn B^{\fk 6}$.
Indeed, the way we choose $\de_1$ and $\de_2$ in Corollary~\ref{Crl:Step6_str_hom} implies $(1,2)\inn I^{\cj}_{\al^\fk 6}$.

\begin{lmm}
	\label{Lm:I^K}
With notation as above,
for every $S'\!\subset\! S_{\al^{\fk 6}}$,
we have $I^{\cj}_{\al^{\fk 6}}\!\subset\!I^{\cj}_{\al^{\fk 6}_{(S')}}$.
\end{lmm}

\begin{proof}
	The statement follows from a direct check, utilizing Theorem~\ref{Thm:tf_smooth}~\ref{Cond:smooth_Z} and Theorem~\ref{Thm:tf_smooth_revert}~\ref{Cond:smooth_Z_revert}.
	We omit further details.
\end{proof}

Let $\cV_x\inn\cV_\al$ be a small chart containing $x\eq(C,D)\inn\fM^{\fk 1}_\al$ that is fixed Corollary~\ref{Crl:Step2_str_hom}.
Recall the locus $\{\ka_{ij}\eq 0\}$ is written as $\wh K_{m;i,j}\!\cap\!\cV_x$ in~\S\ref{Subsec:M2Pnd_loc_eqn}.
We denote by $\wc K_{m;i,j}\!\cap\!\cU^{({\fk 6})}$
the proper transform of $\wh K_{m;i,j}\!\cap\!\cV_x$,
i.e.~$\wc K_{m;i,j}\!\cap\!\cU^{({\fk 6})}\eq\{\wc\ka^{({\fk 6})}_{ij}\eq 0\}$.

\begin{lmm}
	\label{Lm:K_step6}
	Let $\al^{\fk 6}$ 
	and $\cU^{({\fk 6})}$ be as in Corollary~\ref{Crl:Step6_str_hom}.
	Then, one of the following holds.
	\begin{itemize}[leftmargin=*]
		\item Either $\bigcap_{(i,j)\in I^\cj_{\al^{\fk 6}}}
		\wc K_{m;i,j}
		\!\cap\!\fM^{\fk 6}_{\al^{\fk 6}}
		\!\cap\!\cU^{({\fk 6})}
		\eq \emptyset$,
		
		\item or $\bigcap_{(i,j)\in I^\cj_{\al^{\fk 6}}}
		\wc K_{m;i,j}
		\!\cap\!\fM^{\fk 6}_{\al^{\fk 6}}
		\!\cap\!\cU^{({\fk 6})}
		\eq 
		\wc K_{m;h,\ell}\!\cap  \!
		\fM^{\fk 6}_{\al^{\fk 6}}
		\!\cap\!\cU^{({\fk 6})}
		$
		for any $(h,\ell)\inn I^\cj_{\al^\fk 6}$.
	\end{itemize}
\end{lmm}

\begin{proof}
	If $\wc\ka_{ij}^{({\fk 6})}$ is invertible for some $(i,j)\inn I^\cj_{\al^{\fk 6}}$,
	then obviously $\bigcap_{(i,j)\in I^\cj_{\al^{\fk 6}}}
	\wc K_{m;i,j}
	\!\cap\!\fM^{\fk 6}_{\al^{\fk 6}}
	\!\cap\!\cU^{({\fk 6})}
	\eq \emptyset$.

    If for every $(i,j)\in I^\cj_{\al^{\fk 6}}$, $\wc\ka_{ij}^{({\fk 6})}$ vanishes along  $\fM^{\fk 6}_{\al^{\fk 6}}
    \!\cap\!\cU^{({\fk 6})}$,
    then for every $(h,\ell)\inn  I^\cj_{\al^{\fk 6}}$,
    by applying (\ref{Eqn:Step6_ka_diff}) to $(1,\ell)$ then to $(h,\ell)$,
    we see $\wp_{\al;2}\!\cap\!\wp_{\al;\ell}\!\ne\!\emptyset$ and $\wp_{\al;1}\!\cap\!\wp_{\al;h}\!\ne\!\emptyset$, 
    hence $\wc K_{m;h,\ell}\!\cap  \!
    \fM^{\fk 6}_{\al^{\fk 6}}
    \!\cap\!\cU^{({\fk 6})}
    \eq \wc K_{m;1,2}\!\cap  \!
    \fM^{\fk 6}_{\al^{\fk 6}}
    \!\cap\!\cU^{({\fk 6})}$.
\end{proof}

\begin{crl}
	\label{Crl:K_step6}
	There exists a (locally closed) substack $K\!\subset\! \fM^{\fk 6}$ such that with $\al^{\fk 6}$ and
	and $\cU^{({\fk 6})}$ as in Corollary~\ref{Crl:Step6_str_hom}, we have
	\begin{align*}
		K\cap \cU^{({\fk 6})}=
		\bigcup_{(\al^{\fk 6})'\in A^{\fk 6}}
		\big(
		\bigcap_{(i,j)\in I^\cj_{(\al^{\fk 6})'}}\!\!\!\!
		\wc K_{m;i,j}
		\cap\fM^{\fk 6}_{(\al^{\fk 6})'}
		\cap\cU^{({\fk 6})}
		\big)\,.
	\end{align*}
	Moreover, such $K$ satisfies~(\ref{Eqn:K_local}).
\end{crl}

\begin{proof}
	The statement follows directly from Lemmas~\ref{Lm:K_step6} and~\ref{Lm:I^K}
	and Proposition~\ref{Prp:M_strata_local}.
\end{proof}

By Corollary~\ref{Crl:K_step6},
we obtain the grafted stratification of $\fM^{\fk 6}$ with respect to $K$ as per Lemma~\ref{Lm:grafted_strat}.
By Proposition~\ref{Prp:grafted_treelike},
each treelike structure $\wc\Up^{\fk h}$,
$\fk 6\!\le\!\fk h\!\le\!\fk 9$,
then determines its  grafted treelike structure $\wc\Up^{\fk h}_{\ex}$ respect to  $K$ accordingly.
On $\cU^{(\fk 6)}$,
notice
the grafted edge exactly corresponds to $\wc\ka_{12}^{({\fk 6})}$ because $(1,2)\inn I^\cj_{\al^{\fk 6}}$,
so we construct 
$\fM^{\fk{7}}/\fM^{\fk 6},$ $\fM^{\fk{8}}/\fM^{\fk 7},$ $\fM^{\fk{9}}/\fM^{\fk 8},$ and $\fM^{\fk{10}}/\fM^{\fk 9}$
by successively adding the RL twisted fields  with respect to the proper transforms of 
$\wc\Up^{\fk 6}_{\ex}$, $\wc\Up^{\fk 7}_{\ex}$, $\wc\Up^{\fk 8}_{\ex}$, and $\wc\Up^{\fk 9}_{\ex}$.

\begin{prp}
	\label{Prp:Step10_str_hom}
$\fM^{\fk{10}}$ is a smooth algebraic stack,
and the forgetful morphism 
$\fM^{\fk{10}}\!\lra\!\fM^{\fk 6}$ is proper and birational.
Moreover,
locally on $\fM^{\fk{10}}$,
under suitable trivialization of (\ref{Eqn:str_hom_0}) and after re-ordering $\de_i$'s if necessary, the pullback of the structural homomorphism $\varphi$ takes the form
\begin{align*}
	\left[
	\begin{matrix}\,
		z_{11} & 0\\
		0 & z_{22}
		\,
	\end{matrix}
	\right]
	\!\cdot\!
	\left[
	\begin{matrix}\,
		1 & 0 & 0 & \cdots & 0\\
		0 & 1 & 0 & \cdots &
		0\,
	\end{matrix}
	\right]\,,
\end{align*}	
where $z_{11}$ and $z_{22}$ (not necessarily square-free) monomials whose variables are the twisted parameters on $\fM^{\fk{10}}$,
satisfying $z_{11}\mid z_{22}$.
\end{prp}

\begin{proof}
	The matrix of the pullback of $\varphi$ on $\fM^{\fk{10}}$ is obtained from Corollary~\ref{Crl:Step6_str_hom},
	Lemma~\ref{Lm:Step6-9}, and Corollary~\ref{Crl:K_step6}.
	Particularly,
	$z_{11}$ is the pullback of $b_{11}$ of Corollary~\ref{Crl:Step6_str_hom},
	whereas $z_{22}$ is the product of the pullback of $b_{22}$ of Corollary~\ref{Crl:Step6_str_hom} and certain twisted parameters obtained in Steps~$\fk 6$-$\fk{10}$.
	Since $b_{11}\mid b_{22}$,
	we obtain $z_{11}\mid z_{22}$.
\end{proof}

\subsection{Proof of Theorem~\ref{Thm:Main}}\label{Subsec:Proof_Main}

We take the proposed $\ti\fD_2^{\tn{tf}}\!\lra\!\fD_2$ of Theorem~\ref{Thm:Main} to be the composite of $\fM^{\fk h+\fk 1}\!\lra\!\fM^{\fk h}$, $\fk 1\!\le\!\fk h\!\le\!\fk 9$.
The properties of $\ti\fD_2^{\tn{tf}}/\fD_2$ in Theorem~\ref{Thm:Main} follow from
Corollaries~\ref{Crl:G1-tf}, \ref{Crl:G2-tf}, \ref{Crl:G3-tf}, \ref{Crl:G4-tf}, and \ref{Crl:G5-tf}, as well as Proposition~\ref{Prp:Step10_str_hom}.

To construct $\ti M_2^{\tn{tf}}(\P^n,d)$,
we first consider the open covering $\fU$ of $\ov M_2(\P^n,d)$ 
given by
\begin{align*}
	\fU=
	\big\{\,
	\tn{open}~U\!\subset\!\ov M_2(\P^n,d):\,
	\bH_U\!\ne\!\emptyset\,
	\big\}\,,
\end{align*}
where $\bH_U$ is as in (\ref{Eqn:H}).
In \S\ref{Subsec:Step1},
$\fM^{\fk 1}\eq\fD_2$ is endowed with an LES.
The modular parameters are the node-smoothing parameters,
whose pullbacks via (\ref{Eqn:Mg_tauto}) satisfy (\ref{Eqn:f_U_assump}). 
It is thus a direct check that for each $U\inn\fU$,
(\ref{Eqn:Mg_tauto}) gives rise to a strongly admissible family of morphisms $U\!\lra\!\ov M_2(\P^n,d)$ as per Definition~\ref{Dfn:weak_and_strong_adm}.
In addition, these families satisfy (\ref{Eqn:f_U_assump2}).
So, by Corollary~\ref{Crl:local_tf_isom},
we obtain $\ov M^{\fk 2}_2(\P^n,d)/\ov M_2(\P^n,d)$
that is locally given by 
$$
U^{\fk 2}_H:=U\times_{f_{U,H};\fM^{\fk 1}}\fM^{\fk 2}\lra U,\qquad U\inn\fU,~H\inn\bH_U.
$$
Moreover,
by Corollary~\ref{Crl:f_U_induction},
the natural morphisms $U^{\fk 2}_H\!\lra\!\fM^{\fk 2}$ determines a strongly admissible family of morphisms.
So, by repeating the above argument, we apply Corollaries~\ref{Crl:local_tf_isom} and~\ref{Crl:f_U_induction} recursively to obtain $\ov M^{r+\fk 1}_2(\P^n,d)/\ov M^{r}_2(\P^n,d)$,
$\fk 2\!\le\!r\!\le\!\fk 5$.
The strongly admissible family of morphisms obtained in each step is denoted by $(U^{r}_H\!\lra\!\fM^{r})$, $\fk 2\!\le\!r\!\le\!\fk 6$.

Starting from $r\!=\! \fk 6$,
there is an extra modular parameter on $\fM^r$ given by $\wc\ka^{(\fk 6)}_{ij}$,
the proper transform of $\ka_{ij}$; c.f.~(\ref{Eqn:ka_ij}).
For every point on the boundary $\De^{\fk 6}_{\ex}$ (c.f.~(\ref{Eqn:boundary})) of the grafted stratification of $\fM^{\fk 6}$ (with respect to the substack $K$ of $\fM^{\fk 6}$ as in Corollary~\ref{Crl:K_step6}),
we write its underlying point on $\fM^{\fk 1}$ as $(C,D\eq\de_1\!+\!\cdots\!+\!\de_m)$,
and let $\al\inn B^{\fk 1}$ be such that $(C,D)\inn\fM^{\fk 1}_\al$.
Then, by Lemma~\ref{Lm:K_step6}, there are
{\it at most two} points among $\de_1,\ldots,\de_m$ that lie on the core of $C$,
for otherwise the right-hand side of the identity of Corollary~\ref{Crl:K_step6} would be empty, which would contradict (\ref{Eqn:grafted_boundary}).
Given $x\eq (C,\bu)\inn U\inn\fU$ satisfying $f_{U,H}(x)\inn\fM^{\fk 1}_\al$ for some (and thus all) $H\inn\bH_U$,
notice that for every $H_1,H_2\inn\bH_U$ and irreducible component $C'$ of $C$,
the divisors $\bu^{-1}(H_1)\!\cap\!C'$ and $\bu^{-1}(H_2)\!\cap\!C'$ are linearly equivalent.
Therefore,
\begin{align*}
	\ka_{ij}\big(f_{U,H_1}(x)\big)= 0
	\qquad\Longleftrightarrow\qquad
	\ka_{ij}\big(f_{U,H_2}(x)\big)= 0\,.
\end{align*}
We conclude that  with respect to the grafted LES of $\fM^{\fk 6}$ and the treelike structure $\Up^{\fk 6}_{\ex}$,
the family $(U^{\fk 6}_H\!\lra\!\fM^{\fk 6})$ is admissible,
so Corollary~\ref{Crl:local_tf_isom} can be applied to obtain $\ov M^{\fk 7}_2(\P^n,d)/\ov M^{\fk 6}_2(\P^n,d)$.
The resulting family of morphisms is denoted by 
$(U^{\fk 7}_H\!\lra\!\fM^{\fk 7})$.
Although $(U^{\fk 7}_H\!\lra\!\fM^{\fk 7})$ is not strongly admissible,
(\ref{Eqn:f_U_assump}) is still valid for all the modular parameters of $\fM^{\fk 7}$ except $\wc\ka_{ij}^{(\fk 7)}$,
and the above argument on $\wc\ka_{ij}^{(\fk 6)}$ applies to $\wc\ka_{ij}^{(\fk 7)}$ likewise.
Therefore,
we mimic the above construction of $\ov M^{\fk 7}_2(\P^n,d)/\ov M^{\fk 6}_2(\P^n,d)$ and obtain $\ov M^{r+\fk 1}_2(\P^n,d)/\ov M^{r}_2(\P^n,d)$, $\fk 7\!\le\!r\!\le\!\fk 9$,
analogously.
Particularly, we set $\ti M_2^{\tn{tf}}(\P^n,d)\!:=\!\ov M^{\fk{10}}_2(\P^n,d)$.

The existence of the isomorphisms 
\begin{align*}
	\big(U\times_{f_{U,H};\fD_2} \ti{\fD}^{\tn{tf}}_2\big)\big/ U  \lra \big(U\times_{\ov M_2(\P^n,d)}\ti M_2^{\tn{tf}}(\P^n,d)\big)\big/ U
\end{align*}
follows from Corollary~\ref{Crl:local_tf_isom}.

Since the properness is preserved under base changes, 
the properness of $\ti\fD_2^{\tn{tf}}/\fD_2$, inherited from that of each $\fM^{\fk h+\fk 1}/\fM^{\fk h}$, 
along with the properness of $\ov M_2(\P^n,d)$,
implies the properness of $\ti M_2^{\tn{tf}}(\P^n,d)\!\lra\!\ov M_2(\P^n,d)$.
Moreover,
$\ti\fD_2^{\tn{tf}}\!\lra\!\fD_2$ restricts to the identity morphism outside the  preimage of the boundary of $\fD_2$, which is closed in $\ti\fD_2^{\tn{tf}}$.
This establishes Theorem~\ref{Thm:Main}~\ref{Cond:MainBirational}.

Theorem~\ref{Thm:Main}~\ref{Cond:MainSmooth} and~\ref{Cond:MainMainComponent} follow from the $m\eq d$ case of (\ref{Eqn:str_hom_0}) and Proposition~\ref{Prp:Step10_str_hom}.
Similarly, Theorem~\ref{Thm:Main}~\ref{Cond:MainLocallyFree}  follows from the $m\eq kd$ case of (\ref{Eqn:str_hom_0}), Proposition~\ref{Prp:Step10_str_hom},
and an argument parallel to the proof of \cite[Proposition~5.5 \& Corollary 5.6]{HLN}.
The proof of Theorem~\ref{Thm:Main} is  complete.

\section{Appendix: auxiliary examples}\label{Sec:eg}
Included in the Appendix are some auxiliary examples related to the operations on the treelike structures in \S\ref{Sec:Induced_and_derived_tf_stacks}.

The first example below shows the order matters when adding (the proper transforms of) two treelike structures successively to a stack.
The example is presented under the RL setting; similar conclusion applies to the LR case as well.

\begin{eg}\label{Eg:PT_simple}
	Consider $\fM\eq\A^2$. 
	Let $\{\ze_b,\ze_c\}$ be a system of local parameters on $\fM$ and $A$ be the power set of $\{b,c\}$.
	For each $\al\inn A$, we set
	\begin{align*}
		S_\al:=\al,\qquad
		\fM_{\al}:=
		\big\{\,\ze_e\eq 0\ \forall\ e\inn\al\,,\ 
		\ze_e\!\ne\!0\ \forall\ e\inn\{b,c\}\bsl\al\,\big\}.
	\end{align*}
	In this way, $\fM$ is endowed with an LES $\fM\eq\bigsqcup_{\al'\in A}\fM_{\al'}$.
	The open stratum $\fM^\mn\eq\fM_\emptyset$,
	whereas the boundary $\De\eq\fM_{\{b,c\}}\!\sqcup\!\fM_{\{b\}}\!\sqcup\!\fM_{\{c\}}$.
	(In other words, the index $\emptyset\inn A$ plays the role of the index $0\inn A$ of Definition~\ref{Dfn:G-adim_fixture}, and $B\eq A\bsl\{\emptyset\}$.)
	
	\begin{figure}[htp]
		\begin{center}
			\begin{tikzpicture}
				\filldraw
				(0,0) circle (1pt)
				(.35,.66) circle (1pt)
				(.7,0) circle (1pt);
				\draw
				(0,0)--(.35,.66)--(.7,0);
				\draw
				(.35,.66) node[above] {\tiny{$o$}}
				(.175,.33) node[left] {\tiny{$b$}}
				(.525,.33) node[right] {\tiny{$c$}}
				;
				
				\filldraw[xshift=3cm]
				(.35,.33) circle (1pt);
				\draw[xshift=3cm]
				(.35,.33) node[above] {\tiny{$o$}}
				;
				
				\filldraw[xshift=7cm]
				(0,0) circle (1pt)
				(0,.33) circle (1pt)
				(0,.66) circle (1pt);
				\draw[xshift=7cm]
				(0,.66)--(0,0);
				\draw[xshift=7cm]
				(0,.165) node[right] {\tiny{$c$}}
				(0,.495) node[right] {\tiny{$b$}}
				(0,.66) node[above] {\tiny{$o$}}
				;
				
				\filldraw[xshift=8.5cm]
				(0,0) circle (1pt)
				(0,.66) circle (1pt);
				\draw[xshift=8.5cm]
				(0,.66)--(0,0);
				\draw[xshift=8.5cm]
				(0,.33) node[right] {\tiny{$c$}}
				(0,.66) node[above] {\tiny{$o$}}
				;
				
				\filldraw[xshift=10cm]
				(0,0) circle (1pt)
				(0,.66) circle (1pt);
				\draw[xshift=10cm]
				(0,.66)--(0,0);
				\draw[xshift=10cm]
				(0,.33) node[right] {\tiny{$b$}}
				(0,.66) node[above] {\tiny{$o$}}
				; 
				
				\filldraw[xshift=11.5cm]
				(0,.33) circle (1pt);
				\draw[xshift=11.5cm]
				(0,.33) node[above] {\tiny{$o$}}
				;
				
				\draw[dashed]
				(-.5,-1) rectangle (5,1.2);
				
				\draw[dashed]
				(6.4,-1) rectangle (12.2,1.2);
				
				\draw[yshift=-.5cm]
				(.35,0) node {\tiny{$\tau_{\{b,c\}}$}}
				(3.35,0) node {\tiny{$\tau_{\{c\}}\eq \tau_{\{b\}}\eq \tau_{\emptyset}$}}
				(7,0) node {\tiny{$\tau^\dag_{\{b,c\}}$}}
				(8.5,0) node {\tiny{$\tau^\dag_{\{c\}}$}}
				(10,0) node {\tiny{$\tau^\dag_{\{b\}}$}}
				(11.5,0) node {\tiny{$\tau^\dag_{\emptyset}$}}
				(2.25,-1) node {\tiny{$\La$}}
				(9.3,-1) node {\tiny{$\La^\dag$}};
			\end{tikzpicture}	
		\end{center}
		\caption{The treelike structures in Example~\ref{Eg:PT_simple}}\label{Fig:PT_simple}	
	\end{figure}
	
	Consider the rooted trees $\tau_\al$ and $\tau_\al^\dag$, $\al\inn A$, given in Figure~\ref{Fig:PT_simple}.
	To each $\tau_{\al}$ (resp.~$\tau^\dag_{\al}$), $\al\inn A$, we assign the obvious inclusion $\be_{\al}\!: \tau_{\al}\!\hookrightarrow\!S_{\al}$  (resp.~$\be^\dag_{\al}\!: \tau^\dag_{\al}\!\hookrightarrow\!S_{\al}$).
	It is straightforward that
	$$
	\La=(\tau_\al,\be_\al)_{\al\in A}
	\qquad\tn{and}\qquad
	\La^\dag=(\tau^\dag_\al,\be^\dag_\al)_{\al\in A}
	$$
	are both treelike structures on $\fM$.
	For $\La$, we observe that
	\begin{align*}
		&\ov\fT_{\{b,c\}}=\big\{\,\ov\lt_{1}\,,\,\ov\lt_{2}\,,\,\ov\lt_{3}\,\big\},\qquad
		\tn{where}\qquad
		\ov\lt_{1}\!:=\!\big(\tau_{\{b,c\}},\{\{b,c\}\},\{b\}\big),
		\\
		&\hspace{1.75in}
		\ov\lt_{2}\!:=\!\big(\tau_{\{b,c\}},\{\{b,c\}\},\{b,c\}\big),\qquad\tn{and}\qquad
		\ov\lt_{3}\eq\big(\tau_{\{b,c\}},\{\{b,c\}\},\{c\}\big),\\
		&
		\ov\fT_{\{c\}}=\ov\fT_{\{b\}}=\ov\fT_{\emptyset}=\{\ov\lt_\bullet\};
	\end{align*}
	c.f.~Figure~\ref{Fig:Treelik_Level} regarding $\ov\fT_{\{b,c\}}$.
	For $\La^\dag$, we similarly have
	\begin{align*}
		&\ov\fT^\dag_{\{b,c\}}=\big\{\,\ov\lt_{4}\!:=\!\big(\tau_{\{b,c\}}^\dag,\{\{b\},\{c\}\},\{b,c\}\big)\,\big\},
		\\
		&
		\ov\fT_{\{c\}}^\dag=\big\{\,\ov\lt_{5}\!:=\!\big(\tau_{\{c\}}^\dag,\{\{c\}\},\{c\}\big)\,\big\},\qquad \ov\fT_{\{b\}}^\dag=\big\{\,\ov\lt_{6}\!:=\!\big(\tau_{\{b\}}^\dag,\{\{b\}\},\{b\}\big)\,\big\},\qquad
		\ov\fT_{\emptyset}^\dag=\{\ov\lt_\bullet\}.
	\end{align*}
	
	By~(\ref{Eqn:beta_ti_al}), 
	The rooted trees of the proper transform $\ov\PT_{\La}(\La^\dag)$ are as follows:
	\begin{align*}
		\tau^\dag_{(\{b,c\},\ov\lt_1)}\eq\{c\},\quad
		\tau^\dag_{(\{b,c\},\ov\lt_2)}\eq\tau_\bullet,\quad
		\tau^\dag_{(\{b,c\},\ov\lt_3)}\eq\{b\},\quad
		\tau^\dag_{(\{c\},\ov\lt_\bullet)}\eq\{c\},\quad
		\tau^\dag_{(\{b\},\ov\lt_\bullet)}\eq\{b\},\quad
		\tau^\dag_{(\emptyset,\ov\lt_\bullet)}\eq\tau_\bullet,
	\end{align*}
	along with the obvious tree orders,
	whereas those of the proper transform $\ov\PT_{\La^\dag}(\La)$ are given by:
	\begin{align*}
		\tau^\dag_{(\{b,c\},\ov\lt_4)}=
		\tau^\dag_{(\{c\},\ov\lt_5)}=
		\tau^\dag_{(\{b\},\ov\lt_6)}=
		\tau^\dag_{(\emptyset,\ov\lt_\bullet)}=\tau_\bullet.
	\end{align*}
	From Theorem~\ref{Thm:tf_smooth}~\ref{Cond:smooth_tf}
	we conclude that $\big(\fM^{\tf}_\La\big)^{\tf}_{\ov\PT_{\La}(\La^\dag)}\big/\fM^{\tf}_\La$ is isomorphic to the identity morphism $\fM^{\tf}_\La\big/\fM^{\tf}_\La$,
	while $\fM^{\tf}_\La\big/\fM$ is not isomorphic to the identity morphism $\fM/\fM$,
	so the composite morphism $\big(\fM^{\tf}_\La\big)^{\tf}_{\ov\PT_{\La}(\La^\dag)}\big/\fM$ is not isomorphic to the identity morphism $\fM/\fM$.
	However,
	following the same argument we can see $\big(\fM^{\tf}_{\La^\dag}\big)^{\tf}_{\ov\PT_{\La^\dag}(\La)}\big/\fM$ is isomorphic to the identity morphism $\fM/\fM$.
	Therefore, $\big(\fM^{\tf}_\La\big)^{\tf}_{\ov\PT_{\La}(\La^\dag)}\big/\fM$ and $\big(\fM^{\tf}_{\La^\dag}\big)^{\tf}_{\ov\PT_{\La^\dag}(\La)}\big/\fM$ are not isomorphic.
	
	In fact, as shown in Theorem~\ref{Thm:tf_smooth}~\ref{Cond:smooth_local_blowup},
	adding twisted fields is equivalent to a sequence of local blowups.
	In this example, on the one hand, $\fM^{\tf}_\La\big/\fM$ is isomorphic to blowing up $\fM$ along the locus $\{\ze_b\eq\ze_c\eq 0\}$, so particularly, $\big(\fM^{\tf}_\La\big)^{\tf}_{\ov\PT_{\La}(\La^\dag)}\big/\fM$ is not isomorphic to the trivial blowup.
	On the other hand,
	$\fM^{\tf}_{\La^\dag}\big/\fM$ is isomorphic to blowing up $\fM$ along the proper transforms of the Cartier divisors $\{\ze_b\eq 0\}$ and $\{\ze_c\eq 0\}$ successively,
	and $\big(\fM^{\tf}_{\La^\dag}\big)^{\tf}_{\ov\PT_{\La^\dag}(\La)}\big/\fM^{\tf}_{\La^\dag}$ is the trivial blowup,
	so $\big(\fM^{\tf}_{\La^\dag}\big)^{\tf}_{\ov\PT_{\La^\dag}(\La)}\big/\fM$ is isomorphic to the trivial blowup.
\end{eg}

The following two examples indicate in general, the derived trees $\dg{\ov\lt}$ do not form a treelike structure, yet the modified derived trees $\varrho_{\ov\lt}$ do.

\begin{eg}\label{Eg:Derived_CE}
	We continue with the setting of Example~\ref{Eg:PT_simple}.
	By Theorem~\ref{Thm:tf_smooth}~\ref{Cond:smooth_Z},
	we observe that
	$$
	\big(\,\{b,c\}\,,\,\ov\lt_1\,\big)_{\{c\}}=\big(\,\{b,c\}\,,\,\ov\lt_2\,\big)\,.
	$$
	If (\ref{Eqn:false_derived_tree_str}) were a treelike structure on $\fM^{\tf}_\La$,
	then $\dg{\ov\lt_2}$ would be isomorphic to $\dg{\ov\lt_1}\!\bsl\, c^\wedge$.
	
	By direct calculation,
	we find that
	\begin{align*}
		\dg{\ov\lt_1}=\big(\,\big\{\{b,c\},c\big\}\,,\,\preceq'\big),\qquad
		\tn{where}\ \preceq'\ \tn{satisfies}\ 
		\{b,c\}\!\not\prec'\!c,\ c\!\not\prec'\!\{b,c\};
	\end{align*}
	i.e.~$\dg{\ov\lt_1}$ is given by the leftmost graph of Figure~\ref{Fig:PT_simple} with the label $b$ replaced by $\{b,c\}$.
	Therefore, 
	$
	\dg{\ov\lt_1}\!\bsl\,c^\wedge\eq\tau_\bullet.
	$
	However,
	also by direct calculation we have
	$
	\bbI({\ov\lt_2})\eq\big\{\{b,c\}\big\},
	$
	hence 
	$
	\dg{\ov\lt_2}\!\ne\!\tau_\bullet.
	$
	Consequently, (\ref{Eqn:false_derived_tree_str}) cannot be a treelike structure.
\end{eg}

\begin{eg}
	In Example~\ref{Eg:Derived_CE},
	observe that
	$
	\ny_2(\ov\lt_1)\eq\emptyset$ and $
	\ny_2(\ov\lt_1)\eq\{b,c\}$.
	Therefore, we have
	$\vr_{\ov\lt_1}\eq\dg{\ov\lt_1}$,
	and $\vr_{\ov\lt_2}\eq\tau_\bullet\eq \vr_{\ov\lt_1}\!\bsl\,c^\wedge$.	
\end{eg}

The following example illustrates the notions of the derived, the doubly derived, and the second-order derived trees.
\begin{eg}
	\label{Eg:2nd_der_treelike}
	Let $\fM$, $\La$ and $\fM^{\tf}$ be as in Theorem~\ref{Thm:tf_smooth}.
	Assume there exits $(\al,\ov\lt)\inn \ov\La$ such that $\ov\lt\eq\big(\tau_\al,\ov\bE,\Dm(\ov\lt)\big)$ is given by the leftmost graph of Figure~\ref{Fig:2nd_derived}.
	Particularly,
	\begin{align*}
		&
		\fE_1=\{a,y\},\ \ 
		\fE_2=\{g,b,y\},\ \ 
		\fE_3=\{h,x,b,y\},\ \ 
		\fE_4=\{h,x,d,c,f,w,z\};\\
		&
		\ov\bE=\{\fE_1,\fE_2,\fE_3,\fE_4\};\quad
		\dot\fE_1=\{a\},\ \ 
		\dot\fE_2=\{g\},\ \ 
		\dot\fE_3=\{b,y\},\ \ 
		\dot\fE_4=\{d,c\}.
	\end{align*}
	We further assume $S_\al\eq \tau_\al$ and $\be_\al\eq\tn{Id}$.
	
	\begin{figure}[htp]
		\begin{center}
			\begin{tikzpicture}
				\filldraw
				(0,.1) circle (1.2pt)
				(0,.4) circle (1.2pt)
				(0,1) circle (1.2pt)
				(0,1.3) circle (1.2pt)
				(-.3,.1) circle (1.2pt)
				(.5,-.1) circle (1.2pt)
				(-1.1,-.1) circle (1.2pt)
				(-.8,-.1) circle (1.2pt)
				(-.3,.7) circle (1.2pt)
				(1.4,-.1) circle (1.2pt)
				(.9,.4) circle (1.2pt)
				(.9,-.1) circle (1.2pt)
				(1.2,-.4) circle (1.2pt)
				(.6,-.4) circle (1.2pt);
				
				\draw
				(0,.1)--(0,1)--(-1.1,-.1)
				(-.8,-.1)--(-.3,.7)
				(0,1)--(0,1.3)--(1.4,-.1)
				(.9,.4)--(.9,-.1)--(1.2,-.4)
				(.9,-.1)--(.6,-.4)
				(-.3,.1)--(0,.4)--(.5,-.1);
				
				\draw[dotted]
				(-1.2,.1)--(1.5,.1)
				(-1.2,.4)--(1.5,.4)
				(-1.2,.7)--(1.5,.7)
				(-1.2,1)--(1.5,1)
				(-1.2,1.3)--(1.5,1.3)
				;
				
				\draw[thin,decorate,decoration=brace]
				(1.5,.35)--(1.5,.15);
				\draw[thin,decorate,decoration=brace]
				(1.5,.65)--(1.5,.45);
				\draw[thin,decorate,decoration=brace]
				(1.5,.95)--(1.5,.75);
				\draw[thin,decorate,decoration=brace]
				(1.5,1.25)--(1.5,1.05)
				;
				
				\draw
				(-.1,1.15) node {\tiny{$a$}}
				(-.1,.55) node {\tiny{$b$}}
				(.87,.58) node {\tiny{$y$}}
				(1.25,.2) node {\tiny{$z$}}
				(.78,.2) node {\tiny{$w$}}
				(-.1,.2) node {\tiny{$c$}}
				(-.33,.255) node {\tiny{$d$}}
				(.34,.25) node {\tiny{$f$}}
				(-.27,.88) node {\tiny{$g$}}
				(-.93,.25) node {\tiny{$h$}}
				(-.625,-.02) node {\tiny{$x$}}
				(.81,-.325) node {\tiny{$u$}}
				(1.17,-.24) node {\tiny{$v$}}
				(1.5,1.15) node[right] {\tiny{$\fE_1$}}
				(1.5,.85) node[right] {\tiny{$\fE_2$}}
				(1.5,.55) node[right] {\tiny{$\fE_3$}}
				(1.5,.25) node[right] {\tiny{$\fE_4$}}
				(.15,-.75) node {\tiny{$\ov\lt=(\tau_\al,\ov\bE,\Dm(\ov\lt))$}}
				;
				
				\draw[loosely dashed]
				(-1.4,-1) rectangle (2.1,1.5);
				
				\filldraw[xshift=3.2cm] 
				(0,.1) circle (1.2pt)
				(0,.4) circle (1.2pt)
				(0,.7) circle (1.2pt)
				(0,1) circle (1.2pt) 
				(0,1.3) circle (1.2pt) 
				(.3,.1) circle (1.2pt)
				(.3,.7) circle (1.2pt)
				(1.2,.7) circle (1.2pt)
				(.9,1) circle (1.2pt)
				(.3,1) circle (1.2pt)
				(.6,.7) circle (1.2pt)
				(.9,.7) circle (1.2pt)
				;
				
				\draw[xshift=3.2cm]
				(0,.1)--(0,1.3)
				(.3,1)--(0,1.3)--(.9,1)
				(.3,.7)--(0,1)--(.6,.7)
				(1.2,.7)--(.9,1)--(.9,.7)
				(0,.4)--(.3,.1)
				(.1,.85) node[left] {\tiny{$\fE_2$}}
				(.1,.25) node[left] {\tiny{$\fE_4$}}
				(.1,.55) node[left] {\tiny{$\fE_3$}}
				(.1,1.15) node[left] {\tiny{$\fE_1$}}
				(.32,.28) node {\tiny{$f$}}
				(.8,.85) node {\tiny{$u$}}
				(1.2,.85) node {\tiny{$v$}}
				(.72,1.15) node {\tiny{$w$}}
				(.42,1.04) node {\tiny{$z$}}
				(.52,.85) node {\tiny{$x$}}
				(.23,.62) node {\tiny{$h$}}
				(.8,-.1) node {\tiny{$\dg{\ov\lt}$}}
				; 
				
				\filldraw[xshift=3.2cm,yshift=-2.4cm] 
				(0,.4) circle (1.2pt)
				(0,.7) circle (1.2pt)
				(0,1) circle (1.2pt) 
				(0,1.3) circle (1.2pt) 
				(.3,.7) circle (1.2pt)
				(1.2,.7) circle (1.2pt)
				(.9,1) circle (1.2pt)
				(.3,1) circle (1.2pt)
				(.6,.7) circle (1.2pt)
				(.9,.7) circle (1.2pt)
				;
				
				\draw[xshift=3.2cm,yshift=-2.4cm]
				(0,.4)--(0,1.3)
				(.3,1)--(0,1.3)--(.9,1)
				(.3,.7)--(0,1)--(.6,.7)
				(1.2,.7)--(.9,1)--(.9,.7)
				(.1,.85) node[left] {\tiny{$\fE_2$}}
				(.1,.55) node[left] {\tiny{$\fE_3$}}
				(.1,1.15) node[left] {\tiny{$\fE_1$}}
				(.8,.85) node {\tiny{$u$}}
				(1.2,.85) node {\tiny{$v$}}
				(.72,1.15) node {\tiny{$w$}}
				(.42,1.04) node {\tiny{$z$}}
				(.52,.85) node {\tiny{$x$}}
				(.23,.62) node {\tiny{$h$}}
				(.8,.2) node {\tiny{$\vr_{\ov\lt}$}}
				; 
				
				\filldraw[xshift=6.8cm]
				(.2,.5) circle (1.2pt)
				(-.6,.9) circle (1.2pt)
				(-.8,.5) circle (1.2pt)
				(0,1.3) circle (1.2pt)
				(-.2,.1) circle (1.2pt)
				(.2,.1) circle (1.2pt)
				(-1,.1) circle (1.2pt)
				(-.6,.1) circle (1.2pt)
				(.6,.1) circle (1.2pt)
				(1,.1) circle (1.2pt)
				;
				
				\draw[xshift=6.8cm]
				(0,1.3)--(-.6,.9)--(-1,.1)
				(.2,.5)--(0,1.3)--(1,.1)
				(-.6,.1)--(-.6,.9)--(-.2,.1)
				(.6,.1)--(.2,.5)--(.2,.1)
				;
				
				\draw[ultra thick, xshift=6.8cm]
				(-.8,.5)--(-1,.1)
				(0,1.3)--(-.6,.9)--(-.2,.1)
				(.2,.1)--(.2,.5)
				;
				
				\draw[dotted,xshift=6.8cm]
				(-1.1,.1)--(1.1,.1)
				(-1.1,.5)--(1.1,.5)
				(-1.1,.9)--(1.1,.9)
				(-1.1,1.3)--(1.1,1.3)
				;
				
				\draw[thin,decorate,decoration=brace,xshift=6.8cm]
				(1.1,1.25)--(1.1,.95);
				\draw[thin,decorate,decoration=brace,xshift=6.8cm]
				(1.1,.85)--(1.1,.55);
				\draw[thin,decorate,decoration=brace,xshift=6.8cm]
				(1.1,.45)--(1.1,.15);
				
				\draw[xshift=6.8cm]
				(0,.7) node {\tiny{$w$}}
				(.37,.2) node {\tiny{$v$}}
				(.75,.25) node {\tiny{$z$}}
				(.07,.25) node {\tiny{$u$}}
				(-.41,.25) node {\tiny{$x$}}
				(-.62,1.1) node {\tiny{$\fE_1$}}
				(-1.15,.26) node {\tiny{$\fE_3$}}
				(-.9,.7) node {\tiny{$\fE_2$}}
				(-.7,.25) node {\tiny{$h$}}
				(1.1,1.1) node[right] {\tiny{$\fF_3$}}
				(1.1,.7) node[right] {\tiny{$\fF_2$}}
				(1.1,.3) node[right] {\tiny{$\fF_1$}}
				(0,-.26) node {\tiny{$\ud\ls$}}
				;
				
				\draw[loosely dashed,xshift=6.8cm]
				(-1.5,-.5) rectangle (1.7,1.5);
				
				\filldraw[xshift=6.8cm,yshift=-2.4cm]
				(0,.1) circle (1.2pt)
				(0,.7) circle (1.2pt)
				(0,1.3) circle (1.2pt)
				(-.3,.1) circle (1.2pt)
				(.4,-.1) circle (1.2pt)
				(-1,-.1) circle (1.2pt)
				(-.6,.1) circle (1.2pt)
				(-.3,.7) circle (1.2pt)
				(.7,-.1) circle (1.2pt)
				(.3,.7) circle (1.2pt)
				(1,-.1) circle (1.2pt)
				;
				
				\draw[xshift=6.8cm,yshift=-2.4cm]
				(-.6,.1)--(-.3,.7)--(-1,-.1)
				(.7,-.1)--(.3,.7)--(1,-.1)
				(0,.1)--(0,1.3)
				(-.3,.1)--(0,.7)--(.4,-.1)
				(-.3,.7)--(0,1.3)--(.3,.7)
				;
				
				\draw[dotted,xshift=6.8cm,yshift=-2.4cm]
				(-1.1,.1)--(1.1,.1)
				(-1.1,.7)--(1.1,.7)
				(-1.1,1.3)--(1.1,1.3)
				;
				
				\draw[thin,decorate,decoration=brace,xshift=6.8cm,yshift=-2.4cm]
				(1.1,1.25)--(1.1,.75);
				\draw[thin,decorate,decoration=brace,xshift=6.8cm,yshift=-2.4cm]
				(1.1,.65)--(1.1,.15);
				
				\draw[xshift=6.8cm,yshift=-2.4cm]
				(-.08,.85) node {\tiny{$b$}}
				(.34,.85) node {\tiny{$y$}}
				(.85,.25) node {\tiny{$z$}}
				(.4,.25) node {\tiny{$w$}}
				(-.1,.2) node {\tiny{$c$}}
				(-.33,.255) node {\tiny{$d$}}
				(.25,-.05) node {\tiny{$f$}}
				(-.33,.85) node {\tiny{$g$}}
				(-.88,.25) node {\tiny{$h$}}
				(-.625,-.02) node {\tiny{$x$}}
				(1.1,1) node[right] {\tiny{$\fE_2$}}
				(1.1,.4) node[right] {\tiny{$\fE_4$}}
				(0,-.45) node {\tiny{$\ov\lt_{(J)}$}}
				;
				
				\draw[loosely dashed,xshift=6.8cm,yshift=-2.4cm]
				(-1.5,-.7) rectangle (1.7,1.5);
				
				\filldraw[xshift=9.6cm] 
				(0,.7) circle (1.2pt)
				(0,1) circle (1.2pt) 
				(0,1.3) circle (1.2pt) 
				(.3,.7) circle (1.2pt)
				(.9,1) circle (1.2pt)
				(.3,1) circle (1.2pt)
				(.9,.7) circle (1.2pt)
				;
				
				\draw[xshift=9.6cm]
				(0,.7)--(0,1.3)
				(.3,1)--(0,1.3)--(.9,1)
				(.3,.7)--(0,1)--(.9,.7)
				(.1,.85) node[left] {\tiny{$\fE_4$}}
				(.1,1.15) node[left] {\tiny{$\fE_2$}}
				(.72,1.15) node {\tiny{$w$}}
				(.42,1.04) node {\tiny{$z$}}
				(.72,.6) node {\tiny{$f$}}
				(.23,.62) node {\tiny{$h$}}
				(.8,.1) node {\tiny{$\dg{\ov\lt_{(J)}}$}}
				; 
				
				\filldraw[xshift=9.6cm,yshift=-2.4cm] 
				(0,1) circle (1.2pt) 
				(0,1.3) circle (1.2pt) 
				(.9,1) circle (1.2pt)
				(.3,1) circle (1.2pt)
				;
				
				\draw[xshift=9.6cm,yshift=-2.4cm]
				(0,1)--(0,1.3)
				(.3,1)--(0,1.3)--(.9,1)
				(.1,1.15) node[left] {\tiny{$\fE_2$}}
				(.72,1.15) node {\tiny{$w$}}
				(.42,1.04) node {\tiny{$z$}}
				(.8,.5) node {\tiny{$\tau_{\wh\al}$}}
				; 
				
				\filldraw[xshift=11.9cm] 
				(0,.4) circle (1.2pt)
				(0,.7) circle (1.2pt)
				(0,1) circle (1.2pt) 
				(0,1.3) circle (1.2pt) 
				(.3,.7) circle (1.2pt)
				(1.2,.7) circle (1.2pt)
				(.9,1) circle (1.2pt)
				(.3,1) circle (1.2pt)
				(.6,.7) circle (1.2pt)
				;
				
				\draw[xshift=11.9cm]
				(0,.4)--(0,1.3)
				(.3,1)--(0,1.3)--(.9,1)
				(.3,.7)--(0,1)--(.6,.7)
				(1.2,.7)--(.9,1)
				(.1,.85) node[left] {\tiny{$\fF_2$}}
				(.1,.55) node[left] {\tiny{$\fF_1$}}
				(.1,1.15) node[left] {\tiny{$\fF_3$}}
				(1.2,.85) node {\tiny{$v$}}
				(.72,1.15) node {\tiny{$w$}}
				(.42,1.04) node {\tiny{$z$}}
				(.82,.62) node {\tiny{$\fE_2$}}
				(.23,.62) node {\tiny{$h$}}
				(.8,.1) node {\tiny{$\dg{\ud\ls}$}}
				; 
				
				\filldraw[xshift=11.9cm,yshift=-2.4cm] 
				(0,1) circle (1.2pt) 
				(0,1.3) circle (1.2pt) 
				(.9,1) circle (1.2pt)
				(.3,1) circle (1.2pt)
				;
				
				\draw[xshift=11.9cm,yshift=-2.4cm]
				(0,1)--(0,1.3)
				(.3,1)--(0,1.3)--(.9,1)
				(.1,1.15) node[left] {\tiny{$\fF_3$}}
				(.72,1.15) node {\tiny{$w$}}
				(.42,1.04) node {\tiny{$z$}}
				(.8,.5) node {\tiny{$\vs_{\wh\al}$}}
				; 
			\end{tikzpicture}		
		\end{center}
		\caption{Doubly derived treelike structure and second-order derived treelike structure}\label{Fig:2nd_derived}
	\end{figure}
	
	In this case,
	$\bbI(\ov\lt)\eq\{\fE_1,\fE_2,\fE_3,\fE_4,h,x,f,w,u,v,z\}$,
	and $\Dm^*(\ov\lt)\eq\{a,b,c,d\}$, so  
	$\ny_2(\ov\lt)\eq\{\fE_4\}$.
	The derived tree $\dg{\ov\lt}$ is thus determined by the tree order (\ref{Eqn:derived_order}):
	\begin{align*}
		\fE_4\prec\fE_3\prec\fE_2\prec\fE_1,\qquad
		h,\,x\prec\fE_1,\qquad
		f\prec\fE_3\prec\fE_2\prec\fE_1,\qquad
		u,\,v\prec w,
	\end{align*}
	and $\vr_{\ov\lt}\eq\dg{\ov\lt}\bsl \fE_4^\wedge\eq\dg{\ov\lt}\bsl\{\fE_4,f\}$;
	c.f.~Figure~\ref{Fig:2nd_derived}.
	
	Consider the LRS with dominant edges $\ud\ls\eq\big(\vr_{\ov\lt},\ud\bF,\Dm(\ud\ls)\big)$ given by the upper graph in the middle of Figure~\ref{Fig:2nd_derived}.
	In particular,
	\begin{align*}
		&
		\fF_1=\{\fE_3,h,x,u,v,z\},\quad
		\fF_2=\{\fE_2,h,x,w,z\},\quad
		\fF_3=\{\fE_1,w,z\};\\
		&
		\ud\bF=\{\fF_1,\fF_2,\fF_3\};\quad
		\udt\fF_1=\{\fE_3,u\},\ \ 
		\udt\fF_2=\{x\},\ \
		\udt\fF_3=\{\fE_1\};\quad
		J:=\Dm(\ud\ls)=\{\fE_1,\fE_3,x,u\}.
	\end{align*}
	With the subset $J\!\subset\!\bbI(\ov\lt)$ given above, we obtain the RLS with dominant edges $\ov\lt_{(J)}$ and its derived tree $\dg{\ov\lt_{(J)}}$;
	see Figure~\ref{Fig:2nd_derived}.
	
	For such $\ov\lt_{(J)}$,
	observe that $\Dm^*(\ov\lt_{(J)})\eq\{b,g,c,d,x\}$,
	hence $\ny_3(\ov\lt_{(J)})\eq\{\fE_4\}$.
	Therefore,
	the rooted tree
	$\tau_{\wh\al}$ of Proposition~\ref{Prp:Derived_treelike_2} is $\dg{\ov\lt_{(J)}}\bsl\{\fE_4\}^\wedge\eq 
	\dg{\ov\lt_{(J)}}\bsl\{\fE_4,f,h\}$ in this case;
	c.f.~Figure~\ref{Fig:2nd_derived}.
	
	Meanwhile,
	notice 
	$\bbI(\ud\ls)\eq\{\fF_1,\fF_2,\fF_3,\fE_2,h,w,v,z\}$,
	and $\Dm^*(\ud\ls)\eq\{\fE_1,x\}$,
	so the derived tree $\dg{\ud\ls}$ is determined by the tree order (\ref{Eqn:derived_order_revert}):
	\begin{align*}
		\fE_1\prec\fE_2\prec\fE_3,\qquad
		\fE_2,\,h\prec\fE_1,\qquad
		v\prec w;
	\end{align*}
	c.f.~Figure~\ref{Fig:2nd_derived}.
	Moreover, we have $\ny_2^{\ov\lt}(\ud\ls)\eq\{\fF_1,\fF_2\}$
	and $\nu^{\ov\lt}(\ud\ls)\eq\{v\}$,
	hence
	the rooted tree
	$\vs_{\wh\al}$ of Proposition~\ref{Prp:Derived_treelike_second_order_1} is $\dg{\ud\ls}\bsl\{\fF_1,\fF_2,v\}^\wedge\eq 
	\dg{\ud\ls}\bsl\{\fF_1,\fF_2,\fE_2,h,v\}$;
	c.f. the last graph of~Figure~\ref{Fig:2nd_derived}.
\end{eg}

The next example demonstrates the roles that the derived, doubly derived, and second-order derived trees of Example~\ref{Eg:2nd_der_treelike} play in diagonalizing (in the sense of~\cite[Definition 3.2]{HL11}) the structural homomorphism $\varphi$ of (\ref{Eqn:str_hom}).

\begin{eg}
	\label{Eg:2nd_der_treelike_phi}
Let $\fM^{\fk 1}\eq\fD_2$ be endowed with the LES as in \S\ref{Subsec:Step1}.
Consider the stratum $\fM^{\fk 1}_\al$ consisting of the following $(C,D)$: the dual graph of $C$ is given by the underlying rooted tree $\tau_\al$ of $\ov\lt$ in Figure~\ref{Fig:2nd_derived},
whose root vertex corresponds to $C_\core$ (so particularly, there are two tails in the sense of \cite[Definition 2.1]{HLN} attached to $C_\core$, which is smooth,
at two nodes labeled by $a$ and~$y$),
and the points of $D$ are all on the rational components corresponding to the {\it leaf} vertices of $\tau_\al$.
These rational components are respectively
denoted by $C_e$, $e\eq h,x,d,c,f,u,v,z$,
satisfying $e\inn C_e$.
Instead of labeling the points of $D$ by $1,\ldots,m$, we write
\begin{align*}
	D\cap C_e:=\{\de_{e_1},\cdots,\de_{e_{m_e}}\},\qquad
	e\eq h,x,d,c,f,u,v,z\,,
\end{align*}
in this example for convenience.

In this case, we have $S_\al\eq\tau_\al$ (as a set).
The treelike structure $\La^{\fk 1}$ of Proposition~\ref{Prp:G1-level} assigns the rooted tree $\tau_\al$, along with the identity map $\be_\al\eq\tn{Id}$, to the stratum $\fM^{\fk 1}_\al$.
The $\ov\lt$ of Example~\ref{Eg:2nd_der_treelike} determines $\al^{\fk 2}\eq(\al,\ov\lt)\inn A^{\fk 2}$. 
By Corollary~\ref{Crl:Step2_str_hom},
we have
\begin{align*}
	&
	\ti\varphi^{(\fk 2)}=
	\xi^{(\fk 2)}_{\fE_1}
	\xi^{(\fk 2)}_{\fE_2}
	\xi^{(\fk 2)}_{\fE_3}
	\xi^{(\fk 2)}_{\fE_4}
	\left[
	\begin{matrix}\,
		1 & \mathbf{0} & \mathbf{0} &
		\mathbf{0} & \mathbf{0} & \mathbf{0}
		& \mathbf{0} & \mathbf{0}
		\\
		0 & \bv_c & \bv_d &
		\bv_x & \bv_h & \bv_u 
		& \bv_v & \bv_z
	\,\end{matrix}
	\right]\,,\qquad\tn{where}\quad 
	\bv_z=\big[\ti\ka_{c_1z_i}^{(\fk 2)}
	\xi^{(\fk 2)}_{z}\big]_{1\le i\le m_z},
	\\
	&
	\bv_c=\big[\ti\ka_{c_1c_i}^{(\fk 2)}\xi^{(\fk 2)}_{\fE_1}
	\xi^{(\fk 2)}_{\fE_2}
	\xi^{(\fk 2)}_{\fE_3}
	\xi^{(\fk 2)}_{\fE_4}\big]_{2\le i\le m_c},\quad
	\bv_d=\big[\ti\ka_{c_1d_i}^{(\fk 2)}\xi^{(\fk 2)}_{\fE_1}
	\xi^{(\fk 2)}_{\fE_2}
	\xi^{(\fk 2)}_{\fE_3}\big]_{1\le i\le m_d},\quad
	\bv_x=\big[\ti\ka_{c_1x_i}^{(\fk 2)}\xi^{(\fk 2)}_{\fE_1}
	\xi^{(\fk 2)}_{x}\big]_{1\le i\le m_x},
	\\
	&
	\bv_h=\big[\ti\ka_{c_1h_i}^{(\fk 2)}\xi^{(\fk 2)}_{\fE_1}
	\xi^{(\fk 2)}_{h}\big]_{1\le i\le m_h},\quad
	\bv_u=\big[\ti\ka_{c_1u_i}^{(\fk 2)}\xi^{(\fk 2)}_{w}
	\xi^{(\fk 2)}_{u}\big]_{1\le i\le m_u},\quad
	\bv_v=\big[\ti\ka_{c_1v_i}^{(\fk 2)}\xi^{(\fk 2)}_{w}
	\xi^{(\fk 2)}_{v}\big]_{1\le i\le m_v}.
\end{align*}
In the above equality, 
each $\ti\ka_{c_1e_i}^{(\fk 2)}$ denotes the pullback of the corresponding $\ka_{c_1e_i}$, which is the same as its proper transform.
If $m_c\eq 1$, then the columns containing the entries of $\bv_c$ should be omitted.

Since $C_\core$ is smooth and $\wp_{\al;c_1},\wp_{\al;d_1}\inn\Dm(\ov\lt)$,
the treelike structures $\La^{r}$, $\fk 2\!\le\!r\!\le\!\fk 4$ all assign the trivial tree to this stratum,
i.e.
\begin{align*}
	\al^{\fk 5}\eq \big(((\al,\ov\lt),\ov\lt_\bullet),\cdots,\ov\lt_\bullet\big)\,.
\end{align*}
Therefore,  the proper transform of $\varrho_{\ov\lt}$ (c.f.~Proposition~\ref{Prp:two_treelike_strs}) is  equal to $\varrho_{\ov\lt}$ itself,
which is assigned to the stratum $\fM^{\fk 5}_{\al^{\fk 5}}$ by~$\La^{\fk 5}$.
With $\ud\ls$ as in Example~\ref{Eg:2nd_der_treelike},
consider $\al^{\fk 6}\eq(\al^{\fk 5},\ud\ls).$
By Theorem~\ref{Thm:tf_smooth_revert}~\ref{Cond:smooth_parameters_revert},
locally near $\fM^{\fk 6}_{\al^{\fk 6}}$,
$\ti\varphi^{(\fk 2)}$ pulls back to
\begin{align*}
	&
	\ti\varphi^{(\fk 6)}=
	\xi^{(\fk 6)}_{\fF_1}
	\xi^{(\fk 6)}_{\fF_2}
	\xi^{(\fk 6)}_{\fF_3}
	\xi^{(\fk 6)}_{\fE_2}
	\xi^{(\fk 6)}_{\fE_4}
	\left[
	\begin{matrix}\,
		1 & 0
		\\
		0 & \xi^{(\fk 6)}_{\fF_1}
		\xi^{(\fk 6)}_{\fF_2}
		\xi^{(\fk 6)}_{\fF_3}
		\,\end{matrix}
	\right]
	\left[
	\begin{matrix}\,
		1 & \mathbf{0} & \mathbf{0} &
		\mathbf{0} & \mathbf{0} & \mathbf{0}
		& \mathbf{0} & \mathbf{0}
		\\
		0 & \bv_c' & \bv_d' &
		\bv_x' & \bv_h' & \bv_u' 
		& \bv_v' & \bv_z'
		\,\end{matrix}
	\right]\,,\quad\tn{where}
	\\
	&
	\bv_d'\eq\big[\ti\ka_{c_1d_i}^{(\fk 6)}
	\xi^{(\fk 6)}_{\fE_2}\big]_{1\le i\le m_d},\quad
	\bv_x'\eq\big[\ti\ka_{c_1x_i}^{(\fk 6)}
	\big]_{1\le i\le m_x},
	\quad
	\bv_h'\eq\big[\ti\ka_{c_1h_i}^{(\fk 6)}
	\xi^{(\fk 6)}_{h}\big]_{1\le i\le m_h},
	\quad
	\bv_u'\eq\big[\ti\ka_{c_1u_i}^{(\fk 6)}\xi^{(\fk 6)}_{w}\big]_{1\le i\le m_u},
	\\
	&
	\bv_c'\eq\big[\ti\ka_{c_1c_i}^{(\fk 6)}
	\xi^{(\fk 6)}_{\fE_2}
	\xi^{(\fk 6)}_{\fE_4}\big]_{2\le i\le m_c},\quad
	\bv_v'\eq\big[\ti\ka_{c_1v_i}^{(\fk 6)}\xi^{(\fk 6)}_{w}
	\xi^{(\fk 6)}_{v}\big]_{1\le i\le m_v},
	\quad
	\bv_z\eq\big[\ti\ka_{c_1z_i}^{(\fk 6)}
	\xi^{(\fk 6)}_{z}\big]_{1\le i\le m_z}.
\end{align*}
Here, each pullback $\ti\ka_{c_1e_i}^{(\fk 6)}$ of the corresponding $\ka_{c_1e_i}$ is still equal to its proper transform because $C_\core$ is smooth.
The above expression of $\ti\varphi^{(\fk 6)}$ is consistent with Corollary~\ref{Crl:Step6_str_hom}.

If  $\ti\ka_{c_1e_i}^{(\fk 6)}$ is invertible near a point $\ti y^{(\fk 6)}\inn\fM^{\fk 6}_{\al^{\fk 6}}$,
then $\bv'_x$ contains an invertible entry,
which implies $\ti\varphi^{(\fk 6)}$ is diagonalized near $\ti y^{(\fk 6)}$,
and so is its pullback to $\fM^{\fk{10}}$ near any lift of $\ti y^{(\fk 6)}$.

If $\ti\ka_{c_1e_i}^{(\fk 6)}$ vanishes at a point $\ti y^{(\fk 6)}\inn\fM^{\fk 6}_{\al^{\fk 6}}$,
then it is a direct check that $\ti y^{(\fk 6)}\inn\fM^{\fk 6}_{\ex;\al^{\fk 6}}$;
c.f.~(\ref{Eqn:grafted_boundary'}) and Corollary~\ref{Crl:K_step6} for notation.
This implies
$\ka_{c_1x_1}(C,D)\!=\!0$,
so the node $a$ is at a Weierstrass point of $C_\core$;
c.f.~\cite[Lemma 2.15]{HLN}.
Consequently, each $\ka_{c_1e_i}$, $e\eq u,v,z$, $1\!\le\!i\!\le\!m_e$,
is invertible,
and so is its pullback.
By (\ref{Eqn:local_eqn_remainder1}) and (\ref{Eqn:local_eqn_remainder3}),
we have
\begin{align*}
	&
	\ti\ka_{c_1c_i}^{(\fk 6)}=
	f_{c_i}\ti\ka_{c_1x_1}^{(\fk 6)}
	+ 
	g_{c_i}\xi^{(\fk 6)}_{\fF_3},
	\quad
	2\!\le\!i\!\le\!m_c,
	\qquad
	\ti\ka_{c_1d_i}^{(\fk 6)}=
	f_{d_i}\ti\ka_{c_1x_1}^{(\fk 6)}
	+ 
	g_{d_i}\xi^{(\fk 6)}_{\fF_3},
	\quad
	1\!\le\!i\!\le\!m_d,
	\\
	&
	\ti\ka_{c_1x_i}^{(\fk 6)}=
	f_{x_i}\ti\ka_{c_1x_1}^{(\fk 6)}
	+ 
	g_{x_i}\!\cdot\!\big(\xi^{(\fk 6)}_{\fF_1}\xi^{(\fk 6)}_{\fF_2}\big)^2\xi^{(\fk 6)}_{\fF_3}
	\xi^{(\fk 6)}_{\fE_2}\xi^{(\fk 6)}_{\fE_4},
	\quad
	2\!\le\!i\!\le\!m_x,
	\\	
	&
	\ti\ka_{c_1h_i}^{(\fk 6)}=
	f_{h_i}\ti\ka_{c_1x_1}^{(\fk 6)}
	+ 
	g_{h_i}\xi^{(\fk 6)}_{\fF_2}\xi^{(\fk 6)}_{\fF_3}\xi^{(\fk 6)}_{\fE_2},
	\quad
	1\!\le\!i\!\le\!m_h,
\end{align*}
where $f_{e_i}$ and $g_{e_i}$ are all invertible functions.
So, after suitable elementary column operations,
we can rewrite $\ti\varphi^{(\fk 6)}$ as
\begin{align*}
	&
	\xi^{(\fk 6)}_{\fF_1}
	\xi^{(\fk 6)}_{\fF_2}
	\xi^{(\fk 6)}_{\fF_3}
	\xi^{(\fk 6)}_{\fE_2}
	\xi^{(\fk 6)}_{\fE_4}
	\left[
	\begin{matrix}\,
		1 & 0
		\\
		0 & \xi^{(\fk 6)}_{\fF_1}
		\xi^{(\fk 6)}_{\fF_2}
		\xi^{(\fk 6)}_{\fF_3}
		\,\end{matrix}
	\right]
	\left[
	\begin{matrix}\,
		1 & 0 & 0 &
		0 & 0 & 0 & \cdots & 0
		\\
		0 & \ti\ka_{c_1x_1}^{(\fk 6)} & \xi^{(\fk 6)}_{\fF_3}\xi^{(\fk 6)}_{\fE_2} &
		\xi^{(\fk 6)}_{w} & \xi^{(\fk 6)}_{z}
		& 0 & \cdots & 0
		\,\end{matrix}
	\right]\,.
\end{align*}
The treelike structures $\wc\Up^{\fk 6}$ and $\wc\Up^{\fk 7}$ as in Proposition~\ref{Prp:r3p3_treelike_proto}
assigns to $\fM^{\fk 6}_{\al^{\fk 6}}$ the rooted tree $\tau_{\wh\al}$ and $\varsigma_{\wh\al}$ as in Figure~\ref{Fig:2nd_derived}, respectively.
Therefore,
$\wc\Up^{\fk 6}_\ex$ and $\wc\Up^{\fk 7}_\ex$ assigns $\ex(\tau_{\wh\al})$ and $\ex(\varsigma_{\wh\al})$ to $\fM^{\fk 6}_{\ex;\al^{\fk 6}}$, respectively;
c.f.~Figure~\ref{Fig:grafted}.
The extra edge corresponds to the local parameter $\ti\ka_{c_1x_1}^{(\fk 6)}$.
Consequently,
the pullback of $\ti\varphi^{(\fk 6)}$ to $\fM^{\fk 8}$ (hence to $\fM^{\fk{10}}$) is diagonalized near any lift of $\ti y^{(\fk 6)}$. 
\end{eg}

\end{document}